\documentclass[aps,physrev,reprint,superscriptaddress,twocolumn,longbibliography]{revtex4-2}
\usepackage{amsmath,amssymb,amsthm,mathtools,booktabs,microtype,graphicx}
\usepackage[colorlinks=true,linkcolor=blue,citecolor=blue,urlcolor=blue]{hyperref}
\newtheorem{theorem}{Theorem}

\DeclareMathOperator{\Cov}{Cov}

\newcommand{\R}{\mathbb R}

\newcommand{\E}{\mathbb E}
\newcommand{\Prob}{\mathbb P}
\newcommand{\wh}{\widehat}
\newcommand{\wt}{\widetilde}
\hypersetup{pdftitle={Automatic Threshold Selection Leaves Finite-Degree Bias in Network Power-Law Fits},pdfauthor={Tiandong Wang}}
\makeatletter
\newcommand{\arxiveqref}[1]{\textup{\tagform@{\ref*{#1}}}}
\makeatother

\begin{document}
\title{Automatic Threshold Selection Leaves Finite-Degree Bias in Network Power-Law Fits}
\author{Tiandong Wang}
\email{td\_wang@fudan.edu.cn}
\affiliation{Shanghai Center for Mathematical Sciences, Fudan University, Shanghai 200438, China}
\affiliation{Shanghai Academy of Artificial Intelligence for Science, Shanghai, China}
\date{November 7, 2026}

\begin{abstract}
Minimum-distance fitting selects a degree threshold and a power-law exponent from the same network. We study the procedure that fits a continuous Pareto law with the Hill estimator and minimizes the Kolmogorov--Smirnov distance over retained ranks. For homogeneous linear preferential attachment and Norros--Reittu graphs with Pareto weights, we derive the joint limiting law of the selected rank, retained fraction of the threshold tie, and exponent error. The selector retains of order $n^{2/(\alpha+2)}$ of $n$ vertices, where $\alpha$ is the survival exponent. At this rank, integer spacing and finite-degree departures from a Pareto tail have the same order as sampling fluctuations. The exponent estimate is consistent and retains a shift on the scale used for confidence intervals. For integer degree laws whose tail sums have a first-order inverse-degree expansion, adding the retained tie fraction divided by the threshold to the Hill statistic removes the leading contribution of the tie. Under homogeneous preferential attachment with a fixed known number of edges per new vertex, removing the remaining model-specific term and calibrating for selection gives confidence intervals with asymptotically exact coverage. Observed networks show sensitivity of the fitted exponent and the classification of degree variance to threshold ties. Simulations assess interval coverage and show improved pooled connectivity predictions after correction under both specified graph laws. Changes to the Norros--Reittu weight distribution can reverse the prediction gain.
\end{abstract}
\maketitle

\section{Network degree tails and threshold selection}\label{sec:intro}
Fitting a power law to network degrees requires choosing a threshold and estimating an exponent from the retained vertices~\cite{ClausetEtAl2009,Newman2005}. An automatic rule chooses both from the same degree sequence, so their joint distribution determines the uncertainty of the reported exponent~\cite{DreesEtAl2020}. Integer degrees add ties at the threshold, and network growth creates dependence between vertices. These features enter the sampling distribution needed for confidence intervals after threshold selection.

For a sequence of graphs with degrees $D_1(n),\ldots,D_n(n)$, suppose that the proportion of vertices of degree $j$ converges in probability to $\pi_j$ for every $j\ge0$, with $\sum_{j\ge0}\pi_j=1$. Write $\overline F(x)=\sum_{j>x}\pi_j$. Following Ref.~\cite{VoitalovEtAl2019}, we call this limiting degree distribution scale-free with survival exponent $\alpha>0$ when it is regularly varying~\cite{deHaanFerreira2006},
\begin{equation}\label{eq:regular-variation}
 \lim_{x\to\infty}\frac{\overline F(tx)}{\overline F(x)}=t^{-\alpha},\qquad t>0.
\end{equation}
Because the definition concerns a limiting degree law, it permits departures from an exact Pareto survival function at finite degrees.

The survival exponent determines whether the limiting degree distribution has a finite second moment: it is finite for $\alpha>2$ and infinite for $\alpha<2$. Network science often reports the exponent $\alpha+1$ of the degree probabilities, so the boundary is $3$ in that convention. For uncorrelated random networks with finite mean degree, a divergent second moment gives a vanishing susceptible-infected-susceptible epidemic threshold in the heterogeneous mean-field approximation~\cite{PastorSatorrasVespignani2001} and resilience to random removal of any fixed fraction of vertices smaller than one in the infinite-size percolation limit~\cite{CohenEtAl2000}. The relation between a degree tail and a physical observable depends on the graph ensemble and dynamics~\cite{CastellanoPastorSatorras2020}. Uncertainty in the fitted exponent therefore matters when using a degree distribution to classify these regimes.

For regularly varying degree distributions, Voitalov et al.~\cite{VoitalovEtAl2019} used estimators from extreme value theory that are consistent for independent observations from every regularly varying law. Integer degrees make their estimates unstable, so they add independent uniform noise on $[-1/2,1/2]$ to each degree. Their classification also lacks confidence intervals for the exponent, because the convergence rates are unknown. Consistency leaves open the distribution needed to calibrate uncertainty after choosing a threshold from the degree sequence.

We examine this question for the Hill estimator~\cite{Hill1975}, which fits a continuous Pareto law to the largest observations, and minimum-distance threshold selection~\cite{ClausetEtAl2009,DreesEtAl2020}. Varying the cutoff changes the number of retained observations and the finite-degree departures from the fitted law. Selection also couples the cutoff to the fluctuations of the estimate. To specify this dependence, let $D_{(1)}(n)\ge\cdots\ge D_{(n)}(n)$ be the ordered degrees and define the Hill statistic and its reciprocal by
\begin{equation}\label{eq:hill-stat}
 H_{n,k}=\frac1k\sum_{i=1}^k\log\frac{D_{(i)}(n)}{D_{(k+1)}(n)},\qquad
 \wh\alpha_{n,k}=H_{n,k}^{-1}.
\end{equation}
A rank is valid when $D_{(k+1)}(n)>0$ and $H_{n,k}>0$. The minimum-distance selection procedure (MDSP) of Clauset et al.~\cite{ClausetEtAl2009}, in the rank-based form studied by Drees et al.~\cite{DreesEtAl2020}, selects a valid rank that minimizes the Kolmogorov--Smirnov distance between the empirical survival function of the ratios $D_{(i)}(n)/D_{(k+1)}(n)$, $i\le k$, and the fitted Pareto survival function,
\begin{equation}\label{eq:mdsp-def}
 \mathcal D_k=\sup_{y\ge1}\left|\frac1k\sum_{i=1}^k
 \mathbf1\{D_{(i)}(n)/D_{(k+1)}(n)>y\}-y^{-\wh\alpha_{n,k}}\right|.
\end{equation}
The unrestricted selector $\wh k$ is the smallest minimizer over $2\le k<n$, with $\mathcal D_k=+\infty$ at invalid ranks. If all ranks are invalid, then $\wh k=2$ and no point estimate is reported. Here $k$ counts the retained vertices and $D_{(k+1)}(n)$ is their threshold. A restricted search minimizes $\mathcal D_k$ over a set of ranks fixed in advance. Several ranks can share an integer degree threshold, so this procedure can retain part of the tie at that degree. Broido and Clauset~\cite{BroidoClauset2019} instead fit a discrete Pareto likelihood and search over distinct degree thresholds, retaining each tie whole. Our selection laws concern the rank search in Eq.~\eqref{eq:mdsp-def}. Section~\ref{sec:real} compares it with the discrete fit on the degree sequences used by Broido and Clauset.

For independent continuous observations with regularly varying tails, Bhattacharya et al.~\cite{BhattacharyaEtAl2020} proved consistency of the exponent estimate at the MDSP rank. For iid Pareto observations, Drees et al.~\cite{DreesEtAl2020} derived a joint limit conditional on uniqueness of the limiting minimizer and analyzed the criterion under specified departures from a Pareto tail. Integer data raise a separate problem. Continuous fits to rounded observations are less accurate than discrete fits~\cite{ClausetEtAl2009}, Hill estimates oscillate across ties~\cite{MatsuiEtAl2013}, and integer tails motivate estimators adapted to their support~\cite{BertailEtAl2025}. None of these results specifies the size of the oscillation at the rank that MDSP selects.

For linear preferential attachment (PA), Drees et al.~\cite{DreesEtAl2020} found numerically that MDSP performance depends strongly on the model parameters. Their network analysis did not determine the scale of the selected rank or the contributions of threshold ties and finite-degree departures at that rank. Wang and Resnick~\cite{WangResnick2019} proved Hill consistency at deterministic ranks in a model with one edge per new vertex, and Wang~\cite{Wang2026} obtained asymptotic normality at ranks $k=o(n^{2/(\alpha+2)})$, where integer spacing is negligible relative to sampling fluctuations. Cirkovic, Wang, and Cline~\cite{CirkovicWangCline2025} established asymptotic normality of Hill estimates based on Norros--Reittu degrees at deterministic ranks under conditions on the weight tail and the growth of the rank. They left the theory of data-driven threshold selection to future work. Applying a normal approximation to MDSP therefore requires determining whether its selected rank lies in a regime where finite-degree terms are negligible.

Inference after minimum-distance selection requires the joint law of $\wh k$, the retained fraction of vertices tied at $D_{(\wh k+1)}(n)$, and $\wh\alpha_{n,\wh k}-\alpha$. This law must determine which finite-degree terms remain at the selected rank and how selection changes the distribution used for confidence intervals.

\section{Observed degree sequences}\label{sec:observed}
\subsection{Ties at the selected threshold}\label{sec:observed-ties}
Observed degrees are integers, and ties at the threshold enter the Hill statistic through its denominator. For each degree $d$, let $N_{n,>d}$ and $N_{n,=d}$ be the numbers of vertices with degree greater than $d$ and equal to $d$. When the rank $k$ increases through the tie of vertices with degree $d=D_{(k+1)}(n)$, each additional vertex has $\log\{D_{(i)}(n)/D_{(k+1)}(n)\}=0$, so it enlarges the divisor $k$ without changing the sum in Eq.~\eqref{eq:hill-stat}. Define the fraction of the threshold tie retained at rank $k$ by
\begin{equation}\label{eq:tie-fraction}
 \theta_{n,k}=\frac{k-N_{n,>D_{(k+1)}(n)}}{N_{n,=D_{(k+1)}(n)}}.
\end{equation}
Its denominator is at least one, and $0\le\theta_{n,k}<1$. Both counts are evaluated at the random threshold $D_{(k+1)}(n)$, so $\theta_{n,k}$ is a statistic of the degree sequence.

Across adjacent ranks with positive thresholds, subtraction in Eq.~\eqref{eq:hill-stat} gives
\begin{equation}\label{eq:adjacent-hill}
 H_{n,k+1}-H_{n,k}=-\frac{H_{n,k}}{k+1}+\log\frac{D_{(k+1)}(n)}{D_{(k+2)}(n)}.
\end{equation}
Equation~\eqref{eq:adjacent-hill} produces repeated increases and drops in $\wh\alpha_{n,k}$ as a function of $k$. Across the ranks of the tie at threshold $d$, from $k=N_{n,>d}$ to $k=N_{n,>d}+N_{n,=d}-1$, the logarithmic term vanishes and the numerator of $H_{n,k}$ is constant, so $\wh\alpha_{n,k}$ grows linearly with $k$ and reaches its minimum and maximum on the tie at its first and last ranks. When the threshold steps down from $d$ to $d-1$, the logarithmic term equals $\log\{d/(d-1)\}=d^{-1}+O(d^{-2})$ and $\wh\alpha_{n,k}$ drops sharply. Two ratios at the selected rank $\wh k$ compare these discrete steps with an iid Pareto reference. The relative asymptotic standard deviation of $\wh\alpha_{n,k}$ is $k^{-1/2}$ for $k$ independent Pareto observations because $\sqrt k(\wh\alpha_{n,k}-\alpha)$ then has asymptotic variance $\alpha^2$~\cite[Theorem~3.2.5]{deHaanFerreira2006}. With $d=D_{(\wh k+1)}(n)$, the ratio $\rho=\wh k^{1/2}/d$ compares the relative integer spacing $1/d$ with $\wh k^{-1/2}$, and the second ratio divides the ramp across the selected tie, $\wh\alpha_{n,N_{n,>d}+N_{n,=d}-1}-\wh\alpha_{n,N_{n,>d}}$, by $\wh\alpha_{n,\wh k}\wh k^{-1/2}$.

The tie also enters the selection criterion. At $y=1$, the empirical survival function in Eq.~\eqref{eq:mdsp-def} equals $N_{n,>D_{(k+1)}(n)}/k$ because the retained vertices of degree $D_{(k+1)}(n)$ do not exceed the threshold, and the fitted Pareto survival function equals one. Equation~\eqref{eq:tie-fraction} therefore gives
\begin{equation}\label{eq:dk-tie-bound}
 \mathcal D_k\ge\frac{\theta_{n,k}N_{n,=D_{(k+1)}(n)}}{k}.
\end{equation}
The position of the rank within its tie thus alters both the exponent estimate and the criterion that selects the rank.

\subsection{Threshold ties in observed networks}\label{sec:real}
We applied the unrestricted selector to the $3440$ degree sequences with at least $50$ positive degrees in the Index of Complex Networks (ICON) corpus analyzed by Broido and Clauset~\cite{ICON,BroidoClauset2019,SFAnalysis}, and to four SNAP graphs~\cite{SNAPCAIDA,SNAPEnron,SNAPAstro,SNAPHepTh}. Each fit uses the positive degrees and replaces $n$ by their number $n_+$. For CAIDA we ignore relationship direction and type, merge reciprocal and repeated edges, and remove loops.

\begin{figure*}[t]
\centering
\includegraphics[width=\textwidth]{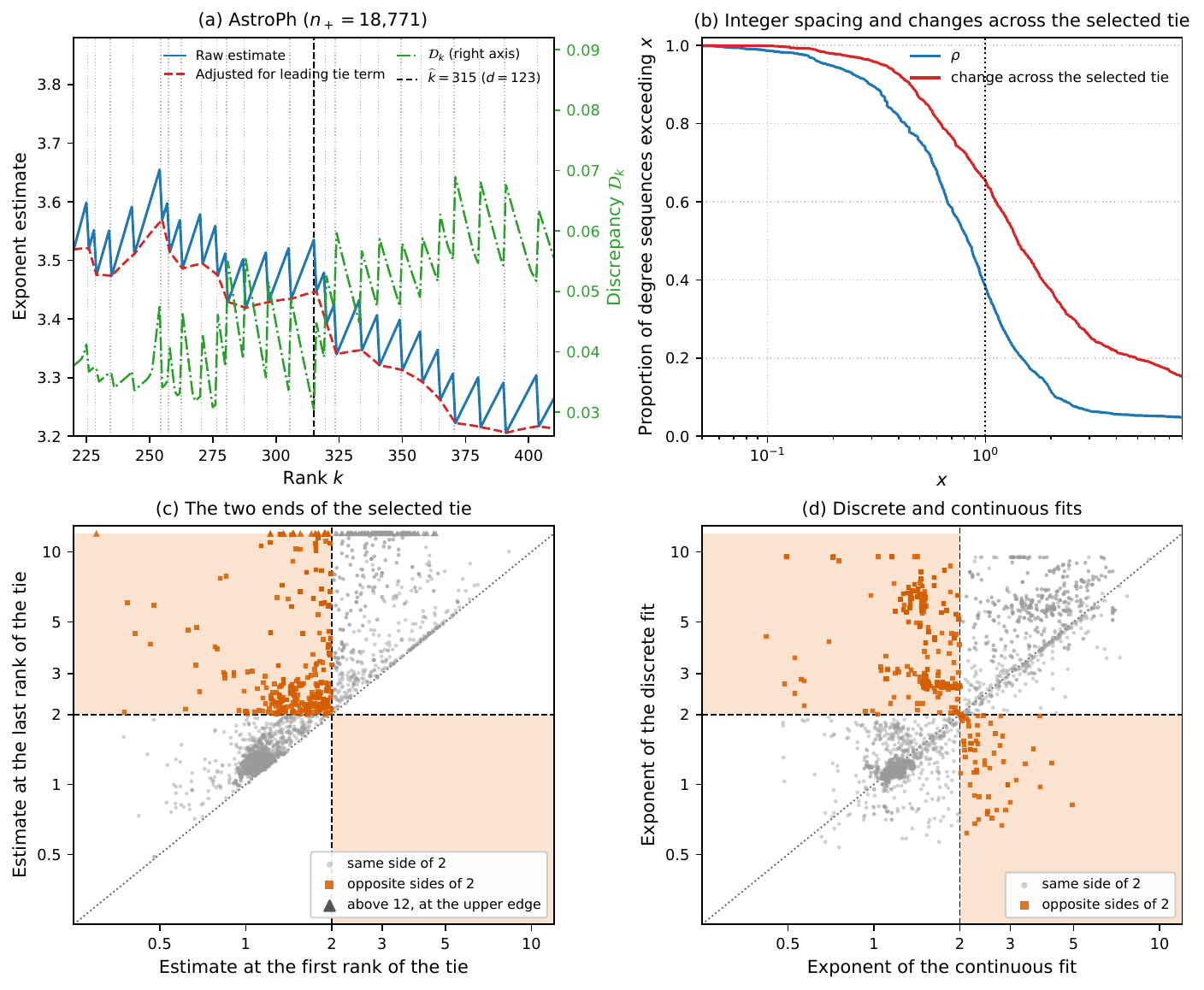}
\caption{Threshold ties and integer spacing in observed networks. (a) The estimates $\wh\alpha_{n_+,k}$ and $(H_{n_+,k}+\theta_{n_+,k}/D_{(k+1)})^{-1}$, the second defined in Sec.~\ref{sec:real}, and the criterion $\mathcal D_k$ against the rank $k$ in the AstroPh collaboration graph ($n_+=18{,}771$). Vertical dotted lines mark changes in $D_{(k+1)}$, and the vertical dashed line marks the unrestricted minimizer $\wh k=315$, at threshold $d=123$. (b) Empirical complementary distribution functions, over the $3440$ ICON degree sequences, of $\rho=\wh k^{1/2}/D_{(\wh k+1)}$ and of $(\wh\alpha_{n_+,N_{n_+,>d}+N_{n_+,=d}-1}-\wh\alpha_{n_+,N_{n_+,>d}})/(\wh\alpha_{n_+,\wh k}\wh k^{-1/2})$, the change of the estimate across the tie at $d=D_{(\wh k+1)}$ in units of $\wh\alpha_{n_+,\wh k}\wh k^{-1/2}$. The vertical dotted line marks one. (c) The estimates $\wh\alpha_{n_+,k}$ at the first rank $k=N_{n_+,>d}$ and the last rank $k=N_{n_+,>d}+N_{n_+,=d}-1$ of the selected tie, for each ICON sequence, with values exceeding 12 drawn at the upper edge. (d) The exponent of the discrete fit of Ref.~\cite{ClausetEtAl2009} against $\wh\alpha_{n_+,\wh k}$ for the same sequences. The row at $9.5$ is the largest exponent on the grid of the discrete fit. In (c) and (d), dashed lines mark estimates equal to $2$, the dotted line is the diagonal, and squares mark the pairs of estimates on opposite sides of $2$, which lie in the shaded quadrants.}
\label{fig:empirical}
\end{figure*}

Figure~\ref{fig:empirical}(a) plots, besides $\wh\alpha_{n_+,k}$, $(H_{n_+,k}+\theta_{n_+,k}/D_{(k+1)})^{-1}$, an adjustment that adds the retained tie fraction divided by the threshold to the Hill statistic. Section~\ref{sec:ties} derives its role in the expansion for integer tails.

Figure~\ref{fig:empirical}(a) shows for the AstroPh collaboration graph that $\wh\alpha_{n_+,k}$ increases within each tie, as Eq.~\eqref{eq:adjacent-hill} predicts, and drops when the threshold decreases. Table~\ref{tab:empirical} reports the four SNAP fits.
\begin{table*}[t]
\caption{Unrestricted MDSP fits to simple undirected projections of four SNAP graphs, with the search over $2\le k<n_+$. The total vertex counts for AstroPh and HepTh are 18,772 and 9,877, and removing loops leaves one and two isolates, respectively. CAIDA and Enron have no isolates in the observed vertex set. The column $n_+$ counts positive degrees, $d=D_{(k+1)}$ is the selected threshold, $\rho=\sqrt k/d$, and $\wh\alpha$ is the estimate $\wh\alpha_{n_+,k}$ of Eq.~\eqref{eq:hill-stat}.}
\label{tab:empirical}
\begin{ruledtabular}
\begin{tabular}{lrrrrrrrr}
Network & $n_+$ & $k$ & $k/n_+$ (\%) & $d$ & $\theta$ & $\rho$ & $\wh\alpha$ & $\mathcal D_k$\\
\hline
CAIDA AS & $26,475$ & $213$ & $0.80$ & $42$ & $0.000$ & $0.347$ & $1.099$ & $0.0320$\\
Enron & $36,692$ & $1164$ & $3.17$ & $57$ & $0.957$ & $0.599$ & $1.478$ & $0.0359$\\
AstroPh & $18,771$ & $315$ & $1.68$ & $123$ & $0.900$ & $0.144$ & $3.535$ & $0.0304$\\
HepTh & $9,875$ & $58$ & $0.59$ & $36$ & $0.714$ & $0.212$ & $6.048$ & $0.0893$\\
\end{tabular}\end{ruledtabular}\end{table*}

The fitted ranks include threshold ties in the email and coauthorship graphs, and Eq.~\eqref{eq:dk-tie-bound} links these ties to the values of $\mathcal D_{\wh k}$ in Table~\ref{tab:empirical}. For HepTh, a search restricted to $162\le k\le646$ selects $k=579$ with $\wh\alpha_{n_+,k}=3.02$, and the unrestricted minimizer $\wh k=58$ gives $\wh\alpha_{n_+,\wh k}=6.05$, so the selected exponent depends strongly on the candidate range.

In the ICON corpus, integer spacing exceeds the iid Pareto reference, $\rho>1$, in $38.0$ percent of fits, and the change across the selected tie exceeds $\wh\alpha_{n_+,\wh k}/\sqrt{\wh k}$ in $64.9$ percent. Figure~\ref{fig:empirical}(b) gives both distributions. Figure~\ref{fig:empirical}(c) shows that the estimates at the two ends of the selected tie lie on opposite sides of $2$, the exponent value that separates finite and infinite limiting degree variance, in $9.9$ percent of fits. These fits can therefore change their classification of the limiting degree variance without changing the degree threshold.

The corpus contains several sequences for some directed and multiplex graphs. Giving each of the $1806$ graphs equal weight, by averaging first within each graph, reduces this proportion to $8.2$ percent. Table~\ref*{S-tab:corpus-weights} of the Supplemental Material gives the two weightings and the proportions by domain.

Broido and Clauset~\cite{BroidoClauset2019} use a discrete Pareto likelihood that retains each tie whole. For each threshold $d$ with an observed degree $d+1$, we fit probabilities proportional to $j^{-\alpha-1}$ on integers $j>d$ and select the threshold by its Kolmogorov--Smirnov distance. The exponent grid is $\alpha\in\{0.01,0.02,\ldots,9.50\}$, matching the code of Ref.~\cite{SFAnalysis}, which reports $\alpha+1$. Figure~\ref{fig:empirical}(d) compares this estimate with the continuous fit. The two estimates lie on opposite sides of $2$ in $10.6$ percent of sequences. The discrete fit retains a median of $3.7$ times as many vertices. The fitted exponent thus depends on both the handling of ties and the tail selected by the fitting rule.

\section{Why selection retains finite-degree bias}\label{sec:pa}
\subsection{The graph model and the degree distribution}
Linear PA supplies a controlled setting in which integer support, finite-degree corrections, and dependence can be separated. Its growth rule fixes the limiting degree probabilities at every finite degree. We use the homogeneous sequential PA model of Ref.~\cite{Wang2026}, with $m\ge1$ edges per new vertex and offset $\delta>-m$, both fixed. Two seed vertices share $m$ parallel edges. At each time $t\ge3$, a new vertex adds $m$ edges successively to older vertices. If $D_v(t,i-1)$ is the degree of vertex $v<t$ immediately before edge $i$, its attachment probability is
\begin{equation}\label{eq:pa-model}
 \frac{D_v(t,i-1)+\delta}{(2m+\delta)(t-1)-2m+i-1}.
\end{equation}
Parallel edges are allowed, loops are excluded, and the new vertex becomes eligible at the next time. The graph at time $n$ therefore has $m(n-1)$ edges, so $m$ is determined by the observed graph. Set
\begin{equation}\label{eq:pa-parameters}
 \alpha=2+\delta/m>1,\quad r=m+\delta,\quad
 c=\left\{\frac{\Gamma(r+\alpha)}{\Gamma(r)}\right\}^{1/\alpha}.
\end{equation}
The limiting degree probabilities and their tail sums satisfy, for integers $j,d\ge m$,
\begin{align}
 \pi_j&=\alpha c^\alpha\frac{\Gamma(j+\delta)}{\Gamma(j+\delta+\alpha+1)},\notag\\
 \sum_{j>d}\pi_j&=c^\alpha\frac{\Gamma(d+1+\delta)}{\Gamma(d+1+\delta+\alpha)}.
 \label{eq:pa-stationary}
\end{align}
Since $\Gamma(d+1+\delta)/\Gamma(d+1+\delta+\alpha)=d^{-\alpha}\{1-\alpha Q/d+O(d^{-2})\}$ with $Q=\delta+(\alpha+1)/2$~\cite[\S5.11(iii)]{DLMF}, the limiting degree distribution satisfies Eq.~\eqref{eq:regular-variation}, and its departure from the exact Pareto tail $c^\alpha d^{-\alpha}$ has relative order $d^{-1}$, the same order as the relative spacing between consecutive integer degrees.

\subsection{The position within a degree tie}\label{sec:ties}
The limiting degree probabilities in Eq.~\eqref{eq:pa-stationary} determine the leading shift of $H_{n,k}$ at a finite threshold $d$. Retaining a fraction $\theta$ of the tie changes this shift through the denominator of the Hill statistic. Replacing the empirical degree proportions in $H_{n,k}$ by their limits $\pi_j$, at a deterministic integer threshold $d\ge m$ with a fraction $\theta\in[0,1]$ of the tie mass $\pi_d$ retained, gives the population ratio on the left of Eq.~\eqref{eq:hill-tie-expansion}. Equation~\eqref{eq:pa-stationary} gives $\pi_d/\sum_{j>d}\pi_j=\alpha/(d+\delta)$, and, with $B=(\delta+\alpha+1/2)/(\alpha+1)$, expanding the Gamma ratios~\cite[\S5.11(iii)]{DLMF} gives, uniformly over $0\le\theta\le1$ as the integer $d\to\infty$,
\begin{equation}\label{eq:hill-tie-expansion}
 \frac{\sum_{j>d}\pi_j\log(j/d)}{\sum_{j>d}\pi_j+\theta\pi_d}
 =\frac1\alpha+\frac{B-\theta}{d}+O(d^{-2}).
\end{equation}
Sec.~\ref*{S-sec:pa-lattice} of the Supplemental Material~\cite{SM} derives this expansion. The term $-\theta/d$ comes from the retained mass $\theta\pi_d$, which enlarges the denominator and adds nothing to the numerator, and $B/d$ combines the Euler--Maclaurin discretization of the sum $\sum_{j>d}\pi_j\log(j/d)$ with the $O(d^{-1})$ relative departure of $\sum_{j>d}\pi_j$ from $c^\alpha d^{-\alpha}$. Only the constant $B$ depends on the PA probabilities. Let $(\pi_j)_{j\ge1}$ be any probability mass function on the positive integers whose tail sums satisfy $\overline F(d)=c^\alpha d^{-\alpha}\{1-\alpha Q/d+O(d^{-1-\eta})\}$ for some $c,\eta>0$ and $Q\in\R$, the expansion assumed in Theorem~\ref*{S-thm:general-lattice} of the Supplemental Material. Then $\pi_d/\overline F(d)=\overline F(d-1)/\overline F(d)-1=\alpha/d+O(d^{-1-\min(\eta,1)})$, and Eq.~\arxiveqref{S-eq:lattice-log-expansion} of the Supplemental Material gives $\sum_{j>d}\pi_j\log(j/d)/\overline F(d)=1/\alpha+O(d^{-1})$. Dividing numerator and denominator by $\overline F(d)$ therefore gives, uniformly over $0\le\theta\le1$,
\begin{equation}\label{eq:tie-term-general}
\begin{aligned}
 \frac{\sum_{j>d}\pi_j\log(j/d)}{\overline F(d)+\theta\pi_d}
 &=\frac{\sum_{j>d}\pi_j\log(j/d)}{\overline F(d)}\\
 &\quad-\frac\theta d+O(d^{-1-\min(\eta,1)}).
\end{aligned}
\end{equation}
Adding $\theta_{n,k}/D_{(k+1)}(n)$ to $H_{n,k}$ thus removes the leading dependence on the position within the tie for every such law, using only the observed degrees and ranks, and leaves the model-dependent term $B/d$. Figure~\ref{fig:empirical}(a) shows for AstroPh that $(H_{n_+,k}+\theta_{n_+,k}/D_{(k+1)})^{-1}$ is nearly flat across each tie and agrees with $\wh\alpha_{n_+,k}$ at the first rank of each tie, where $\theta_{n_+,k}=0$. Across the ICON corpus of Sec.~\ref{sec:real}, the adjusted estimate at the last rank of the selected tie differs from its value at the first rank by more than $\wh\alpha_{n_+,\wh k}/\sqrt{\wh k}$ in $14.7$ percent of fits, and the two ends lie on opposite sides of $2$ in $5.1$ percent.

Taking the reciprocal of Eq.~\eqref{eq:hill-tie-expansion} gives $\alpha-\alpha^2(B-\theta)/d+O(d^{-2})$, so a rank with threshold $d$ and retained tie fraction $\theta$ shifts $\wh\alpha_{n,k}-\alpha$ by
\begin{equation}\label{eq:exponent-shift}
 -\frac{\alpha^2(B-\theta)}{d}.
\end{equation}
For the Barab\'asi--Albert model ($m=1$, $\delta=0$, $\alpha=2$), $B=5/6$, so the shift is negative over most of the tie and changes sign only when $\theta$ crosses $5/6$.

Discrete likelihood, a half-integer threshold shift~\cite{ClausetEtAl2009},
and independent uniform noise on $[-1/2,1/2]$~\cite{VoitalovEtAl2019}
remove dependence on a split tie. Lemma~\ref*{S-lem:noise-expansion}
of the Supplemental Material defines their population fits and proves
that their reciprocal exponents equal
$1/\alpha+(B-1/2)/d+O(d^{-2})$ for the PA degree law.
For a discrete power law, $B=1/2$, so this term vanishes.
For PA, $B-1/2=(\delta+\alpha/2)/(\alpha+1)$ is positive at
$\alpha\ge2$. These integer-data treatments therefore leave the
finite-degree PA bias even after removing tie sensitivity.

\subsection{The rank where bias and fluctuations balance}\label{sec:vanishing}
The threshold selected by $\mathcal D_k$ determines whether the shift in Eq.~\eqref{eq:exponent-shift} is negligible on the scale used for uncertainty. The deterministic departure from a continuous Pareto tail has order $d^{-1}$, and the random fluctuations of the empirical tail have order $k^{-1/2}$. Because $\sum_{j>d}\pi_j\sim c^\alpha d^{-\alpha}$, a rank $k$ corresponds to a degree scale $c(n/k)^{1/\alpha}$. Equating $(k/n)^{1/\alpha}$ and $k^{-1/2}$ suggests the rank and degree sequences
\begin{equation}\label{eq:balance}
 a_n=n^{2/(\alpha+2)},\qquad b_n=c(n/a_n)^{1/\alpha}=c\sqrt{a_n}.
\end{equation}
For fixed $t>0$ and $k=\lfloor a_nt\rfloor$, Proposition~\ref*{S-prop:pa-fixed-rank} of the Supplemental Material gives $D_{(k+1)}(n)/b_n\xrightarrow{p}t^{-1/\alpha}$ in $\R$. In the scaling limit, $t$ measures the retained rank in units of $a_n$, $x=t^{-1/\alpha}$ measures the degree threshold in units of $b_n$, and both the systematic shift $(B-\theta)/d$ and the sampling fluctuations of $H_{n,k}$ have order $a_n^{-1/2}=n^{-1/(\alpha+2)}$. The scale of the systematic term therefore remains comparable to the scale of the sampling fluctuations as both vanish.

The random fluctuations of the tail are governed by the centered counts $a_n^{-1/2}(N_{n,>b_nx}-\E N_{n,>b_nx})$ of vertices with degree greater than $b_nx$. The centered limits of Ref.~\cite{Wang2026} for sums of a function of the degrees, Theorem~\ref*{S-thm:pa-functional}, and Proposition~\ref*{S-prop:pa-beta-covariance} of the Supplemental Material show that, for every $\epsilon>0$, this centered process converges in distribution in $\ell^\infty([\epsilon,\infty))$ to a centered Gaussian process $\mathbb G$ on $(0,\infty)$ with continuous paths and covariance, for $x,y>0$,
\begin{equation}\label{eq:covariance}
\begin{split}
 \Cov\{\mathbb G(x),\mathbb G(y)\}
 &=(x\vee y)^{-\alpha}\\
 &\quad-(1+1/m)\big\{
 x^{-\alpha}\mathrm I_{x/(x+y)}(r+\alpha,r)\\
 &\hspace{20mm}+y^{-\alpha}\mathrm I_{y/(x+y)}(r+\alpha,r)\big\},
\end{split}
\end{equation}
where $\mathrm I_z(a,b)$ is the Beta$(a,b)$ distribution function at $z$ and $r=m+\delta$ is given in Eq.~\eqref{eq:pa-parameters}. In Eq.~\eqref{eq:covariance}, $(x\vee y)^{-\alpha}$ is the covariance for independent draws from the limiting degree distribution, and the subtracted Beta term combines two effects of PA growth, as in the variance of the Hill estimator in Ref.~\cite{Wang2026}: vertices born at different times have different degree laws, which gives the factor $1$ in $1+1/m$, and they compete for each new edge, which gives the factor $1/m$. The integrated process
\[
 \mathbb L(x)=\int_x^\infty\mathbb G(z)\,\frac{\mathrm dz}{z}
\]
is well defined almost surely, as established in Sec.~\ref*{S-app:pa-tail-measure} of the Supplemental Material. It is the limit of the sum of logarithmic degree excesses beyond $b_nx$, centered by its mean and multiplied by $a_n^{-1/2}$.

Fix $t>0$ and $\theta\in[0,1]$, and choose deterministic integer thresholds $d_n$ with $d_n/b_n\to x=t^{-1/\alpha}$. Retaining a fraction tending to $\theta$ of the vertices tied at $d_n$ defines random ranks with $k/a_n\xrightarrow{p}t$ in $\R$, as specified in Theorem~\ref*{S-thm:pa-degree-clt} of the Supplemental Material. At these ranks, the rescaled Hill statistic $\sqrt{a_n}(H_{n,k}-1/\alpha)$ converges in distribution in $\R$ to
\begin{equation}\label{eq:pa-A}
 \mathcal A(t,\theta)
 =\frac{\mathbb L(x)-\mathbb G(x)/\alpha}{t}+\frac{B-\theta}{cx},
\end{equation}
which adds the Gaussian fluctuation of the logarithmic excesses to the shift $(B-\theta)/d$ from Eq.~\eqref{eq:hill-tie-expansion}. Likewise, at a relative threshold $y\ge1$ in Eq.~\eqref{eq:mdsp-def}, the rescaled discrepancy between the empirical and fitted survival functions has the Gaussian part
\begin{equation}\label{eq:pa-Z}
\begin{split}
 \mathcal Z_{\rm lat}(t,y)
 &=\frac{\mathbb G(xy)-y^{-\alpha}\mathbb G(x)}t\\
 &\quad-\frac{\alpha^2y^{-\alpha}\log y}{t}
       \{\mathbb L(x)-\mathbb G(x)/\alpha\}
\end{split}
\end{equation}
together with a deterministic drift from the jumps at integer degrees and the shift of $\wh\alpha_{n,k}$:
\begin{equation}\label{eq:pa-R}
\begin{split}
 \mathcal R(t,\theta;y,\psi)
 &=\mathcal Z_{\rm lat}(t,y)+\frac{y^{-\alpha}}{cx}\bigg[
 \alpha\left\{\frac{\psi-Q}{y}+Q-\theta\right\}\\
 &\hspace{30mm}-\alpha^2(B-\theta)\log y\bigg],
\end{split}
\end{equation}
where $\psi\in\{0,1\}$ indexes the two one-sided limits of the empirical survival function across an integer degree. Taking the supremum over $y\ge1$ and $\psi\in\{0,1\}$ gives the limit criterion
\begin{equation}\label{eq:pa-D}
 \mathcal D_{\rm lat}(t,\theta)
 =\sup_{y\ge1}\max_{\psi\in\{0,1\}}
              |\mathcal R(t,\theta;y,\psi)|.
\end{equation}
Theorem~\ref{thm:vanishing}, proved in Sec.~\ref*{S-sec:unrestricted-proofs} of the Supplemental Material, shows that the unrestricted MDSP rank localizes at scale $a_n$ and converges jointly with its tie fraction and exponent estimate to the minimizer of $\mathcal D_{\rm lat}$.
\begin{theorem}\label{thm:vanishing}
Under Eq.~\eqref{eq:pa-model}, with fixed $m\ge1$ and $\delta>-m$, $\mathcal D_{\rm lat}$ has almost surely a unique global minimizer $(T,\Theta)\in(0,\infty)\times[0,1]$, and the unrestricted minimizer $\wh k$ of Eq.~\eqref{eq:mdsp-def} satisfies
\begin{equation}\label{eq:unrestricted-joint}
\begin{split}
 \left(\frac{\wh k}{a_n},\theta_{n,\wh k},
       \sqrt{a_n}(\wh\alpha_{n,\wh k}-\alpha),
       \sqrt{a_n}\mathcal D_{\wh k}\right)\qquad\\
 \Rightarrow\left(T,\Theta,-\alpha^2\mathcal A(T,\Theta),
                         \mathcal D_{\rm lat}(T,\Theta)\right)
 \quad\text{in }\R^4.
\end{split}
\end{equation}
Jointly, $D_{(\wh k+1)}(n)/b_n\Rightarrow T^{-1/\alpha}$ in $\R$. In particular,
\begin{equation}\label{eq:rank-bounds}
 \lim_{L\to\infty}\limsup_{n\to\infty}
 \Prob\{\wh k/a_n\notin[L^{-1},L]\}=0,
\end{equation}
$\wh k/n\xrightarrow{p}0$ in $\R$, and $\wh\alpha_{n,\wh k}\xrightarrow{p}\alpha$ in $\R$.
\end{theorem}

Theorem~\ref{thm:vanishing} identifies the regime selected by the fitting criterion itself. The positive, finite limit $T$ places $\wh k$ at the boundary of the regime $k=o(a_n)$ where the finite-degree shift is negligible. Joint convergence and $b_n=c\sqrt{a_n}$ also give
\begin{equation}\label{eq:selected-fit-spacing}
 \left(\mathcal D_{\wh k},\frac{\sqrt{\wh k}}{D_{(\wh k+1)}(n)}\right)
 \Rightarrow\left(0,\frac{T^{1/2+1/\alpha}}c\right)
       \quad\text{in }\R^2.
\end{equation}
Thus the minimized discrepancy tends to zero in probability, and relative integer spacing remains comparable to $\wh k^{-1/2}$. The difference $\wh\alpha_{n,\wh k}-\alpha$ has scale $n^{-1/(\alpha+2)}$, with the term $-\alpha^2(B-\Theta)T^{1/\alpha}/c$ in its limit. Its sampling fluctuation is a functional of the Gaussian path that also selects $(T,\Theta)$. A vanishing fitting discrepancy therefore leaves a finite-degree shift on the scale used for confidence intervals.

\subsection{Selection in a graph with latent weights}\label{sec:nr}
To determine whether the rank balance requires PA growth, consider the NR multigraph~\cite{NorrosReittu2006}. Take iid weights with $\Prob(W_i>w)=(w/w_0)^{-\alpha}$ for $w\ge w_0$, where $\alpha,w_0>0$. Conditional on the weights, let the edge multiplicities be independent Poisson variables with means $W_iW_j/\sum_vW_v$ for $i\le j$, and count each self-loop once in the degree.

Cirkovic, Wang, and Cline~\cite{CirkovicWangCline2025} establish normality of tail-index estimators under their tail assumptions at deterministic ranks satisfying $k\log^\alpha n/n^{1/(4\alpha+1)}\to0$. Theorem~\ref*{S-NR-thm:process-localization} of the Supplemental Material proves the centered functional limit at $a_n=n^{2/(\alpha+2)}$ for the Pareto weights specified here, and places the unrestricted MDSP rank at that scale. This scale lies outside the rank regime covered by their normality theorem, since $2/(\alpha+2)>1/(4\alpha+1)$ for every $\alpha>0$. The limit has covariance $(x\vee y)^{-\alpha}$, and the degree-tail expansion gives $c=w_0$, $Q=(1-\alpha)/2$, and $B=1/\{2(\alpha+1)\}$ in Eqs.~\eqref{eq:pa-A}--\eqref{eq:pa-D}. These constants replace the PA constants throughout this subsection.

For the following NR limit statements, set the raw exponent and criterion coordinate to zero when every rank is invalid. With $\mathbb G(x)=\mathbb W(x^{-\alpha})$, where $\mathbb W$ is standard Brownian motion, Theorem~\ref*{S-NR-thm:selected-joint} of the Supplemental Material gives a unique global minimizing pair $(T,\Theta)$ and the joint limit in Eq.~\eqref{eq:unrestricted-joint}, together with $D_{(\wh k+1)}(n)/b_n\Rightarrow T^{-1/\alpha}$ in $\R$. Equation~\eqref{eq:selected-fit-spacing} follows here too, with both coordinates taken as zero if no rank is valid, an event whose probability tends to zero. A vanishing fitting discrepancy coexists with finite-degree terms at the fluctuation scale even with the iid Pareto covariance. The rank balance occurs in both models, and each covariance determines the corresponding calibration. Theorem~\ref*{S-thm:general-lattice} of the Supplemental Material isolates what the balance requires. For positive integer observations with arbitrary dependence within each row, the tail expansion used in Eq.~\eqref{eq:tie-term-general} and the conditions of Assumption~\ref*{S-ass:lattice-inputs} on the expected counts and the centered tail empirical measure at scale $a_n$ give the limit of $\sqrt{a_n}\mathcal D_k$ on every compact range of $k/a_n$, with the drift terms of Eq.~\eqref{eq:pa-R}. The two models enter this limit only through $c$, $Q$, $B$, and the covariance of $\mathbb G$.

\section{Inference after selecting the threshold}\label{sec:inference}
\subsection{First- and second-order corrections and covariance}
Under the PA model with known $m$, $\delta=m(\alpha-2)$, so the constant $B$ is a function $B(\alpha,m)$. Replacing $d$ by the observed threshold $D_{(k+1)}(n)$, $\theta$ by $\theta_{n,k}$, and $\alpha$ by $\wh\alpha_{n,k}$ in the term $(B-\theta)/d$ of Eq.~\eqref{eq:hill-tie-expansion} gives the corrected statistics
\begin{align}
 \wh B_{n,k}&=\frac{m(\wh\alpha_{n,k}-2)+\wh\alpha_{n,k}+1/2}{\wh\alpha_{n,k}+1},\notag\\
 \wt H_{n,k}&=H_{n,k}-\frac{\wh B_{n,k}-\theta_{n,k}}{D_{(k+1)}(n)},\qquad
 \wt\alpha_{n,k}=\wt H_{n,k}^{-1}.
 \label{eq:corrected-estimator}
\end{align}
Set $\wt\alpha_{n,k}=0$, and report the fit as invalid, when the rank is invalid or $\wt H_{n,k}\le0$. At $k=\wh k$ this event has probability tending to zero, because $\wt\alpha_{n,\wh k}\to\alpha$ in probability. The correction changes only the reported exponent, and $\wh k$ still minimizes $\mathcal D_k$ in Eq.~\eqref{eq:mdsp-def}, which uses $\wh\alpha_{n,k}$.

At moderate thresholds, the shift of $\wh\alpha_{n,k}$ in Eq.~\eqref{eq:exponent-shift} carries over to $\wh B_{n,k}=B(\wh\alpha_{n,k},m)$, and the $O(d^{-2})$ remainder in Eq.~\eqref{eq:hill-tie-expansion} is not negligible. Carrying the expansion of the Gamma ratios~\cite[\S5.11(iii)]{DLMF} and the Euler--Maclaurin formula~\cite[Eq.~2.10.1]{DLMFEulerMaclaurin} for the sums in Eq.~\eqref{eq:pa-stationary} one order further gives the coefficient $C_2(\alpha,m,\theta)=B_2(\alpha,m)-\alpha\theta B+\theta(\delta+\alpha\theta)$ of $d^{-2}$ in Eq.~\eqref{eq:hill-tie-expansion}, with $B_2$ given in Eq.~\arxiveqref{S-eq:pa-C2} of the Supplemental Material. Evaluating $B$ and $C_2$ at $\bar\alpha_{n,k}=\wt\alpha_{n,k}$ when $\wt\alpha_{n,k}>1$, and at $\max(\wh\alpha_{n,k},1)$ otherwise, defines the second-order estimator
\begin{align}
 \wt H_{n,k}^{(2)}&=H_{n,k}-\frac{B(\bar\alpha_{n,k},m)-\theta_{n,k}}{D_{(k+1)}(n)}
   -\frac{C_2(\bar\alpha_{n,k},m,\theta_{n,k})}{D_{(k+1)}(n)^2},\notag\\
 \wt\alpha_{n,k}^{(2)}&=\bigl(\wt H_{n,k}^{(2)}\bigr)^{-1},
 \label{eq:second-order-estimator}
\end{align}
with the same rule for invalid fits. For fixed $t>0$, at $k=\lfloor a_nt\rfloor$ and at $k=\wh k$, the threshold $D_{(k+1)}(n)$ has order $b_n=c\sqrt{a_n}$, so $\sqrt{k}/D_{(k+1)}(n)^2\to0$ in probability, and Eq.~\arxiveqref{S-eq:pa-second-order-equivalence} of the Supplemental Material gives $\sqrt{k}(\wt\alpha_{n,k}^{(2)}-\wt\alpha_{n,k})\xrightarrow{p}0$ in $\R$.

Once the $O(d^{-1})$ shift is removed, the fluctuations of $\wt\alpha_{n,k}$ at a deterministic rank $k=\lfloor a_nt\rfloor$ are Gaussian, with a variance strictly smaller than the iid Pareto variance $\alpha^2$~\cite[Theorem~3.2.5]{deHaanFerreira2006}. For iid draws, Corollary~\ref*{S-cor:iid-lattice} of the Supplemental Material gives variance $z^{-\alpha}$ for $\mathbb G(z)$. The subtracted Beta term in Eq.~\eqref{eq:covariance} reduces this variance by the factor $\sigma(\alpha,m)^2/\alpha^2$ at every threshold $z$, where
\begin{equation}\label{eq:pa-factor}
\begin{split}
 \sigma(\alpha,m)^2=\alpha^2\big[1-2(1+1/m)\qquad\qquad\\
 {}\times\mathrm I_{1/2}\{m(\alpha-1)+\alpha,m(\alpha-1)\}\big].
\end{split}
\end{equation}
Because $\mathrm I_{1/2}(a,b)>0$, Propositions~\ref*{S-prop:pa-beta-covariance} and \ref*{S-U-prop:positive-covariance} of the Supplemental Material give $0<\sigma(\alpha,m)^2<\alpha^2$. Reference~\cite{Wang2026} obtained the same variance for the Hill estimate at deterministic ranks $k=o(a_n)$, at which the shift is negligible. For $m=1$ and $\alpha=2$, $\sigma(2,1)^2=2$ is half the iid variance $\alpha^2=4$. Corollary~\ref*{S-cor:pa-corrected} of the Supplemental Material shows that, for fixed $t>0$ and $k=\lfloor a_nt\rfloor$,
\begin{equation}\label{eq:corrected-clt}
 \sqrt{k}(\wt\alpha_{n,k}-\alpha)
 \Rightarrow N\{0,\sigma(\alpha,m)^2\}\quad\text{in }\R.
\end{equation}
For the uncorrected estimate, Proposition~\ref*{S-prop:pa-fixed-rank} of the Supplemental Material gives limits of $\sqrt{a_n}(\wh\alpha_{n,k}-\alpha)$ along subsequences on which $b_nt^{-1/\alpha}-\lfloor b_nt^{-1/\alpha}\rfloor$ converges. These limits retain the integer rounding of the threshold.

\subsection{Limit law of the corrected estimator at the unrestricted cutoff}\label{sec:selected-law}
At the selected cutoff $\wh k$, the scale $T$ and the Gaussian fluctuation in Eq.~\eqref{eq:pa-A} depend on the same path of $\mathbb G$. The deterministic-rank limit in Eq.~\eqref{eq:corrected-clt} therefore does not determine the selected law. Subtracting $(\wh B_{n,\wh k}-\theta_{n,\wh k})/D_{(\wh k+1)}(n)$ in Eq.~\eqref{eq:corrected-estimator} cancels the shift $(B-\Theta)T^{1/\alpha}/c$ in $\mathcal A(T,\Theta)$, and Corollary~\ref*{S-cor:pa-corrected-selected} of the Supplemental Material gives
\begin{equation}\label{eq:selected-convergence}
 \frac{\sqrt{\wh k}(\wt\alpha_{n,\wh k}-\alpha)}{\sigma(\alpha,m)}
 \Rightarrow\mathcal E_{\alpha,m}\quad\text{in }\R,
\end{equation}
where
\begin{equation}\label{eq:selected-limit}
 \mathcal E_{\alpha,m}
 =-\frac{\alpha^2}{\sigma(\alpha,m)\sqrt T}
       \{\mathbb L(T^{-1/\alpha})-\mathbb G(T^{-1/\alpha})/\alpha\}.
\end{equation}
By Slutsky's lemma, the same limit holds with $\wt\alpha_{n,\wh k}^{(2)}$ in place of $\wt\alpha_{n,\wh k}$. For $z\in\R$ and $0<p<1$, define the distribution function and quantiles of $\mathcal E_{\alpha,m}$ by
\begin{align}
 F_{\alpha,m}(z)&=\Prob\{\mathcal E_{\alpha,m}\le z\},\notag\\
 \mathfrak q_p(\alpha,m)&=\inf\{z:F_{\alpha,m}(z)\ge p\}.
 \label{eq:selected-cdf}
\end{align}
Because $T$ is chosen by minimizing $\mathcal D_{\rm lat}$ on the same path of $\mathbb G$ that enters Eq.~\eqref{eq:selected-limit}, the post-selection law $F_{\alpha,m}$ need not be Gaussian or symmetric, and confidence intervals use its lower and upper quantiles separately.

For the Norros--Reittu model of Sec.~\ref{sec:nr}, replace $\wh B_{n,k}$ in Eq.~\eqref{eq:corrected-estimator} by $1/\{2(1+\wh\alpha_{n,k})\}$. Theorem~\ref*{S-NR-thm:selected-joint} of the Supplemental Material gives, jointly with $(\wh k/a_n,\theta_{n,\wh k})\Rightarrow(T,\Theta)$ in $\R^2$,
\begin{equation}\label{eq:nr-corrected-selected}
 \sqrt{\wh k}(\wt\alpha_{n,\wh k}-\alpha)
 \Rightarrow-\frac\alpha{\sqrt T}
    \left\{\int_0^T\frac{\mathbb W(v)}v\,dv-\mathbb W(T)\right\}
       \quad\text{in }\R.
\end{equation}
The correction removes the explicit shift involving the tie fraction. The remaining term and the minimizing time $T$ depend on the same Brownian path, so deterministic-rank normality does not by itself give the law after selection.

\subsection{Confidence intervals at the unrestricted cutoff}\label{sec:intervals}
Under the PA model with fixed known $m$, for a coverage error $\upsilon\in(0,1)$, substitute $\wt\alpha_{n,\wh k}$ for $\alpha$ in both $\sigma(\alpha,m)$ and the quantiles of $\mathcal E_{\alpha,m}$. When $\wt\alpha_{n,\wh k}>1$, define
\begin{equation}\label{eq:unrestricted-interval}
\begin{split}
 C_n(\upsilon)=\bigg[&\wt\alpha_{n,\wh k}
 -\frac{\sigma(\wt\alpha_{n,\wh k},m)}{\sqrt{\wh k}}
       \mathfrak q_{1-\upsilon/2}(\wt\alpha_{n,\wh k},m),\\
 &\wt\alpha_{n,\wh k}
 -\frac{\sigma(\wt\alpha_{n,\wh k},m)}{\sqrt{\wh k}}
       \mathfrak q_{\upsilon/2}(\wt\alpha_{n,\wh k},m)\bigg]\cap(1,\infty).
\end{split}
\end{equation}
If the fit is invalid or $\wt\alpha_{n,\wh k}\le1$, set $C_n(\upsilon)=(1,\infty)$. Since $\wt\alpha_{n,\wh k}\to\alpha>1$ in probability, $\Prob\{\wt\alpha_{n,\wh k}\le1\}\to0$. Substituting $\wt\alpha_{n,\wh k}^{(2)}$ for $\wt\alpha_{n,\wh k}$ in Eq.~\eqref{eq:unrestricted-interval} defines the second-order interval $C_n^{(2)}(\upsilon)$. Both intervals use the single unrestricted minimizer $\wh k$, with no separate pilot or estimate of $a_n$.

The next theorem establishes the asymptotic coverage of intervals formed with the quantiles in Eq.~\eqref{eq:selected-cdf}.
\begin{theorem}\label{thm:coverage}
Assume the homogeneous PA model in Eq.~\eqref{eq:pa-model}, with fixed known $m\ge1$ and $\delta>-m$. For every fixed $\upsilon\in(0,1)$,
\begin{equation}\label{eq:unrestricted-coverage}
 \Prob\{\alpha\in C_n(\upsilon)\}\longrightarrow1-\upsilon,
 \qquad
 \Prob\{\alpha\in C_n^{(2)}(\upsilon)\}\longrightarrow1-\upsilon.
\end{equation}
\end{theorem}
Theorem~\ref*{S-thm:pa-unrestricted-interval} of the Supplemental Material proves this result and its extension to simulated quantiles under the uniform approximation condition in Eq.~\arxiveqref{S-eq:pa-calibration-accuracy}. To compute the quantiles, simulate Gaussian paths with covariance~\eqref{eq:covariance} at the corrected exponent, minimize Eq.~\eqref{eq:pa-D} on each path, and evaluate Eq.~\eqref{eq:selected-limit} at its minimizer. The simulation therefore repeats the threshold selection that determines the observed estimate. Section~\ref{sec:pa-eval} assesses the finite calibration grid and the coverage in finite graphs.

\section{Consequences for uncertainty and connectivity}\label{sec:empirical}

\subsection{The scale selected by the fit and the effect of growth}\label{sec:scaling-evidence}
To assess the scale selected by the fit, we use independent PA graphs with $m=4$, $\alpha=3$, and $n\in\{10^4,10^5,10^6\}$. For each graph, $\wh k$ minimizes Eq.~\eqref{eq:mdsp-def} over all valid ranks. The simulated limit is obtained by minimizing $\mathcal D_{\rm lat}$ on Gaussian paths, using the calibration described in Sec.~\ref*{S-sec:unrestricted-design} of the Supplemental Material. Figure~\ref{fig:mechanisms}(a) compares the distributions of $\wh k/a_n$ across the three sizes with the limit $T$. Their comparable distributions support the rank scale in Theorem~\ref{thm:vanishing}. Panel (b) compares the raw and corrected exponent errors at the selected rank with the corrected limit law. The correction moves the distribution of $\sqrt{\wh k}(\wt\alpha_{n,\wh k}-\alpha)/\sigma(\alpha,m)$ toward $F_{\alpha,m}$, with a residual difference assessed through coverage in Sec.~\ref{sec:pa-eval}.

To separate growth from the limiting degree law, we generated PA graphs and independent observations with the probabilities $\pi_j$ in Eq.~\eqref{eq:pa-stationary}, for each $n\in\{10^4,10^5\}$, $m\in\{1,4\}$, and $\alpha\in\{1.8,2.2\}$. Each sample is evaluated at the deterministic rank $k=\lfloor a_n\rfloor$. The same correction in Eq.~\eqref{eq:corrected-estimator} applies because both samples have the same limiting degree probabilities. Equation~\eqref{eq:corrected-clt} and Corollary~\ref*{S-cor:iid-lattice} of the Supplemental Material give limiting variances $\sigma(\alpha,m)^2/\alpha^2$ and $1$, respectively, for $\sqrt{k}(\wt\alpha_{n,k}-\alpha)/\alpha$. Figure~\ref{fig:mechanisms}(c) compares these fluctuations, and panel (d) compares their empirical variances. The construction and replication counts are given in Sec.~\ref*{S-sec:variance-decisions} of the Supplemental Material.

\begin{figure*}[t]
\centering
\includegraphics[width=\textwidth]{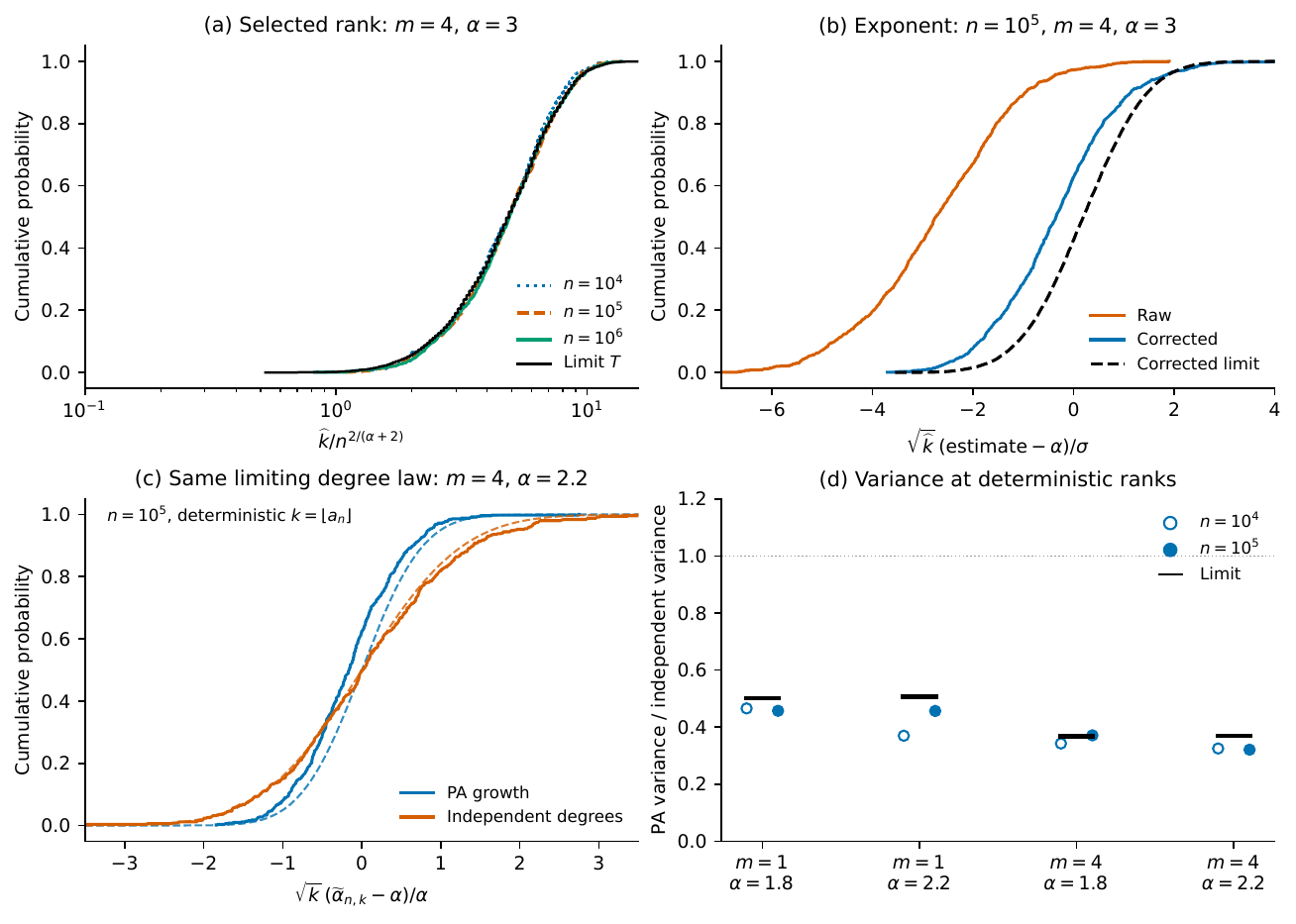}
\caption{Threshold selection and dependence in PA. (a) Empirical distribution functions of $\wh k/a_n$ for $m=4$, $\alpha=3$, and the indicated sizes, with the simulated limit $T$. (b) Distributions of $\sqrt{\wh k}(\wh\alpha_{n,\wh k}-\alpha)/\sigma(\alpha,m)$ and $\sqrt{\wh k}(\wt\alpha_{n,\wh k}-\alpha)/\sigma(\alpha,m)$, with the limit of the latter. Panels (a) and (b) use 1000 graphs per setting and 20,000 draws of the limit. (c) Distributions of $\sqrt{k}(\wt\alpha_{n,k}-\alpha)/\alpha$ at deterministic $k=\lfloor a_n\rfloor$ for PA and independent draws from the same limiting degree law. Dashed curves are the corresponding centered Gaussian limits. (d) Ratios of the empirical PA and iid variances of this quantity, with their predicted limits $\sigma(\alpha,m)^2/\alpha^2$. Panels (c) and (d) use 500 replications for each sampling law and setting.}
\label{fig:mechanisms}
\end{figure*}

\subsection{Coverage at the selected rank}\label{sec:pa-eval}
We evaluated $C_n(\upsilon)$ and $C_n^{(2)}(\upsilon)$ at $\upsilon=0.05$ on simulated PA graphs of sizes matching observed networks: $1000$ independent graphs from Eq.~\eqref{eq:pa-model} for each combination of $n\in\{10^4,10^5\}$, $m\in\{1,4\}$, and $\alpha\in\{1.5,2,2.5,3\}$, with $\wh k$ selected over all $2\le k<n$. The quantiles use $20{,}000$ simulated Gaussian paths per
calibration setting, with linear interpolation in $\alpha$.
Sec.~\ref*{S-sec:unrestricted-design} specifies the grids and
refinement checks. In its six checked settings, doubling the time
range and refining the grids changes either endpoint quantile by at
most $0.036$, including Monte Carlo variation from new paths.

\begin{figure*}[t]
\centering
\includegraphics[width=\textwidth]{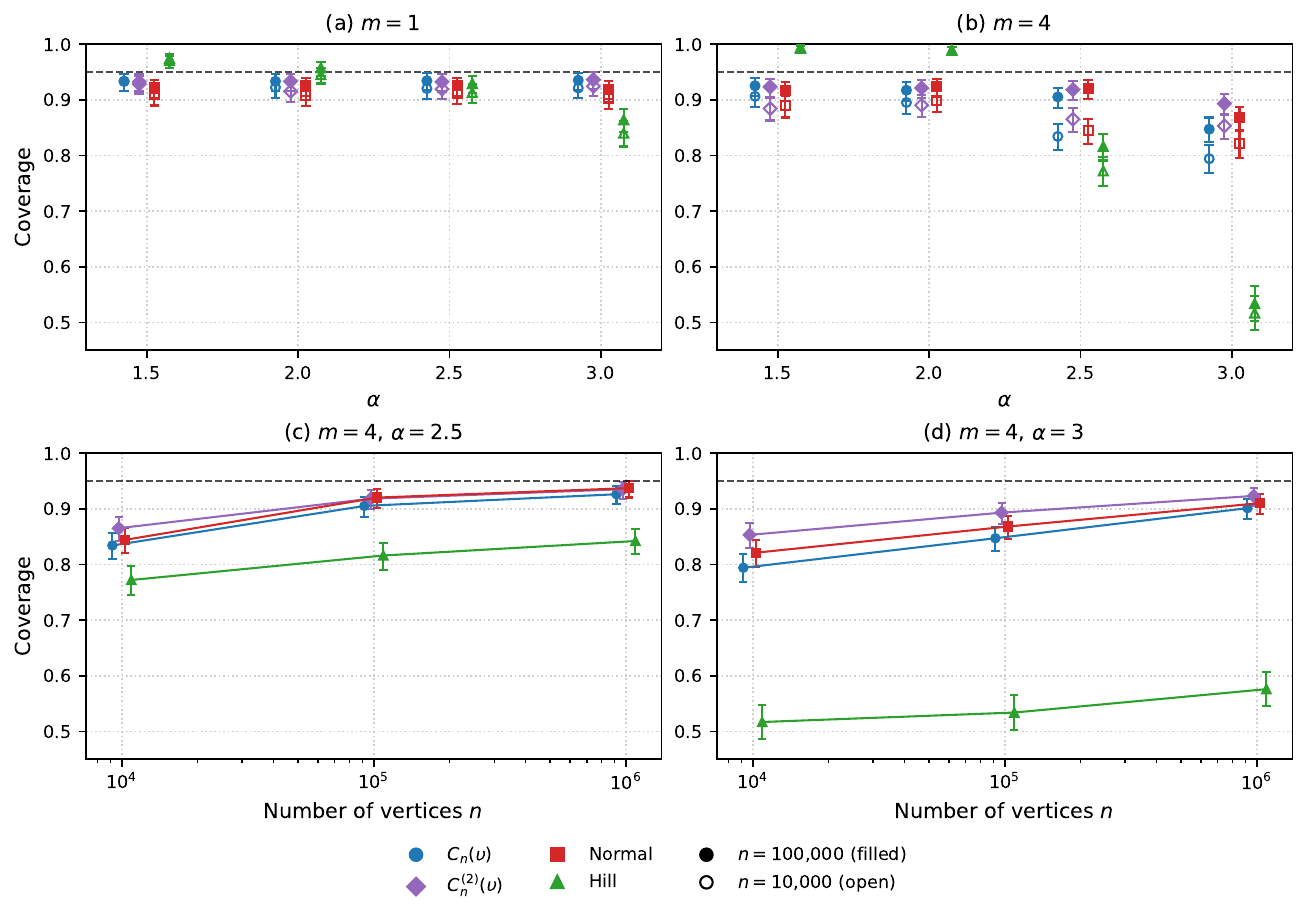}
\caption{Finite-size coverage at the unrestricted selected rank. Panels (a) and (b) compare $C_n(\upsilon)$, $C_n^{(2)}(\upsilon)$, the corrected normal interval, and the Hill interval across $\alpha$, with filled markers for $n=10^5$ and open markers for $n=10^4$. Panels (c) and (d) follow the two settings with the largest coverage deficits to $n=10^6$. Horizontal offsets separate the methods. Each point uses 1000 PA graphs, and bars are 95\% Wilson intervals~\cite{Wilson1927}. The dashed line marks the nominal coverage $1-\upsilon=0.95$.}
\label{fig:coverage}
\end{figure*}

To isolate the roles of the bias correction, the PA variance reduction, and the post-selection law $\mathcal E_{\alpha,m}$, we compare $C_n(\upsilon)$ and $C_n^{(2)}(\upsilon)$ with two reference intervals at the same rank $\wh k$: the normal interval $\wt\alpha_{n,\wh k}\pm1.96\,\sigma(\wt\alpha_{n,\wh k},m)\wh k^{-1/2}$, which uses the corrected estimate and the PA variance with standard normal quantiles, and the naive Hill interval $\wh\alpha_{n,\wh k}\pm1.96\,\wh\alpha_{n,\wh k}\wh k^{-1/2}$, which uses the uncorrected estimate and the iid Pareto variance. For intervals $C^{(b)}$, $b=1,\ldots,R$, the coverage frequency is $R^{-1}\sum_{b=1}^R\mathbf1\{\alpha\in C^{(b)}\}$. For estimates $\wh\alpha^{(b)}$, $b=1,\ldots,R$, the empirical root mean squared error (RMSE) is $\{R^{-1}\sum_{b=1}^R(\wh\alpha^{(b)}-\alpha)^2\}^{1/2}$. Table~\ref*{S-tab:unrestricted} of the Supplemental Material reports the coverage, median interval width, and RMSE across the original settings and the larger graphs, and Fig.~\ref{fig:coverage} plots the coverage frequencies with 95\% Wilson intervals~\cite{Wilson1927}.

Figure~\ref{fig:coverage} shows that increasing network size can narrow the Hill interval without restoring its nominal coverage. At $m=4$ and $\alpha=3$, coverage is $51.7$, $53.4$, and $57.6$ percent for $n=10^4$, $10^5$, and $10^6$, respectively, against the nominal $95$ percent. Over these sizes its median width decreases from $0.745$ to $0.321$. The direction of miscalibration depends on the parameters. For $m=4$ and $\alpha\le2$, the Hill interval covers more than $98$ percent of graphs.

For $m=1$ across all tested exponents, and for $m=4$ at $\alpha\le2$, the first correction gives $89.5$ to $93.5$ percent coverage and reduces the median width by $12$ to $39$ percent. The corrected normal interval is within $2.7$ percentage points of $C_n(0.05)$ in every original setting. Both intervals already use the finite-degree correction and PA variance, so this comparison isolates the additional effect of replacing Gaussian quantiles by those in Eq.~\eqref{eq:selected-cdf}. Most of the observed gain over the Hill interval comes from removing the shift and using the PA variance. The law after selection supplies the asymptotic calibration.

The largest remaining deficits occur for $m=4$ and $\alpha\ge2.5$, where moderate selected thresholds leave appreciable terms of order $d^{-2}$. The second correction reduces these terms. Inverting the exact population ratio in Eq.~\eqref{eq:hill-tie-expansion} changes its coverage by at most $1.6$ percentage points across the original settings, as reported in Sec.~\ref*{S-app:unrestricted-experiments} of the Supplemental Material. At $n=10^6$, $C_n^{(2)}(0.05)$ reaches $93.5$ percent coverage for $\alpha=2.5$ and $92.3$ percent for $\alpha=3$. Table~\ref*{S-tab:unrestricted} gives the remaining finite-size deficits, widths, and estimation errors.

\subsection{Connectivity prediction and model error}\label{sec:percolation}\label{sec:nr-percolation}
We use the fitted exponent to predict connectivity after random edge
loss~\cite{CallawayEtAl2000}. Give each edge
an independent uniform label $U_e$ and retain it when $U_e\le p$.
If $L_{n,1}(p)$ is the largest component size, define
\begin{equation}\label{eq:percolation-observable}
 \mathcal Q_{n,f}=\inf\{p\in[0,1]:L_{n,1}(p)/n\ge f\},
\end{equation}
with value one if the target is unreached. We predict the survival
probability needed to connect $f=0.10$ of the vertices of an independent
test graph from the degrees of a training graph with the same parameters.
$\mathcal Q_{n,f}$ measures the crossing of a fixed component fraction, distinct from the
asymptotic percolation threshold~\cite{HazraEtAl2023}.

Both PA and Pareto NR use $n\in\{10^4,10^5\}$,
$\alpha\in\{1.8,2.2,2.6,3.0\}$, and known $m\in\{2,4\}$ or
$w_0\in\{2,4\}$, with $300$ independent training/test pairs per setting.
An independent calibration table converts each exponent into a mean
crossing under its graph law. All methods predict the same test outcome.
For predictions $\wh{\mathcal Q}_{n,f}^{(b,j)}$ from method $j$, the loss is
\begin{equation}\label{eq:percolation-mse}
 \operatorname{MSE}_j(f)=\frac1{300}\sum_{b=1}^{300}
 (\wh{\mathcal Q}_{n,f}^{(b,j)}-\mathcal Q_{n,f}^{(b)})^2.
\end{equation}
Pooled MSE weights settings equally. Bootstrap intervals resample both
calibration graphs and prediction pairs, preserving the pairing of
methods~\cite{Efron1979}. Sec.~\ref*{S-sec:percolation-experiment}
specifies the protocols and reports every setting and secondary fraction.

\begin{figure*}[t]
\centering
\includegraphics[width=\textwidth]{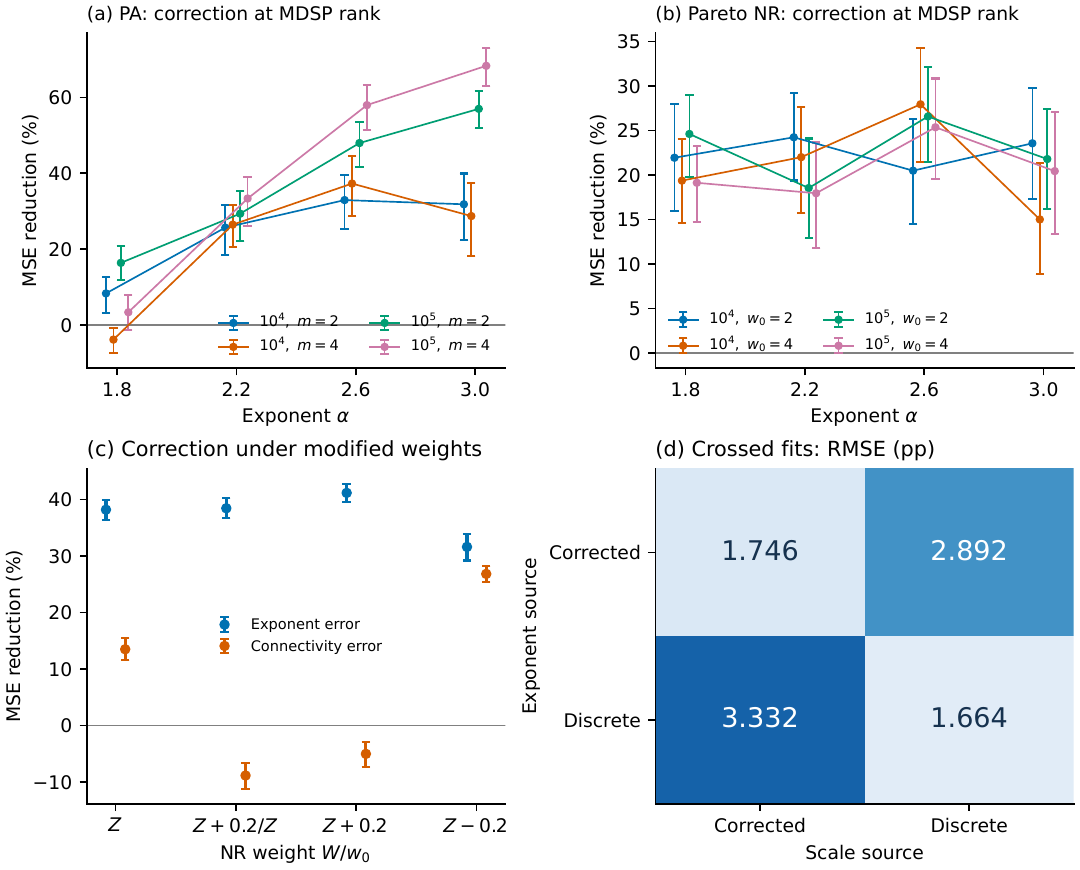}
\caption{When correcting the selected exponent improves connectivity
prediction. (a,b) Percentage MSE reduction relative to raw Hill for every
PA and Pareto NR setting, with known $m$ or $w_0$. Legends give $n$ and
the model parameter. (c) Percentage reductions in exponent and prediction
MSE under four NR weight laws, pooling each law's 16 settings equally.
Prediction estimates scale and retains Pareto calibration. Bars in
(a)--(c) are pointwise 95\% paired bootstrap intervals; prediction
intervals include calibration uncertainty. Each setting uses 300 pairs.
(d) Prediction RMSE in percentage points for crossed corrected and
discrete fits, pooled over the 48 non-Pareto settings. This exploratory
comparison reuses the samples in (c). All prediction panels use a
10-percent component target. Negative reductions indicate increased MSE.}
\label{fig:connectivity}
\end{figure*}

Figure~\ref{fig:connectivity}(a,b) shows the change in prediction MSE
when the model correction retains the MDSP rank. The pooled reduction is
$27.6$ percent for PA, with a $95$ percent interval of $[24.8,30.2]$,
and $22.5$ percent for NR, with an interval of $[20.3,24.6]$.
PA has one adverse setting: at $n=10^4$, $m=4$, $\alpha=1.8$, MSE
increases by $3.8$ percent, with an interval of $[0.8,7.3]$.
NR improves in all settings; halving its rank instead raises pooled MSE
by $24.3$ percent. The full PA likelihood~\cite{GaoVanDerVaart2017}
uses all degrees and reduces pooled prediction MSE by $55.3$ percent.
The correction therefore improves the selected tail fit under these
specified laws, with the full model providing a stronger PA predictor.

To assess model error, a separate NR experiment retains the Poisson
edge law and changes the weights to $w_0Z$, $w_0(Z+0.2/Z)$,
$w_0(Z+0.2)$, or $w_0(Z-0.2)$, where $\Prob(Z>z)=z^{-\alpha}$
for $z\ge1$. Each law uses the same parameter grid and replication
count. Prediction still uses Pareto calibration and the Pareto correction,
but estimates scale as
$\wh w_0(a)=(a-1)a^{-1}n^{-1}\sum_iD_i(n)$ for fitted exponent $a$.
This rule matches the Pareto mean degree to the observed mean.
Sec.~\ref*{S-sec:nr-robustness-experiment} gives the full design.

The correction reduces pooled exponent MSE under every weight law,
as Fig.~\ref{fig:connectivity}(c) shows. Its connectivity effect depends
on the law: prediction MSE increases by $8.9$ and $5.0$ percent for
$w_0(Z+0.2/Z)$ and $w_0(Z+0.2)$, and decreases by $26.8$ percent for
$w_0(Z-0.2)$. Each law-specific interval excludes zero. Across the
non-Pareto settings, the discrete tail likelihood of
Clauset et al.~\cite{ClausetEtAl2009} gives $9.2$ percent lower prediction
MSE than the correction, with an interval of $[4.2,14.2]$, despite
larger exponent MSE. In an exploratory analysis of these saved samples,
we cross the corrected and discrete exponents with their scale estimates.
Figure~\ref{fig:connectivity}(d) shows that replacing either parameter
alone worsens prediction, and replacing both improves it. Compensation
between exponent and scale in the predictor explains this ordering;
connectivity loss does not rank exponent accuracy by itself.

\section{Scope and implications}\label{sec:discussion}
\subsection{Inference for the selected exponent}\label{sec:testing}
The joint law in Eq.~\eqref{eq:unrestricted-joint} specifies the
finite-degree shift in $\wh\alpha_{n,\wh k}$ and its dependence on the
selected rank and retained tie fraction. The correction in
Eq.~\eqref{eq:corrected-estimator} removes this shift, and the law in
Eq.~\eqref{eq:selected-convergence} retains the dependence introduced by
selection. Together they give the PA intervals in
Theorem~\ref{thm:coverage}.

Equation~\eqref{eq:regular-variation} defines $\alpha$ through regular
variation of the limiting degree distribution.
Broido and Clauset~\cite{BroidoClauset2019} assess an exact discrete
power law beyond a fitted cutoff. Voitalov et al.~\cite{VoitalovEtAl2019}
use regular variation, which permits departures from an exact power law
at finite degrees. Their Section~V explains why finite data cannot
distinguish unrestricted asymptotic tail classes.

The selected order $a_n$ also balances the squared magnitudes of the
leading deterministic-rank bias and fluctuation terms,
$(k/n)^{2/\alpha}$ and $k^{-1}$, when the leading bias coefficient is
nonzero. Selection at this order therefore requires bias correction and
an uncertainty law that includes selection. Model-based likelihood
uses additional information when the graph law is correct, as in the PA
comparison. Tail methods permit weaker structural assumptions, with
advantages under some network misspecifications~\cite{WanWangDavisResnick2020}.

The finite-degree coefficient can also be estimated from degrees.
Sec.~\ref*{S-sec:degree-bias-estimation} regresses Hill statistics over strict exceedances against inverse degree at thresholds of order $\log n$.
Proposition~\ref*{S-prop:degree-bias-consistency} proves consistency under
explicit tail and sampling bounds, verified for homogeneous PA and
Pareto NR, and preserves their corrected selected limits. A fresh
validation of one fixed implementation improves PA estimation, has
little effect in NR, and usually withholds the correction in Zipf
samples. It fails the criterion fixed for applying it to the corpus.
The supplement retains the full comparison, including a discrete
likelihood and a direct adaptation of the extended Pareto estimator
of Ref.~\cite{BeirlantJoossensSegers2009}.

\subsection{Limitations}\label{sec:limitations}
A small $\mathcal D_{\wh k}$ alone supplies neither a calibrated
goodness-of-fit test nor a comparison with competing tails.
The PA coverage theorem assumes homogeneous growth with fixed known
$m$ and $\delta>-m$. Pareto-weight NR has a selected joint law here,
without an interval coverage theorem. Both require model-specific
finite-degree constants and covariance; regular variation alone fixes
neither. The observed-network analysis establishes sensitivity to ties,
without identifying a generating law or measuring exponent accuracy.
Convergence of normalized bias and MSE to moments of the limiting laws
would additionally require uniform integrability. Numerical moment
comparisons concern the stated finite samples.

Simulated quantiles require the uniform approximation in
Eq.~\arxiveqref{S-eq:pa-calibration-accuracy}. Grid refinement and coverage
assess the implementation, which retains finite-size deficits.
Exact inversion of the limiting population ratio does not remove the
difference between expected finite-graph degree frequencies and $\pi_j$.
The tie adjustment removes only the term in $\theta$; the remaining
$B/d$ needs a graph-specific coefficient or an accurate estimate.
Consistency of the coefficient estimator in
Proposition~\ref*{S-prop:degree-bias-consistency} does not establish an
effective correction across unknown graph families.

Connectivity predictions concern independent test graphs and independent
edge removal. They do not condition on the observed adjacency matrix.
The NR weight changes retain the conditional Poisson edge law, and the
crossed-fit comparison is exploratory. Their losses assess the fitted
exponent and scale within the assumed prediction family. The finite
comparisons establish neither an optimal selector nor a universally
preferable tail estimator.

\subsection{Reporting a selected degree exponent}\label{sec:conclusion}
Report $\wh\alpha_{n,\wh k}$ with the threshold $D_{(\wh k+1)}(n)$,
rank $\wh k$, retained tie fraction $\theta_{n,\wh k}$, and ratio
$\sqrt{\wh k}/D_{(\wh k+1)}(n)$. This ratio compares relative integer
spacing with nominal relative sampling error. Reporting
$(H_{n,\wh k}+\theta_{n,\wh k}/D_{(\wh k+1)}(n))^{-1}$ alongside
$\wh\alpha_{n,\wh k}$ identifies sensitivity to the retained tie.
For homogeneous PA with known $m$, retain $\wh k$, apply the correction
in Eq.~\eqref{eq:corrected-estimator}, and use the quantiles of
$F_{\alpha,m}$ as specified in Eq.~\eqref{eq:unrestricted-interval}.
Theorem~\ref{thm:coverage} gives asymptotically exact coverage for these
intervals. Figure~\ref{fig:coverage} reports their coverage in finite
graphs, and Fig.~\ref{fig:connectivity} assesses the consequence of
exponent correction for connectivity prediction.

\section*{Data and code availability}
The empirical degree sequences are available from the corpus repository~\cite{SFAnalysis} and the SNAP sources in Refs.~\cite{SNAPCAIDA,SNAPEnron,SNAPAstro,SNAPHepTh}. Simulation code, protocols, outputs, and input hashes are available from the corresponding author upon reasonable request. These materials have not yet been deposited in a public archive.

\begin{acknowledgments}
The author used OpenAI Codex (ChatGPT 6.0 Astra) to assist with verification of theoretical arguments, language editing, and routine code refactoring. The author reviewed and edited the resulting material and takes full responsibility for the content of this article.
\end{acknowledgments}
\nocite{BoucheronLugosiMassart2013,DuppalaEtAl2025,JoagDevProschan1983,Massart1990,Mourrat2013,Phelps1978,RevuzYor1999,vanDerVaart1998,Vershynin2026}
\bibliography{prx-references}
\end{document}


\maketitle
\tableofcontents
\bigskip

\section{Conventions and notation}\label{sec:setup}
Equation and theorem numbers without an S prefix refer to the main article.
All probabilities refer to the stated graph or sampling law, including simulation draws when specified.

For ordered observations $V_{(1)}\ge\cdots\ge V_{(n)}\ge0$, use the Hill statistic~\cite{Hill1975}
\begin{equation}\label{eq:statistics}
 H_{n,k}=\frac1k\sum_{i=1}^k\log\frac{V_{(i)}}{V_{(k+1)}},
 \qquad \wh\alpha_{n,k}=H_{n,k}^{-1}.
\end{equation}
The rank is valid when $V_{(k+1)}>0$ and $H_{n,k}>0$.
At invalid ranks set $\wh\alpha_{n,k}=0$, and also set $H_{n,k}=0$
when $V_{(k+1)}=0$. These assignments make the random variables in
the limit statements real-valued on every sample.
The rank version of minimum-distance selection~\cite{ClausetEtAl2009,DreesEtAl2020} uses
\begin{equation}\label{eq:Dk}
 \mathcal D_k=\sup_{y\ge1}\left|\frac1k\sum_{i=1}^k
 \mathbf1\{V_{(i)}/V_{(k+1)}>y\}-y^{-\wh\alpha_{n,k}}\right|,
 \qquad \wh k=\min\argmin_{2\le k<n}\mathcal D_k.
\end{equation}
Set $\mathcal D_k=+\infty$ at invalid ranks. If every rank is invalid, set
$\wh k=2$ and report no point estimate. For degrees, $V_v=D_v(n)$,
with ties retained. The reference ratio
\begin{equation}\label{eq:lattice-ratio}
 \rho_{n,k}=\frac{\sqrt{k}}{V_{(k+1)}}
\end{equation}
compares relative integer spacing with the iid Pareto relative standard
error $k^{-1/2}$~\cite{deHaanFerreira2006}.
For any exponent estimate, bias and root mean squared error are
\begin{equation}\label{eq:mdsp-estimation-error}
 \E(\wh\alpha-\alpha),\qquad
 \{\E[(\wh\alpha-\alpha)^2]\}^{1/2},
\end{equation}
when these moments exist. Their empirical versions replace expectations
by averages over independent replications.

For PA degrees, use the model, constants $\alpha,r,c$, and sequences
$a_n,b_n$ in Eqs.~\arxiveqref{M-eq:pa-model}, \arxiveqref{M-eq:pa-parameters}, and
\arxiveqref{M-eq:balance}, and put $\gamma=1/\alpha$.
The notation in Wang~\cite{Wang2026} is related by
$p=\alpha$, $\alpha=\gamma$, and $\beta=c$.
Its centered degree limits require only an intermediate rank; its
centered normal limit for $\sqrt{k}(H_{n,k}-1/\alpha)$ additionally requires $k=o(a_n)$.
The coefficients of the inverse-degree terms are
\begin{equation}\label{eq:pa-bias-constants}
 Q=\delta+\frac{\alpha+1}{2},\qquad
 B=\frac{\delta+\alpha+1/2}{\alpha+1}.
\end{equation}
The limiting probabilities $\pi_j$ are given in Eq.~\arxiveqref{M-eq:pa-stationary}.

\subsection{Counts and threshold ties}\label{sec:pa-interpretation}
For integer observations, write
\[
 N_{n,>d}=\sum_{v=1}^n\mathbf1\{V_v>d\},\qquad
 N_{n,=d}=\sum_{v=1}^n\mathbf1\{V_v=d\}.
\]
The fraction of the threshold tie retained at rank $k$ is
\begin{equation}\label{eq:pa-tie-fraction}
 \theta_{n,k}=\frac{k-N_{n,>V_{(k+1)}}}{N_{n,=V_{(k+1)}}}\in[0,1).
\end{equation}
At $y=1$, Eq.~\eqref{eq:Dk} gives
\begin{equation}\label{eq:pa-criterion-tie}
 \mathcal D_k\ge \theta_{n,k}N_{n,=V_{(k+1)}}/k.
\end{equation}
In PA, $V_v=D_v(n)$. The main article defines the continuous Gaussian
process $\mathbb G$, its integral $\mathbb L$, and the functionals
$\mathcal A,\mathcal Z_{\rm lat},\mathcal R,\mathcal D_{\rm lat}$
in Eqs.~\arxiveqref{M-eq:covariance} and \arxiveqref{M-eq:pa-A}--\arxiveqref{M-eq:pa-D}.
Since $\mathcal R$ is affine in $\psi$, the identity
$\max(|z-a|,|z+a|)=|z|+a$ for $a\ge0$ gives
\begin{equation}\label{eq:pa-D-midpoint}
 \mathcal D_{\rm lat}(t,\theta)=\sup_{y\ge1}
 \left\{|\mathcal R(t,\theta;y,1/2)|+
 \frac{\alpha y^{-\alpha-1}}{2ct^{-1/\alpha}}\right\}.
\end{equation}

\section{Limit processes for preferential attachment}\label{app:pa-limits}
The proof of Theorem~\ref*{M-thm:vanishing} uses the following centered limits, tail expansions, and bounds for the homogeneous PA model
\arxiveqref{M-eq:pa-model} with fixed $m\ge1$ and $\delta>-m$. The limit of
the criterion on compact ranges of ranks, Theorem~\ref{thm:general-lattice},
holds for triangular arrays of positive integers with arbitrary
dependence within rows, and Corollary~\ref{cor:pa-conditional} applies
it to PA degrees. The section uses four
results of Wang~\cite{Wang2026}: the central limit theorems for
centered sums of a function of the
degrees~\cite[Proposition~3.7 and Section~3.4]{Wang2026}, the
variance bound~\cite[Lemma~3.3]{Wang2026}, the bound on the expected
degree profile~\cite[Lemma~3.1]{Wang2026}, and the asymptotic
variance of the Hill statistic~\cite[Theorem~2.1 and
Lemma~3.9]{Wang2026}.

\subsection{The tail empirical measure and its Gaussian limit}\label{app:pa-tail-measure}
For an intermediate sequence $K_n\to\infty$, $K_n/n\to0$, put
$b_n=c(n/K_n)^{1/\alpha}$, which is $c\sqrt{a_n}$ when $K_n=a_n$, and
define measures on $(0,\infty)$ and a centered process by
\begin{equation}\label{eq:pa-measure}
 \nu_n=K_n^{-1}\sum_{v=1}^n\delta_{D_v(n)/b_n},\qquad
 \nu(\dd z)=\alpha z^{-\alpha-1}\dd z,\qquad
 S_n(f)=\sqrt{K_n}\{\nu_nf-\E\nu_nf\}.
\end{equation}
For $a>0$ let $G_a$ have the Gamma distribution with shape $a$ and rate
one, and for $s>0$ write $X_s=c^{-1}s^{-1/\alpha}G_r$. For test functions $f$, define
\begin{align}
 \omega_f(s)&=r\E\{f(c^{-1}s^{-1/\alpha}G_{r+1})-f(X_s)\},\notag\\
 \Omega_f(u)&=\int_0^u\omega_f(s)\dd s,\notag\\
 V(f,g)&=\int_0^\infty\Cov(f(X_s),g(X_s))\dd s
 -\frac1{m\alpha^2}\int_0^\infty
                 \frac{\Omega_f(u)\Omega_g(u)}{u^2}\dd u.
 \label{eq:pa-covariance}
\end{align}
For $x>0$ let $f_x(z)=\mathbf1\{z>x\}$ and
$h_x(z)=\{\log(z/x)-1/\alpha\}\mathbf1\{z>x\}$.
Proposition~3.7 and Section~3.4 of~\cite{Wang2026} give
$S_n(f)\wto N\{0,V(f,f)\}$ in $\R$ for bounded piecewise Lipschitz
functions $f$ with compact support in $(0,\infty)$, and equations
(3.26)--(3.27) there give it for $h_1$. The same truncation at infinity,
with $x$ in place of $1$, gives it for $h_x$, $x>0$. The indicators $f_x$
follow by truncation at infinity, using the variance
bound~\cite[Lemma~3.3]{Wang2026} and
$\limsup_n\E\nu_n(M,\infty)\le CM^{-\alpha}$ for $M>0$, by
\arxiveqref{M-eq:pa-stationary} and~\cite[Lemma~3.1]{Wang2026}. The
Cram\'er--Wold device, applied to linear combinations, gives joint
convergence of any finite collection of $S_n(f)$ in Euclidean space.
Only $K_n\to\infty$ and $K_n/n\to0$ are needed for these centered limits.

Let $\mathbb G$ be the centered Gaussian process indexed by $x>0$, a
degree divided by $b_n$, with covariance
\begin{equation}\label{eq:pa-G}
 \Cov(\mathbb G(x),\mathbb G(y))=V(f_x,f_y),
 \qquad
 \mathbb L(x)=\int_x^\infty\mathbb G(z)\frac{\dd z}{z}.
\end{equation}
This process has a continuous version, and $\mathbb L(x)$ exists almost
surely. Indeed, the covariance is positive semidefinite by the
finite-dimensional CLTs, and
\begin{equation}\label{eq:pa-G-bound}
 V(f_x-f_y,f_x-f_y)\le\nu\{(f_x-f_y)^2\}=|x^{-\alpha}-y^{-\alpha}|.
\end{equation}
The inequality follows from \eqref{eq:pa-covariance} with
$f=g=f_x-f_y$: Fubini's theorem and $c^\alpha=\E G_r^\alpha$ give
$\int_0^\infty\E f(X_s)^2\dd s=\nu(f^2)$, and the other two terms are
nonpositive. Gaussian moment bounds and the Kolmogorov continuity
theorem~\cite[Chapter~I, Theorem~2.1]{RevuzYor1999} apply in the time
coordinate $v=x^{-\alpha}$.
Also $\E|\mathbb G(z)|\le z^{-\alpha/2}$, so Fubini's theorem gives
$\E\int_x^\infty|\mathbb G(z)|\dd z/z<\infty$. Since
$\log(z/x)\mathbf1\{z>x\}=\int_x^\infty f_y(z)\dd y/y$ and
$\int_x^\infty\sqrt{V(f_y,f_y)}\dd y/y<\infty$, this integral commutes
with the Gaussian isometry. In particular,
$S_n(h_x)\wto\mathbb L(x)-\mathbb G(x)/\alpha$ in $\R$.

Two further properties of $V$ and $\mathbb G$ enter the proofs that the
selected rank has order $a_n$ and that $\mathcal D_{\rm lat}$ has a
unique minimizer. Strict positivity of $V$ at nonzero functions
bounds the probability that the criterion is small at small times in
Proposition~\ref{prop:pa-attainment} and at small ranks in
Lemma~\ref{U-lem:small-discrepancy}. Nondegeneracy of the law of
$\mathbb G$ on a separable Banach space of paths lets the proof of
Proposition~\ref{U-prop:unique-global} apply the theorem of
Phelps~\cite[Theorem~6]{Phelps1978} on G\^ateaux differentiability
outside Gaussian null sets. Set $Y(v)=\mathbb G(v^{-1/\alpha})$, which
indexes $\mathbb G$ by the time coordinate $v=x^{-\alpha}$.
Its covariance is
\begin{equation}\label{U-eq:Y-covariance}
 \Cov\{Y(v),Y(u)\}=V(f_{v^{-1/\alpha}},f_{u^{-1/\alpha}}),\qquad Y(0)=0.
\end{equation}
For $T>0$ and $0<\iota<1/2$, define
\begin{equation}\label{U-eq:weighted-path-space}
 E_\iota(T)=
 \left\{y\in C[0,T]:y(0)=0,\ \lim_{v\downarrow0}v^{-\iota}y(v)=0\right\},
 \qquad
 \|y\|_\iota=\sup_{0<v\le T}v^{-\iota}|y(v)|.
\end{equation}
Multiplication by $v^{-\iota}$ maps $E_\iota(T)$ isometrically onto the
continuous functions on $[0,T]$ that vanish at zero, so $E_\iota(T)$ is
a separable Banach space.

The next result proves $V(f,f)>0$ for every nonzero $f\in L^2(\nu)$ and
nondegeneracy of the law of $Y$ on $E_\iota(T)$.
\begin{proposition}\label{U-prop:positive-covariance}
The covariance $V$ extends continuously to $L^2(\nu)$ and
\begin{equation}\label{U-eq:strict-covariance}
 V(f,f)>0\qquad\text{for every nonzero }f\in L^2(\nu).
\end{equation}
For every $T>0$ and $0<\iota<1/2$, the process $Y$ has a version in
$E_\iota(T)$ whose law is a nondegenerate Gaussian measure, in the sense
that $\Var\ell(Y)>0$ for every nonzero continuous linear functional
$\ell$ on $E_\iota(T)$.
\end{proposition}
\begin{proof}
We first prove \eqref{U-eq:covariance-combined}, which combines the two
subtracted terms of \eqref{eq:pa-covariance} into one integral.
For $f,g\in C_c^1((0,\infty))$, integration by parts in the Gamma
density gives
\begin{equation}\label{U-eq:gamma-derivative}
 \omega_f(s)=\E\{X_s f'(X_s)\}=-\alpha s\frac{d}{ds}\E f(X_s),
 \qquad
 \Omega_f(u)=\alpha\left\{\int_0^u\E f(X_s)\,ds-u\,\E f(X_u)\right\}.
\end{equation}
If the support of $f$ is contained in $[A_f,B_f]\subset(0,\infty)$,
then $\E f(X_s)=O(s^{r/\alpha})$ as $s\downarrow0$, and $|\E f(X_s)|$ is
bounded by a polynomial in $s$ times $\exp(-cA_f s^{1/\alpha})$ as
$s\to\infty$. These bounds also hold with $g$ in place of $f$, so the product
$u^{-1}\int_0^u\E f(X_s)\,ds\int_0^u\E g(X_s)\,ds$ tends to zero as
$u\downarrow0$ and as $u\to\infty$. Differentiating this product,
substituting \eqref{U-eq:gamma-derivative}, and integrating over
$(0,\infty)$ yield
\[
 \int_0^\infty\frac{\Omega_f(u)\Omega_g(u)}{u^2}\,du
   =\alpha^2\int_0^\infty\E f(X_s)\,\E g(X_s)\,ds.
\]
Fubini's theorem and $c^\alpha=\E G_r^\alpha$ give
$\int_0^\infty\E\{f(X_s)g(X_s)\}\,ds=\nu(fg)$.
Consequently
\begin{equation}\label{U-eq:covariance-combined}
 V(f,g)=\nu(fg)
             -(1+1/m)\int_0^\infty\E f(X_s)\,\E g(X_s)\,ds.
\end{equation}

The logarithmic change of variable $z=e^u$ makes the integral in
\eqref{U-eq:covariance-combined} a Fourier multiplier. Define
\begin{align}
 (\mathsf U_\alpha f)(u)&=\sqrt\alpha\,e^{-\alpha u/2}f(e^u),\nonumber\\
 (\mathsf T\mathsf U_\alpha f)(v)
   &=\sqrt\alpha\,c^{-\alpha/2}e^{-\alpha v/2}\E f(e^vG_r)\nonumber\\
   &=c^{-\alpha/2}\E\{G_r^{\alpha/2}(\mathsf U_\alpha f)(v+\log G_r)\}.
 \label{U-eq:covariance-transform}
\end{align}
Then $\|\mathsf U_\alpha f\|_2^2=\nu(f^2)$, so $\mathsf U_\alpha$ extends
to a unitary map from $L^2(\nu)$ onto $L^2(\R)$, and
$\|\mathsf T\mathsf U_\alpha f\|_2^2=\int_0^\infty\{\E f(X_s)\}^2\,ds$.
For the unitary Fourier transform with kernel $e^{-i\xi v}$,
the multiplier of $\mathsf T$ is
$c^{-\alpha/2}\Gamma(r+\alpha/2+i\xi)/\Gamma(r)$.
Its squared modulus is
\begin{equation}\label{U-eq:gamma-multiplier}
 \Psi(\xi)=
 \frac{|\Gamma(r+\alpha/2+i\xi)|^2}
      {\Gamma(r)\Gamma(r+\alpha)}.
\end{equation}
Plancherel's identity therefore gives
\begin{equation}\label{U-eq:covariance-Fourier}
 V(f,f)=\int_\R
       \{1-(1+1/m)\Psi(\xi)\}|\widehat{\mathsf U_\alpha f}(\xi)|^2\,d\xi.
\end{equation}
The multiplier $1-(1+1/m)\Psi(\xi)$ is bounded and continuous.
Positive semidefiniteness of $V$, together with the density of
$\{\mathsf U_\alpha f:f\in C_c^1((0,\infty))\}$ in $L^2(\R)$, forces
$1-(1+1/m)\Psi(\xi)\ge0$ for every $\xi$.
Indeed, a negative value would persist on an interval and contradict
\eqref{U-eq:covariance-Fourier} for a function whose Fourier transform
is supported there, approximated in $L^2$ by these functions.
For $\xi\ne0$, the triangle inequality
$|\E\{G_r^{\alpha/2}e^{i\xi\log G_r}\}|\le\E G_r^{\alpha/2}$ is strict,
because $G_r$ has a positive density on $(0,\infty)$, so that
$\xi\log G_r$ is not almost surely constant modulo $2\pi$. Hence
$\Psi(\xi)<\Psi(0)$, and the multiplier in
\eqref{U-eq:covariance-Fourier} is positive at every $\xi\ne0$.
Continuity and density extend $V$ and \eqref{U-eq:covariance-Fourier}
to $L^2(\nu)$, which proves \eqref{U-eq:strict-covariance}.

For the path assertion, \eqref{eq:pa-G-bound} with $x=v^{-1/\alpha}$ and
$y=u^{-1/\alpha}$ gives
\begin{equation}\label{U-eq:Y-increment}
 \E|Y(v)-Y(u)|^2\le\nu\{(f_{v^{-1/\alpha}}-f_{u^{-1/\alpha}})^2\}=|v-u|.
\end{equation}
Gaussian moment bounds and the Kolmogorov continuity theorem give a
version of $Y$ on $[0,T]$ that is H\"older continuous with any fixed
exponent smaller than $1/2$. Taking that exponent larger than $\iota$
proves membership in $E_\iota(T)$ almost surely.

Let $\ell$ be a nonzero continuous linear functional on $E_\iota(T)$.
By the Riesz representation theorem there is a nonzero finite signed
measure $\mu$ on $(0,T]$ such that
\[
 \ell(y)=\int_{(0,T]}v^{-\iota}y(v)\,\mu(dv).
\]
An atom at zero contributes nothing, since $v^{-\iota}y(v)\to0$, and is
omitted. Define
\[
 f_\mu(z)=\int_{(0,T]}v^{-\iota}
                       \mathbf1\{z>v^{-1/\alpha}\}\,\mu(dv).
\]
This function vanishes for $z\le T^{-1/\alpha}$ and satisfies
$|f_\mu(z)|\le\|\mu\|z^{\alpha\iota}$, so it belongs to
$L^2(\nu)$ because $2\alpha\iota<\alpha$. The integral against $\mu$
commutes with the Gaussian isometry because
\[
 \int_{(0,T]}v^{-\iota}\|f_{v^{-1/\alpha}}\|_{L^2(\nu)}\,|\mu|(dv)
   =\int_{(0,T]}v^{1/2-\iota}\,|\mu|(dv)<\infty.
\]
Thus $\Var\ell(Y)=V(f_\mu,f_\mu)$. The map $\mu\mapsto f_\mu$ is
injective, because $f_\mu(z)=\int_{(z^{-\alpha},T]}v^{-\iota}\mu(dv)$
and one-sided limits of these values in $z$ determine the measure
$v^{-\iota}\mu(dv)$ on every $[\varepsilon,T]$. If $f_\mu=0$ almost
everywhere, those limits force $\mu=0$, a contradiction.
Equation~\eqref{U-eq:strict-covariance} therefore gives
$\Var\ell(Y)>0$.
\end{proof}

The substitution $s\mapsto x^{-\alpha}s$ in \eqref{eq:pa-covariance} gives
\begin{equation}\label{eq:pa-homogeneity}
 V(f_x,f_y)=x^{-\alpha}V(f_1,f_{y/x}),\qquad
 V(h_x,h_x)=x^{-\alpha}V(h_1,h_1).
\end{equation}
In the second integral of \eqref{eq:pa-covariance} the substitution acts
through $\Omega_{f_x}(u)=x^{-\alpha}\Omega_{f_1}(x^\alpha u)$.

The covariances of $\mathbb G(x)$ and of $\mathbb L(x)-\mathbb G(x)/\alpha$
are the values $V(f_x,f_y)$ and $V(h_x,h_y)$ of the integrals in
\eqref{eq:pa-covariance}. In particular, $V(h_1,h_1)$ is the asymptotic
variance of $\sqrt k(H_{n,k}-1/\alpha)$ at ranks $k=o(a_n)$, which
Theorem~2.1 of~\cite{Wang2026} writes as an integral of the square of the
upper incomplete Gamma function and Lemma~3.9 there evaluates with the
Beta distribution function $\mathrm I_{1/2}(a,b)$.
Proposition~\ref{prop:pa-beta-covariance} extends this evaluation to
$V(h_x,h_y)$ and $V(f_x,f_y)$ at all $x,y>0$, with $\mathrm I_z(a,b)$ in
place of $\mathrm I_{1/2}(a,b)$. This gives the closed forms
\arxiveqref{M-eq:covariance} and \arxiveqref{M-eq:pa-factor}, which the simulation
and the calibration use, and shows how they depend on $\alpha$.
For $\alpha>1$, write $V_{\alpha,m}$ for $V$ at $\delta=m(\alpha-2)$, so
that $r=m(\alpha-1)$ and $c$ is given by \arxiveqref{M-eq:pa-parameters}.
\begin{proposition}\label{prop:pa-beta-covariance}
For $\alpha>1$, $r=m(\alpha-1)$, and $x,y>0$,
\begin{align}
 V_{\alpha,m}(h_x,h_y)
 =\frac1{\alpha^2}\bigg[&(x\vee y)^{-\alpha}
 -(1+1/m)\big\{
 x^{-\alpha}\mathrm I_{x/(x+y)}(r+\alpha,r)\notag\\
 &\hspace{35mm}+y^{-\alpha}\mathrm I_{y/(x+y)}(r+\alpha,r)
                         \big\}\bigg].\label{eq:pa-beta-covariance}
\end{align}
The covariance $V_{\alpha,m}(f_x,f_y)$ is \arxiveqref{M-eq:covariance}, and
$\sigma(\alpha,m)^2$ in \arxiveqref{M-eq:pa-factor} satisfies
$\sigma(\alpha,m)^2=\alpha^4V_{\alpha,m}(h_1,h_1)>0$.
\end{proposition}
\begin{proof}
With $X_u=c^{-1}u^{-1/\alpha}G_r$ as before \eqref{eq:pa-covariance}, put
$p_x(u)=\Prob\{X_u>x\}=\Prob\{G_r>cxu^{1/\alpha}\}$.
Since $\log(X_u/x)\mathbf1\{X_u>x\}=\int_x^\infty\mathbf1\{X_u>z\}\dd z/z$,
Fubini's theorem and the substitution $z=x(v/u)^{1/\alpha}$ give
\begin{equation}\label{eq:pa-log-tail-mean}
 \E h_x(X_u)
 =\frac1\alpha\left\{\int_u^\infty p_x(v)\frac{\dd v}{v}-p_x(u)\right\}.
\end{equation}
For $x,y>0$, integration over $u$ first, on
$\{(u,v,w):0<u<\min(v,w)\}$, yields
\[
 \int_0^\infty
 \left\{\int_u^\infty p_x(v)\frac{\dd v}{v}-p_x(u)\right\}
 \left\{\int_u^\infty p_y(w)\frac{\dd w}{w}-p_y(u)\right\}\dd u
 =\int_0^\infty p_x(u)p_y(u)\dd u.
\]
Indeed, the integral over $u$ of the product of the two inner integrals
equals the sum of the two cross terms, because $\min(v,w)/(vw)$ equals
$1/w$ on $\{v<w\}$ and $1/v$ on $\{w<v\}$. These integrals are finite.
As $u\to\infty$, $p_x(u)$ decays faster than any power of $u$, and as
$u\downarrow0$, $\int_u^\infty p_x(v)\dd v/v\le\log(1/u)+C$, which is
square integrable near zero.

For $x>0$, integration by parts in the Gamma density gives
\[
 \omega_{f_x}(s)
 =r\bigl[\Prob\{G_{r+1}>cxs^{1/\alpha}\}-\Prob\{G_r>cxs^{1/\alpha}\}\bigr]
 =(cxs^{1/\alpha})^re^{-cxs^{1/\alpha}}/\Gamma(r).
\]
Integrating over $0<s<u$ and using
$c^\alpha=\Gamma(r+\alpha)/\Gamma(r)$ give
\begin{equation}\label{eq:pa-indicator-Omega}
 \Omega_{f_x}(u)=\alpha x^{-\alpha}
                     \Prob\{G_{r+\alpha}\le cxu^{1/\alpha}\}.
\end{equation}
Since $h_x=\int_x^\infty f_z\,\dd z/z-f_x/\alpha$ and $f\mapsto\Omega_f$
is linear, integrating \eqref{eq:pa-indicator-Omega} over thresholds
$z\ge x$ against $\dd z/z$, subtracting $\Omega_{f_x}(u)/\alpha$, and
applying Fubini's theorem give
\[
 \Omega_{h_x}(u)
 =\int_x^\infty\alpha z^{-\alpha}
      \Prob\{G_{r+\alpha}\le czu^{1/\alpha}\}\frac{\dd z}{z}
       -x^{-\alpha}\Prob\{G_{r+\alpha}\le cxu^{1/\alpha}\}
 =u p_x(u).
\]
For the last equality, integrate by parts in $z$ and use
$c^\alpha=\Gamma(r+\alpha)/\Gamma(r)$.
Since $\E G_r^\alpha=\Gamma(r+\alpha)/\Gamma(r)=c^\alpha$, Fubini's
theorem gives $\int_0^\infty p_z(u)\dd u=z^{-\alpha}$ for every $z>0$.
Thus the law of $X_u$ integrated over $u$ is $\nu$, and
$\int_0^\infty\E\{h_x(X_u)h_y(X_u)\}\dd u=\nu(h_xh_y)$.
Also $\nu(h_xh_y)=(x\vee y)^{-\alpha}/\alpha^2$:
if $x\le y$, then $h_x=h_y+\log(y/x)$ on $(y,\infty)$,
and $\nu h_y=0$, $\nu h_y^2=y^{-\alpha}/\alpha^2$.
Substituting these identities, \eqref{eq:pa-log-tail-mean}, and
$\Omega_{h_x}(u)=up_x(u)$ into \eqref{eq:pa-covariance} gives
\[
 V_{\alpha,m}(h_x,h_y)=\frac1{\alpha^2}
  \left\{(x\vee y)^{-\alpha}
          -(1+1/m)\int_0^\infty p_x(u)p_y(u)\dd u\right\}.
\]
For independent $G_r,G'_r$, the integral $\int_0^\infty p_x(u)p_y(u)\dd u$
equals $c^{-\alpha}\E\min(G_r/x,G'_r/y)^\alpha$.
Split this expectation according to which argument is smaller.
Weighting the density of $G_r$ by $G_r^\alpha/c^\alpha$ changes
its shape to $r+\alpha$. The ratio of independent Gamma variables
of shapes $r+\alpha$ and $r$ to their sum has the Beta$(r+\alpha,r)$
distribution. Hence the two terms are
$x^{-\alpha}\mathrm I_{x/(x+y)}(r+\alpha,r)$ and
$y^{-\alpha}\mathrm I_{y/(x+y)}(r+\alpha,r)$.
This proves \eqref{eq:pa-beta-covariance}. Setting $x=y=1$ and
comparing with \arxiveqref{M-eq:pa-factor} gives
$\alpha^4V_{\alpha,m}(h_1,h_1)=\sigma(\alpha,m)^2$.
The term subtracted in \arxiveqref{M-eq:pa-factor} is positive because
$\mathrm I_{1/2}(r+\alpha,r)>0$, so $\sigma(\alpha,m)^2<\alpha^2$. Also
$\nu(h_1^2)=\alpha^{-2}>0$, so
Proposition~\ref{U-prop:positive-covariance} gives
$\sigma(\alpha,m)^2=\alpha^4V_{\alpha,m}(h_1,h_1)>0$.
For the indicators, \eqref{U-eq:covariance-combined} extends from
$C_c^1((0,\infty))$ to $L^2(\nu)$: Proposition~\ref{U-prop:positive-covariance}
extends $V$ continuously, and the integral in
\eqref{U-eq:covariance-combined} is the inner product of
$\mathsf T\mathsf U_\alpha f$ and $\mathsf T\mathsf U_\alpha g$, which is
continuous because $\Psi$ in \eqref{U-eq:gamma-multiplier}, the squared
modulus of the multiplier of $\mathsf T$, is bounded. With
$\E f_x(X_u)=p_x(u)$, it gives
$V(f_x,f_y)=\nu(f_xf_y)-(1+1/m)\int_0^\infty p_x(u)p_y(u)\dd u$.
Since $\nu(f_xf_y)=(x\vee y)^{-\alpha}$, the Beta evaluation of this
integral proves \arxiveqref{M-eq:covariance}.
\end{proof}

Setting $x=y$ in \arxiveqref{M-eq:covariance} and $x=y=1$ in
\eqref{eq:pa-beta-covariance}, and using $\sigma(\alpha,m)^2>0$, gives
the variance of $\mathbb G(x)$,
\begin{equation}\label{eq:pa-indicator-variance}
 V(f_x,f_x)=\alpha^2x^{-\alpha}V(h_1,h_1)>0.
\end{equation}

\subsection{Expansions of the stationary tail sums}\label{sec:pa-lattice}
At ranks $k$ of order $a_n$ the threshold $d=D_{(k+1)}(n)$ has order
$b_n=c\sqrt{a_n}$, so terms of relative order $d^{-1}$ in the mean degree
profile are of the same order $a_n^{-1/2}$ as the sampling error of
$H_{n,k}$. The centering in
Theorem~\ref{thm:pa-degree-clt} therefore needs the stationary
probabilities \arxiveqref{M-eq:pa-stationary} to second order.
Sec.~B.2 of~\cite{Wang2026} derives the following expansions, as
the integer $d\to\infty$, from the expansion of ratios of Gamma
functions~\cite[Eq.~5.11.13]{DLMF} and the Euler--Maclaurin
formula~\cite[Eq.~2.10.1]{DLMFEulerMaclaurin}, with $Q$ and $B$ as in
\eqref{eq:pa-bias-constants}:
\begin{align}
 \sum_{j>d}\pi_j
 &=c^\alpha d^{-\alpha}\{1-\alpha Q/d+O(d^{-2})\},
 \label{eq:pa-strict-expansion}\\
 \pi_d&=\alpha c^\alpha d^{-\alpha-1}\{1+O(d^{-1})\},
 \label{eq:pa-atom}\\
 \sum_{j>d}\pi_j\log(j/d)
 &=c^\alpha d^{-\alpha}\left\{\frac1\alpha
 -\frac{\alpha(\delta+\alpha/2)}{(\alpha+1)d}+O(d^{-2})\right\}.
 \label{eq:pa-log-expansion}
\end{align}
Equation~\arxiveqref{M-eq:pa-stationary} gives the exact ratio
$\pi_d/\sum_{j>d}\pi_j=\alpha/(d+\delta)$ for every integer $d\ge m$, and
this ratio and \eqref{eq:pa-strict-expansion} give \eqref{eq:pa-atom}.
Dividing \eqref{eq:pa-log-expansion} by \eqref{eq:pa-strict-expansion} gives
\begin{equation}\label{eq:pa-population-hill}
 \frac{\sum_{j>d}\pi_j\log(j/d)}{\sum_{j>d}\pi_j}
       =\frac1\alpha+\frac B d+O(d^{-2}).
\end{equation}

The ingredients of $C_2(\alpha,m,\theta)$ are the coefficients $q_2$ and $\ell_2$ of $d^{-2}$ in the
expansions of $c^{-\alpha}d^\alpha\sum_{j>d}\pi_j$ and
$c^{-\alpha}d^\alpha\sum_{j>d}\pi_j\log(j/d)$, and its value $B_2$ at
$\theta=0$:
\begin{align}
 q_2&=\frac{\alpha(\alpha+1)}{2}Q^2+\frac{\alpha(\alpha^2-1)}{24},
 \qquad
 \ell_2=\frac{q_2+\alpha Q/2+1/3}{\alpha+2}
        -\frac{\alpha Q+1/2}{2}+\frac{\alpha+1}{12},\notag\\
 B_2&=\ell_2+\alpha QB-\frac{q_2}{\alpha},\qquad
 C_2(\alpha,m,\theta)=B_2-\alpha\theta B+\theta(\delta+\alpha\theta).
 \label{eq:pa-C2}
\end{align}

Carrying~\cite[Eq.~5.11.13]{DLMF} to $O(d^{-3})$ gives
$\sum_{j>d}\pi_j=c^\alpha d^{-\alpha}\{1-\alpha Q/d+q_2/d^2+O(d^{-3})\}$
with $q_2$ in \eqref{eq:pa-C2}. Summation by parts writes
$\sum_{j>d}\pi_j\log(j/d)=\sum_{i\ge d}(\sum_{j>i}\pi_j)\log(1+1/i)$, and
Euler--Maclaurin summation~\cite[Eq.~2.10.1]{DLMFEulerMaclaurin} gives its
$d^{-\alpha-2}$ coefficient
$c^\alpha\ell_2$, with $\ell_2$ in \eqref{eq:pa-C2}. Division by
$\sum_{j>d}\pi_j+\theta\pi_d=\{1+\alpha\theta/(d+\delta)\}\sum_{j>d}\pi_j$
yields
\[
 \frac{\sum_{j>d}\pi_j\log(j/d)}{\sum_{j>d}\pi_j+\theta\pi_d}
 =\frac1\alpha+\frac{B-\theta}{d}
   +\frac{C_2(\alpha,m,\theta)}{d^2}+O(d^{-3})
\]
uniformly for $0\le\theta\le1$, with $C_2(\alpha,m,\theta)$ as in
\eqref{eq:pa-C2}.

For a real threshold $u\to\infty$, such as $D_{(k+1)}(n)y$ in
\eqref{eq:Dk}, $\sum_{j>u}\pi_j=\sum_{j>\lfloor u\rfloor}\pi_j$, and
substituting
$\lfloor u\rfloor^{-\alpha}=u^{-\alpha}
\{1+\alpha(u-\lfloor u\rfloor)/u+O(u^{-2})\}$
into \eqref{eq:pa-strict-expansion} gives
\begin{equation}\label{eq:pa-real-tail}
 \sum_{j>u}\pi_j=c^\alpha u^{-\alpha}
 \{1+\alpha(u-\lfloor u\rfloor-Q)/u+O(u^{-2})\}.
\end{equation}
The logarithmic sum $\sum_{j>u}\pi_j\log(j/u)$ satisfies
\eqref{eq:pa-log-expansion} with $d$ replaced by $u$, with no
first-order term in $u-\lfloor u\rfloor$. To see this, write
$d=\lfloor u\rfloor$ and add $\log(d/u)\sum_{j>d}\pi_j$ to
$\sum_{j>d}\pi_j\log(j/d)$. The first-order change of
$c^\alpha d^{-\alpha}/\alpha$ under $d\mapsto u$ cancels the added term.
Sec.~B.2 of~\cite{Wang2026} gives the corresponding first-order
coefficient $B-(u-\lfloor u\rfloor)$ of
$\sum_{j>u}\pi_j\log(j/u)/\sum_{j>u}\pi_j$.

Take $K_n=a_n$ and $b_n=c\sqrt{a_n}$ as in
\arxiveqref{M-eq:balance}. Lemma~3.1 of~\cite{Wang2026} bounds
the departure $e_j(n)=\E N_{n,=j}-n\pi_j$ of the expected degree
profile from its stationary limit:
\begin{equation}\label{eq:pa-profile-error}
 \sum_j|e_j(n)|\le C,\qquad \sum_jj|e_j(n)|\le Cn^{1/\alpha}.
\end{equation}
Consequently $|\E N_{n,>d}-n\sum_{j>d}\pi_j|\le C$ for every $d$, and,
for $d$ comparable to $b_n$, Jensen's inequality gives
$\sum_{j>d}|e_j(n)|\log(j/d)\le C(1+\log a_n)$. Thus, at thresholds of
order $b_n$, the stationary expansions give $\E N_{n,>d}$ and
$\sum_{j>d}\E N_{n,=j}\log(j/d)$ up to $o(\sqrt{a_n})$, smaller than the
order $\sqrt{a_n}$ of their centered fluctuations.

\subsection{A Gaussian limit at an integer threshold}\label{sec:pa-degree-clt}
Recall $\mathcal A(t,\theta)$ from \arxiveqref{M-eq:pa-A}, where $x=t^{-1/\alpha}$
is the limit of $D_{(k+1)}(n)/b_n$ when $k/a_n\to t$, because
$nc^\alpha b_n^{-\alpha}=a_n$. The next theorem gives the limit of $\sqrt{a_n}(H_{n,k}-1/\alpha)$ at a
deterministic integer threshold.

\begin{theorem}\label{thm:pa-degree-clt}
Fix $t>0$ and $\theta\in[0,1]$, and let deterministic integers $d_n$ satisfy
$d_n/b_n\to x=t^{-1/\alpha}$ in $\R$. Choose the rank
\[
 k=N_{n,>d_n}+\min\{\lfloor\theta N_{n,=d_n}\rfloor,N_{n,=d_n}-1\}
\]
when $N_{n,=d_n}>0$ and this expression is at least two, and set $k=2$
otherwise. The minimum specifies the number of retained observations
tied at $d_n$, leaving at least one at the threshold. Then
$\Prob\{D_{(k+1)}(n)=d_n\}\to1$, $k/a_n\xrightarrow{p}t$ in $\R$, and
\begin{equation}\label{eq:pa-degree-limit}
 \sqrt{a_n}(H_{n,k}-1/\alpha)\wto\mathcal A(t,\theta)
 \quad\hbox{in }\R.
\end{equation}
In particular,
\begin{equation}\label{eq:pa-degree-normal}
 \sqrt{k}(\wh\alpha_{n,k}-\alpha)
 \wto N\left(-\frac{\alpha^2}{c}t^{1/2+1/\alpha}(B-\theta),
                   \alpha^4V(h_1,h_1)\right)\quad\hbox{in }\R.
\end{equation}
The convergence \eqref{eq:pa-degree-limit} holds jointly for finitely
many pairs $(t,\theta)$, with all limits defined from one process
$\mathbb G$.
\end{theorem}
\begin{proof}
Throughout, take $K_n=a_n$ in \eqref{eq:pa-measure}, and write
$\log_+u=\max(\log u,0)$ for $u>0$.
We expand $k$ and $\sum_{v=1}^n\log_+(D_v(n)/d_n)$ jointly and take
their ratio, which is $H_{n,k}$ on $\{D_{(k+1)}(n)=d_n\}$.
The variance bound of~\cite[Lemma~3.3]{Wang2026} is
\begin{equation}\label{eq:pa-variance-domination}
 \E S_n(f)^2\le C\E\nu_nf^2
\end{equation}
for bounded measurable $f$, with $C$ depending only on $m,\delta$.
Truncation and the tail bound of~\cite[Lemma~3.2]{Wang2026} extend this bound to the
functions $z\mapsto\log_+(z/y)$, $y>0$. Since $d_n/b_n\to x>0$,
\eqref{eq:pa-real-tail} and \eqref{eq:pa-profile-error} give
\[
 \E\nu_n\bigl(\min\{d_n/b_n,x\},\max\{d_n/b_n,x\}\bigr]\longrightarrow0.
\]
Thus \eqref{eq:pa-variance-domination} implies
$S_n(f_{d_n/b_n})-S_n(f_x)\xrightarrow{p}0$ in $\R$.
For $z>0$,
\[
 \left|\log_+(zb_n/d_n)-\log_+(z/x)\right|
 \le\left|\log\frac{d_n}{b_nx}\right|
             \mathbf1\{z>\min(d_n/b_n,x)\}.
\]
For large $n$, the $\nu_n$-expectation of the squared right side is at
most $\log^2\{d_n/(b_nx)\}\,\E\nu_n(x/2,\infty)\to0$. Applying
\eqref{eq:pa-variance-domination} and the joint centered limits of
Sec.~\ref{app:pa-tail-measure} therefore gives
\begin{equation}\label{eq:pa-moving-centered}
 \left(
 \frac{N_{n,>d_n}-\E N_{n,>d_n}}{\sqrt{a_n}},
 \frac{\sum_{v=1}^n\log_+(D_v(n)/d_n)
       -\E\sum_{v=1}^n\log_+(D_v(n)/d_n)}{\sqrt{a_n}}
 \right)
 \wto(\mathbb G(x),\mathbb L(x))
 \quad\hbox{in }\R^2.
\end{equation}
Here $\mathbb L(x)$, the integral in \eqref{eq:pa-G}, is the second limit
because $\log_+(z/x)=h_x(z)+f_x(z)/\alpha$ and
$S_n(h_x)\wto\mathbb L(x)-\mathbb G(x)/\alpha$.

The number $N_{n,=d_n}$ of vertices of degree $d_n$ diverges in
probability. Indeed, \eqref{eq:pa-atom}, \eqref{eq:pa-profile-error},
and \eqref{eq:pa-variance-domination} give
\[
 \frac{\E N_{n,=d_n}}{\sqrt{a_n}}
 =\frac{\alpha}{c}(d_n/b_n)^{-\alpha-1}+O(a_n^{-1/2})
 \longrightarrow\frac{\alpha}{c}x^{-\alpha-1},
 \qquad
 \Var\left(\frac{N_{n,=d_n}}{\sqrt{a_n}}\right)
 \le\frac{C\E N_{n,=d_n}}{a_n}\longrightarrow0.
\]
Chebyshev's inequality proves
\begin{equation}\label{eq:pa-tie-size}
 \frac{N_{n,=d_n}}{\sqrt{a_n}}\xrightarrow{p}
       \frac{\alpha}{c}x^{-\alpha-1}>0
 \quad\hbox{in }\R.
\end{equation}
Equations \eqref{eq:pa-strict-expansion}, \eqref{eq:pa-profile-error},
and \eqref{eq:pa-moving-centered} also
give $N_{n,>d_n}/a_n\xrightarrow{p}t>0$ in $\R$.
Consequently, with probability tending to one, $k$ is given by the
expression in Theorem~\ref{thm:pa-degree-clt}, $D_{(k+1)}(n)=d_n$, and
$\theta_{n,k}=(k-N_{n,>d_n})/N_{n,=d_n}$, so that
\[
 \left|\frac{k-N_{n,>d_n}}{N_{n,=d_n}}-\theta\right|
 \le\frac1{N_{n,=d_n}}.
\]
Equation~\eqref{eq:pa-tie-size} makes the right side tend to zero in
probability, so $\theta_{n,k}\xrightarrow{p}\theta$ in $\R$.

The expectations $\E N_{n,>d_n}$ and
$\E\sum_{v=1}^n\log_+(D_v(n)/d_n)$ come from
\eqref{eq:pa-strict-expansion}--\eqref{eq:pa-log-expansion}.
Since $nc^\alpha d_n^{-\alpha}=a_n(d_n/b_n)^{-\alpha}$, those
expansions and \eqref{eq:pa-profile-error} give
\begin{align}
 \frac{\E N_{n,>d_n}-nc^\alpha d_n^{-\alpha}}{\sqrt{a_n}}
 &\longrightarrow-\frac{\alpha Q}{c}x^{-\alpha-1},
 \notag\\
 \frac{\E\sum_{v=1}^n\log_+(D_v(n)/d_n)
              -nc^\alpha d_n^{-\alpha}/\alpha}{\sqrt{a_n}}
 &\longrightarrow
 -\frac{\alpha(\delta+\alpha/2)}{c(\alpha+1)}x^{-\alpha-1}.
 \label{eq:pa-moving-means}
\end{align}
Before taking $d_n/b_n\to x$, the two first-order terms are proportional
to $(d_n/b_n)^{-\alpha-1}$ and the remainders are bounded by
$C(1+\log a_n)/\sqrt{a_n}$. Adding
\eqref{eq:pa-moving-centered}, \eqref{eq:pa-tie-size}, and
\eqref{eq:pa-moving-means} proves the following limits jointly in $\R^2$:
\begin{align}
 \frac{k-nc^\alpha d_n^{-\alpha}}{\sqrt{a_n}}
 &\wto\mathbb G(x)+\frac{\alpha}{c}(\theta-Q)x^{-\alpha-1},
 \label{eq:pa-k-expansion}\\
 \frac{\sum_{v=1}^n\log_+(D_v(n)/d_n)
                         -nc^\alpha d_n^{-\alpha}/\alpha}{\sqrt{a_n}}
 &\wto\mathbb L(x)-
 \frac{\alpha(\delta+\alpha/2)}{c(\alpha+1)}x^{-\alpha-1}.
 \label{eq:pa-sum-expansion}
\end{align}
In particular, \eqref{eq:pa-k-expansion} implies
$k/a_n\xrightarrow{p}x^{-\alpha}=t$ in $\R$.

On $\{D_{(k+1)}(n)=d_n\}$, the Hill statistic satisfies the exact identity
\begin{equation}\label{eq:pa-H-ratio}
 \sqrt{a_n}(H_{n,k}-1/\alpha)
 =\frac{a_n}{k}\left[
 \frac{\sum_{v=1}^n\log_+(D_v(n)/d_n)
                     -nc^\alpha d_n^{-\alpha}/\alpha}{\sqrt{a_n}}
 -\frac{k-nc^\alpha d_n^{-\alpha}}{\alpha\sqrt{a_n}}\right].
\end{equation}
Substituting \eqref{eq:pa-k-expansion}--\eqref{eq:pa-sum-expansion}
and $k/a_n\xrightarrow{p}t$ into \eqref{eq:pa-H-ratio} gives the weak limit
\[
 \frac{\mathbb L(x)-\mathbb G(x)/\alpha}{t}
 +\frac{Q-\theta-\alpha(\delta+\alpha/2)/(\alpha+1)}{cx}
 =\mathcal A(t,\theta),
\]
because $Q-\alpha(\delta+\alpha/2)/(\alpha+1)=B$.
This proves \eqref{eq:pa-degree-limit}. The identity
\[
 \sqrt{k}(\wh\alpha_{n,k}-\alpha)
 =-\frac{\alpha}{H_{n,k}}\sqrt{\frac{k}{a_n}}\,
                  \sqrt{a_n}(H_{n,k}-1/\alpha)
\]
and $H_{n,k}\xrightarrow{p}1/\alpha$ in $\R$ yield
\eqref{eq:pa-degree-normal}, because $\mathbb L(x)-\mathbb G(x)/\alpha$
has variance $V(h_x,h_x)=tV(h_1,h_1)$ by \eqref{eq:pa-homogeneity}.
The centered limits of Sec.~\ref{app:pa-tail-measure} hold
jointly for finitely many functions, so \eqref{eq:pa-moving-centered}
holds jointly for finitely many thresholds, which proves the joint
assertion.
\end{proof}

Theorem~2.1 of~\cite{Wang2026} gives a centered normal limit for
$\sqrt k(H_{n,k}-1/\alpha)$, and hence by the delta method for
$\sqrt k(\wh\alpha_{n,k}-\alpha)$, when
$\sqrt k(k/n)^{1/\alpha}\to0$, equivalently $k=o(a_n)$.
Sec.~B.2 of~\cite{Wang2026} computes the term $B/d$ in
\eqref{eq:pa-population-hill}, and Lemma~3.2 there shows that the mean
bias is negligible at those ranks. At
$k$ of order $a_n$ this term persists, and the mean of the normal limit
in \eqref{eq:pa-degree-normal} depends on $\theta$ through $B-\theta$.

\subsection{Limits at a deterministic rank of order \texorpdfstring{$a_n$}{a(n)}}\label{sec:pa-rank}
Theorem~\ref{thm:pa-degree-clt} fixes the threshold $d_n$ and the tie
fraction $\theta$, and lets the rank $k$ be random. Fixing instead $t>0$
and the rank $k=\lfloor a_nt\rfloor$, with $x=t^{-1/\alpha}$, makes the
threshold $D_{(k+1)}(n)$ and the tie fraction $\theta_{n,k}$ random. The tie
fraction enters $\mathcal A$ through $B-\theta$, so identifying the limit of
$\sqrt{a_n}(H_{n,k}-1/\alpha)$ requires the joint law of the Gaussian
terms and the tie fraction. We obtain this law by locating
$D_{(k+1)}(n)$ among finitely many integers near $b_nx$, then letting
the number of integers increase.

By \eqref{eq:pa-atom} and \eqref{eq:pa-profile-error},
$\E N_{n,=d}/\sqrt{a_n}\to w_x$ for integers $d$ with $d/b_n\to x$, where
\[
 w_x=\frac{\alpha}{c}x^{-\alpha-1}.
\]
Suppose a subsequence satisfies
$b_nx-\lfloor b_nx\rfloor\to\zeta\in[0,1]$, and define
\begin{equation}\label{eq:pa-rounding}
 U_\zeta=\frac{\mathbb G(x)}{w_x}+\zeta-Q,
 \qquad J_\zeta=\lceil U_\zeta\rceil,
 \qquad \Theta_\zeta=J_\zeta-U_\zeta.
\end{equation}
Sec.~\ref{sec:pa-rank} shows that, along this subsequence
and jointly over any finite set of integers $j$,
\[
 \frac{N_{n,>\lfloor b_nx\rfloor+j}-k}{w_x\sqrt{a_n}}\wto U_\zeta-j
 \quad\hbox{in }\R.
\]
Thus $U_\zeta$ is the limiting excess of $N_{n,>\lfloor b_nx\rfloor}$ over
$k$, in units of $w_x\sqrt{a_n}$.
The threshold $D_{(k+1)}(n)$ equals $\lfloor b_nx\rfloor+j$ precisely when
$N_{n,>\lfloor b_nx\rfloor+j}\le k<N_{n,>\lfloor b_nx\rfloor+j-1}$.
In the limit this reads $j-1<U_\zeta\le j$, so
$D_{(k+1)}(n)-\lfloor b_nx\rfloor$ and $\theta_{n,k}$ converge in
distribution to $J_\zeta$ and $\Theta_\zeta$.
Equation~\eqref{eq:pa-indicator-variance} gives
$\Var U_\zeta=V(f_x,f_x)/w_x^2>0$, so $u\mapsto\lceil u\rceil$ is almost
surely continuous at $U_\zeta$.
\begin{proposition}\label{prop:pa-fixed-rank}
Fix $t>0$, and put $x=t^{-1/\alpha}$ and $k=\lfloor a_nt\rfloor$. Along
every subsequence on which $b_nx-\lfloor b_nx\rfloor\to\zeta\in[0,1]$,
\begin{equation}\label{eq:pa-fixed-rank-law}
 \left(D_{(k+1)}(n)-\lfloor b_nx\rfloor,
       \theta_{n,k},\sqrt{a_n}(H_{n,k}-1/\alpha)\right)
 \wto\left(J_\zeta,\Theta_\zeta,\mathcal A(t,\Theta_\zeta)\right)
 \quad\hbox{in }\R^3.
\end{equation}
\end{proposition}
\begin{proof}
We first locate $D_{(k+1)}(n)$ among the integers near $b_nx$, with
$x=t^{-1/\alpha}$, then evaluate
$\theta_{n,k}$ and $\sum_{v=1}^n\log_+(D_v(n)/D_{(k+1)}(n))$.
For a fixed integer $j$, take $d=\lfloor b_nx\rfloor+j$.
Taylor's formula gives
\[
 \frac{nc^\alpha d^{-\alpha}-a_nt}{\sqrt{a_n}}
 =-w_x\{j-(b_nx-\lfloor b_nx\rfloor)\}+O(a_n^{-1/2}).
\]
Along the subsequence on which $b_nx-\lfloor b_nx\rfloor\to\zeta$,
combining this expansion with
\eqref{eq:pa-moving-centered}, \eqref{eq:pa-moving-means}, and
$|k-a_nt|\le1$ yields, jointly for every finite collection of integers $j$,
\begin{align}
 \frac{N_{n,>d}-k}{\sqrt{a_n}}
 &\wto w_x(U_\zeta-j),\label{eq:pa-integer-survival}\\
 \frac{N_{n,=d}}{\sqrt{a_n}}
 &\xrightarrow{p}w_x,\label{eq:pa-integer-mass}\\
 \frac{\sum_{v=1}^n\log_+(D_v(n)/d)-k/\alpha}{\sqrt{a_n}}
 &\wto\mathbb L(x)-\frac{t}{cx}
       \left\{j-\zeta+\frac{\alpha(\delta+\alpha/2)}{\alpha+1}\right\}.
 \label{eq:pa-integer-log}
\end{align}
The weak convergence is in the corresponding finite-dimensional
Euclidean space, and \eqref{eq:pa-integer-mass} holds in $\R$.
The Gaussian terms are common to all these $j$ because $d/b_n\to x$ for
each $j$, and \eqref{eq:pa-tie-size} with $d_n=d$ gives
\eqref{eq:pa-integer-mass}.

For every integer $d$, $D_{(k+1)}(n)>d$ if and only if $N_{n,>d}>k$.
Hence, for every positive integer $M$,
\begin{align*}
 \Prob\{D_{(k+1)}(n)-\lfloor b_nx\rfloor>M\}
 &=\Prob\{N_{n,>\lfloor b_nx\rfloor+M}>k\},\\
 \Prob\{D_{(k+1)}(n)-\lfloor b_nx\rfloor<-M\}
 &=\Prob\{N_{n,>\lfloor b_nx\rfloor-M-1}\le k\}.
\end{align*}
Since $U_\zeta$ is normal with positive variance by
\eqref{eq:pa-indicator-variance}, the portmanteau theorem and
\eqref{eq:pa-integer-survival} give
\begin{equation}\label{eq:pa-integer-tightness}
 \limsup_n\Prob\{|D_{(k+1)}(n)-\lfloor b_nx\rfloor|>M\}
 \le\Prob\{U_\zeta>M\}+\Prob\{U_\zeta\le-M-1\}
 \longrightarrow0
\end{equation}
as $M\to\infty$.

For $|j|\le M$, the event $D_{(k+1)}(n)=\lfloor b_nx\rfloor+j$ is
\[
 \{N_{n,>\lfloor b_nx\rfloor+j}\le k
            <N_{n,>\lfloor b_nx\rfloor+j-1}\}.
\]
In the limit \eqref{eq:pa-integer-survival} these inequalities become
$j-1<U_\zeta\le j$, whose boundaries have probability zero. By the
continuous mapping theorem, the indicator of this event therefore
converges jointly with \eqref{eq:pa-integer-survival} to
$\mathbf1\{J_\zeta=j\}$, where $J_\zeta=\lceil U_\zeta\rceil$.
On this event the exact identities
\begin{equation}\label{eq:pa-selected-integer}
 \theta_{n,k}=\frac{k-N_{n,>d}}{N_{n,=d}},
 \qquad
 \sqrt{a_n}(H_{n,k}-1/\alpha)
 =\frac{a_n}{k}
       \frac{\sum_{v=1}^n\log_+(D_v(n)/d)-k/\alpha}{\sqrt{a_n}}
\end{equation}
combine with \eqref{eq:pa-integer-survival}--\eqref{eq:pa-integer-log}
to give the joint limit on $\{J_\zeta=j\}$
\[
 \left(J_\zeta,\ J_\zeta-U_\zeta,\
 \frac{\mathbb L(x)}t-\frac{J_\zeta-\zeta}{cx}
       -\frac{\alpha(\delta+\alpha/2)}{c(\alpha+1)x}\right).
\]
By \eqref{eq:pa-integer-tightness}, the events with $|j|>M$ have
probability tending to zero as $M\to\infty$, for the sample and for
$J_\zeta$. Finally,
$w_x=\alpha t/(cx)$ and \eqref{eq:pa-rounding} give
\[
 \mathcal A(t,\Theta_\zeta)
 =\frac{\mathbb L(x)}t+
   \frac{B-J_\zeta+\zeta-Q}{cx}
 =\frac{\mathbb L(x)}t-\frac{J_\zeta-\zeta}{cx}
       -\frac{\alpha(\delta+\alpha/2)}{c(\alpha+1)x}.
\]
The third coordinate is therefore $\mathcal A(t,\Theta_\zeta)$, which
proves \eqref{eq:pa-fixed-rank-law}.
\end{proof}

Since $b_n\to\infty$, the first coordinate of
\eqref{eq:pa-fixed-rank-law} implies
$D_{(k+1)}(n)/b_n\xrightarrow{p}x$ in $\R$ along every such subsequence,
hence along the full sequence, because every subsequence has a further
subsequence on which $b_nx-\lfloor b_nx\rfloor$ converges. The
proposition also gives the limit law of the integer
$D_{(k+1)}(n)-\lfloor b_nx\rfloor$, on which $\theta_{n,k}$ depends.
The third coordinate of the limit contains the rounding
$J_\zeta=\lceil U_\zeta\rceil$ of a Gaussian variable, and its law could
depend on $\zeta$. Retaining all or none of the vertices tied at
$D_{(k+1)}(n)$, instead of the first $k$ ranks, could change the limit of
$\sqrt{a_n}(H_{n,k}-1/\alpha)$.

\subsection{Exponential bounds}\label{sec:pa-functional-proof-input}
The variance bound of~\cite[Lemma~3.3]{Wang2026} controls one
function at a time. Chaining over thresholds needs an exponential
bound, and Lemma~\ref{lem:pa-exponential} supplies one by comparing
the moment generating function of a sum of functions of PA degrees
with a product over vertices. The functions of the degrees in
Lemma~\ref{lem:pa-exponential} need not be monotone.
\begin{lemma}\label{lem:pa-exponential}
For bounded functions $g_v$ on the possible degrees, let
$Z_v=g_v(D_v(n))-\E g_v(D_v(n))$. For every real $\lambda$,
\begin{equation}\label{eq:pa-product-mgf}
 \E\exp\left(\lambda\sum_{v=1}^n Z_v\right)
 \le\prod_{v=1}^n\{\E e^{2\lambda Z_v}\}^{1/2}.
\end{equation}
In particular, for an interval $A\subset\R$, a positive normalizer $K$, and
$b=c(n/K)^{1/\alpha}$, put
\[
 \mu_n(A)=K^{-1}\sum_v\Prob\{D_v(n)/b\in A\},\qquad
 S_n(A)=K^{-1/2}\sum_v
   [\mathbf1\{D_v(n)/b\in A\}-\Prob\{D_v(n)/b\in A\}].
\]
For every $u>0$,
\begin{equation}\label{eq:pa-Bernstein}
 \Prob\{|S_n(A)|>2\sqrt{\mu_n(A)u}+2u/\sqrt K\}\le2e^{-u}.
\end{equation}
\end{lemma}
\begin{proof}
We bound the conditional expectation of a product over vertices across
each added edge and iterate the bound. First let $(I_v)$ be a finite vector
of indicators with $\sum_v I_v=1$ almost surely. For positive functions $q_v:\{0,1\}\to(0,\infty)$,
we claim
\begin{equation}\label{eq:categorical-product}
 \E\prod_v q_v(I_v)^{1/2}
       \le\prod_v\{\E q_v(I_v)\}^{1/2}.
\end{equation}
Put $\mathcal V_+=\{v:q_v(1)\ge q_v(0)\}$ and
$\mathcal V_-=\{v:q_v(1)<q_v(0)\}$. The vector $(I_v)$ has a
multinomial distribution with one trial, so its coordinates are
negatively associated~\cite{JoagDevProschan1983}. Since the $q_v$ are
nondecreasing on $\mathcal V_+$ and nonincreasing on $\mathcal V_-$,
negative association gives
$\E\prod_{v\in\mathcal V_\pm}q_v(I_v)\le\prod_{v\in\mathcal V_\pm}\E q_v(I_v)$.
Cauchy--Schwarz
then gives
\[
 \E\prod_v q_v(I_v)^{1/2}
 \le
 \left\{\E\prod_{v\in\mathcal V_+}q_v(I_v)\right\}^{1/2}
 \left\{\E\prod_{v\in\mathcal V_-}q_v(I_v)\right\}^{1/2}
 \le\prod_v\{\E q_v(I_v)\}^{1/2}.
\]

We now apply \eqref{eq:categorical-product} to each added edge
$e=(t,i)$, the $i$th edge sent by the vertex arriving at time $t$ in
\arxiveqref{M-eq:pa-model}. Let $\mathcal B_e$ be the
sigma-field generated by the graph through that edge, and use $e-$
for the instant before it. The denominator in \arxiveqref{M-eq:pa-model} for
this edge is the deterministic number $W_e=(2m+\delta)(t-1)-2m+i-1$.
For an eligible vertex $v<t$, set
\[
 q_{v,e}(d)=\E[e^{2\lambda Z_v}\mid D_v(e)=d].
\]
Given $\mathcal B_{e-}$, vertex $v$ receives edge $e$ with probability
$(D_v(e-)+\delta)/W_e$, so the degree of $v$ is a Markov chain with
respect to $(\mathcal B_e)$, as used in the proof of~\cite[Lemma~3.3]{Wang2026}.
Hence $\E[e^{2\lambda Z_v}\mid\mathcal B_e]=q_{v,e}(D_v(e))$, and the
backward recursion is
\begin{equation}\label{eq:pa-exponential-recursion}
 q_{v,e-}(d)=
 \left(1-\frac{d+\delta}{W_e}\right)q_{v,e}(d)
       +\frac{d+\delta}{W_e}q_{v,e}(d+1).
\end{equation}
Until vertex $v$ becomes eligible, interpret
$q_{v,e}(D_v(e))$ as the constant $\E e^{2\lambda Z_v}$.
This convention is valid because $v$ has the deterministic degree $m$
at eligibility and all later $W_e$ are deterministic.

Conditionally on $\mathcal B_{e-}$, apply \eqref{eq:categorical-product}
with $I_v$ the indicator that the eligible vertex $v<t$ receives $e$,
$q_v(0)=q_{v,e}(D_v(e-))$, and $q_v(1)=q_{v,e}(D_v(e-)+1)$.
The factors for vertices not yet eligible are constant. Equation
\eqref{eq:pa-exponential-recursion} therefore gives
\begin{equation}\label{eq:pa-exponential-conditional}
 \E\left[\prod_{v=1}^n q_{v,e}(D_v(e))^{1/2}
                      \,\middle|\,\mathcal B_{e-}\right]
 \le\prod_{v=1}^n q_{v,e-}(D_v(e-))^{1/2}.
\end{equation}
When a vertex becomes eligible, its constant factor equals
$q_{v,e}(m)$, so the product on the left of
\eqref{eq:pa-exponential-conditional} is unchanged at that instant.
At the deterministic seed graph this product equals
$\prod_{v=1}^n(\E e^{2\lambda Z_v})^{1/2}$.
After the last edge it equals $e^{\lambda\sum_{v=1}^n Z_v}$.
Iterating expectations in \eqref{eq:pa-exponential-conditional}
proves \eqref{eq:pa-product-mgf}.

For indicator summands $Z_v=\mathbf1\{D_v(n)/b\in A\}-\Prob\{D_v(n)/b\in A\}$,
we have $|Z_v|\le1$ and
$\sum_v\E Z_v^2\le K\mu_n(A)$. Bernstein's inequality for bounded
variables~\cite[Chapter~2]{BoucheronLugosiMassart2013} gives
$\log\E e^{sZ_v}\le s^2\E Z_v^2/[2(1-|s|/3)]$ for $|s|<3$.
Apply this inequality with $s=2\lambda/\sqrt K$ in
\eqref{eq:pa-product-mgf} to obtain
\begin{equation}\label{eq:pa-mgf-Bernstein}
 \log\E e^{\lambda S_n(A)}
 \le\frac{\lambda^2\mu_n(A)}{1-2|\lambda|/(3\sqrt K)},
 \qquad |\lambda|<3\sqrt K/2.
\end{equation}
The exponential Markov inequality applied to $S_n(A)$ and $-S_n(A)$, with
\eqref{eq:pa-mgf-Bernstein} optimized over $\lambda$ as for sub-gamma
variables in~\cite[Chapter~2]{BoucheronLugosiMassart2013}, gives
\[
 \Prob\left\{|S_n(A)|>
       2\sqrt{\mu_n(A)u}+\frac{2u}{3\sqrt K}\right\}\le2e^{-u}.
\]
Increasing $2u/(3\sqrt K)$ to $2u/\sqrt K$ proves
\eqref{eq:pa-Bernstein}.
\end{proof}

The product bound \eqref{eq:pa-product-mgf} is where the PA setting
differs from independent observations. Given
\eqref{eq:pa-product-mgf}, the tail bound \eqref{eq:pa-Bernstein} is
the sub-gamma tail bound of Boucheron, Lugosi, and
Massart~\cite[Chapter~2]{BoucheronLugosiMassart2013}.

\subsection{A maximal inequality for nested intervals}\label{app:nested-maximal}
Two maximal inequalities enter the proof of
Theorem~\ref{thm:pa-functional}. It bounds the oscillation of
$\mathbb G_n$ on short intervals by \eqref{eq:nested-local-maximal}, and
$\sup_{x\ge M}|\mathbb G_n(x)|$ by \eqref{eq:nested-maximal} with
$A=(M,\infty)$. Both follow from \eqref{eq:pa-Bernstein} alone, so they
also apply to independent observations satisfying that inequality.
Write $S_n(A)=\sqrt K\{\nu_n(A)-\mu_n(A)\}$, with
$\nu_n(A)=K^{-1}\sum_v\mathbf1\{D_v(n)/b\in A\}$ and $\mu_n=\E\nu_n$ from
Lemma~\ref{lem:pa-exponential}, and assume
\eqref{eq:pa-Bernstein} for every interval $A$.
For a fixed interval $A$ with $\mu_n(A)\le C_0$, put
$A_x=A\cap(x,\infty)$. There is a constant
$C$, depending only on $C_0$, such that, for $K\ge2$,
\begin{equation}\label{eq:nested-maximal}
 \E\sup_x|S_n(A_x)|
 \le C\left\{\sqrt{\mu_n(A)}+
                      \frac{\log^2(K+2)}{\sqrt K}\right\}.
\end{equation}
For intervals $A_1,\ldots,A_N$ with
$\max_{1\le i\le N}\mu_n(A_i)\le\varrho\le C_0$, the corresponding bound is
\begin{equation}\label{eq:nested-local-maximal}
 \E\max_{1\le i\le N}\sup_x|S_n(A_i\cap(x,\infty))|
 \le C\left\{\sqrt{\varrho\log(2N)}+
 \frac{\log(K+2)\{\log(K+2)+\log(2N)\}}{\sqrt K}\right\}.
\end{equation}

\begin{proof}
We chain over finite families of sets that approximate the sets $A_x$
in $\mu_n$-mass. For $j\ge0$, let $\mathcal C_j$ contain
$A\cap[z,\infty)$ and $A\cap(z,\infty)$ for every atom $z$ of $\mu_n$
with $\mu_n\{z\}>\mu_n(A)2^{-j}$, and the sets $A\cap(y,\infty)$ at the
endpoints $y$ of a partition of the remaining mass of $A$ into
consecutive intervals of $\mu_n$-mass at most $2\mu_n(A)2^{-j}$.
No $A_x$ lies strictly between $A\cap(z,\infty)$ and $A\cap[z,\infty)$,
so each $A_x$ satisfies $A''\subset A_x\subset A'$ for some
$A',A''\in\mathcal C_j$ with $\mu_n(A'\setminus A'')\le2\mu_n(A)2^{-j}$,
and $C2^j$ sets suffice. Replacing $\mathcal C_j$ by
$\mathcal C_0\cup\cdots\cup\mathcal C_j$ keeps this size bound and makes
the families increase in $j$. For $j\ge1$, each $A'\in\mathcal C_j$ is
nested with some $A''\in\mathcal C_{j-1}$ with
$\mu_n(A'\triangle A'')\le C\mu_n(A)2^{-j}$. A union bound in
\eqref{eq:pa-Bernstein} over these pairs and integration of the
resulting tail bound give
\[
 \E\max_{A'\in\mathcal C_j}|S_n(A')-S_n(A'')|
 \le C\sqrt{\mu_n(A)(j+1)2^{-j}}+C(j+1)/\sqrt K.
\]
Sum these bounds over $1\le j\le J=\lceil2\log_2(K+2)\rceil$, and bound
the members of $\mathcal C_0$ by \eqref{eq:pa-Bernstein} directly.
If $A''\subset A_x\subset A'$ with $A',A''\in\mathcal C_J$, monotonicity
of $\nu_n$ and $\mu_n$ gives
\[
 |S_n(A_x)-S_n(A')|\le|S_n(A')-S_n(A'')|+2\sqrt K\mu_n(A'\setminus A'').
\]
A final application of \eqref{eq:pa-Bernstein} bounds the expected
supremum over $x$ of the left side by
\[
 C\sqrt{\mu_n(A)(J+1)2^{-J}}+C(J+1)/\sqrt K
                    +C\sqrt K\mu_n(A)2^{-J}.
\]
Since $2^{-J}\le(K+2)^{-2}$, this proves \eqref{eq:nested-maximal}.
Taking each union bound over all $N$ intervals $A_i$ replaces $j+1$ by
$j+1+\log(2N)$ and $\mu_n(A)$ by $\varrho$, which proves
\eqref{eq:nested-local-maximal}.
\end{proof}

\subsection{Functional convergence}\label{sec:pa-functional-result}
Sec.~\ref{app:pa-tail-measure} derives the joint limits of
$S_n(f_x)$ at finitely many fixed $x$ from Proposition~3.7 and
Section~3.4 of~\cite{Wang2026}. Lemma~3.8 there also controls
$\nu_n(x,\infty)$ for $x$ within $O(K_n^{-1/2})$ of $1$
under the additional condition $\sqrt{K_n}/b_n\to0$, which Section~4
there shows is necessary for that local approximation and which fails at
$K_n=a_n$, where $\sqrt{a_n}/b_n=1/c$.
The criterion $\mathcal D_k$ in \eqref{eq:Dk} takes a supremum over $y\ge1$ and
is then minimized over ranks, so its limit requires convergence of
$S_n(f_x)$ uniformly in $x\in[\epsilon,\infty)$ without that additional
condition.
Theorem~\ref{thm:pa-functional} establishes this convergence.
Its proof uses \eqref{eq:nested-local-maximal} to control increments
on compact intervals and \eqref{eq:nested-maximal} to control
$\sup_{x\ge M}|S_n(f_x)|$ as $M\to\infty$.
Sec.~\ref{app:nested-maximal} derives both inequalities from
\eqref{eq:pa-Bernstein}.
\begin{theorem}\label{thm:pa-functional}
For every deterministic $K_n\to\infty$ with $K_n/n\to0$, take
$b_n=c(n/K_n)^{1/\alpha}$ and $\nu_n$ as in \eqref{eq:pa-measure}.
For every $\epsilon>0$,
\begin{equation}\label{eq:pa-functional}
 \mathbb G_n(x)=\sqrt{K_n}\{\nu_n(x,\infty)-\E\nu_n(x,\infty)\}
 \wto\mathbb G(x)
 \quad\hbox{in }\ell^\infty([\epsilon,\infty)).
\end{equation}
The limit is the continuous Gaussian process in \eqref{eq:pa-G}, with
$\mathbb G(\infty)=0$. Jointly with \eqref{eq:pa-functional}, the
integrals $\int_x^\infty\mathbb G_n(z)\dd z/z$ converge in distribution
to $\mathbb L(x)$ in $\ell^\infty$ of every compact subinterval of
$(0,\infty)$.
\end{theorem}
\begin{proof}
Take $K=K_n$ and $b=b_n$ in Lemma~\ref{lem:pa-exponential}, so that
$\mu_n=\E\nu_n$ and $\mathbb G_n(x)=S_n((x,\infty))$.
The joint centered limits of Sec.~\ref{app:pa-tail-measure}
give the finite-dimensional convergence, with covariance
\eqref{eq:pa-covariance}. We prove asymptotic tightness by controlling
increments on compact intervals and the supremum over $[M,\infty)$.
Equations \eqref{eq:pa-real-tail}, \eqref{eq:pa-atom}, and
\eqref{eq:pa-profile-error} give, for every $0<e<E<\infty$,
\begin{equation}\label{eq:pa-mean-uniform}
 \sup_{e\le x\le E}|\mu_n(x,\infty)-x^{-\alpha}|\longrightarrow0,
 \qquad
 \sup_{eb_n\le d\le Eb_n}\mu_n\{d/b_n\}\longrightarrow0.
\end{equation}
Fix $0<\varrho<1$. For all sufficiently large $n$, every atom of $\mu_n$
in $[e,E]$ has mass at most $\varrho/4$. Starting at $e$, close each
interval as soon as its $\mu_n$-mass reaches $\varrho/2$. Except for the
last interval, its mass is then between $\varrho/2$ and $3\varrho/4$.
Equation~\eqref{eq:pa-mean-uniform} bounds the number of intervals by
$C/\varrho$. Every pair $x,y$ with
$\mu_n((x\wedge y,x\vee y])\le\varrho/2$ lies in the union of at most
three such intervals, and \eqref{eq:nested-local-maximal} applied to
these unions gives
\[
 \limsup_n\E\sup_{\substack{x,y\in[e,E]\\
                   \mu_n((x\wedge y,x\vee y])\le\varrho/2}}
       |\mathbb G_n(x)-\mathbb G_n(y)|
 \le C\sqrt{\varrho\log(2/\varrho)}.
\]
For $e\le x<y\le E$, the mean value theorem gives
\[
 \mu_n((x,y])\le\alpha e^{-\alpha-1}|y-x|
       +2\sup_{e\le z\le E}|\mu_n(z,\infty)-z^{-\alpha}|.
\]
Hence $\alpha e^{-\alpha-1}|x-y|\le\varrho/4$ implies
$\mu_n((x\wedge y,x\vee y])\le\varrho/2$ for large $n$, and Markov's
inequality gives, for every $\eta>0$,
\begin{equation}\label{eq:pa-equicontinuity}
 \limsup_n\Prob\left\{
   \sup_{\substack{x,y\in[e,E]\\\alpha e^{-\alpha-1}|x-y|\le\varrho/4}}
       |\mathbb G_n(x)-\mathbb G_n(y)|>\eta\right\}
 \le\frac C\eta\sqrt{\varrho\log(2/\varrho)}.
\end{equation}
Letting $\varrho\downarrow0$ proves asymptotic tightness in
$\ell^\infty([e,E])$, with continuous limits.

For $A=(M,\infty)$, \eqref{eq:nested-maximal} and
$\mu_n(M,\infty)\to M^{-\alpha}$ give
\begin{equation}\label{eq:pa-uniform-tail-maximal}
 \limsup_n\E\sup_{x\ge M}|\mathbb G_n(x)|\le CM^{-\alpha/2}.
\end{equation}
For every finite set of thresholds in $[M,\infty)$, the portmanteau
theorem transfers \eqref{eq:pa-uniform-tail-maximal} to the expected
maximum of $|\mathbb G|$. Increasing these sets to a countable dense
subset, continuity and monotone convergence give
$\E\sup_{x\ge M}|\mathbb G(x)|\le CM^{-\alpha/2}$, so $\mathbb G(x)\to0$
as $x\to\infty$. By Markov's inequality, both suprema over $[M,\infty)$
are small in probability for large $M$, uniformly in large $n$, which
extends the convergence to $\ell^\infty([\epsilon,\infty))$.

For $M>\epsilon$, the map $y\mapsto\int_\cdot^My(z)\dd z/z$ is Lipschitz
on $\ell^\infty([\epsilon,M])$, so the continuous mapping theorem
handles the integrals over $[x,M]$. By Fubini's theorem,
$\int_M^\infty\mathbb G_n(z)\dd z/z=S_n(\log_+(\cdot/M))$. The variance
bound \eqref{eq:pa-variance-domination}, in which $K_n$ cancels, and the
bound $\mu_n(z,\infty)\le C\{z^{-\alpha}+K_n^{1/\alpha-1}z^{-1}\}$ from
\arxiveqref{M-eq:pa-stationary} and \eqref{eq:pa-profile-error} yield
\begin{equation}\label{eq:pa-log-tail-finite}
 \E\left|\int_M^\infty\mathbb G_n(z)\frac{\dd z}{z}\right|^2
 \le C\{M^{-\alpha}+K_n^{1/\alpha-1}M^{-1}\},
\end{equation}
since $\E\nu_n\log_+^2(\cdot/M)=\int_M^\infty2\log(z/M)\mu_n(z,\infty)\dd z/z$.
Since $\alpha>1$ and $K_n\to\infty$, letting $n\to\infty$ and then
$M\to\infty$ makes this bound vanish. For $\mathbb G$, Minkowski's
integral inequality and $\Var\mathbb G(z)\le z^{-\alpha}$, which is
\eqref{eq:pa-G-bound} as $y\to\infty$, give
$\E|\int_M^\infty\mathbb G(z)\dd z/z|^2\le4\alpha^{-2}M^{-\alpha}$.
Splitting the integrals at $M$ therefore proves the asserted joint
uniform convergence of the logarithmic integrals.
\end{proof}

In the rest of this section take $K_n=a_n$ as in
\arxiveqref{M-eq:balance}.

\subsection{The criterion on a compact range of ranks}\label{sec:general-lattice}
The proof of Theorem~\ref*{M-thm:vanishing} first restricts the
search to ranks $k$ with $k/a_n$ in a compact interval
$I\subset(0,\infty)$. On these ranks, the minimum of $\mathcal D_k$ requires
convergence of $\sqrt{a_n}\mathcal D_k$ to $\mathcal D_{\rm lat}$ in
\arxiveqref{M-eq:pa-D} uniformly in $k$, with the terms from integer degrees
retained. This convergence uses two properties of PA degrees, the
expansions of the stationary tail sums in Sec.~\ref{sec:pa-lattice}
and the functional limit in Theorem~\ref{thm:pa-functional}.
Theorem~\ref{thm:general-lattice} takes their analogues as hypotheses,
for triangular arrays of positive integers with arbitrary dependence
within rows. Corollary~\ref{cor:pa-conditional} verifies these hypotheses
for PA degrees, and Corollary~\ref{cor:iid-lattice} verifies them for iid
integer-valued observations used in the matched comparison of
Section~\ref{sec:variance-decisions}.

Let $V_{n,1},\ldots,V_{n,n}$ be positive integer-valued random variables,
with arbitrary dependence within each row. In this subsection, let $(\pi_j)_{j\ge1}$ be
a probability mass function satisfying, for some $\alpha,c,\eta>0$
and $Q\in\R$,
\begin{equation}\label{eq:lattice-reference-tail}
 \sum_{j>d}\pi_j=c^\alpha d^{-\alpha}
        \{1-\alpha Q/d+O(d^{-1-\eta})\},\qquad d\in\N, d\to\infty.
\end{equation}
Use the normalizations $a_n=n^{2/(\alpha+2)}$ and $b_n=c\sqrt{a_n}$ of
\arxiveqref{M-eq:balance}, now with the $\alpha$ and $c$ of
\eqref{eq:lattice-reference-tail}. They satisfy
$n\sum_{j>b_nx}\pi_j\sim a_nx^{-\alpha}$ for every $x>0$, so the limits
in Theorem~\ref{thm:general-lattice} use two coordinates: the time
coordinate $t$ of a rank
$k=\lfloor a_nt\rfloor$ and the threshold coordinate $x=t^{-1/\alpha}$ of
an observation divided by $b_n$. Put
\begin{equation}\label{eq:lattice-mean-error}
 \nu_n=a_n^{-1}\sum_{v=1}^n\delta_{V_{n,v}/b_n},\qquad
 e_j(n)=\E\sum_v\mathbf1\{V_{n,v}=j\}-n\pi_j,
 \qquad B=\frac{Q+\alpha/2}{\alpha+1}.
\end{equation}
Here $e_j(n)$ is the difference between the expected number of
observations equal to $j$ and $n\pi_j$. The constant $B$ is the
coefficient of $1/d$ in the ratio
$\sum_{j>d}\pi_j\log(j/d)/\sum_{j>d}\pi_j$, which is $H_{n,k}$ at
threshold $d$ with the proportions of observations replaced by $\pi_j$:
dividing \eqref{eq:lattice-log-expansion} in
Sec.~\ref{sec:general-lattice} by \eqref{eq:lattice-reference-tail}
gives $1/\alpha+B/d+O(d^{-1-\min(\eta,1)})$ for this ratio.

\begin{assumption}\label{ass:lattice-inputs}
The differences $e_j(n)$ satisfy
\begin{equation}\label{eq:lattice-mean-condition}
 \frac1{\sqrt{a_n}}\sum_{j\ge1}
       \{1+\log(1+j/b_n)\}|e_j(n)|\longrightarrow0.
\end{equation}
For every $\epsilon>0$, the centered process
$\mathbb G_n(x)=\sqrt{a_n}\{\nu_n(x,\infty)-\E\nu_n(x,\infty)\}$
converges in distribution in $\ell^\infty([\epsilon,\infty))$ to a
centered Gaussian process $\mathbb G$ with continuous paths on
$(0,\infty)$ and $\mathbb G(x)\to0$ as $x\to\infty$. For every $x>0$ the
integral $\mathbb L(x)=\int_x^\infty\mathbb G(z)\dd z/z$ converges almost
surely, and, for every $u>0$,
\begin{equation}\label{eq:lattice-log-condition}
 \lim_{M\to\infty}\limsup_n
 \Prob\left\{\left|\int_M^\infty\mathbb G_n(z)\frac{\dd z}{z}\right|>u\right\}=0.
\end{equation}
\end{assumption}
The conditions in Assumption~\ref{ass:lattice-inputs} are cumulative, and
every conclusion of Theorem~\ref{thm:general-lattice} uses all of them:
\eqref{eq:lattice-mean-condition} controls
$\E\sum_v\mathbf1\{V_{n,v}>d\}$ and $\E\sum_v\log_+(V_{n,v}/d)$, the
convergence of $\mathbb G_n$ controls the empirical survival function in
the KS supremum, and the convergence of $\mathbb L$ together with
\eqref{eq:lattice-log-condition} controls the Hill statistic, which enters
$\mathcal D_k$ through $\wh\alpha_{n,k}$. Almost sure convergence of $\mathbb L(1)$
gives $\int_M^\infty\mathbb G(z)\dd z/z\to0$ almost surely as
$M\to\infty$, the analogue of \eqref{eq:lattice-log-condition} for the
limit process.

Define $N_{n,>d}$, $N_{n,=d}$, and $\theta_{n,k}$ as in
Section~\ref{sec:pa-interpretation}, with $V_{n,v}$ in place of $D_v(n)$,
and $H_{n,k}$ and $\mathcal D_k$ by \eqref{eq:statistics}--\eqref{eq:Dk} with
$V_v=V_{n,v}$. Use $\mathcal A$, $\mathcal R$, and
$\mathcal D_{\rm lat}$ from
\arxiveqref{M-eq:pa-A} and \arxiveqref{M-eq:pa-Z}--\arxiveqref{M-eq:pa-D}, with the
constants and Gaussian process specified in this subsection. Fix
$0<t_-<t_+<\infty$, put $I=[t_-,t_+]$, and let $\wh k_I$ minimize $\mathcal D_k$ over
$\lceil a_nt_-\rceil\le k\le\lfloor a_nt_+\rfloor$.
Theorem~\ref{thm:general-lattice} gives the limits of the minimum
$\sqrt{a_n}\mathcal D_{\wh k_I}$, of the selected rank $\wh k_I$ with its tie
fraction, Hill statistic, and fitted exponent, and of the corrected
exponent along $I$.

\begin{theorem}\label{thm:general-lattice}
Under \eqref{eq:lattice-reference-tail} and
Assumption~\ref{ass:lattice-inputs},
\begin{equation}\label{eq:lattice-min-limit}
 \sqrt{a_n}\mathcal D_{\wh k_I}\wto
       \min_{(t,\theta)\in I\times[0,1]}\mathcal D_{\rm lat}(t,\theta)
       \quad\hbox{in }\R.
\end{equation}
Every subsequential limit in distribution of
\[
 \bigl(\wh k_I/a_n,\theta_{n,\wh k_I},
 \sqrt{a_n}(H_{n,\wh k_I}-1/\alpha)\bigr)
\]
has the form
$(T_I,\Theta_I,\mathcal A(T_I,\Theta_I))$, where
$(T_I,\Theta_I)$ minimizes $\mathcal D_{\rm lat}$ on $I\times[0,1]$
almost surely. If $\mathcal D_{\rm lat}$ has an almost surely unique
minimizer $(T_I,\Theta_I)$ on $I\times[0,1]$, then
\begin{equation}\label{eq:lattice-selected-law}
 \left(\frac{\wh k_I}{a_n},\theta_{n,\wh k_I},
 \sqrt{a_n}(\wh\alpha_{n,\wh k_I}-\alpha),\sqrt{a_n}\mathcal D_{\wh k_I}\right)
 \wto\left(T_I,\Theta_I,-\alpha^2\mathcal A(T_I,\Theta_I),
                    \mathcal D_{\rm lat}(T_I,\Theta_I)\right)
 \quad\hbox{in }\R^4.
\end{equation}

Suppose also, regardless of whether $\mathcal D_{\rm lat}$ has a unique
minimizer, that statistics $\wh B_{n,k}$ satisfy
$\sup_{\lfloor a_nt_-\rfloor\le k\le\lceil a_nt_+\rceil}
|\wh B_{n,k}-B|\xrightarrow{p}0$ in $\R$.
Define $\wt H_{n,k}=H_{n,k}-(\wh B_{n,k}-\theta_{n,k})/V_{(k+1)}$
and $\wt\alpha_{n,k}=\wt H_{n,k}^{-1}$, with $\wt\alpha_{n,k}=0$ when
$\wt H_{n,k}\le0$ or $H_{n,k}=0$, as in
\arxiveqref{M-eq:corrected-estimator}. Then
\begin{equation}\label{eq:lattice-corrected-process}
 \left\{\sqrt{\lfloor a_nt\rfloor}
       (\wt\alpha_{n,\lfloor a_nt\rfloor}-\alpha):t\in I\right\}
 \wto
 \left\{-\frac{\alpha^2}{\sqrt t}
       [\mathbb L(t^{-1/\alpha})-\mathbb G(t^{-1/\alpha})/\alpha]:t\in I\right\}
 \quad\hbox{in }\ell^\infty(I).
\end{equation}
\end{theorem}
\begin{proof}
We first expand the tail sums of $(\pi_j)$, then expand the rank $k$ and
the Hill sum uniformly over the candidate ranks, and finally pass to the
limit in the KS supremum and in the minimum over ranks. Put
$\eta_0=\min(\eta,1)$. Since $\pi_d=\sum_{j>d-1}\pi_j-\sum_{j>d}\pi_j$,
\eqref{eq:lattice-reference-tail} gives
$\pi_d=\alpha c^\alpha d^{-\alpha-1}\{1+O(d^{-\eta_0})\}$.
For real $u\to\infty$, replacing $u$ by $\lfloor u\rfloor$ gives
\begin{equation}\label{eq:lattice-real-expansion}
 \sum_{j>u}\pi_j=c^\alpha u^{-\alpha}
 \left\{1+\frac{\alpha(u-\lfloor u\rfloor-Q)}u+O(u^{-1-\eta_0})\right\}.
\end{equation}
For integer $d$, Fubini's theorem gives
$\sum_{j>d}\pi_j\log(j/d)=\int_d^\infty\sum_{j>u}\pi_j\dd u/u$, and
integrating \eqref{eq:lattice-real-expansion} yields
\begin{equation}\label{eq:lattice-log-expansion}
 \sum_{j>d}\pi_j\log(j/d)
 =c^\alpha d^{-\alpha}
 \left\{\frac1\alpha+
   \frac{\alpha(1/2-Q)}{(\alpha+1)d}+O(d^{-1-\eta_0})\right\}.
\end{equation}
The factor $1/2$ is the mean of $u-\lfloor u\rfloor$ over a unit
interval: the function $u-\lfloor u\rfloor-1/2$ has mean zero on each unit
interval and a bounded primitive, so integration by parts gives
$\int_d^\infty(u-\lfloor u\rfloor-1/2)u^{-\alpha-2}\dd u=O(d^{-\alpha-2})$.
Fix $0<e<E<\infty$. All remainders in
\eqref{eq:lattice-real-expansion}--\eqref{eq:lattice-log-expansion}
are uniform for $u,d\ge eb_n$.

For $d\ge eb_n$,
$\log_+(j/d)\le C_e\{1+\log(1+j/b_n)\}$, so
\eqref{eq:lattice-mean-condition} gives
\[
 \sup_{d\ge eb_n}\left\{\Bigl|\E N_{n,>d}-n\sum_{j>d}\pi_j\Bigr|
 +\Bigl|\E\sum_v\log_+\frac{V_{n,v}}d-n\sum_{j>d}\pi_j\log\frac jd\Bigr|\right\}
 =o(\sqrt{a_n}).
\]
The centered count at one integer $d$ is
$\{N_{n,=d}-\E N_{n,=d}\}/\sqrt{a_n}
=\mathbb G_n((d-1)/b_n)-\mathbb G_n(d/b_n)$, which tends to zero
uniformly over $eb_n\le d\le Eb_n$ by asymptotic equicontinuity of
$\mathbb G_n$. Since $\E N_{n,=d}=n\pi_d+e_d(n)$, the expansion of
$\pi_d$ and \eqref{eq:lattice-mean-condition} give
\begin{equation}\label{eq:lattice-atoms}
 \sup_{eb_n\le d\le Eb_n}
 \left|\frac{N_{n,=d}}{\sqrt{a_n}}-
               \frac{\alpha}{c}(d/b_n)^{-\alpha-1}\right|
 \xrightarrow{p}0\quad\hbox{in }\R.
\end{equation}
Hence, with probability tending to one, $N_{n,=d}\ge1$ for every integer
$d\in[eb_n,Eb_n]$, and the spacing $1/N_{n,=d}$ of the values of
$\theta_{n,k}$ at threshold $d$ tends to zero uniformly in probability.

Put $x=d/b_n$ and $t=x^{-\alpha}$, and write
$k=N_{n,>d}+\theta N_{n,=d}$ with $\theta\in[0,1]$, so that
$\theta=\theta_{n,k}$ when $k$ is a candidate rank with $V_{(k+1)}=d$.
Fubini's theorem gives
$\sum_v\log_+(V_{n,v}/d)=a_n\int_x^\infty\nu_n(z,\infty)\dd z/z$.
Uniformly over integers $d$ with $x\in[e,E]$ and over $\theta\in[0,1]$,
\eqref{eq:lattice-real-expansion}, \eqref{eq:lattice-log-expansion}, and
the bound on the expectations give
\begin{align}
 \frac{k-a_nt}{\sqrt{a_n}}
 &=\mathbb G_n(x)+\frac{\alpha}{c}(\theta-Q)x^{-\alpha-1}+R_{n,1}(x,\theta),\notag\\
 \frac{\sum_v\log_+(V_{n,v}/d)-a_nt/\alpha}{\sqrt{a_n}}
 &=\int_x^\infty\mathbb G_n(z)\frac{\dd z}{z}
     +\frac{\alpha(1/2-Q)}{c(\alpha+1)}x^{-\alpha-1}+R_{n,2}(x).
 \label{eq:lattice-uniform-expansions}
\end{align}
Here $\sup_{x\in[e,E],\theta\in[0,1]}|R_{n,1}(x,\theta)|\xrightarrow{p}0$
and $\sup_{x\in[e,E]}|R_{n,2}(x)|\xrightarrow{p}0$ in $\R$. The
deterministic parts of $R_{n,1}$ and $R_{n,2}$ are bounded by
$C\sqrt{a_n}b_n^{-1-\eta_0}$ plus $C_e$ times the left side of
\eqref{eq:lattice-mean-condition}, and the random part of $R_{n,1}$ is
$\theta\{N_{n,=d}/\sqrt{a_n}-(\alpha/c)x^{-\alpha-1}\}$, which tends to
zero uniformly by \eqref{eq:lattice-atoms}. On $[e,M]$, integration
against $\dd z/z$ is continuous in the uniform norm, so the convergence of
$\mathbb G_n$, \eqref{eq:lattice-log-condition}, and the almost sure
convergence of $\mathbb L$ give convergence in
distribution of $x\mapsto\int_x^\infty\mathbb G_n(z)\dd z/z$ to
$\mathbb L$ in $\ell^\infty([e,E])$, jointly with $\mathbb G_n$.
From now on take $e<t_+^{-1/\alpha}$ and $E>t_-^{-1/\alpha}$. Since
$\sup_{x\in[e,E]}|\nu_n(x,\infty)-x^{-\alpha}|\to0$ in probability, by
\eqref{eq:lattice-real-expansion}, \eqref{eq:lattice-mean-condition}, and
the convergence of $\mathbb G_n$, every candidate with $k/a_n\in I$ has
the threshold $V_{(k+1)}$ in $[eb_n,Eb_n]$ with probability tending to
one.

The first line of \eqref{eq:lattice-uniform-expansions} gives
$k/a_n-t\to0$ uniformly in probability. Subtract $1/\alpha$ times that
line from the second and divide by $k/a_n$. Since
$Q+\alpha(1/2-Q)/(\alpha+1)=B$ and
$H_{n,k}=k^{-1}\sum_v\log_+(V_{n,v}/d)$ at every candidate rank $k$ with
$V_{(k+1)}=d$, this gives, uniformly over the candidate ranks with
$d/b_n\in[e,E]$ and $\theta=\theta_{n,k}$,
\begin{equation}\label{eq:lattice-hill-expansion}
 \sqrt{a_n}\Bigl(H_{n,k}-\frac1\alpha\Bigr)
 =\frac1t\left\{\int_x^\infty\mathbb G_n(z)\frac{\dd z}{z}
       -\frac{\mathbb G_n(x)}\alpha\right\}
 +\frac{B-\theta}{cx}
\end{equation}
plus a remainder whose supremum tends to zero in probability. The right
side is $\mathcal A(t,\theta)$ in \arxiveqref{M-eq:pa-A} with
$\mathbb G_n$ and $\int_x^\infty\mathbb G_n(z)\dd z/z$ in place of
$\mathbb G$ and $\mathbb L$. 

The Taylor expansion of the fitted survival function is uniform on the
unbounded range $y\ge1$. For every $\vartheta\in[\alpha/2,3\alpha/2]$,
\begin{equation}\label{eq:lattice-pareto-remainder}
 \left|y^{-\vartheta}-y^{-\alpha}
             +(\vartheta-\alpha)y^{-\alpha}\log y\right|
 \le \tfrac12(\vartheta-\alpha)^2 y^{-\alpha/2}(\log y)^2,
 \qquad y\ge1.
\end{equation}
This is Taylor's formula in the exponent. The dominating function on
the right is bounded and tends to zero as $y\to\infty$.
Equation~\eqref{eq:lattice-hill-expansion} and inversion at
$1/\alpha$ imply that, for every $\varepsilon>0$, some finite $C_\varepsilon$ satisfies
$\sqrt{a_n}|\wh\alpha_{n,k}-\alpha|\le C_\varepsilon$ simultaneously over the
candidate ranks with probability at least $1-\varepsilon$ for all large $n$.
On this event, \eqref{eq:lattice-pareto-remainder} with
$\vartheta=\wh\alpha_{n,k}$ has supremum at most
$C C_\varepsilon^2/a_n$.
Multiplication by $\sqrt{a_n}$ therefore makes this remainder vanish.

For $y\ge1$, the empirical survival
function in \eqref{eq:Dk} equals $N_{n,>dy}/k$, where
$N_{n,>u}=\sum_v\mathbf1\{V_{n,v}>u\}$ for real $u$. Expanding
$N_{n,>dy}/k$ by \eqref{eq:lattice-real-expansion} and the first line of
\eqref{eq:lattice-uniform-expansions}, and the fitted survival
$y^{-\wh\alpha_{n,k}}$ by \eqref{eq:lattice-hill-expansion}, shows that
$\sqrt{a_n}(N_{n,>dy}/k-y^{-\wh\alpha_{n,k}})$ is
$\mathcal R(t,\theta;y,dy-\lfloor dy\rfloor)$, with $\mathbb G_n$ and
$\int_x^\infty\mathbb G_n(z)\dd z/z$ in place of $\mathbb G$ and
$\mathbb L$, up to a remainder whose supremum tends to zero in probability.

For $y\in[1,M]$, the fractional part $dy-\lfloor dy\rfloor$ runs through
$[0,1)$ on every complete interval between consecutive multiples of
$1/d$. The deterministic coefficients of $\mathcal R$ are uniformly
continuous on compact sets. The only term that jumps as $y$ varies is
$\mathbb G_n(xy)$, and asymptotic equicontinuity gives
\begin{equation}\label{eq:lattice-local-modulus}
 \sup_{\substack{u,v\in[e,EM+1]\\|u-v|\le b_n^{-1}}}
 |\mathbb G_n(u)-\mathbb G_n(v)|\xrightarrow{p}0
 \quad\hbox{in }\R.
\end{equation}
Changing $y$ by at most $1/d$ changes $xy$ by at most $1/b_n$.
Thus \eqref{eq:lattice-local-modulus} permits choosing the fractional
part near either endpoint without changing the other terms in the
limit. Consequently, replacing the fractional part by a free
$\psi\in[0,1]$ changes the supremum over $y\in[1,M]$ by a term tending
to zero in probability, uniformly over the candidate ranks.
Since $\mathcal R$ is affine in $\psi$, the maximum of $|\mathcal R|$ over
$\psi\in[0,1]$ is attained at $\psi=0$ or $\psi=1$, which gives the
maximum in \arxiveqref{M-eq:pa-D}.
For $y>M$, the terms of $\mathcal R$ without $\mathbb G_n$ or its
logarithmic integral are bounded by
$Cy^{-\alpha}(1+\log y)$, and the other terms tend to zero
uniformly in probability as $n\to\infty$ and then $M\to\infty$, by the
convergence of $\mathbb G_n$ to a limit with $\mathbb G(z)\to0$ as
$z\to\infty$ and by \eqref{eq:lattice-log-condition}. Hence, uniformly
over the candidate ranks $k$, $\sqrt{a_n}\mathcal D_k$ differs from
$\mathcal D_{\rm lat}(t,\theta)$, computed with $\mathbb G_n$ and
$\int_x^\infty\mathbb G_n(z)\dd z/z$ in place of $\mathbb G$ and
$\mathbb L$ and evaluated
at the attainable pair
$(t,\theta)=((V_{(k+1)}/b_n)^{-\alpha},\theta_{n,k})$, by a remainder that
tends to zero in probability.
Work on a representation on which the convergence in
Assumption~\ref{ass:lattice-inputs} and of the logarithmic integrals holds
almost surely. On it, $\mathbb G_n\to\mathbb G$ uniformly on $[e,\infty)$
and $\int_x^\infty\mathbb G_n(z)\dd z/z\to\mathbb L(x)$ uniformly for
$x\in[e,E]$. The function $\mathcal R$ in \arxiveqref{M-eq:pa-R} is linear in
$\mathbb G(xy)$, $\mathbb G(x)$, and $\mathbb L(x)$, with coefficients
bounded uniformly for $x\in[e,E]$ and $y\ge1$, so on this representation
\begin{equation}\label{eq:lattice-uniform-criterion}
 \max_{\lceil a_nt_-\rceil\le k\le\lfloor a_nt_+\rfloor}
 \Bigl|\sqrt{a_n}\mathcal D_k
   -\mathcal D_{\rm lat}\bigl((V_{(k+1)}/b_n)^{-\alpha},\theta_{n,k}\bigr)\Bigr|
 \xrightarrow{p}0\quad\hbox{in }\R.
\end{equation}

Adjacent integer thresholds in $[eb_n,Eb_n]$ give values of
$t=(d/b_n)^{-\alpha}$ that differ by at most $C/b_n$, and the values of
$\theta_{n,k}$ at each threshold have spacing tending to zero. The rank
restriction removes values of $\theta_{n,k}$ only at the two thresholds
whose candidate ranks $N_{n,>d}\le k<N_{n,>d}+N_{n,=d}$ include
$\lceil a_nt_-\rceil$ or $\lfloor a_nt_+\rfloor$. Hence the largest
distance from a point of $I\times[0,1]$ to the nearest attainable pair
tends to zero in probability. The criterion $\mathcal D_{\rm lat}$ is
continuous on the compact set $I\times[0,1]$, because $\mathbb G$ and
$\mathbb L$ are continuous and $\mathbb G(z)\to0$ as $z\to\infty$. A
finite mesh of $I\times[0,1]$ therefore approximates its minimum to within
its modulus of continuity, and the density of the attainable pairs and
\eqref{eq:lattice-uniform-criterion} give \eqref{eq:lattice-min-limit}.
For every $\varepsilon>0$, compactness and continuity give a strictly
positive gap between this minimum and the infimum of $\mathcal D_{\rm lat}$
over the pairs at distance at least $\varepsilon$ from its set of
minimizers, whenever such pairs exist.
Equation~\eqref{eq:lattice-uniform-criterion}, together with
$k/a_n-t\to0$, excludes these pairs as limits of
$(\wh k_I/a_n,\theta_{n,\wh k_I})$.
Every subsequential limit of $(\wh k_I/a_n,\theta_{n,\wh k_I})$ therefore
minimizes $\mathcal D_{\rm lat}$, and \eqref{eq:lattice-hill-expansion}
gives the third coordinate. If the minimizer is unique, every
subsequential limit equals it. On the representation,
$\sqrt{a_n}\mathcal D_{\wh k_I}$ then converges in probability to the minimum
$\mathcal D_{\rm lat}(T_I,\Theta_I)$, jointly with the other coordinates.
Since $H_{n,\wh k_I}\to1/\alpha$ in probability, the identity
$\sqrt{a_n}(\wh\alpha_{n,k}-\alpha)=-\alpha\sqrt{a_n}(H_{n,k}-1/\alpha)/H_{n,k}$
at $k=\wh k_I$ proves \eqref{eq:lattice-selected-law}.

For $\wt H_{n,k}$, subtract
$\sqrt{a_n}(\wh B_{n,k}-\theta_{n,k})/d=(\wh B_{n,k}-\theta_{n,k})/(cx)$
from \eqref{eq:lattice-hill-expansion}. The term $(B-\theta)/(cx)$
becomes $(B-\wh B_{n,k})/(cx)$, which tends to zero uniformly in
probability by the hypothesis on $\wh B_{n,k}$. Since
$k/a_n-x^{-\alpha}\to0$ uniformly in probability, asymptotic
equicontinuity of $\mathbb G_n$ and of
$x\mapsto\int_x^\infty\mathbb G_n(z)\dd z/z$ permits replacing $x$ by
$(k/a_n)^{-1/\alpha}$, and $(\lfloor a_nt\rfloor/a_n)^{-1/\alpha}\to
t^{-1/\alpha}$ uniformly for $t\in I$. Finally,
$\wt H_{n,k}\to1/\alpha$ in probability uniformly over these ranks, so the
derivative $-\alpha^2$ of $h\mapsto1/h$ at $1/\alpha$ linearizes the
inversion uniformly. Multiplication by $\sqrt{k/a_n}$ gives
\eqref{eq:lattice-corrected-process}.
\end{proof}

\subsection{Verification for PA degrees}\label{sec:pa-verification}
The next corollary applies the compact-range result to PA degrees.
\begin{corollary}\label{cor:pa-conditional}
In the homogeneous PA model \arxiveqref{M-eq:pa-model} with fixed $m\ge1$ and
$\delta>-m$, the degrees $V_{n,v}=D_v(n)$ satisfy
\eqref{eq:lattice-reference-tail} and
Assumption~\ref{ass:lattice-inputs}, with $\pi_j$ in
\arxiveqref{M-eq:pa-stationary}, $\eta=1$, $Q$ in \eqref{eq:pa-bias-constants},
and $\mathbb G$ in \eqref{eq:pa-G}. The constant $B$ in
\eqref{eq:lattice-mean-error} is then $B$ in
\eqref{eq:pa-bias-constants}, and all conclusions of
Theorem~\ref{thm:general-lattice} hold.
\end{corollary}
\begin{proof}
The expansion \eqref{eq:pa-strict-expansion} is
\eqref{eq:lattice-reference-tail} with $\eta=1$ and $Q$ in
\eqref{eq:pa-bias-constants}, and substituting $Q=\delta+(\alpha+1)/2$
into $B=(Q+\alpha/2)/(\alpha+1)$ gives $B$ in
\eqref{eq:pa-bias-constants}. By \eqref{eq:pa-profile-error} and Jensen's
inequality, the left side of \eqref{eq:lattice-mean-condition} is at most
$C(1+\log a_n)/\sqrt{a_n}$. Theorem~\ref{thm:pa-functional} with
$K_n=a_n$ gives the convergence of $\mathbb G_n$ to $\mathbb G$ in
\eqref{eq:pa-G}, with $\mathbb G(x)\to0$ as $x\to\infty$, and the bound
$\E|\mathbb G(z)|\le z^{-\alpha/2}$ in Sec.~\ref{app:pa-limits} gives
the almost sure convergence of $\mathbb L$. Finally,
\eqref{eq:pa-log-tail-finite} with $K_n=a_n$ gives
\[
 \limsup_{n\to\infty}\E\left|
      \int_M^\infty\mathbb G_n(z)\frac{\dd z}{z}\right|^2
 \le CM^{-\alpha},
\]
because $K_n^{1/\alpha-1}\to0$ for $\alpha>1$, and Chebyshev's
inequality gives \eqref{eq:lattice-log-condition}.
\end{proof}

\subsection{Verification for independent integer observations}\label{sec:iid-verification}
For iid observations, the limit process is a time-changed Brownian motion.

\begin{corollary}\label{cor:iid-lattice}
Suppose that, for each $n$, $V_{n,1},\ldots,V_{n,n}$ are iid with
$\Prob(V_{n,1}=j)=\pi_j$ for $j\ge1$, where $(\pi_j)$ satisfies
\eqref{eq:lattice-reference-tail}. Then
Assumption~\ref{ass:lattice-inputs} holds with
$\mathbb G(x)=\mathbb W(x^{-\alpha})$, where $\mathbb W$ is standard
Brownian motion, and all conclusions of Theorem~\ref{thm:general-lattice}
hold.
\end{corollary}
\begin{proof}
For iid observations, $e_j(n)=0$, so \eqref{eq:lattice-mean-condition}
holds. For a finite linear combination of the indicators
$\mathbf1\{V_{n,v}>b_nx\}$, the centered summands of $\mathbb G_n$ are
bounded by a constant times $a_n^{-1/2}\to0$, so they satisfy the
Lindeberg condition, and the Lindeberg--Feller central limit
theorem~\cite[Proposition~2.27]{vanDerVaart1998} gives
finite-dimensional convergence with limiting covariance
\[
 \frac{n}{a_n}\{\Prob(V_{n,1}>b_n(x\vee y))
       -\Prob(V_{n,1}>b_nx)\Prob(V_{n,1}>b_ny)\}
       \longrightarrow(x\vee y)^{-\alpha}.
\]
This is the covariance of $\mathbb W(x^{-\alpha})$, which has continuous
paths and tends to zero as $x\to\infty$. By Bernstein's
inequality~\cite[Chapter~2]{BoucheronLugosiMassart2013}, sums of
independent indicators satisfy \eqref{eq:pa-Bernstein}, and
\eqref{eq:lattice-reference-tail} gives \eqref{eq:pa-mean-uniform} for
$\mu_n=\E\nu_n$. From these two facts, through \eqref{eq:nested-maximal}
and \eqref{eq:nested-local-maximal}, the proof of
Theorem~\ref{thm:pa-functional} in Sec.~\ref{sec:pa-functional-result}
derives \eqref{eq:pa-equicontinuity} and
\eqref{eq:pa-uniform-tail-maximal}, which give tightness in
$\ell^\infty([\epsilon,\infty))$ and
\[
 \limsup_n\E\sup_{x\ge M}|\mathbb G_n(x)|\le CM^{-\alpha/2}.
\]
The variance
of $\int_M^\infty\mathbb G_n(z)\dd z/z$ is at most
$\E\nu_n[\log^2(z/M)\mathbf1\{z>M\}]\le CM^{-\alpha}$,
by integrating the tail bound $\E\nu_n(z,\infty)\le Cz^{-\alpha}$, and
Chebyshev's inequality gives
\eqref{eq:lattice-log-condition}. For the limit process,
$\E|\mathbb W(z^{-\alpha})|\le z^{-\alpha/2}$, so Fubini's theorem gives
the almost sure convergence of $\mathbb L(x)$ for every $x>0$. These facts
verify every part of Assumption~\ref{ass:lattice-inputs}.
\end{proof}

In both cases the terms of $\mathcal A$ and $\mathcal R$ that contain
neither $\mathbb G$ nor $\mathbb L$ depend on the observations only through
$\alpha$, $c$, $Q$, and the tie fraction $\theta$. Dependence among the
observations enters only through the covariance of $\mathbb G$, which is
that of $\mathbb W(x^{-\alpha})$ for iid observations and
\arxiveqref{M-eq:covariance} for PA. The correction removes the term
$(B-\theta)/(cx)$, the only term of $\mathcal A$ that depends on the tie
fraction, so the limit in \eqref{eq:lattice-corrected-process} is a
centered Gaussian process with continuous paths.
Theorem~\ref{thm:general-lattice} restricts the search to
$k/a_n\in I$. For PA, Theorem~\ref*{M-thm:vanishing} removes this
restriction through the localization \arxiveqref{M-eq:rank-bounds} and gives
the joint limit \arxiveqref{M-eq:unrestricted-joint} of the unrestricted
selector.

\subsection{Attainment of the global minimum}\label{sec:pa-attainment}
For the PA process $\mathbb G$ in \eqref{eq:pa-G}, a deterministic lower
bound on $\mathcal D_{\rm lat}$ comes from $y=1$,
where $\mathcal Z_{\rm lat}(t,1)=0$ and
$\mathcal R(t,\theta;1,\psi)=\alpha(\psi-\theta)t^{1/\alpha}/c$.
Hence
\begin{equation}\label{eq:pa-large-t}
 \mathcal D_{\rm lat}(t,\theta)\ge
 \frac{\alpha}{c}t^{1/\alpha}\max(\theta,1-\theta)
 \ge\frac{\alpha}{2c}t^{1/\alpha}.
\end{equation}
Thus $\inf_{0\le\theta\le1}\mathcal D_{\rm lat}(t,\theta)\to\infty$ as
$t\to\infty$, even at $\theta=B\in[0,1]$, where the term
$\alpha^2(B-\theta)\log y$ in \arxiveqref{M-eq:pa-R} vanishes. For $\mathcal D_k$,
Lemma~\ref{lem:pa-upper-localization} excludes the ranks $\wh k>La_n$
for large $L$.

The infimum over $\theta$ also diverges as $t\downarrow0$.
\begin{proposition}\label{prop:pa-attainment}
Almost surely,
\[
 \lim_{t\downarrow0}\inf_{0\le\theta\le1}\mathcal D_{\rm lat}(t,\theta)
 =\lim_{t\to\infty}\inf_{0\le\theta\le1}\mathcal D_{\rm lat}(t,\theta)
 =\infty.
\]
Consequently, almost surely, $\mathcal D_{\rm lat}$ attains its global
minimum on $(0,\infty)\times[0,1]$, and its set of minimizers is a
nonempty compact set.
\end{proposition}
\begin{proof}
The limit as $t\to\infty$ follows from \eqref{eq:pa-large-t}. For
$t\downarrow0$, we combine \arxiveqref{M-eq:pa-R} at three values of $y$ so as
to cancel the terms proportional to $y^{-\alpha}$ and to
$y^{-\alpha}\log y$, which come from threshold normalization and
exponent fitting.
Fix $1=r_1<r_2<r_3$ and choose coefficients
$e_1=1,e_2,e_3$ satisfying
\begin{equation}\label{eq:pa-contrast-coefficients}
 \sum_{j=1}^3e_jr_j^{-\alpha}=0,\qquad
 \sum_{j=1}^3e_jr_j^{-\alpha}\log r_j=0.
\end{equation}
The determinant of the two equations for $e_2,e_3$ is
$r_2^{-\alpha}r_3^{-\alpha}\log(r_3/r_2)>0$.
The function $\sum_j e_j f_{r_j}$ equals one on $(r_1,r_2)$, so it is
nonzero in $L^2(\nu)$, and Proposition~\ref{U-prop:positive-covariance} gives
\begin{equation}\label{eq:pa-contrast-variance}
 \sigma_e^2=V\left(\sum_j e_j f_{r_j},\sum_j e_j f_{r_j}\right)>0.
\end{equation}

For $b\le t\le2b$ in the time coordinate, put $x=t^{-1/\alpha}$ and
$z_j=b^{-1/\alpha}r_j$ in the threshold coordinate, so that
$y=z_j/x\ge1$. Multiply \arxiveqref{M-eq:pa-R} at $y=z_j/x$ and $\psi=0$
by $t$, and sum with coefficients $e_j$. Since
$\sum_je_jz_j^{-\alpha}=b\sum_je_jr_j^{-\alpha}$ and
$\sum_je_jz_j^{-\alpha}\log z_j
=b\sum_je_jr_j^{-\alpha}(\log r_j-\alpha^{-1}\log b)$,
\eqref{eq:pa-contrast-coefficients} cancels every term proportional to
$z_j^{-\alpha}$ or $z_j^{-\alpha}\log z_j$. As
$|\mathcal R(t,\theta;y,0)|\le\mathcal D_{\rm lat}(t,\theta)$ and
$t\le2b$,
\[
 \left|\sum_j e_j\left\{\mathbb G(z_j)
                  -\frac{\alpha Q}{c}z_j^{-\alpha-1}\right\}\right|
 \le 2b\,\mathcal D_{\rm lat}(t,\theta)\sum_j|e_j|
 \qquad(\theta\in[0,1]).
\]
The left side does not depend on $(t,\theta)$, and by
\eqref{eq:pa-homogeneity} $\sum_je_j\mathbb G(z_j)$ has variance
$b\sigma_e^2$. A normal density with this variance is at most
$(2\pi b\sigma_e^2)^{-1/2}$ whatever its mean, so
\begin{equation}\label{eq:pa-small-t}
 \Prob\left\{\inf_{\substack{b\le t\le2b\\0\le\theta\le1}}
             \mathcal D_{\rm lat}(t,\theta)\le z\right\}
 \le C_e z\sqrt b,\qquad
 C_e=\frac{4\sum_j|e_j|}{\sqrt{2\pi\sigma_e^2}},\qquad b,z>0.
\end{equation}
For each positive integer $z$, these bounds at $b=2^{-i}$ are summable
in $i$, so by the Borel--Cantelli lemma, almost surely,
$\mathcal D_{\rm lat}>z$ on $[2^{-i},2^{1-i}]\times[0,1]$ for all large
$i$. Intersecting over $z$ proves the limit as $t\downarrow0$. The
continuity of $\mathbb G$, with $\mathbb G(x)\to0$ as $x\to\infty$ by
Theorem~\ref{thm:pa-functional}, and the decay in $y$ of the
deterministic terms of \arxiveqref{M-eq:pa-R} make $\mathcal D_{\rm lat}$
continuous on compact subsets of $(0,\infty)\times[0,1]$. With the two
limits, $\mathcal D_{\rm lat}$ attains its minimum, and its set of
minimizers is closed and bounded, hence compact.
\end{proof}

\subsection{Added noise, a shifted threshold, and the discrete likelihood}\label{sec:integer-treatments}
Clauset et al.~\cite{ClausetEtAl2009} treat integer observations at a
threshold $d$ in two ways, both using the observations at least $d$:
the maximum likelihood fit of the probabilities proportional to $j^{-a-1}$
on the integers $j\ge d$, and its approximation by the reciprocal of the
mean of $\log\{x/(d-1/2)\}$ over these observations. Voitalov
et al.~\cite{VoitalovEtAl2019} instead add to each observation an
independent random variable, uniform on $[-1/2,1/2]$, before computing the
Hill statistic. The next lemma gives the population versions of the three
estimates, in which the proportions of observations are replaced by their
limits, for every law whose tail sums satisfy
\eqref{eq:lattice-reference-tail}. Compared with
\arxiveqref{M-eq:hill-tie-expansion}, all three remove the dependence on
the position within the tie and keep a term of order $d^{-1}$.
\begin{lemma}\label{lem:noise-expansion}
Let $X$ be a random variable on the nonnegative integers whose
probabilities $\pi_j=\Prob(X=j)$ have tail sums satisfying
\eqref{eq:lattice-reference-tail} for some $\alpha,c,\eta>0$ and
$Q\in\R$, put $B=(Q+\alpha/2)/(\alpha+1)$ and $\eta_0=\min(\eta,1)$, and
let $\check X$ be the sum of $X$ and an independent random variable
uniform on $[-1/2,1/2]$. For integers $d\ge1$ and real $a>0$, put
\[
 r_d=\E\{\log(X/d)\mid X\ge d\},\qquad
 g_d(a)=\frac{\sum_{j\ge d}j^{-a-1}\log(j/d)}{\sum_{j\ge d}j^{-a-1}}.
\]
For all large $d$, the equation $g_d(a)=r_d$ has exactly one solution
$a=a_d>0$. As the integer $d\to\infty$,
\begin{equation}\label{eq:noise-population-hill}
 \frac{\E[\log(\check X/y)\mathbf1\{\check X>y\}]}{\Prob(\check X>y)}
 =\frac1\alpha+\frac{B-1/2}d+O(d^{-1-\eta_0})
\end{equation}
uniformly over $y\in[d-1/2,d+1/2)$, and both
$r_d+\log\{d/(d-1/2)\}=\E[\log\{X/(d-1/2)\}\mid X\ge d]$ and $1/a_d$ equal
the right side of \eqref{eq:noise-population-hill}.
\end{lemma}

\begin{proof}
The proof first expands the ratio in \arxiveqref{M-eq:hill-tie-expansion}
for these probabilities and then treats the three quantities in turn.

We first show that, uniformly over $\theta\in[0,1]$,
\begin{equation}\label{eq:lattice-tie-expansion}
 \frac{\sum_{j>d}\pi_j\log(j/d)}{\sum_{j>d}\pi_j+\theta\pi_d}
 =\frac1\alpha+\frac{B-\theta}d+O(d^{-1-\eta_0}).
\end{equation}
The proof of Theorem~\ref{thm:general-lattice} derives from
\eqref{eq:lattice-reference-tail} alone the expansions
$\pi_d=\alpha c^\alpha d^{-\alpha-1}\{1+O(d^{-\eta_0})\}$ and
\eqref{eq:lattice-log-expansion}. Dividing \eqref{eq:lattice-log-expansion}
by \eqref{eq:lattice-reference-tail} gives
$\sum_{j>d}\pi_j\log(j/d)/\sum_{j>d}\pi_j=1/\alpha+B/d+O(d^{-1-\eta_0})$,
and the two expansions give
$\pi_d/\sum_{j>d}\pi_j=\alpha/d+O(d^{-1-\eta_0})$. Dividing the first ratio
by $1+\theta\pi_d/\sum_{j>d}\pi_j$ proves \eqref{eq:lattice-tie-expansion}.
In particular, $\pi_d\le Cd^{-1}\sum_{j>d}\pi_j$ for all large $d$.

For the noise, put $\theta=d+1/2-y$, which lies in $(0,1]$. Since
$|\check X-X|\le1/2$, the event $\{\check X>y\}$ contains
$\{X\ge d+1\}$, has probability zero on $\{X\le d-1\}$, and on $\{X=d\}$
equals $\{\check X-X>y-d\}$, whose conditional probability is $\theta$.
Hence
\begin{equation}\label{eq:noise-denominator}
 \Prob(\check X>y)=\sum_{j>d}\pi_j+\theta\pi_d.
\end{equation}
For $j\ge d+1$ and $|u|\le1/2$,
$\log\{(j+u)/y\}=\log(j/d)+\log(1+u/j)-\log(y/d)$. The inequality
$|\log(1+x)-x|\le x^2$ for $|x|\le1/2$, applied with $x=(\check X-X)/j$
and with $x=(y-d)/d$, together with $\E(\check X-X)=0$ and
$\E(\check X-X)^2=1/12$, gives $|\E\log\{1+(\check X-X)/j\}|\le1/(12d^2)$
and $|\log(y/d)-(y-d)/d|\le1/(4d^2)$. On $\{X=d,\ \check X>y\}$,
$0<\log(\check X/y)\le\log\{(d+1/2)/(d-1/2)\}\le2/d$. Therefore
\begin{equation}\label{eq:noise-numerator}
 \E[\log(\check X/y)\mathbf1\{\check X>y\}]
 =\sum_{j>d}\pi_j\log(j/d)-\frac{y-d}d\sum_{j>d}\pi_j+R_d(y),
\end{equation}
with $|R_d(y)|\le\{1/(3d^2)\}\sum_{j>d}\pi_j+(2/d)\pi_d
\le Cd^{-2}\sum_{j>d}\pi_j$ for every $y\in[d-1/2,d+1/2)$. The right side
of \eqref{eq:noise-denominator} is at least $\sum_{j>d}\pi_j$ and at most
$(1+Cd^{-1})\sum_{j>d}\pi_j$, and $|y-d|\le1/2$, so dividing
\eqref{eq:noise-numerator} by \eqref{eq:noise-denominator} and applying
\eqref{eq:lattice-tie-expansion} gives
\[
 \frac{\E[\log(\check X/y)\mathbf1\{\check X>y\}]}{\Prob(\check X>y)}
 =\frac1\alpha+\frac{B-\theta}d-\frac{y-d}d+O(d^{-1-\eta_0})
\]
uniformly in $y$. Since $\theta+(y-d)=1/2$, this is
\eqref{eq:noise-population-hill}.

For the shifted threshold, the observations equal to $d$ contribute zero
to $r_d$, so $r_d$ is the left side of \eqref{eq:lattice-tie-expansion} at
$\theta=1$, which equals $1/\alpha+(B-1)/d+O(d^{-1-\eta_0})$. Adding
$\log\{d/(d-1/2)\}=1/(2d)+O(d^{-2})$ gives the right side of
\eqref{eq:noise-population-hill}.

For the discrete likelihood, the Euler--Maclaurin
formula~\cite[Eq.~2.10.1]{DLMFEulerMaclaurin} with one correction term,
applied on $[d,\infty)$ to $x\mapsto x^{-a-1}$ and to
$x\mapsto x^{-a-1}\log(x/d)$, whose second derivatives are bounded by
constant multiples of $x^{-a-3}\{1+\log(x/d)\}$ with constants continuous
in $a$, gives
\[
 \sum_{j\ge d}j^{-a-1}=\frac{d^{-a}}a\left\{1+\frac a{2d}+O(d^{-2})\right\},
 \qquad
 \sum_{j\ge d}j^{-a-1}\log\frac jd=\frac{d^{-a}}{a^2}\{1+O(d^{-2})\},
\]
uniformly for $a$ in compact subsets of $(0,\infty)$. Hence
$g_d(a)=1/a-1/(2d)+O(d^{-2})$ uniformly over $a\in[\alpha/2,2\alpha]$.
The function $g_d(a)$ is the mean of $\log(j/d)$ under the probabilities
proportional to $j^{-a-1}$ on the integers $j\ge d$, and differentiation
under the sums shows that $g_d'(a)$ equals minus the variance of $\log j$
under these probabilities, so that $g_d$ is strictly decreasing. Since $r_d\to1/\alpha$, $g_d(\alpha/2)\to2/\alpha$, and
$g_d(2\alpha)\to1/(2\alpha)$, for all large $d$ the solution $a_d$ exists,
lies in $(\alpha/2,2\alpha)$, and is unique. Then
$1/a_d=g_d(a_d)+1/(2d)+O(d^{-2})=r_d+1/(2d)+O(d^{-2})$, which is the
expansion of the shifted threshold.
\end{proof}

For the discrete power law, with $\pi_j$ proportional to $j^{-a-1}$ for
$j\ge1$, the first expansion in the proof, less the term $j=d$, gives
\eqref{eq:lattice-reference-tail} with $\alpha=a$ and $Q=1/2$, so $B=1/2$
and the three quantities in Lemma~\ref{lem:noise-expansion} equal
$1/a+O(d^{-2})$. For PA degrees, \eqref{eq:pa-strict-expansion} gives
\eqref{eq:lattice-reference-tail} with $\eta=1$ and $Q$ in
\eqref{eq:pa-bias-constants}, and $B-1/2=(\delta+\alpha/2)/(\alpha+1)$,
which is positive for $\alpha\ge2$ and vanishes only when
$\delta=-\alpha/2$. At these exponents the term of order $d^{-1}$ makes
the reciprocals of the three quantities smaller than $\alpha$ for all
large $d$.

\section{Unrestricted threshold selection in preferential attachment}\label{sec:unrestricted-proofs}
On each compact range of ranks, Corollary~\ref{cor:pa-conditional} and
Theorem~\ref{thm:general-lattice} give the limits of $\sqrt{a_n}\mathcal D_k$ and
$\sqrt{a_n}(H_{n,k}-1/\alpha)$. The proof of
Theorem~\ref*{M-thm:vanishing} adds the localization
\arxiveqref{M-eq:rank-bounds} and the uniqueness of the global minimizer
$(T,\Theta)$. Lemma~\ref{lem:pa-upper-localization} excludes ranks large
relative to $a_n$, because the jumps of the empirical survival function
at integer degrees make $\mathcal D_k$ large there.
Proposition~\ref{prop:pa-attainment} shows that $\mathcal D_{\rm lat}$
attains its minimum, but attainment does not give
$\lim_{L\to\infty}\limsup_n\Prob\{\wh k<a_n/L\}=0$.
Lemma~\ref{U-lem:small-discrepancy} supplies a quantitative
bound instead. For $2\le h\le a_n$ with $2h<n$ and $0<u\le1$, it bounds
$\Prob\{\min_{h\le k\le2h}\sqrt h\,\mathcal D_k\le u\}$ by $Cu$ plus error terms
that are summable over the dyadic intervals $[h,2h]$ covering the ranks
smaller than $a_n/L$, and Lemma~\ref{U-lem:lower-localization} sums
these bounds. The proof of Lemma~\ref{U-lem:small-discrepancy} uses a
smooth function whose integral against the Pareto measure
$\vartheta z^{-\vartheta-1}\dd z$ vanishes, together with its derivative
in $\vartheta$, at $\vartheta=\alpha$, followed by the quantitative
martingale central limit theorem of
Proposition~\ref{U-prop:quantitative-normal}. These two results rest on
the centered martingale and comparison estimates of
Wang~\cite[Proposition~3.7, Lemmas~3.1--3.6, and Lemmas~A.2 and~A.5]{Wang2026}.
The proof of uniqueness in Proposition~\ref{U-prop:unique-global}
accounts separately for the term $\alpha y^{-\alpha-1}/(2cx)$ in
\eqref{eq:pa-D-midpoint}, which comes from the two sides of each integer
jump, and for the minimization over the tie fraction.
Proposition~\ref{U-prop:unique-global} gives uniqueness only for the
global minimizer, so Sec.~\ref{U-sec:joint-limit} completes the
proof by comparing $\wh k$ with ranks at rational pairs $(t,\theta)$.
\subsection{Exclusion of large ranks}\label{app:pa-upper-localization}
Both localization proofs compare $\mathcal D_k$ with its value at the rank
$\lfloor a_n\rfloor$, which satisfies
\begin{equation}\label{U-eq:benchmark}
 \lim_{M\to\infty}\limsup_n
 \Prob\{\sqrt{a_n}\mathcal D_{\lfloor a_n\rfloor}>M\}=0.
\end{equation}
Indeed, apply \eqref{eq:lattice-uniform-criterion}, which holds for PA
degrees by Corollary~\ref{cor:pa-conditional}, with $I=[1/2,2]$. On the
representation used there, $\sqrt{a_n}\mathcal D_{\lfloor a_n\rfloor}$ differs by a
term tending to zero in probability from
$\mathcal D_{\rm lat}(t,\theta_{n,\lfloor a_n\rfloor})$, where
$t=\{D_{(\lfloor a_n\rfloor+1)}(n)/b_n\}^{-\alpha}\to1$ in probability by
Proposition~\ref{prop:pa-fixed-rank}, and $\mathcal D_{\rm lat}$ is
continuous, hence almost surely bounded, on $I\times[0,1]$.

The next lemma proves the part of the localization
\arxiveqref{M-eq:rank-bounds} in Theorem~\ref*{M-thm:vanishing} that
excludes ranks $\wh k>La_n$.
\begin{lemma}\label{lem:pa-upper-localization}
Under the homogeneous PA model \arxiveqref{M-eq:pa-model} with fixed $m\ge1$
and $\delta>-m$, the rank $\wh k$ defined in \eqref{eq:Dk} satisfies
\begin{equation}\label{eq:pa-upper-localization}
 \lim_{L\to\infty}\limsup_n\Prob\{\wh k>La_n\}=0.
\end{equation}
\end{lemma}
\begin{proof}
We compare $\mathcal D_{\lfloor a_n\rfloor}$ with ranks
whose threshold is much smaller than $b_n$.
Write $N_{n,\ge d}=N_{n,>d}+N_{n,=d}$. For every integer $m\le d\le b_n$,
\arxiveqref{M-eq:pa-stationary} and \eqref{eq:pa-profile-error} give constants
$c_1,C_1>0$, depending only on $m$ and $\delta$, such that
\[
 \E N_{n,=d+1}\ge c_1 n(d+1)^{-\alpha-1},\qquad
 \E N_{n,\ge d}\le C_1 n d^{-\alpha}
\]
for all sufficiently large $n$. The lower bound is at least a constant
times $\sqrt{a_n}$, because $nb_n^{-\alpha-1}=c^{-\alpha-1}\sqrt{a_n}$.
Apply \eqref{eq:pa-Bernstein} with $K=1$ to these counts, with
deviations equal to fixed fractions, at most one half, of their means.
A union bound over at most $b_n$ integers gives, with probability at
least $1-Cb_n\exp(-\sqrt{a_n}/C)\to1$ and a constant $c_2>0$,
\begin{equation}\label{eq:pa-uniform-large-ranks}
 N_{n,=d+1}\ge c_2 n(d+1)^{-\alpha-1},\qquad
 N_{n,\ge d}\le 2C_1 n d^{-\alpha}
 \quad(m\le d\le b_n).
\end{equation}
If a valid candidate $k$ has $D_{(k+1)}(n)=d$ in this range, then
$k\le N_{n,\ge d}$. In \eqref{eq:Dk}, the empirical survival
$N_{n,>\lfloor dy\rfloor}/k$ jumps by $N_{n,=d+1}/k$ at $y=(d+1)/d$,
whereas $y^{-\wh\alpha_{n,k}}$ is continuous, so one side of the jump
differs from it by at least half that size.
Equation~\eqref{eq:pa-uniform-large-ranks} and $d/(d+1)\ge1/2$ yield
\begin{equation}\label{eq:pa-large-rank-distance}
 \mathcal D_k\ge\frac{N_{n,=d+1}}{2k}\ge\frac{c_3}{d+1},
 \qquad c_3=2^{-\alpha-2}c_2/C_1.
\end{equation}

Fix a large $L$ and set $d_{n,L}=\lceil b_n(2/L)^{1/\alpha}\rceil$, so that
$d_{n,L}/b_n$ tends to the threshold coordinate $(2/L)^{1/\alpha}$ of the
time $t=L/2$. By \eqref{eq:pa-strict-expansion} and
\eqref{eq:pa-profile-error}, $\E N_{n,>d_{n,L}}/a_n\to L/2$, and
\eqref{eq:pa-Bernstein} implies $\Prob(N_{n,>d_{n,L}}\ge La_n)\to0$.
Hence, for $L>2$ and with probability tending to one, every rank
$k\ge La_n$ has $m\le D_{(k+1)}(n)\le d_{n,L}\le b_n$. Since
$d_{n,L}+1\le2b_n(2/L)^{1/\alpha}$ for large $n$,
\eqref{eq:pa-large-rank-distance} gives
\[
 \min_{k\ge La_n}\sqrt{a_n}\mathcal D_k\ge c_4 L^{1/\alpha},
 \qquad c_4=2^{-1-1/\alpha}c_3/c,
\]
with probability tending to one. Invalid candidates have infinite
criterion and cannot improve the minimum. Given $\eta>0$, choose $M$ by
\eqref{U-eq:benchmark} with
$\limsup_n\Prob\{\sqrt{a_n}\mathcal D_{\lfloor a_n\rfloor}>M\}\le\eta$. If $L>2$
and $c_4L^{1/\alpha}>M$, then $\limsup_n\Prob\{\wh k>La_n\}\le\eta$,
which proves \eqref{eq:pa-upper-localization}.
\end{proof}

\subsection{A quantitative normal approximation}
Lemma~\ref{U-lem:small-discrepancy} in Sec.~\ref{U-sec:small-ranks}
bounds the probability that $\mathcal D_k$ is small at some rank $k$ in a dyadic
range $[h,2h]$, and these bounds must be summable over the ranges with
$h\le a_n/L$. It therefore needs a normal approximation in Kolmogorov
distance with an explicit rate in $h$ and $h/n$.
Retain the homogeneous model \arxiveqref{M-eq:pa-model}, its parameters, and
$V$ in \eqref{eq:pa-covariance}. For an integer $h\ge2$, take
$K_n=h$ in \eqref{eq:pa-measure} and write
\begin{equation}\label{U-eq:smooth-statistic}
 \nu_{n,h}=h^{-1}\sum_{v=1}^n
       \delta_{D_v(n)/\{c(n/h)^{1/\alpha}\}},\qquad
 S_{n,h}(f)=\sqrt h\{\nu_{n,h}f-\E\nu_{n,h}f\}.
\end{equation}
Let $\Phi$ denote the standard normal distribution function.

The next proposition bounds the Kolmogorov distance between
$S_{n,h}(f)$ and $N\{0,V(f,f)\}$ for each fixed smooth $f$.
\begin{proposition}\label{U-prop:quantitative-normal}
Fix $f\in C_c^\infty((0,\infty))$ with $V(f,f)>0$.
There are constants $C<\infty$ and $\kappa>0$, depending on
$m,\delta,f$, such that, for every $n$ and every integer $h$ with
$2\le h\le n/2$,
\begin{equation}\label{U-eq:quantitative-normal}
 \sup_{z\in\R}\left|
 \Prob\{S_{n,h}(f)\le z\}-\Phi\{z/\sqrt{V(f,f)}\}\right|
 \le C\{h^{-\kappa}+(h/n)^\kappa\}.
\end{equation}
\end{proposition}

\begin{proof}
The martingale argument of
\cite[Proposition~3.7 and Lemmas~3.3--3.6]{Wang2026} gives convergence
without a rate. The increments are already uniformly small, so the
remaining task is to bound the expected absolute error of the
predictable variance. We keep the dependence on the vertex and time
cutoffs explicit, then choose them as powers of $\min(h,n/h)$.
Constants denoted by $C$ depend only on $m,\delta,f$. Write
$b=c(n/h)^{1/\alpha}$ and choose $x_0>0$ such that $f=0$ on $(0,x_0)$.
Let $\mathcal B_e$ and $\mathcal B_{e-}$ be the sigma-fields of the graph
just after and before edge $e=(t,i)$, and let $J_e$ be the vertex
receiving $e$. By \arxiveqref{M-eq:pa-model}, $J_e=v$ has conditional
probability $\rho_{v,e}=(D_v(t,i-1)+\delta)/W_e$ given $\mathcal B_{e-}$,
with deterministic $W_e=m\alpha(t-1)-2m+i-1$. Hence
$\E\{f(D_v(n)/b)\mid\mathcal B_e\}$ depends only on the degree of $v$
after $e$. Let $\xi_{v,e}$ be the difference between these conditional
terminal expectations when the degree just after $e$ is
$D_v(t,i-1)+1$ and when it is $D_v(t,i-1)$.
The Doob martingale $\E\{S_{n,h}(f)\mid\mathcal B_e\}$ has increments
\begin{equation}\label{U-eq:increments}
 \Delta_e=h^{-1/2}\sum_{v<t}
       \{\mathbf1\{J_e=v\}-\rho_{v,e}\}\xi_{v,e},
       \qquad |\Delta_e|\le4\|f\|_\infty h^{-1/2}.
\end{equation}
Write
\[
 \bar\xi_e=\sum_{v<t}\rho_{v,e}\xi_{v,e},\qquad
 \Xi_v=\sum_e\rho_{v,e}(1-\rho_{v,e})\xi_{v,e}^2,
\]
where sums over edges include only those added after the birth of $v$.
Then $\E\Xi_v=\Var f(D_v(n)/b)$ by
\cite[Lemma~3.3, Eq.~(3.14)]{Wang2026}. The predictable variance is
\begin{equation}\label{U-eq:variance-decomposition}
 \begin{split}
 V_{n,h}&=\sum_e\E(\Delta_e^2\mid\mathcal B_{e-})
 =h^{-1}\sum_v\Xi_v
       +h^{-1}\sum_{e,v}\rho_{v,e}^2\xi_{v,e}^2-K_{n,h},\\
 K_{n,h}&=h^{-1}\sum_e\bar\xi_e^2.
 \end{split}
\end{equation}
The first term in \eqref{U-eq:variance-decomposition} contributes the
integral of marginal variances in $V(f,f)$. The middle term vanishes,
and $K_{n,h}$ contributes the term due to competition for edges.

With $X_s$ as in \eqref{eq:pa-covariance} and $3\le v\le n$, the
comparison inequality \cite[Lemma~A.5, Eq.~(A.9)]{Wang2026} and the product
estimate \cite[Lemma~A.2, Eq.~(A.4)]{Wang2026} imply
\begin{equation}\label{U-eq:marginal-rate}
 \left|\E\varphi(D_v(n)/b)-\E\varphi(X_{v/h})\right|
 \le C_\varphi\left\{(h/v)^\gamma
       [v^{-1}+(v/n)^{\gamma/2}]+(h/n)^\gamma\right\}
\end{equation}
for every fixed bounded Lipschitz $\varphi$.
To prove \eqref{U-eq:marginal-rate}, put
$p_{v,n}=\prod_{v<e\le n}(1+W_e^{-1})^{-1}$, a product over the edges
added after the birth of $v$.
For the pure-birth comparison in \cite[Lemma~A.1]{Wang2026}, the weight
multiplied by $p_{v,n}$ has mean $r$ and variance $r(1-p_{v,n})$.
Its $L^2$ martingale limit has law $G_r$ and variance $r$. Martingale
orthogonality therefore gives mean-square approximation error
$rp_{v,n}$, and Cauchy--Schwarz gives mean absolute error at most
$C\sqrt{p_{v,n}}$.
There are three contributions after division by $b$:
\begin{enumerate}
\item Lemma~A.5 of~\cite{Wang2026} bounds the difference between the
PA chain and the comparison chain by $C(h/v)^\gamma/v$.
\item Since $p_{v,n}^{-1}=(n/v)^\gamma\{1+O(v^{-1})\}$ by
\cite[Lemma~A.2]{Wang2026}, the Gamma approximation contributes at most
$C(h/v)^\gamma\{v^{-1}+(v/n)^{\gamma/2}\}$.
\item Subtracting the offset $\delta$ from the weight contributes at
most $C/b$.
\end{enumerate}
These bounds prove \eqref{U-eq:marginal-rate}. Starting the chain at
$m+1$ replaces $r$ by $r+1$ in its moments and replaces $X_{v/h}$ by
$c^{-1}(v/h)^{-\gamma}G_{r+1}$ in that estimate.

To bound the variance error, retain vertex labels in
$[\varepsilon h,Lh]$ and edge times in $[h/L,Lh]$. The longer time
range ending at $Th$ is used only to compare two degree paths. Choose
$L>1$, $0<\varepsilon<L^{-1}$, and $T>L$, with
$\varepsilon h\ge3$ and $Th\le n$.
For vertices $v_1\ne v_2$ with labels at least $\varepsilon h$, the
coupling inequality \cite[Lemma~A.5, Eq.~(A.10)]{Wang2026}
couples their degree paths up to time $Th$ with independent paths and
bounds the probability that the coupling fails by
\begin{equation}\label{U-eq:pair-coupling-rate}
 C\sum_{v\in\{v_1,v_2\}}\sum_{v<t\le Th}
                    v^{-2\gamma}t^{2\gamma-2}
 \le C\varepsilon^{-2\gamma-2}T^{2\gamma+1}/h.
\end{equation}
The second inequality uses at most $Th$ summands and
$t\ge\varepsilon h$, and \eqref{U-eq:pair-coupling-rate} also holds with
$L$ in place of $T$.

The proof of \cite[Lemma~3.4]{Wang2026} gives
\begin{equation}\label{U-eq:diagonal-remainder-rate}
 h^{-1}\sum_{e,v}\E(\rho_{v,e}^2\xi_{v,e}^2)\le C\{\varepsilon+h^{-1}
                                +\varepsilon^{-2\gamma}h^{-1}\}.
\end{equation}
In $h^{-1}\sum_v\Xi_v$, labels $v<\varepsilon h$
contribute at most $C(\varepsilon+h^{-1})$ in expectation.
For an integer $q$ with $q\gamma\ge2$, moments of the degrees and the
fact that $f$ vanishes on $(0,x_0)$ bound the contribution of labels
$v>Lh$ by $CL^{1-q\gamma}\le C/L$.
Boundedness of $f$ and the Gamma moments give these bounds for the
omitted parts of $\int_0^\infty\Var f(X_s)\,ds$ as well.
For $\varepsilon h\le v\le Lh$, \eqref{U-eq:marginal-rate} applied to
$f$ and $f^2$ bounds $|\Var f(D_v(n)/b)-\Var f(X_{v/h})|$.
Since $zf'(z)$ is bounded, differentiation of $f(c^{-1}s^{-\gamma}G_r)$
gives $|\partial_s\Var f(X_s)|\le C/s$ by dominated convergence.
Thus
$h^{-1}\sum_{\varepsilon h\le v\le Lh}\Var f(X_{v/h})$ differs from
$\int_\varepsilon^L\Var f(X_s)\,ds$ by at most $CL/(\varepsilon h)$.

To control $h^{-1}\sum_v\Xi_v$, let $\Xi_v^T$ be the sum defining
$\Xi_v$ restricted to edges at times $t\le Th$. Moments of the degrees
and the bound $|\xi_{v,e}|\le C(h/t)^\gamma$ give
\begin{equation}\label{U-eq:xi-moment-rate}
 \|\Xi_v\|_4\le C(h/v)^{2\gamma},\qquad
 \|\Xi_v-\Xi_v^T\|_4
       \le C(h/v)^\gamma T^{-\gamma}
       \quad(v\ge\varepsilon h).
\end{equation}
Indeed, Minkowski's inequality bounds these norms by the respective
sums of $C(t/v)^\gamma t^{-1}(h/t)^{2\gamma}$.
Since $\Xi_v^T$ is a deterministic function of the degree path of $v$
up to time $Th$, evaluate that function on the independent paths in
\eqref{U-eq:pair-coupling-rate}, setting it to zero at degrees the
chain cannot reach. The results also obey \eqref{U-eq:xi-moment-rate}.
Cauchy--Schwarz on the event that the coupling fails, followed by the
second bound in \eqref{U-eq:xi-moment-rate}, yields
\begin{equation}\label{U-eq:xi-covariance-rate}
 |\Cov(\Xi_v,\Xi_w)|
 \le C\{\varepsilon^{-5\gamma-1}T^{\gamma+1/2}h^{-1/2}
                 +\varepsilon^{-3\gamma}T^{-\gamma}\},
 \qquad v\ne w\in[\varepsilon h,Lh].
\end{equation}
Summing the bound on $\Var\Xi_v$ from \eqref{U-eq:xi-moment-rate} and
\eqref{U-eq:xi-covariance-rate} shows that, for a finite constant $P$
independent of the cutoffs,
\begin{align}
 \E\left|h^{-1}\sum_v\Xi_v-
                   \int_0^\infty\Var f(X_s)\,ds\right|
 &\le C\{\varepsilon+L^{-1}
              +L\varepsilon^{-P}T^{-\gamma/2}\}\nonumber\\
 &\quad+C(\varepsilon^{-1}+L+T)^P
                   \{h^{-1/4}+(h/n)^{\gamma/2}\}.
 \label{U-eq:average-variance-rate}
\end{align}
Later bounds reuse $P$, enlarged if necessary, with
$P\ge\max(1,2\gamma)$. Each factor $(\cdot)^P$ bounds a finite product
of powers of the cutoffs, so no bound depends exponentially on a cutoff.

We next estimate $K_{n,h}$ by approximating each $\bar\xi_e$.
Set $\zeta_{v,e}=(D_v(t,i-1)+\delta)\xi_{v,e}$, so that
$\bar\xi_e=W_e^{-1}\sum_{v<t}\zeta_{v,e}$.
By \cite[Lemma~3.6, Eq.~(3.24)]{Wang2026}, $\E\zeta_{v,e}$ is $r$ times
the difference of the expectations of $f(\,\cdot\,/b)$ at time $n$
under degree chains started at $m+1$ and $m$ at time $v$, both omitting
edge $e$.
Under the coupling by common uniforms, omitting that edge changes a
chain with probability at most $C(t/v)^\gamma/t$.
For $\varepsilon h\le v<t$ and $h/L\le t\le Lh$, this is at most
$C(\varepsilon^{-1}+L)^P/h$. Applying \eqref{U-eq:marginal-rate} to
the two initial degrees, and using the definition of $\omega_f$ in
\eqref{eq:pa-covariance}, gives, for every $i=1,\ldots,m$,
\begin{equation}\label{U-eq:zeta-mean-rate}
 \sup_{\substack{\varepsilon h\le v<t\\h/L\le t\le Lh}}
 |\E\zeta_{v,(t,i)}-\omega_f(v/h)|
 \le C(\varepsilon^{-1}+L)^P
                  \{h^{-1}+(h/n)^{\gamma/2}\}.
\end{equation}
The moment bounds \cite[Lemma~A.3, Eq.~(A.6)]{Wang2026}, as applied in
the proof of \cite[Lemma~3.6]{Wang2026}, give
\[
 \|\zeta_{v,e}\|_4\le C(h/\max(v,2))^\gamma.
\]
Consequently $W_e^{-1}\sum_{v<\varepsilon h}\zeta_{v,e}$
has $L^2$ norm at most
\begin{equation}\label{U-eq:zeta-small-labels}
 CL\{\varepsilon^{1-\gamma}+h^{\gamma-1}\}.
\end{equation}
For distinct labels $v,w\ge\varepsilon h$,
\eqref{U-eq:pair-coupling-rate} with $T=L$ bounds
$|\Cov(\zeta_{v,e},\zeta_{w,e})|$ by
$C\varepsilon^{-3\gamma-1}L^{\gamma+1/2}h^{-1/2}$.
Summing these covariances and the variances and dividing by $W_e^2$
bounds the standard deviation of
$W_e^{-1}\sum_{\varepsilon h\le v<t}\zeta_{v,e}$ by
$C(\varepsilon^{-1}+L)^Ph^{-1/4}$.
Since $\partial_sf(c^{-1}s^{-1/\alpha}G)=-(\alpha s)^{-1}zf'(z)$ at
$z=c^{-1}s^{-1/\alpha}G$ and $zf'(z)$ is bounded, the derivative of
$\omega_f$ on $[\varepsilon,L]$ is at most $C/\varepsilon$. Hence,
by \eqref{U-eq:zeta-mean-rate}, the Riemann sum
$W_e^{-1}\sum_{\varepsilon h\le v<t}\E\zeta_{v,e}$ differs from
$\Omega_f(t/h)/(m\alpha t/h)$ by at most
$C(\varepsilon^{-1}+L)^P\{h^{-1}+(h/n)^{\gamma/2}\}
+CL\varepsilon$.
It follows that
\begin{align}
 \sup_{\substack{h/L\le t\le Lh\\1\le i\le m}}
 \left\|\bar\xi_{(t,i)}-
         \frac{\Omega_f(t/h)}{m\alpha t/h}\right\|_2
 &\le CL\varepsilon^{1-\gamma}\nonumber\\
 &\quad+C(\varepsilon^{-1}+L)^P
       \{h^{-1/4}+(h/n)^{\gamma/2}+h^{\gamma-1}\}.
 \label{U-eq:weighted-rate}
\end{align}
Both $\bar\xi_{(t,i)}$ and $\Omega_f(t/h)/(m\alpha t/h)$ are bounded,
because $f$ is bounded and $|\Omega_f(u)|\le Cu$, so the expected
absolute difference of their squares is at most a constant times the
norm in \eqref{U-eq:weighted-rate}.
The estimate $\E\bar\xi_{(t,i)}^2\le C\min\{1,(h/t)^2\}$
of \cite[Lemma~3.6]{Wang2026} bounds the contributions of $t<h/L$ and
$t>Lh$ to $K_{n,h}$ by $C/L$, and
$|\Omega_f(u)|\le C\min(u,1)$ gives this bound for the integrals over
$u<L^{-1}$ and $u>L$.
On $[L^{-1},L]$, $\Omega_f'=\omega_f$ and $|\omega_f|\le C$ bound
the derivative of $\{\Omega_f(u)/u\}^2$ by $CL$. Its Riemann sum
therefore differs from the integral by at most $CL^2/h$.
Thus
\begin{align}
 \E\left|K_{n,h}-\frac1{m\alpha^2}
                      \int_0^\infty\Omega_f(u)^2u^{-2}\,du\right|
 &\le C\{L^{-1}+L^2\varepsilon^{1-\gamma}\}\nonumber\\
 &\quad+C(\varepsilon^{-1}+L)^P
       \{h^{-1/4}+(h/n)^{\gamma/2}+h^{\gamma-1}\}.
 \label{U-eq:K-rate}
\end{align}

To obtain a rate, choose the cutoffs as powers of $\min(h,n/h)$. Take
$A>\max\{1,3/(1-\gamma)\}$ and
$D>\max\{A,2(2+AP)/\gamma\}$, and put
\begin{equation}\label{U-eq:cutoffs}
 R=\min(h,n/h)^\lambda,\qquad
 \varepsilon=R^{-A},\qquad L=R,\qquad T=R^D.
\end{equation}
Choose $\lambda>0$ so small that
$\lambda D<1/2$ and
$\lambda DP<\min\{1/4,\gamma/2,1-\gamma\}/2$.
The conditions imposed on the cutoffs before
\eqref{U-eq:pair-coupling-rate} then hold once $\min(h,n/h)$ is large.
The choice of $A$ gives $L^2\varepsilon^{1-\gamma}\le R^{-1}$,
and the choice of $D$ gives $L\varepsilon^{-P}T^{-\gamma/2}\le R^{-1}$.
The choice of $\lambda$ makes each polynomial factor times
$h^{-1/4}$, $(h/n)^{\gamma/2}$, or $h^{\gamma-1}$ at most $CR^{-1}$.
Thus \eqref{U-eq:diagonal-remainder-rate},
\eqref{U-eq:average-variance-rate}, and \eqref{U-eq:K-rate} give
\begin{equation}\label{U-eq:qv-L1-rate}
 \E|V_{n,h}-V(f,f)|\le CR^{-1}.
\end{equation}

The martingale bound stated in \cite[Theorem~1.2]{Mourrat2013}, with
$p=1$ and $p'=2$, now applies directly. Its two inputs are the expected
absolute error of the normalized predictable variance and the sum of
fourth moments of the normalized increments. Equation~\eqref{U-eq:qv-L1-rate}
gives the first, and \eqref{U-eq:increments} gives
\[
 \sum_e\E|\Delta_e|^4
 \le Ch^{-1}\sum_e\E\Delta_e^2
 =Ch^{-1}\E V_{n,h}\le C/h.
\]
Moreover $|\E V_{n,h}-V(f,f)|\le CR^{-1}$, so $\E V_{n,h}$ is bounded
away from zero once $\min(h,n/h)$ is large. Normalizing by
$\sqrt{\E V_{n,h}}$ in the cited theorem yields
\begin{equation}\label{U-eq:normal-variance-bound}
 \sup_z\left|\Prob\{S_{n,h}(f)\le z\}
       -\Phi\{z/\sqrt{\E V_{n,h}}\}\right|
 \le C\{R^{-1/3}+h^{-1/5}\}.
\end{equation}
Changing the normal variance to $V(f,f)$ adds at most $CR^{-1}$.
Since $\lambda<1/8$ and $R=\min(h,n/h)^\lambda$,
\eqref{U-eq:normal-variance-bound} proves
\eqref{U-eq:quantitative-normal} with $\kappa=\lambda/3$.
Increasing $C$ covers the bounded values of $\min(h,n/h)$.
\end{proof}

\subsection{Exclusion of small ranks}\label{U-sec:small-ranks}
Lemma~\ref{lem:pa-upper-localization} excludes ranks large relative to
$a_n$, and this subsection excludes ranks small relative to $a_n$. Use
$\mathcal D_k$, $H_{n,k}$, and the smallest unrestricted minimizer $\wh k$ from
\eqref{eq:Dk}, with $\mathcal D_k=+\infty$ for invalid candidates.
The next lemma bounds the probability that $\sqrt h\,\mathcal D_k$ is small at
some rank in a dyadic interval $[h,2h]$ with $h\le a_n$.
\begin{lemma}\label{U-lem:small-discrepancy}
There are $C<\infty$ and $\kappa>0$, depending on $m,\delta$, such
that, for every $n$, every integer $h$ with $2\le h\le a_n$ and $2h<n$,
and every $0<u\le1$,
\begin{equation}\label{U-eq:small-discrepancy}
 \Prob\left\{\min_{h\le k\le2h}\sqrt h\,\mathcal D_{k}\le u\right\}
 \le Cu+C\{h^{-\kappa}+(h/n)^\kappa\}.
\end{equation}
\end{lemma}

\begin{proof}
For $\vartheta>0$, write
$\nu_\vartheta(\dd z)=\vartheta z^{-\vartheta-1}\dd z$, so that
$\nu_\alpha=\nu$ in \eqref{eq:pa-measure}. The proof applies
Proposition~\ref{U-prop:quantitative-normal} to an $f$ with
$\nu_\vartheta f=O\{(\vartheta-\alpha)^2\}$.
Choose a nonzero $\eta\in C_c^\infty((\log5,\log6))$ and set
\begin{equation}\label{U-eq:orthogonal-function}
 f(z)=z^\alpha\eta''(\log z).
\end{equation}
Integration by parts gives
$\nu_\alpha f=\alpha\int\eta''(w)\,dw=0$ and
$\left.\partial_\vartheta\nu_\vartheta f\right|_{\vartheta=\alpha}
=\int\eta''(w)\,dw-\alpha\int w\eta''(w)\,dw=0$.
Consequently
\begin{equation}\label{U-eq:orthogonal-bound}
 |\nu_\vartheta f|\le C(\vartheta-\alpha)^2,
 \qquad \alpha/2\le\vartheta\le2\alpha.
\end{equation}
Since $\eta$ is nonzero with compact support, $\eta''$ and $f$ are
nonzero, and Proposition~\ref{U-prop:positive-covariance} gives
$V(f,f)>0$.

With $b=c(n/h)^{1/\alpha}$, \eqref{eq:pa-real-tail} and
\eqref{eq:pa-profile-error} give, at every fixed $z>0$,
\begin{equation}\label{U-eq:mean-tail-rate}
 |\E\nu_{n,h}(z,\infty)-z^{-\alpha}|
                                      \le C_z(b^{-1}+h^{-1}).
\end{equation}
Apply \eqref{eq:pa-Bernstein} with $K=h$ and $A=(z,\infty)$ for
$z\in\{4^{-1/\alpha},2,3,4\}$.
For all sufficiently large $h$ and $n/h$, with probability at least
$1-Ch^{-3}$ the following two conditions hold:
\begin{enumerate}
\item $4^{-1/\alpha}\le D_{(k+1)}(n)/b\le2$ for every $h\le k\le2h$,
so that every such threshold $D_{(k+1)}(n)$ is less than $3b$.
Indeed, $\nu_{n,h}(4^{-1/\alpha},\infty)>2$ and
$\nu_{n,h}(2,\infty)<1$ except with probability at most $Ce^{-h/C}$,
because their limiting means are $4$ and $2^{-\alpha}<1$.
\item For $z=3,4$,
$|\nu_{n,h}(z,\infty)-z^{-\alpha}|
 \le C\{b^{-1}+\sqrt{\log h/h}\}$, which locates the fitted exponent.
By \eqref{U-eq:mean-tail-rate} and \eqref{eq:pa-Bernstein} with
$3\log h$ in place of its parameter $u$, this fails with probability
at most $Ch^{-3}$.
\end{enumerate}
For bounded $h$ or $n/h$, \eqref{U-eq:small-discrepancy} holds after
enlarging $C$, because its right side then exceeds one.

Suppose these conditions hold and some $k\in[h,2h]$ satisfies
$\mathcal D_{k}\le u/\sqrt h$.
For $z\ge3$, \eqref{eq:Dk} and condition (i) imply
\begin{equation}\label{U-eq:fitted-tail}
 \left|\nu_{n,h}(z,\infty)
        -(k/h)(D_{(k+1)}(n)/b)^{\wh\alpha_{n,k}}z^{-\wh\alpha_{n,k}}\right|
 \le2\mathcal D_{k}.
\end{equation}
By condition (ii), the second term in \eqref{U-eq:fitted-tail} is
bounded away from zero at $z=3,4$. Taking the ratio of these two
values shows that
\begin{equation}\label{U-eq:fitted-exponent}
 |\wh\alpha_{n,k}-\alpha|\le C\{\mathcal D_{k}+b^{-1}+\sqrt{\log h/h}\},
 \qquad \alpha/2\le\wh\alpha_{n,k}\le2\alpha,
 \qquad (k/h)(D_{(k+1)}(n)/b)^{\wh\alpha_{n,k}}\le C.
\end{equation}
For large $h$ these statements follow from the bounded derivative
of the logarithm on a fixed neighborhood of $3^{-\alpha}$ and
$4^{-\alpha}$.

Since $f$ is supported in $(5,6)$, integration of
\eqref{U-eq:fitted-tail} against $f'(z)$ gives
\[
 |\nu_{n,h}f-(k/h)(D_{(k+1)}(n)/b)^{\wh\alpha_{n,k}}
   \nu_{\wh\alpha_{n,k}}f|
       \le2\|f'\|_1\mathcal D_{k}.
\]
Using \eqref{U-eq:orthogonal-bound} and \eqref{U-eq:fitted-exponent}
therefore yields
\begin{equation}\label{U-eq:contrast-reduction}
 \sqrt h\,|\nu_{n,h}f|
 \le C\{u+\sqrt h/b^2+\log h/\sqrt h\}
 \le C\{u+\log h/\sqrt h\}.
\end{equation}
The last inequality uses $h\le a_n$, which gives $b\ge c\sqrt h$.
Thus $S_{n,h}(f)$ lies within $C\{u+\log h/\sqrt h\}$ of the
deterministic point $-\sqrt h\,\E\nu_{n,h}f$.
Proposition~\ref{U-prop:quantitative-normal} bounds the probability
that $S_{n,h}(f)$ lies in an interval by $C$ times its length plus
$C\{h^{-\kappa}+(h/n)^\kappa\}$, whatever the center.
Adding the exceptional probability $Ch^{-3}$, and replacing $\kappa$
by $\min(\kappa,1/4)$ to absorb $h^{-3}+\log h/\sqrt h$, proves
\eqref{U-eq:small-discrepancy}.
\end{proof}

The next lemma shows that $\wh k/a_n$ is bounded away from zero in
probability.
\begin{lemma}\label{U-lem:lower-localization}
Under the homogeneous PA model \arxiveqref{M-eq:pa-model} with fixed $m\ge1$
and $\delta>-m$,
\begin{equation}\label{U-eq:lower-localization}
 \lim_{L\to\infty}\limsup_{n\to\infty}
       \Prob\{\wh k<a_n/L\}=0.
\end{equation}
\end{lemma}

\begin{proof}
The empirical survival function in \eqref{eq:Dk} has a jump of size at
least $1/k$, whereas the fitted Pareto survival function is
continuous. Hence $\mathcal D_k\ge1/(2k)$, also for invalid candidates.
On the event $\{\sqrt{a_n}\mathcal D_{\lfloor a_n\rfloor}\le M\}$, the bound
$\mathcal D_{\wh k}\le \mathcal D_{\lfloor a_n\rfloor}$ gives $\wh k\ge\sqrt{a_n}/(2M)$.

Cover $[\sqrt{a_n}/(2M),a_n/L]$ by consecutive intervals
$[h,2h]$, starting with an integer $h$ between
$\sqrt{a_n}/(4M)$ and $\sqrt{a_n}/(2M)$.
At most one last interval crosses $a_n/L$.
Put $u=M\sqrt{h/a_n}$, which is at most one when $L\ge M^2$.
On $\{\sqrt{a_n}\mathcal D_{\lfloor a_n\rfloor}\le M\}$, $\wh k\in[h,2h]$ implies
$\min_{h\le k\le2h}\sqrt h\,\mathcal D_{k}\le u$.
Summing \eqref{U-eq:small-discrepancy} over these intervals gives
\begin{align}
 \Prob\{\wh k<a_n/L\}
 &\le\Prob\{\sqrt{a_n}\mathcal D_{\lfloor a_n\rfloor}>M\}
       +CM/\sqrt L\nonumber\\
 &\quad+C_Ma_n^{-\kappa/2}+C(a_n/n)^\kappa,
 \label{U-eq:lower-localization-bound}
\end{align}
because the sums of $Cu$, $Ch^{-\kappa}$, and $C(h/n)^\kappa$ over the
dyadic values of $h$ are geometric.
Take the limit superior in $n$, then let $L$ and $M$ tend to infinity
in that order. Equations \eqref{U-eq:benchmark} and
\eqref{U-eq:lower-localization-bound} prove \eqref{U-eq:lower-localization}.
\end{proof}

\subsection{Uniqueness of the global minimizing pair}
Theorem~\ref{thm:general-lattice}, which applies to PA degrees by
Corollary~\ref{cor:pa-conditional}, identifies the limit of the selected
pair only when $\mathcal D_{\rm lat}$ has a unique minimizer. To prove
uniqueness, write $\mathcal D_{\rm lat}$ in the time coordinate as a
maximum of an absolute residual, affine in the Gaussian path, plus a
strictly increasing deterministic function. This argument uses two
properties of the Gaussian path, which we state explicitly so that
it also applies to the general lattice criterion.
For arbitrary $\alpha,c>0$, $Q\in\R$, and
$B=(Q+\alpha/2)/(\alpha+1)$, let $Y$ be a continuous centered Gaussian
process with $Y(0)=0$, and put
\begin{equation}\label{U-eq:shifted-path}
 X(v)=Y(v)+\frac{\alpha}{c}(1/2-Q)v^{1+1/\alpha}.
\end{equation}
In PA, $\alpha,c$ are given by \arxiveqref{M-eq:pa-parameters}, $Q,B$ by
\eqref{eq:pa-bias-constants}, and $Y(v)=\mathbb G(v^{-1/\alpha})$
has covariance \eqref{U-eq:Y-covariance}.
For $t>0$ and $0<\iota<1/2$, define the continuous linear operator
$\mathcal L_t$ on $E_\iota(t)$ in \eqref{U-eq:weighted-path-space} by
\begin{equation}\label{eq:mdsp-linear-operator}
 (\mathcal L_t x)(v)=x(v)-\frac vt x(t)
 +\frac vt\log\frac vt\left\{\int_0^t\frac{x(u)}u\,du-x(t)\right\},
 \qquad 0\le v\le t.
\end{equation}
The integral converges because $|x(u)|\le u^\iota\|x\|_\iota$, and
$\mathcal L_tx$ vanishes at $v=0$ and $v=t$. For $0\le v\le t$, define
\begin{align}
 a(t,\theta)&=\frac{\alpha}{c}(1/2-\theta)t^{1/\alpha},
 &k_t(v)&=v\{1+\log(v/t)\},
 &w(v)&=\frac{\alpha}{2c}v^{1+1/\alpha}.
 \label{U-eq:criterion-ingredients}
\end{align}
Substituting $v=ty^{-\alpha}$ for the relative excess $y\ge1$ in
\eqref{eq:pa-D-midpoint} gives
\begin{equation}\label{U-eq:criterion-transformed}
 \mathcal D_{\rm lat}(t,\theta;X)=\frac1t\max_{0\le v\le t}
       \left\{|(\mathcal L_tX)(v)+a(t,\theta)k_t(v)|+w(v)\right\}.
\end{equation}
For a path $x\in E_\iota(t)$, $\mathcal D_{\rm lat}(t,\theta;x)$ denotes
the right side with $X$ replaced by $x$, and the expression inside the
absolute value is the residual. The identity
\[
 \mathcal L_t(X-Y)(v)
 =\frac{\alpha}{c}(1/2-Q)
 \left\{v^{1+1/\alpha}-t^{1/\alpha}v
       -\frac{t^{1/\alpha}v\log(v/t)}{\alpha+1}\right\}
\]
shows that adding $a(t,\theta)k_t(v)$ gives the coefficients $Q-\theta$
and $B-\theta$ in \arxiveqref{M-eq:pa-R}. Evaluating the maximum at $v=t$, where
$(\mathcal L_tX)(t)=0$ and $k_t(t)=t$, gives
\begin{equation}\label{U-eq:criterion-endpoint}
 t\mathcal D_{\rm lat}(t,\theta;X)\ge |a(t,\theta)|t+w(t),\qquad
 |a(t,\theta)|\le w(t)/t.
\end{equation}

The next lemma shows that, almost surely, $Y$ plus any deterministic
continuously differentiable function is Lipschitz from the left at no
$t>0$. Proposition~\ref{U-prop:unique-global} uses it to exclude a maximizer
at $v=t$ in \eqref{U-eq:criterion-transformed} and an interval of
maximizers.
\begin{lemma}\label{U-lem:Gaussian-roughness}
For the PA process $Y$ in \eqref{U-eq:Y-covariance}, almost surely,
for every $t>0$,
there are no finite $M$ and positive $\varepsilon$ such that
\begin{equation}\label{U-eq:no-left-Lipschitz}
 |Y(t-h)-Y(t)|\le Mh\qquad(0<h<\varepsilon).
\end{equation}
This assertion also holds for Brownian motion in place of $Y$, and
for either process after adding any deterministic function that is
continuously differentiable on $(0,\infty)$. Moreover, there is
a coupling $Y=\mathbb W-H$, where $\mathbb W$ is Brownian motion and $H$
has continuously differentiable paths on every compact subinterval of
$(0,\infty)$.
\end{lemma}

\begin{proof}
With $X_s$ as in \eqref{eq:pa-covariance} and
$p_v(s)=\Prob\{G_r>c(s/v)^{1/\alpha}\}=\E f_{v^{-1/\alpha}}(X_s)$
for $s,v>0$, \eqref{U-eq:covariance-combined} gives
\[
 \Cov\{Y(v),Y(u)\}
    =\min(v,u)-(1+1/m)\int_0^\infty p_v(s)p_u(s)\,ds.
\]
Take an isonormal process $\mathcal W$ on $L^2((0,\infty),ds)$,
independent of $Y$, and put
$H(v)=\sqrt{1+1/m}\,\mathcal W(p_v)$, $H(0)=0$.
Then $\mathbb W=Y+H$ has Brownian covariance. Also
$\E|H(v)-H(u)|^2\le|v-u|$, because independence gives
$\E|H(v)-H(u)|^2+\E|Y(v)-Y(u)|^2=\E|\mathbb W(v)-\mathbb W(u)|^2=|v-u|$.
Gaussian moment bounds and the Kolmogorov continuity theorem give
continuous versions of $Y$ and $H$, hence of $\mathbb W$, so
$\mathbb W$ is Brownian motion.

Writing $z=c(s/v)^{1/\alpha}$, direct differentiation gives
\[
 \partial_v p_v(s)=\frac{z^re^{-z}}{\alpha v\Gamma(r)},\qquad
 \partial_v^2p_v(s)
   =\frac{z^re^{-z}}{\alpha v^2\Gamma(r)}
                         \left(-1-\frac r\alpha+\frac z\alpha\right).
\]
On every compact interval of $v$ in $(0,\infty)$, the two derivatives
have common square-integrable dominating functions in $s$. Near $s=0$
they are bounded by $Cs^{r/\alpha}$, and as $s\to\infty$ they are
bounded by a polynomial in $s$ times $\exp(-C s^{1/\alpha})$.
On such an interval the second derivative bounds
$\|\partial_vp_v-\partial_vp_u\|_2$ by $C|v-u|$.
Gaussian moment bounds and the Kolmogorov continuity theorem give a
continuous version of $\sqrt{1+1/m}\,\mathcal W(\partial_vp_v)$.
Since $\int_u^v\|\partial_yp_y\|_2\,dy<\infty$ for $0<u<v$, the Gaussian
isometry commutes with $p_v-p_u=\int_u^v\partial_yp_y\,dy$, so
$H(v)-H(u)=\sqrt{1+1/m}\int_u^v\mathcal W(\partial_yp_y)\,dy$ almost
surely, first for rational $u,v$ and then for all by continuity. Hence
$H'=\sqrt{1+1/m}\,\mathcal W(\partial_vp_v)$ is continuous.

We first prove \eqref{U-eq:no-left-Lipschitz} for $\mathbb W$ in place
of $Y$, simultaneously for all $t>0$. Fix an integer $M$ and a bounded
interval of positive $t$. If a point $t$ of this interval satisfied
\eqref{U-eq:no-left-Lipschitz} for $\mathbb W$, then for every
sufficiently small dyadic $\varepsilon$ a dyadic interval of length
$\varepsilon$ inside $[t-2\varepsilon,t]$ would have oscillation at most
$4M\varepsilon$. Divide a fixed such interval into four equal parts.
The oscillation bound forces four independent centered Gaussian
increments, each of variance $\varepsilon/4$, to have absolute value
at most $4M\varepsilon$. The probability of this event is at most
$CM^4\varepsilon^2$. There are at most $C/\varepsilon$ candidate
intervals in a fixed bounded interval. The resulting probabilities
are summable over dyadic $\varepsilon$. The Borel--Cantelli lemma and
countable unions over $M$ and the bounded intervals exclude every such
$t$. Adding a locally continuously differentiable function preserves
this property, since such a function is locally Lipschitz. Applying it
to $Y=\mathbb W-H$, with or without an added deterministic continuously
differentiable function, proves the lemma.
\end{proof}

The next proposition separates the Gaussian hypotheses used for
uniqueness from the proof that the minimum is attained.
\begin{proposition}\label{U-prop:unique-global}
Fix $\alpha,c>0$, $Q\in\R$, and $0<\iota<1/2$.
Suppose that, for every $T_0>0$, the continuous centered Gaussian
process $Y$ has a nondegenerate Gaussian law on $E_\iota(T_0)$,
meaning $\Var\ell(Y)>0$ for every nonzero continuous linear
functional $\ell$ on that space. Suppose also that, almost surely,
no $t>0$ admits the bound \eqref{U-eq:no-left-Lipschitz}, and
$\mathcal D_{\rm lat}$ defined by
\eqref{U-eq:shifted-path}--\eqref{U-eq:criterion-transformed}
attains its global minimum.
Then, almost surely, it has exactly one global minimizing pair
$(T,\Theta)\in(0,\infty)\times[0,1]$.
For every deterministic compact interval $I\subset(0,\infty)$, almost
surely all minimizers of $\mathcal D_{\rm lat}$ on $I\times[0,1]$ share
one time coordinate.
\end{proposition}
For PA, Proposition~\ref{U-prop:positive-covariance} verifies
nondegeneracy, Lemma~\ref{U-lem:Gaussian-roughness} excludes
\eqref{U-eq:no-left-Lipschitz}, and
Proposition~\ref{prop:pa-attainment} gives attainment.

\begin{proof}
We first prove that all minimizing pairs on a compact interval have
one time coordinate. We then rule out two tie fractions at a global
minimizing time. Call $v$ active when it attains the maximum in
\eqref{U-eq:criterion-transformed}.

We will need the residual to be nonzero at each active point.
For $v<t$ this follows from $w(v)<w(t)$ and
\eqref{U-eq:criterion-endpoint}. At $v=t$ it follows from
$a(t,\theta)\ne0$, except when $\theta=1/2$.
In that case, if $t$ were active, the maximum would be $w(t)$ and
\[
 |(\mathcal L_tX)(v)|\le w(t)-w(v)\le C(t-v)
\]
for $v<t$ close to $t$. The difference between
$(\mathcal L_tX)(v)$ and $X(v)-X(t)$ is bounded by $C(t-v)$ there.
This would make $X$ Lipschitz from the left at $t$, contrary to
the hypothesis on $Y$, because $X-Y$ is continuously
differentiable. Thus, on one event of probability one,
\begin{equation}\label{U-eq:midpoint-endpoint}
 \theta=1/2\quad\Longrightarrow\quad
 t\text{ is not active, for every }t>0.
\end{equation}
Consequently every active point has a unique residual sign, and zero
is never active because the maximum is positive.

To compare minimizing times, fix $[u,T_0]\subset(0,\infty)$ and
use the exponent $\iota$ in the hypothesis.
On $E_\iota(T_0)$ define the minimum value
\[
 \mathcal F(x)=\min_{\substack{u\le t\le T_0\\0\le\theta\le1}}
                       \mathcal D_{\rm lat}(t,\theta;x).
\]
The bounds
\[
 \left|\int_0^t x(v)\frac{dv}{v}\right|
 \le t^\iota\|x\|_\iota/\iota,
 \qquad \sup_{0\le s\le1}|s\log s|<\infty
\]
show that $\mathcal D_{\rm lat}$ is Lipschitz in the path, uniformly in
$(t,\theta)\in[u,T_0]\times[0,1]$. Substituting $v=ts$ in its maximum
also shows joint continuity in $(t,\theta,x)$. Hence $\mathcal F$ is
Lipschitz. Theorem~6 of Phelps~\cite{Phelps1978} makes $\mathcal F$
G\^ateaux differentiable outside a set null for every nondegenerate
Gaussian law on this separable Banach space. By hypothesis, $Y$ has
such a law. Its deterministic translate $X$ does too, because
$v\mapsto v^{1+1/\alpha}$ belongs to $E_\iota(T_0)$.
Almost surely, therefore, $\mathcal F$ has a continuous linear
derivative $\ell$ at $X$.

This derivative gives one supporting inequality for every minimizing
pair. Fix such a pair $(t,\theta)$ and put $z=\mathcal F(X)>0$.
The criterion at fixed $t$ is convex jointly in the path and tie
fraction. For $0<\varepsilon<1$, $h\in E_\iota(T_0)$, and
$\theta'\in[0,1]$, convexity gives
\[
 \mathcal F(X+\varepsilon h)
 \le\mathcal D_{\rm lat}(t,(1-\varepsilon)\theta+\varepsilon\theta';X+\varepsilon h)
 \le(1-\varepsilon)z+\varepsilon\mathcal D_{\rm lat}(t,\theta';X+h).
\]
Subtract $z$, divide by $\varepsilon$, and let $\varepsilon\downarrow0$.
We obtain
\begin{equation}\label{U-eq:supporting-bound}
 \mathcal D_{\rm lat}(t,\theta';X+h)\ge z+\ell(h),
 \qquad \theta'\in[0,1],\ h\in E_\iota(T_0).
\end{equation}

We express $\ell$ using the points that attain the maximum.
Writing the absolute value in \eqref{U-eq:criterion-transformed} as
the maximum over its two signs makes the criterion a maximum of
affine functions of $(x,\theta)$. Their active gradients form a
weak-star compact set. Its closed convex hull consists of averages
of those gradients against probability measures. The unique sign
at every active point turns such an average into a signed measure
$\sigma_t$ on $[0,t]$, supported on the active points and having their
residual signs. Equation~\eqref{U-eq:supporting-bound} implies that
one can choose this measure so that
\begin{align}
 t\ell(h)&=\int(\mathcal L_th)(v)\,\sigma_t(dv),\nonumber\\
 (\theta'-\theta)\int k_t(v)\,\sigma_t(dv)&\le0
                         \qquad(0\le\theta'\le1).
 \label{U-eq:joint-subgradient}
\end{align}
Indeed, otherwise separation from this compact convex set would give
a feasible direction $(h,s)$ in which every active affine function
has derivative less than $\ell(h)$. The directional derivative of the
maximum is the maximum of the active derivatives, contradicting
\eqref{U-eq:supporting-bound}. The sign in the second line of
\eqref{U-eq:joint-subgradient} uses
$\partial_\theta a(t,\theta)=-\alpha t^{1/\alpha}/c<0$.
Since the signs are unique, the averaging measure has no cancellation,
so
\begin{equation}\label{U-eq:active-variation}
 |\sigma_t|([0,t])=1.
\end{equation}

Some of this mass must lie before $t$.
If all of it were at $t$, then $\theta\ne1/2$ by
\eqref{U-eq:midpoint-endpoint}, and
$\int k_t(v)\,\sigma_t(dv)=t\operatorname{sign}(1/2-\theta)$.
Taking $\theta'=1/2$ in \eqref{U-eq:joint-subgradient} would give a
strictly positive left side. Hence
\begin{equation}\label{U-eq:active-interior-mass}
 |\sigma_t|([0,t))>0.
\end{equation}

Suppose now that $(t_1,\theta_1)$ and $(t_2,\theta_2)$ both minimize,
with $t_1<t_2$. The representation at $t_1$ gives $\ell(h)=0$ for every
$h$ vanishing on $[0,t_1]$. To use this at $t_2$, let $\sigma^\circ$
be the restriction of $\sigma_{t_2}$ to $[0,t_2)$ and put
\[
 d_0=\int\frac v{t_2}\log(v/t_2)\,\sigma^\circ(dv).
\]
The endpoint contributes nothing to $\mathcal L_{t_2}h$. Expanding
\eqref{eq:mdsp-linear-operator} in \eqref{U-eq:joint-subgradient} gives
\begin{equation}\label{U-eq:dual-supported}
 t_2\ell(h)=\int h(v)\,\sigma^\circ(dv)
 +d_0\int_0^{t_2}h(v)\frac{dv}{v}
 -\left\{d_0+\int\frac v{t_2}\,\sigma^\circ(dv)\right\}h(t_2).
\end{equation}
For continuous $h$ supported in $(t_1,t_2)$ this implies
$\sigma^\circ(dv)=-d_0\,dv/v$ on that interval.
If $d_0\ne0$, every point there is active, with residual sign
$-\operatorname{sign}(d_0)$. The equality defining the maximum would
then express $X(v)$ as a linear combination of $1,v,v\log v,w(v)$
on $(t_1,t_2)$, contradicting the hypothesis on $Y$ because $X-Y$
is continuously differentiable there.
Thus $d_0=0$ and $\sigma^\circ$ gives no mass to $(t_1,t_2)$.
Taking $h=0$ on $[0,t_1]$ and $h(t_2)=1$ in
\eqref{U-eq:dual-supported} now gives $\int v\,\sigma^\circ(dv)=0$.
Together with $d_0=0$, this proves
\begin{equation}\label{U-eq:dual-moments}
 \operatorname{supp}(\sigma^\circ)\subset[0,t_1],\qquad
 \int v\,\sigma^\circ(dv)=\int v\log v\,\sigma^\circ(dv)=0.
\end{equation}

These two vanishing moments contradict $t_1<t_2$.
The residual at either pair differs from $X(v)$ by a linear
combination of $v$ and $v\log v$. Integrating the equality at the
active points of $t_2$ against $|\sigma^\circ|$, then using the
criterion bound at $t_1$, therefore gives
\begin{equation}\label{U-eq:dual-contradiction}
 \begin{split}
 t_2z\,|\sigma^\circ|([0,t_1])
 &=\int X(v)\,\sigma^\circ(dv)+\int w(v)\,|\sigma^\circ|(dv)\\
 &\le t_1z\,|\sigma^\circ|([0,t_1]).
 \end{split}
\end{equation}
The mass is positive by \eqref{U-eq:active-interior-mass}, and $z>0$.
Thus all minimizers on $[u,T_0]\times[0,1]$ have one time coordinate.

It remains to show that two tie fractions cannot minimize globally
at the same time. Suppose $\theta_1<\theta_2$ both give the global
minimum $z$ at time $t$. Convexity makes every intermediate $\theta$
a minimizer. At such a $\theta$, each active affine function has zero
slope in $\theta$, since a nonzero slope would make its value exceed
$z$ at one of the endpoints $\theta_1,\theta_2$.
Its slope is proportional to $k_t(v)$, whose only positive zero is
$v_*=t/e$. Hence the maximum has one active point $v_*$ and one sign.

To decrease the criterion, we vary time while keeping one coefficient
of the residual fixed. This avoids differentiating the Gaussian path.
Choose $\theta\in(\theta_1,\theta_2)$, write
$I_s=\int_0^sX(v)\,dv/v$, and set
$\beta_*=(I_t-X(t))/t+a(t,\theta)$.
The residual at time $t$ becomes
\begin{equation}\label{U-eq:residual-fixed-beta}
 X(v)-\frac vtI_t+\beta_*v\{1+\log(v/t)\}.
\end{equation}
At time $s$, the possible values of $(I_s-X(s))/s+a(s,\theta')$ form
\[
 \left[\frac{I_s-X(s)}s-\frac{w(s)}s,
       \frac{I_s-X(s)}s+\frac{w(s)}s\right].
\]
The fixed value $\beta_*$ is in the interior at $s=t$. Continuity
therefore gives a continuous choice $\theta(s)\in(0,1)$ for $s$ near
$t$ that keeps this value fixed. Since $I_s'=X(s)/s$, the derivative
of \eqref{U-eq:residual-fixed-beta} with respect to time, at $s=t$
and fixed $v$, is $-a(t,\theta)v/t$.

The unique maximizer $v_*=t/e$ is strictly inside $(0,t)$.
Continuity preserves the gap at the moving upper endpoint, so the
maximum for $s$ near $t$ can be taken over a fixed compact interval.
Every maximizing point tends to $v_*$ as $s\to t$, by uniqueness and
compactness. On a neighborhood of $v_*$ the residual has constant
sign and a continuous time derivative. Comparing the maxima at $s$
and $t$ at their respective maximizing points therefore differentiates
the maximum at $v_*$. If $R_*$ is the nonzero residual there, this gives
\begin{equation}\label{U-eq:criterion-derivative}
 \left.\frac{d}{ds}\mathcal D_{\rm lat}(s,\theta(s);X)\right|_{s=t}
 =\frac{-\operatorname{sign}(R_*)a(t,\theta)/e-z}{t}<0.
\end{equation}
The inequality follows from $|a(t,\theta)|\le w(t)/t\le z$ and $z>0$.
Increasing $s$ would thus lower the global minimum, a contradiction.

Finally, intersect the events of probability one for
$[u,T_0]=[1/j,j]$, $j\ge2$, with the assumed event of global attainment.
Any two global minimizing times
would lie in one of these intervals and contradict
\eqref{U-eq:dual-contradiction}. Two minimizing tie fractions at their
common time contradict \eqref{U-eq:criterion-derivative}. The global
minimizing pair is therefore unique almost surely.
\end{proof}

\subsection{Completion of the proof}\label{U-sec:joint-limit}
Lemmas~\ref{lem:pa-upper-localization} and~\ref{U-lem:lower-localization}
give the localization \arxiveqref{M-eq:rank-bounds}, and
Proposition~\ref{U-prop:unique-global} gives the unique global minimizer
$(T,\Theta)$ of $\mathcal D_{\rm lat}$. It remains to prove
\arxiveqref{M-eq:unrestricted-joint} jointly with the limit of
$D_{(\wh k+1)}(n)/b_n$, and the two probability limits of
Theorem~\ref*{M-thm:vanishing}.

\begin{proof}
Place the rank $k$ at the time $t=(D_{(k+1)}(n)/b_n)^{-\alpha}$ of its
threshold, with $b_n=c\sqrt{a_n}$, and at the tie fraction
$\theta=\theta_{n,k}$. For each deterministic $I=[L^{-1},L]$,
\eqref{eq:lattice-uniform-criterion} and \eqref{eq:lattice-hill-expansion}
in the proof of Theorem~\ref{thm:general-lattice}, which applies to PA
degrees by Corollary~\ref{cor:pa-conditional}, show that
$\sqrt{a_n}\mathcal D_k$ and $\sqrt{a_n}(H_{n,k}-1/\alpha)$ converge to
$\mathcal D_{\rm lat}(t,\theta)$ and $\mathcal A(t,\theta)$ uniformly over
ranks with $k/a_n\in I$. That proof also shows that these pairs
$(t,\theta)$ become dense in $I\times[0,1]$ and that $k/a_n-t\to0$
uniformly in probability.

With Theorem~\ref{thm:pa-functional} and these uniform limits, the
localization \arxiveqref{M-eq:rank-bounds} implies joint tightness of the
four coordinates in \arxiveqref{M-eq:unrestricted-joint}, of
$D_{(\wh k+1)}(n)/b_n$ and $\sqrt{a_n}(H_{n,\wh k}-1/\alpha)$, of
$\mathbb G_n$ in \eqref{eq:pa-functional} with $K_n=a_n$, and of its
logarithmic integrals. In every joint subsequential limit, the last two
have continuous versions with the laws of $\mathbb G$ and $\mathbb L$ in
\eqref{eq:pa-G}. Let $(T',\Theta')$ be the limit of
$(\wh k/a_n,\theta_{n,\wh k})$ along the subsequence; then
$0<T'<\infty$ almost surely by \arxiveqref{M-eq:rank-bounds}, and
$\sqrt{a_n}\mathcal D_{\wh k}$ and $\sqrt{a_n}(H_{n,\wh k}-1/\alpha)$ converge to
$\mathcal D_{\rm lat}(T',\Theta';X)$ and $\mathcal A(T',\Theta')$.

Proposition~\ref{U-prop:unique-global} gives uniqueness of the pair only
for the global minimizer, whereas \eqref{eq:lattice-selected-law} in
Theorem~\ref{thm:general-lattice} assumes a unique minimizer on
$I\times[0,1]$. We therefore compare $\wh k$ with ranks at rational
pairs. Fix $s\in\mathbb Q_{>0}$ and $\phi\in\mathbb Q\cap[0,1]$,
and let $d_n(s)=\lfloor b_ns^{-1/\alpha}\rfloor$, so that $d_n(s)/b_n$
tends to the threshold coordinate $s^{-1/\alpha}$ of the time $s$.
On the event $N_{n,=d_n(s)}>0$, define, as in
Theorem~\ref{thm:pa-degree-clt},
\begin{equation}\label{U-eq:rational-candidate}
 k_n(s,\phi)=N_{n,>d_n(s)}
 +\min\{\lfloor\phi N_{n,=d_n(s)}\rfloor,N_{n,=d_n(s)}-1\}.
\end{equation}
Equation~\eqref{eq:lattice-atoms}, which holds for PA degrees by
Corollary~\ref{cor:pa-conditional}, and
Theorem~\ref{thm:pa-degree-clt} with $(t,\theta)=(s,\phi)$ give
\[
 \frac{N_{n,=d_n(s)}}{\sqrt{a_n}}
       \xrightarrow{p}\frac\alpha c s^{1+1/\alpha}>0,
 \qquad
 \frac{k_n(s,\phi)}{a_n}\xrightarrow{p}s,
 \qquad \theta_{n,k_n(s,\phi)}\xrightarrow{p}\phi
 \quad\text{in }\R.
\]
Hence $2\le k_n(s,\phi)<n$ and
$H_{n,k_n(s,\phi)}>0$ with probability tending to one.
On the complement use rank $n-1$, which is valid for $n\ge3$:
the newest vertex has degree $m$ and the total degree
$2m(n-1)>nm$ forces a larger degree. This choice changes
$k_n(s,\phi)$ only on an event of probability tending to zero.

By definition of $\wh k$, $\mathcal D_{\wh k}\le \mathcal D_{k_n(s,\phi)}$.
Along the subsequence, $\sqrt{a_n}\mathcal D_{k_n(s,\phi)}$ converges to
$\mathcal D_{\rm lat}(s,\phi;X)$ jointly with the selected coordinates, so
\[
 \mathcal D_{\rm lat}(T',\Theta';X)\le \mathcal D_{\rm lat}(s,\phi;X)\quad\text{almost surely}.
\]
The rational pairs are countable, so these inequalities hold on one
event of probability one. Continuity of $\mathcal D_{\rm lat}$ extends them to every
$(s,\phi)\in(0,\infty)\times[0,1]$. Thus $(T',\Theta')$ is a global
minimizer, and Proposition~\ref{U-prop:unique-global} gives
$(T',\Theta')=(T,\Theta)$. All subsequential limits therefore have one
law. This proves the convergence of $\wh k/a_n$, $\theta_{n,\wh k}$,
$\sqrt{a_n}\mathcal D_{\wh k}$, and $D_{(\wh k+1)}(n)/b_n=t^{-1/\alpha}$ at
$k=\wh k$, whose limit is $T^{-1/\alpha}$ because $k/a_n-t\to0$, together
with
\[
 \sqrt{a_n}(H_{n,\wh k}-1/\alpha)
           \Rightarrow\mathcal A(T,\Theta)\quad\text{in }\R.
\]
Since $H_{n,\wh k}\to1/\alpha$ in probability, the identity
$\sqrt{a_n}(\wh\alpha_{n,k}-\alpha)=-\alpha\sqrt{a_n}(H_{n,k}-1/\alpha)/H_{n,k}$
at $k=\wh k$ gives the third coordinate of
\arxiveqref{M-eq:unrestricted-joint}.
All these operations use one Gaussian process and one minimizing pair,
so the convergence is joint.

By \eqref{eq:pa-large-t},
$\mathcal D_{\rm lat}(t,\theta;X)\ge\alpha t^{1/\alpha}/(2c)$.
Attainment and $T\in(0,\infty)$ imply $0<\mathcal D_{\rm lat}(T,\Theta)<\infty$
almost surely, so the last limit in \arxiveqref{M-eq:unrestricted-joint} is
positive and finite.
Finally $a_n/n\to0$ and $\sqrt{a_n}\to\infty$ give
$\wh k/n\to0$ and $\wh\alpha_{n,\wh k}\to\alpha$ in probability.
\end{proof}

\section{Corrections and confidence intervals}\label{app:unrestricted-inference}
\subsection{Correction at a deterministic rank}\label{sec:pa-correction}
Use the corrected estimators in Eqs.~\arxiveqref{M-eq:corrected-estimator}
and \arxiveqref{M-eq:second-order-estimator}, with the stated invalid-fit conventions.
The coefficient $C_2$ is given in Eq.~\eqref{eq:pa-C2}.
Proposition~\ref{prop:pa-beta-covariance} identifies
$\sigma(\alpha,m)^2$ in Eq.~\arxiveqref{M-eq:pa-factor} as
$\alpha^4V(h_1,h_1)>0$. Subtracting the inverse-degree term removes
the rounding contribution at every fixed normalized rank.
\begin{corollary}\label{cor:pa-corrected}
In the PA model \arxiveqref{M-eq:pa-model} with fixed known $m$ and $\delta>-m$,
fix $t>0$ and put $k=\lfloor a_nt\rfloor$, where $a_n$ is given by
\arxiveqref{M-eq:balance}. Along the full sequence,
\begin{align}
 \sqrt{k}(\wt H_{n,k}-1/\alpha)
       &\wto N\{0,\sigma(\alpha,m)^2/\alpha^4\}\quad\hbox{in }\R,\notag\\
 \sqrt{k}(\wt\alpha_{n,k}-\alpha)
       &\wto N\{0,\sigma(\alpha,m)^2\}\quad\hbox{in }\R.
 \label{eq:pa-corrected-clt}
\end{align}
\end{corollary}
\begin{proof}
Subtracting $(\wh B_{n,k}-\theta_{n,k})/D_{(k+1)}(n)$ in
\arxiveqref{M-eq:corrected-estimator} cancels the term
$(B-\Theta_\zeta)/(cx)$, $x=t^{-1/\alpha}$, of $\mathcal A(t,\Theta_\zeta)$
in \eqref{eq:pa-fixed-rank-law}. To verify this along the full sequence,
start with any subsequence and use compactness of $[0,1]$ to select a
further subsequence on which
$b_nt^{-1/\alpha}-\lfloor b_nt^{-1/\alpha}\rfloor\to\zeta$ in $\R$.
Proposition~\ref{prop:pa-fixed-rank} implies, along that subsequence,
\[
 \frac{D_{(k+1)}(n)}{b_n}\xrightarrow{p}t^{-1/\alpha},
 \qquad H_{n,k}\xrightarrow{p}1/\alpha,
 \qquad \wh B_{n,k}\xrightarrow{p}B
 \quad\hbox{in }\R,
\]
where the last convergence uses the continuity of $B(\cdot,m)$ in
\arxiveqref{M-eq:corrected-estimator}.
Since $b_n=c\sqrt{a_n}$, Slutsky's theorem and
\eqref{eq:pa-fixed-rank-law} yield
\begin{equation}\label{eq:pa-correction-joint}
 \left(\sqrt{a_n}(H_{n,k}-1/\alpha),\
       \frac{\sqrt{a_n}(\wh B_{n,k}-\theta_{n,k})}{D_{(k+1)}(n)}\right)
 \wto\left(\mathcal A(t,\Theta_\zeta),\
               \frac{t^{1/\alpha}}c(B-\Theta_\zeta)\right)
 \quad\hbox{in }\R^2.
\end{equation}
Subtracting the second coordinate from the first, as prescribed in
\arxiveqref{M-eq:corrected-estimator}, gives
\begin{equation}\label{eq:pa-correction-limit}
 \sqrt{a_n}(\wt H_{n,k}-1/\alpha)
 \wto
 \frac{\mathbb L(t^{-1/\alpha})-\mathbb G(t^{-1/\alpha})/\alpha}{t}
 \quad\hbox{in }\R.
\end{equation}
The right side does not depend on $\zeta$.
Every subsequence therefore has a further subsequence with this
limit, which proves \eqref{eq:pa-correction-limit} along the full sequence.
Its variance is
$V(h_{t^{-1/\alpha}},h_{t^{-1/\alpha}})/t^2=V(h_1,h_1)/t$
by \eqref{eq:pa-homogeneity}, and $\alpha^4V(h_1,h_1)=\sigma(\alpha,m)^2$
by \eqref{eq:pa-beta-covariance} at $x=y=1$ and \arxiveqref{M-eq:pa-factor}.
Multiplying by $\sqrt{k/a_n}\to\sqrt t$ in $\R$ proves the first
limit in \eqref{eq:pa-corrected-clt}.
In particular $\wt H_{n,k}$, $H_{n,k}$, and $D_{(k+1)}(n)$ are positive
with probability tending to one, and then
\[
 \sqrt{k}(\wt\alpha_{n,k}-\alpha)
 =-\frac{\alpha}{\wt H_{n,k}}\sqrt{k}(\wt H_{n,k}-1/\alpha).
\]
The factor $-\alpha/\wt H_{n,k}$ converges in probability to $-\alpha^2$, proving the
second limit in \eqref{eq:pa-corrected-clt}.
\end{proof}

\subsection{Equivalence of the second-order correction}
Fix $t>0$. At the deterministic rank $k=\lfloor a_nt\rfloor$, and at
$k=\wh k$, Proposition~\ref{prop:pa-fixed-rank} and
Theorem~\ref*{M-thm:vanishing} give $\wh\alpha_{n,k}\to\alpha$ in
probability in $\R$ and convergence in distribution of $D_{(k+1)}(n)/b_n$
to a positive limit. Since $b_n=c\sqrt{a_n}\to\infty$, $\wh B_{n,k}\to B$,
and $0\le\theta_{n,k}<1$, it follows that $\wt\alpha_{n,k}\to\alpha$ and
$\bar\alpha_{n,k}\to\alpha$ in probability, and that $\sqrt k/D_{(k+1)}(n)$
converges in distribution.
Hence $\sqrt k/D_{(k+1)}(n)^2\xrightarrow{p}0$ in $\R$, and continuity of
$B$ and $C_2$ makes $\sqrt k$ times each term of
\[
 \wt H_{n,k}^{(2)}-\wt H_{n,k}
 =\frac{B(\wh\alpha_{n,k},m)-B(\bar\alpha_{n,k},m)}{D_{(k+1)}(n)}
 -\frac{C_2(\bar\alpha_{n,k},m,\theta_{n,k})}{D_{(k+1)}(n)^2}
\]
converge to zero in probability in $\R$. Because
$\wt H_{n,k}\to1/\alpha$ in probability, the identity
$\wt\alpha_{n,k}^{(2)}-\wt\alpha_{n,k}
=(\wt H_{n,k}-\wt H_{n,k}^{(2)})/(\wt H_{n,k}\wt H_{n,k}^{(2)})$,
valid when both denominators are positive, gives
\begin{equation}\label{eq:pa-second-order-equivalence}
 \sqrt{k}(\wt\alpha_{n,k}^{(2)}-\wt\alpha_{n,k})\xrightarrow{p}0
 \quad\hbox{in }\R
\end{equation}
at both ranks.

\subsection{Correction at the selected rank}\label{sec:pa-unrestricted-inference}
The joint law in Theorem~\ref*{M-thm:vanishing} also determines the corrected exponent at the selected rank.
\begin{corollary}\label{cor:pa-corrected-selected}
Under the hypotheses of Theorem~\ref*{M-thm:vanishing}, with $m$ known
and $(T,\Theta)$ the minimizer of $\mathcal D_{\rm lat}$ given there,
\begin{equation}\label{eq:pa-corrected-selected}
 \sqrt{\wh k}(\wt\alpha_{n,\wh k}-\alpha)
 \wto-\frac{\alpha^2}{\sqrt T}
       \{\mathbb L(T^{-1/\alpha})-\mathbb G(T^{-1/\alpha})/\alpha\}
       \quad\hbox{in }\R.
\end{equation}
\end{corollary}
\begin{proof}
Theorem~\ref*{M-thm:vanishing} gives
$D_{(\wh k+1)}(n)/b_n\wto T^{-1/\alpha}$ jointly with
\arxiveqref{M-eq:unrestricted-joint}. Inverting the third coordinate of
\arxiveqref{M-eq:unrestricted-joint} by the delta method for $h\mapsto h^{-1}$
gives $\sqrt{a_n}(H_{n,\wh k}-1/\alpha)\wto\mathcal A(T,\Theta)$ jointly.
Consistency of $\wh\alpha_{n,\wh k}$ and continuity of $B(\cdot,m)$ give
$\wh B_{n,\wh k}\to B$ in probability. Since $b_n=c\sqrt{a_n}$, the
continuous mapping theorem gives
\[
 \frac{\sqrt{a_n}(\wh B_{n,\wh k}-\theta_{n,\wh k})}{D_{(\wh k+1)}(n)}
 \wto\frac{(B-\Theta)T^{1/\alpha}}{c}
\]
jointly, and this limit is the second term of $\mathcal A(T,\Theta)$ in
\arxiveqref{M-eq:pa-A}. Subtracting it as prescribed in
\arxiveqref{M-eq:corrected-estimator} gives
\[
 \sqrt{a_n}(\wt H_{n,\wh k}-1/\alpha)\wto
 \frac{\mathbb L(T^{-1/\alpha})-\mathbb G(T^{-1/\alpha})/\alpha}{T}
 \quad\hbox{in }\R.
\]
In particular $\wt H_{n,\wh k}\to1/\alpha>0$ in probability, so
$\wt\alpha_{n,\wh k}=\wt H_{n,\wh k}^{-1}$ with probability tending to one.
The map $h\mapsto h^{-1}$ has derivative $-\alpha^2$ at $1/\alpha$, and
$\sqrt{\wh k/a_n}\wto\sqrt T$ jointly, which proves
\eqref{eq:pa-corrected-selected}.
\end{proof}

Dividing this limit by $\sigma(\alpha,m)$ gives
$\mathcal E_{\alpha,m}$ in Eq.~\arxiveqref{M-eq:selected-limit}.
Use its distribution function and quantiles in Eq.~\arxiveqref{M-eq:selected-cdf}
and the interval $C_n(\upsilon)$ in Eq.~\arxiveqref{M-eq:unrestricted-interval}.
The next theorem also permits simulated quantiles, with probability
taken jointly over the graph and simulation draws.
\begin{theorem}\label{thm:pa-unrestricted-interval}
Assume the homogeneous PA model \arxiveqref{M-eq:pa-model} with fixed known
$m\ge1$ and $\delta>-m$. For every fixed $\upsilon\in(0,1)$,
\begin{equation}\label{eq:pa-unrestricted-coverage}
 \Prob\{\alpha\in C_n(\upsilon)\}\longrightarrow1-\upsilon.
\end{equation}
Let $\wh F_n$ be a random distribution function, defined on a
probability space that carries the graph and the simulation draws, such
that, with the left side taken as zero when $\wt\alpha_{n,\wh k}\le1$,
\begin{equation}\label{eq:pa-calibration-accuracy}
 \mathbf1\{\wt\alpha_{n,\wh k}>1\}\sup_{z\in\R}
 |\wh F_n(z)-F_{\wt\alpha_{n,\wh k},m}(z)|
       \xrightarrow{p}0\quad\hbox{in }\R.
\end{equation}
Then \eqref{eq:pa-unrestricted-coverage} also holds, with $\Prob$ the
joint law of the graph and the draws, when the quantiles in
\arxiveqref{M-eq:unrestricted-interval} are replaced by
$\inf\{z:\wh F_n(z)\ge p\}$ at $p=\upsilon/2$ and $p=1-\upsilon/2$.
\end{theorem}
The proof follows the continuity result in this section. Its main input,
Proposition~\ref{prop:pa-selected-distribution}, shows that
$F_{\alpha,m}$ is continuous and depends continuously on $\alpha$ in
the supremum norm. That proof uses $\wt\alpha_{n,\wh k}$ only through its
consistency and the weak limit \eqref{eq:pa-corrected-selected}. By
\eqref{eq:pa-second-order-equivalence} and Slutsky's lemma,
$\wt\alpha_{n,k}^{(2)}$ has both properties, so
Corollaries~\ref{cor:pa-corrected} and \ref{cor:pa-corrected-selected} and
Theorem~\ref{thm:pa-unrestricted-interval} also hold with
$\wt\alpha_{n,k}^{(2)}$ in place of $\wt\alpha_{n,k}$.

Section~\ref{sec:pa-unrestricted-inference} writes $\mathcal E_{\alpha,m}$
in the time coordinate, in which $Y(t)=\mathbb G(t^{-1/\alpha})$.
Replacing $T$ there by $t>0$ defines
\begin{equation}\label{eq:pa-error-time}
 \mathcal E_{\alpha,m}(t)=-\frac{\alpha}{\sigma(\alpha,m)\sqrt t}
     \left\{\int_0^tY(v)\frac{\dd v}{v}-Y(t)\right\},
\end{equation}
so that $\mathcal E_{\alpha,m}=\mathcal E_{\alpha,m}(T)$.

\subsection{An independent component of \texorpdfstring{$Y$ along $v\log v$}{Y along v log v}}
Adding $v\log v$ to a path $y$ lowers $\int_0^ty(v)\,\dd v/v-y(t)$ by
$t$, and so shifts $\mathcal E_{\alpha,m}(t)$, but leaves $\mathcal L_ty$
in \eqref{eq:mdsp-linear-operator}, and hence the criterion
\eqref{U-eq:criterion-transformed}, unchanged.
Lemma~\ref{lem:pa-Gaussian-direction} shows that, on every interval
$[0,L]$, $Y$ has the law of a Gaussian process plus an independent normal
multiple of $v\log v$.

\begin{lemma}\label{lem:pa-Gaussian-direction}
For every $L>0$, there exist a centered Gaussian process $Y^\circ$ on
$[0,L]$ and an independent centered normal variable $Z_L$ with
strictly positive variance such that, with $v\log v=0$ at $v=0$,
\begin{equation}\label{eq:pa-Gaussian-direction}
 \{Y(v):0\le v\le L\}\ \stackrel{d}{=}\
 \{Y^\circ(v)+Z_L v\log v:0\le v\le L\}.
\end{equation}
For every $0<\iota<1/2$, both $Y$ and $Y^\circ$ have versions with
paths in $E_\iota(L)$.
\end{lemma}

\begin{proof}
We obtain $Z_L$ from the operator $\mathsf A_\alpha$ on $L^2(\nu)$
defined by $V(f,g)=\langle f,\mathsf A_\alpha g\rangle$, which is bounded
and positive because Proposition~\ref{U-prop:positive-covariance} extends
$V$ continuously to $L^2(\nu)$. By \eqref{U-eq:covariance-Fourier}, the
unitary map $\mathsf U_\alpha$ in \eqref{U-eq:covariance-transform},
followed by the Fourier transform, carries $\mathsf A_\alpha$ to
multiplication by
\begin{equation}\label{eq:pa-covariance-multiplier}
 a_\alpha(\xi)=1-(1+1/m)
 \frac{|\Gamma(r+\alpha/2+i\xi)|^2}{\Gamma(r)\Gamma(r+\alpha)}.
\end{equation}
This multiplier satisfies
\begin{equation}\label{eq:pa-multiplier-lower}
 \inf_{\xi\ne0}\frac{a_\alpha(\xi)}{\min(\xi^2,1)}>0.
\end{equation}
Indeed, the proof of Proposition~\ref{U-prop:positive-covariance}
gives $a_\alpha(0)\ge0$ and $a_\alpha(\xi)>0$ for $\xi\ne0$.
The ratio $\Gamma(r+\alpha/2+i\xi)/\Gamma(r+\alpha/2)$ is the
characteristic function of $\log G_{r+\alpha/2}$, which has finite,
strictly positive variance, so expanding the squared modulus as
$\xi\to0$ gives $a_\alpha(\xi)=a_\alpha(0)+c'\xi^2+o(\xi^2)$ with
$c'>0$. Since $\log G_{r+\alpha/2}$ has a density, the Riemann--Lebesgue
lemma gives $a_\alpha(\xi)\to1$ as $|\xi|\to\infty$. Continuity and
strict positivity of $a_\alpha$ for $\xi\ne0$ complete the proof of
\eqref{eq:pa-multiplier-lower}.

Choose a smooth function $\varphi$ on $(0,\infty)$ equal to $v\log v$
on $(0,L]$, zero for sufficiently large $v$, and satisfying
\begin{equation}\label{eq:pa-direction-moment}
 \int_0^\infty v^{-1/2}\varphi'(v)\dd v=0.
\end{equation}
To construct it, first multiply $v\log v$ by a smooth cutoff equal to
one on $(0,L]$. Then add a multiple of a nonnegative smooth function,
not identically zero, with compact support in $(L,\infty)$. Integration
by parts shows that the left side of \eqref{eq:pa-direction-moment}, with
this last function in place of $\varphi$, is half its integral against
$v^{-3/2}\dd v$, which is strictly positive. A suitable multiple therefore
enforces \eqref{eq:pa-direction-moment}.
Put $g_\varphi(z)=\varphi'(z^{-\alpha})$.
The change of variable $v=z^{-\alpha}$ gives
\[
 \|g_\varphi\|_{L^2(\nu)}^2
     =\int_0^\infty|\varphi'(v)|^2\dd v<\infty,
 \qquad
 \langle f_{t^{-1/\alpha}},g_\varphi\rangle
     =\int_0^t\varphi'(v)\dd v=\varphi(t).
\]
Moreover $\mathsf U_\alpha g_\varphi$ is integrable with a finite first
absolute moment: after the change of variable $v=e^{-\alpha u}$, the
integrand of $\int_\R(1+|u|)|\mathsf U_\alpha g_\varphi(u)|\,du$ is
bounded near $v=0$ by $Cv^{-1/2}(1+|\log v|)^2$ and vanishes for large
$v$. By \eqref{eq:pa-direction-moment},
$\int_\R\mathsf U_\alpha g_\varphi(u)\,du=0$, so the Fourier transform
of $\mathsf U_\alpha g_\varphi$ is $O(|\xi|)$ as $\xi\to0$. This bound,
\eqref{eq:pa-multiplier-lower}, and Plancherel's identity imply
\[
 \int_\R
 \frac{|\widehat{\mathsf U_\alpha g_\varphi}(\xi)|^2}
      {a_\alpha(\xi)}\dd\xi<\infty.
\]
Thus $\mathsf A_\alpha^{-1/2}g_\varphi$ is a nonzero element of
$L^2(\nu)$.

Let $\mathbb W_\nu$ be an isonormal Gaussian process on $L^2(\nu)$ and
realize $Y(t)=\mathbb W_\nu(\mathsf A_\alpha^{1/2}f_{t^{-1/\alpha}})$,
which has covariance \eqref{U-eq:Y-covariance}. Define
\[
 Z_L=\frac{\mathbb W_\nu(\mathsf A_\alpha^{-1/2}g_\varphi)}
          {\|\mathsf A_\alpha^{-1/2}g_\varphi\|^2},
 \qquad Y^\circ(t)=Y(t)-Z_L\varphi(t).
\]
Then $\Var Z_L=\|\mathsf A_\alpha^{-1/2}g_\varphi\|^{-2}>0$ and
$\Cov\{Y(t),Z_L\}=\langle f_{t^{-1/\alpha}},g_\varphi\rangle\Var Z_L
=\varphi(t)\Var Z_L$, so the jointly Gaussian $Y^\circ$ and $Z_L$ are
uncorrelated and hence independent.
On $[0,L]$, $\varphi(t)=t\log t$, which proves
\eqref{eq:pa-Gaussian-direction}. The path assertion follows from
Proposition~\ref{U-prop:positive-covariance} and
$t^{1-\iota}\log t\to0$ as $t\downarrow0$.
\end{proof}

\subsection{Continuity of \texorpdfstring{$F_{\alpha,m}$}{F(alpha,m)}}
Lemma~\ref{lem:pa-Gaussian-direction} makes $\mathcal E_{\alpha,m}(t)$,
evaluated at the minimizing time of $\mathcal D_{\rm lat}$ over a compact
interval of $t$, conditionally normal, and bounds on
$\Prob\{T\notin[L^{-1},L]\}$ extend absolute continuity to
$\mathcal E_{\alpha,m}=\mathcal E_{\alpha,m}(T)$.

\begin{proposition}\label{prop:pa-selected-distribution}
For every fixed $m\ge1$ and $\alpha>1$, the law of
$\mathcal E_{\alpha,m}$ is absolutely continuous with respect to Lebesgue
measure. For every deterministic sequence $\alpha_j\to\alpha>1$,
\begin{equation}\label{eq:pa-selected-cdf-continuity}
 \sup_{z\in\R}|F_{\alpha_j,m}(z)-F_{\alpha,m}(z)|\longrightarrow0.
\end{equation}
\end{proposition}

\begin{proof}
First restrict the minimization to $I=[u,L]\subset(0,\infty)$, and let
$T_I$ be the smallest minimizing time of $\mathcal D_{\rm lat}$ on
$I\times[0,1]$. For the function $\varphi$ in the proof of
Lemma~\ref{lem:pa-Gaussian-direction}, which equals $v\log v$ on
$(0,L]$, direct integration gives
\begin{equation}\label{eq:pa-direction-identities}
 \mathcal L_t\varphi(v)=0\quad(0\le v\le t\le L),\qquad
 \int_0^t\varphi(v)\frac{\dd v}{v}-\varphi(t)=-t\quad(0<t\le L).
\end{equation}
The term added to $Y$ in \eqref{U-eq:shifted-path} to form $X$ is
deterministic. Hence, under the decomposition
\eqref{eq:pa-Gaussian-direction}, \eqref{U-eq:criterion-transformed} and
\eqref{eq:pa-direction-identities} show that the criterion on
$I\times[0,1]$, and therefore $T_I$, depend only on $Y^\circ$.
By \eqref{eq:pa-error-time} and \eqref{eq:pa-direction-identities},
$\mathcal E_{\alpha,m}(T_I)$ is a function of $Y^\circ$ plus
$\alpha Z_L\sqrt{T_I}/\sigma(\alpha,m)$. Conditional on $Y^\circ$, it is
therefore normal, with variance
\begin{equation}\label{eq:pa-conditional-variance}
 \frac{\alpha^2T_I}{\sigma(\alpha,m)^2}\Var Z_L>0,
\end{equation}
and so its law assigns probability zero to every Lebesgue-null set.
On the event $T\in I$, the unique global minimizer is also the
minimizer on $I\times[0,1]$, so $T_I=T$ and
$\mathcal E_{\alpha,m}(T_I)=\mathcal E_{\alpha,m}$. For a Lebesgue-null
set $A$, this gives
$\Prob\{\mathcal E_{\alpha,m}\in A\}\le\Prob\{T\notin I\}$.
Taking $I=[L^{-1},L]$ and then $L\to\infty$ proves absolute continuity.

We next prove continuity in $\alpha$, with a bound on
$\Prob\{T\notin[L^{-1},L]\}$ uniform in $\alpha$. Fix
$1<\alpha_-<\alpha_+<\infty$ and let $\alpha$ range over
$[\alpha_-,\alpha_+]$. The covariance of $Y$ is continuous in $\alpha$,
by \arxiveqref{M-eq:covariance}, and \eqref{U-eq:Y-increment} and Gaussian
moments give
\[
 \E|Y(v)-Y(u)|^p\le C_p|v-u|^{p/2}
\]
for every even integer $p\ge2$, uniformly in $\alpha$.
For fixed $0<\iota<1/2$, choose $\iota<\iota'<1/2$ and then choose $p$
large enough that $\iota'<1/2-1/p$. The moment form of the Kolmogorov
continuity theorem gives
\begin{equation}\label{eq:pa-uniform-holder}
 \sup_{\alpha\in[\alpha_-,\alpha_+]}\Prob\left\{
    \sup_{0\le u<v\le L}
       \frac{|Y(v)-Y(u)|}{|v-u|^{\iota'}}>R\right\}
 \le C_{p,\iota',L}R^{-p},\qquad R>0,
\end{equation}
with a constant uniform in $\alpha$ because $C_p$ is. By the
Arzel\`a--Ascoli theorem, the functions on $[0,L]$ that vanish at zero
and have the inner supremum in \eqref{eq:pa-uniform-holder} at most $R$
form a relatively compact subset of $E_\iota(L)$. The laws of $Y$ are
therefore tight in that space, and, their finite-dimensional laws being
Gaussian with covariances continuous in $\alpha$, the law of $Y$ on
$E_\iota(L)$ is weakly continuous in $\alpha$. The maps
$y\mapsto\int_0^ty(v)\,\dd v/v$, $0<t\le L$, are continuous on
$E_\iota(L)$, uniformly in $t$, since
\begin{equation}\label{eq:pa-integral-continuity}
 \sup_{0<t\le L}\left|\int_0^t y(v)\frac{\dd v}{v}\right|
 \le L^\iota\|y\|_\iota/\iota.
\end{equation}

The bounds \eqref{eq:pa-large-t} and \eqref{eq:pa-small-t} hold
uniformly for $\alpha\in[\alpha_-,\alpha_+]$. For $t\ge1$,
\eqref{eq:pa-large-t} gives
$\mathcal D_{\rm lat}(t,\theta)\ge c_*t^{1/\alpha_+}$, where $c_*>0$ is
the minimum of $\alpha/(2c)$ over $[\alpha_-,\alpha_+]$. The coefficients
$e_2,e_3$ in \eqref{eq:pa-contrast-coefficients} and the variance
$\sigma_e^2$ in \eqref{eq:pa-contrast-variance} depend continuously on
$\alpha$, and $\sigma_e^2>0$ for every $\alpha$ by
Proposition~\ref{U-prop:positive-covariance}, because
$\sum_je_jf_{r_j}$ equals one on $(r_1,r_2)$. Compactness therefore
bounds $C_e$ in \eqref{eq:pa-small-t} by one constant $C$. Summing
\eqref{eq:pa-small-t} over $b=L^{-1}/2,L^{-1}/4,\ldots$ and enlarging
$C$ bounds by $CzL^{-1/2}$ the probability that
$\inf_{0\le\theta\le1}\mathcal D_{\rm lat}(t,\theta)\le z$ for some
$t\le L^{-1}$. Since
$\mathcal D_{\rm lat}(T,\Theta)\le\mathcal D_{\rm lat}(1,1/2)$, for
$L\ge1$ and $z>0$,
\begin{equation}\label{eq:pa-parameter-localization}
\begin{split}
 \sup_{\alpha\in[\alpha_-,\alpha_+]}\Prob\{T\notin[L^{-1},L]\}
 &\le \sup_{\alpha\in[\alpha_-,\alpha_+]}
      \Prob\{\mathcal D_{\rm lat}(1,1/2)>z\}\\
 &\quad+CzL^{-1/2}
      +\mathbf1\{c_*L^{1/\alpha_+}\le z\}.
\end{split}
\end{equation}
The laws of $\mathcal D_{\rm lat}(1,1/2)$, $\alpha\in[\alpha_-,\alpha_+]$,
are tight by \eqref{eq:pa-uniform-holder},
\eqref{eq:pa-integral-continuity}, and continuity of $c,Q,B$ in $\alpha$,
so taking $L\to\infty$ and then $z\to\infty$ in
\eqref{eq:pa-parameter-localization} gives
\begin{equation}\label{eq:pa-uniform-localization}
 \lim_{L\to\infty}\sup_{\alpha\in[\alpha_-,\alpha_+]}
 \Prob\{T\notin[L^{-1},L]\}=0.
\end{equation}

Now let $\alpha_j\to\alpha$ be a deterministic sequence in
$[\alpha_-,\alpha_+]$. The path convergence just proved holds jointly
on every interval $[0,L]$, using the norms of $E_\iota(L)$ for integer
$L\ge1$. Together with \eqref{eq:pa-uniform-localization}, it gives
tightness of the paths and the selected pairs jointly. On a Skorokhod
representation of any convergent subsequence, the criteria converge
uniformly on each compact subset of $(0,\infty)\times[0,1]$.
The limiting selected time lies in $(0,\infty)$ almost surely by
\eqref{eq:pa-uniform-localization}.

For each fixed rational pair $(s,\theta)$ with $s>0$ and
$0\le\theta\le1$, the selected criterion is no larger than the criterion
at $(s,\theta)$. Compact uniform convergence passes this inequality
to the limit. The rational pairs are countable, so their inequalities
hold simultaneously almost surely. Continuity extends them to every
pair, and Proposition~\ref{U-prop:unique-global} identifies the limiting
selected pair as the unique global minimizer at $\alpha$.
Equation~\eqref{eq:pa-integral-continuity} then passes
\eqref{eq:pa-error-time} to the limit at the selected time. This
identifies every subsequential limit, so
$\mathcal E_{\alpha_j,m}\wto\mathcal E_{\alpha,m}$ in $\R$.
Since $F_{\alpha,m}$ is continuous by the first part of the proof,
P\'olya's theorem~\cite[Lemma~2.11]{vanDerVaart1998} gives
\eqref{eq:pa-selected-cdf-continuity}.
\end{proof}

\subsection{Proof of Theorem~\ref{thm:pa-unrestricted-interval}}
\begin{proof}
Proposition~\ref{prop:pa-selected-distribution} gives continuity of
$F_{\alpha,m}$ and \eqref{eq:pa-selected-cdf-continuity}. By
Corollary~\ref{cor:pa-corrected-selected} and $\wh k\to\infty$ in
probability, $\wt\alpha_{n,\wh k}\to\alpha>1$ in
probability, so the event $\{\wt\alpha_{n,\wh k}\le1\}$, which contains
every invalid fit, has probability tending to zero. On this event, set
$\sqrt{\wh k}(\wt\alpha_{n,\wh k}-\alpha)/\sigma(\wt\alpha_{n,\wh k},m)$
to zero and $F_{\wt\alpha_{n,\wh k},m}$ to $\Phi$.
Corollary~\ref{cor:pa-corrected-selected}, continuity and strict
positivity of $\sigma(\cdot,m)$ at $\alpha$, and Slutsky's lemma give
\begin{equation}\label{eq:pa-selected-standardization}
 \frac{\sqrt{\wh k}(\wt\alpha_{n,\wh k}-\alpha)}
      {\sigma(\wt\alpha_{n,\wh k},m)}
       \wto\mathcal E_{\alpha,m}\quad\hbox{in }\R.
\end{equation}

Because \eqref{eq:pa-selected-cdf-continuity} holds along every
deterministic sequence $\alpha_j\to\alpha$, consistency of
$\wt\alpha_{n,\wh k}$ gives
$\sup_z|F_{\wt\alpha_{n,\wh k},m}(z)-F_{\alpha,m}(z)|\xrightarrow{p}0$
in $\R$, and, since $\Prob\{\wt\alpha_{n,\wh k}>1\}\to1$, condition
\eqref{eq:pa-calibration-accuracy} then gives
$\sup_z|\wh F_n(z)-F_{\alpha,m}(z)|\xrightarrow{p}0$ in $\R$.
Exact calibration is the case $\wh F_n=F_{\wt\alpha_{n,\wh k},m}$.
Let $\wh q_{p,n}=\inf\{z:\wh F_n(z)\ge p\}$.
Fix $p\in(0,1)$ and $0<2\varepsilon<\min(p,1-p)$.
On the event $\sup_z|\wh F_n(z)-F_{\alpha,m}(z)|<\varepsilon$, whose
probability tends to one,
\begin{equation}\label{eq:pa-quantile-comparison}
 \mathfrak q_{p-2\varepsilon}(\alpha,m)
 \le\wh q_{p,n}\le
 \mathfrak q_{p+2\varepsilon}(\alpha,m).
\end{equation}
On $\{\wt\alpha_{n,\wh k}>1\}$, the interval
\arxiveqref{M-eq:unrestricted-interval} with the quantiles $\wh q_{p,n}$
contains $\alpha>1$ if and only if the left side of
\eqref{eq:pa-selected-standardization} lies in
$[\wh q_{\upsilon/2,n},\wh q_{1-\upsilon/2,n}]$. Continuity of
$F_{\alpha,m}$ gives
$F_{\alpha,m}(\mathfrak q_{p\pm2\varepsilon}(\alpha,m))=p\pm2\varepsilon$.
Applying \eqref{eq:pa-quantile-comparison} at $p=\upsilon/2$ and
$p=1-\upsilon/2$, and then \eqref{eq:pa-selected-standardization}, shows
that $\liminf_n\Prob\{\alpha\in C_n(\upsilon)\}\ge1-\upsilon-4\varepsilon$
and $\limsup_n\Prob\{\alpha\in C_n(\upsilon)\}\le1-\upsilon+4\varepsilon$.
Letting $\varepsilon\downarrow0$ proves \eqref{eq:pa-unrestricted-coverage}.
The argument uses only continuity of $F_{\alpha,m}$, which need not be
strictly increasing.
\end{proof}

\subsection{Numerical approximation of \texorpdfstring{$F_{\alpha,m}$}{F(alpha,m)}}
Fix $\alpha>1$ and $m\ge1$. For $L>1$, simulate the values of $Y$ at the
points of a partition of $[0,L]$ that contains zero and has maximum
spacing at most $\Delta$. Use the piecewise linear interpolant of these
values, equal to zero at zero, and integrate it exactly against
$\dd v/v$. Include the deterministic drift of $X$ in
\eqref{U-eq:shifted-path} analytically. Minimize
\eqref{U-eq:criterion-transformed} over grids in
$[L^{-1},L]\times[0,1]$ with maximum spacings at most $\Delta$,
including their endpoints. Minimization over $0\le\theta\le1$ may
alternatively be exact. Evaluate the inner maximum on a grid in
$[0,t]$ that includes zero and $t$ and has maximum spacing at most
$\Delta$. The selected pair need only attain the minimum of the
discretized criterion within $\Delta$. For each of $M$ independent
simulated paths, evaluate $\mathcal E_{\alpha,m}(t)$ in
\eqref{eq:pa-error-time}, with the interpolated path in place of $Y$,
at the selected time $t$. Let $\wh F_{M,L,\Delta}(\,\cdot\,;\alpha,m)$
be the empirical distribution function of the $M$ resulting values.

\begin{proposition}\label{prop:pa-simulation}
For fixed $m\ge1$, $1<\alpha_-<\alpha_+<\infty$, and every
$\varepsilon>0$,
\begin{equation}\label{eq:pa-simulation-consistency}
 \lim_{L\to\infty}\lim_{\Delta\downarrow0}\lim_{M\to\infty}
 \sup_{\alpha\in[\alpha_-,\alpha_+]}
 \Prob\left\{\sup_{z\in\R}
  |\wh F_{M,L,\Delta}(z;\alpha,m)-F_{\alpha,m}(z)|
       >\varepsilon\right\}=0.
\end{equation}
There are deterministic sequences of truncations $L\to\infty$, meshes
$\Delta\downarrow0$, and simulation sizes $M\to\infty$ along which
$\sup_{z\in\R}|\wh F_{M,L,\Delta}(z;\alpha,m)-F_{\alpha,m}(z)|\to0$ in
probability, uniformly in $\alpha$ on each compact subset of
$(1,\infty)$. If, on a probability space carrying the graph and the
simulation draws, the paths are simulated conditionally on the graph at
$\alpha=\wt\alpha_{n,\wh k}$ whenever $\wt\alpha_{n,\wh k}>1$, with
$(L,\Delta,M)$ along such sequences as $n\to\infty$, then
$\wh F_n=\wh F_{M,L,\Delta}(\,\cdot\,;\wt\alpha_{n,\wh k},m)$ satisfies
\eqref{eq:pa-calibration-accuracy}.
\end{proposition}

\begin{proof}
Fix $L$, realize each simulated Gaussian vector as the values of $Y$ at
the points of its partition, and take $\iota'$ from
\eqref{eq:pa-uniform-holder}. For $v\le\Delta$, the path and its linear
interpolant are bounded in absolute value by a constant times
$\Delta^{\iota'-\iota}v^\iota$ times the inner supremum in
\eqref{eq:pa-uniform-holder}, and for $v>\Delta$, $v^{-\iota}$ times
their difference is at most a constant times $\Delta^{\iota'-\iota}$
times that supremum. Hence \eqref{eq:pa-uniform-holder} makes the
interpolants converge to $Y$ in $E_\iota(L)$ in probability as
$\Delta\downarrow0$, uniformly in $\alpha\in[\alpha_-,\alpha_+]$. By
\eqref{eq:pa-integral-continuity} and the joint continuity of the
criterion in $t$, $\theta$, and the path, shown in the proof of
Proposition~\ref{U-prop:unique-global} by setting $v=ts$, the discretized
criteria converge uniformly on $[L^{-1},L]\times[0,1]$ in probability.

Let $T_I$ be the minimizing time of $\mathcal D_{\rm lat}$ on
$I\times[0,1]$ for the deterministic interval $I=[L^{-1},L]$, almost
surely unique by Proposition~\ref{U-prop:unique-global}. The time selected
from the discretized criterion therefore converges in probability to
$T_I$, and since $\mathcal E_{\alpha,m}(t)$ depends on the selected pair
only through $t$, the computed value converges in probability to
$\mathcal E_{\alpha,m}(T_I)$ even if the minimizing $\theta$ is not
unique. By \eqref{eq:pa-conditional-variance} and P\'olya's
theorem~\cite[Lemma~2.11]{vanDerVaart1998}, the distribution function of
the computed value converges to that of $\mathcal E_{\alpha,m}(T_I)$
uniformly in $z$. This convergence is uniform in
$\alpha\in[\alpha_-,\alpha_+]$: otherwise choose parameter values along
which the uniform distance stays bounded below, pass to a convergent
subsequence, and apply at its limit the convergence of the interpolants,
the weak continuity in $\alpha$ of the law of $Y$ from the proof of
Proposition~\ref{prop:pa-selected-distribution},
Proposition~\ref{U-prop:unique-global}, and \eqref{eq:pa-conditional-variance}.

On $\{T\in I\}$, $T_I=T$, so the distribution functions of
$\mathcal E_{\alpha,m}(T_I)$ and $\mathcal E_{\alpha,m}$ differ by at
most $\Prob\{T\notin I\}$, uniformly in $z$, which vanishes as
$L\to\infty$ uniformly in $\alpha$ by \eqref{eq:pa-uniform-localization}.
Finally, the $M$ computed values are independent with a common
distribution function, and the Dvoretzky--Kiefer--Wolfowitz inequality
with Massart's constant~\cite{Massart1990} bounds by
$2e^{-2M\varepsilon^2}$, whatever this distribution function, the
probability that $\wh F_{M,L,\Delta}(\,\cdot\,;\alpha,m)$ differs from
it by more than $\varepsilon$ in the supremum norm. These three bounds
and the triangle inequality prove \eqref{eq:pa-simulation-consistency}.

For each $j\ge2$, \eqref{eq:pa-simulation-consistency} with
$[\alpha_-,\alpha_+]=[1+1/j,j]$ and $\varepsilon=1/j$ provides
$(L,\Delta,M)$ with $L\ge j$, $\Delta\le1/j$, and $M\ge j$ at which the
probability in \eqref{eq:pa-simulation-consistency} is at most $1/j$ for
every $\alpha\in[1+1/j,j]$, which proves the second assertion. For the
third, take $[\alpha_-,\alpha_+]\subset(1,\infty)$ containing $\alpha$
in its interior. Conditionally on the graph, $\wh F_n$ is distributed as
$\wh F_{M,L,\Delta}(\,\cdot\,;\alpha,m)$ with $\wt\alpha_{n,\wh k}$
substituted for $\alpha$, so the bound, uniform over
$[\alpha_-,\alpha_+]$, applies on
$\{\wt\alpha_{n,\wh k}\in[\alpha_-,\alpha_+]\}$, whose probability tends
to one by consistency. On $\{\wt\alpha_{n,\wh k}\le1\}$ the left side of
\eqref{eq:pa-calibration-accuracy} is zero.
\end{proof}

\section{Tail fitting in Norros--Reittu graphs}
\label{app:nr}

Cirkovic, Wang, and Cline~\cite[Theorem~2]{CirkovicWangCline2025}
prove asymptotic normality of the Hill and Pickands estimators in
Norros--Reittu graphs with weight tail index $\alpha>0$.
Their theorem assumes regular variation,
a second-order condition, a lower bound on the weight quantile ratios,
and deterministic ranks satisfying
$k\log^\alpha n/n^{1/(4\alpha+1)}\to0$.
Their comparison of large degrees with latent weights also guides the
argument here. For Pareto weights, we compare centered sums of functions
of degrees and weights. This comparison permits the larger rank
$a_n=n^{2/(\alpha+2)}$, where the integer-degree terms enter the MDSP
criterion, and gives bounds uniform over the unrestricted search.

Use the Norros--Reittu multigraph~\cite{NorrosReittu2006} with the
degree convention of~\cite[Section~2]{CirkovicWangCline2025}.
Fix $\alpha,w_0>0$, and let $W_1,\ldots,W_n$ be iid with
$\Prob(W_i>w)=(w/w_0)^{-\alpha}$ for $w\ge w_0$.
Given these weights, take independent
$A_{ij}\sim\operatorname{Poisson}(W_iW_j/\sum_vW_v)$ for $i\le j$.
Set $A_{ji}=A_{ij}$ and $D_i(n)=\sum_jA_{ij}$, so each self-loop
contributes once to the degree. In particular,
\begin{equation}\label{NR-eq:conditional-degree}
 D_i(n)\mid W_1,\ldots,W_n\sim\operatorname{Poisson}(W_i).
\end{equation}
All probabilities in this section average over both the weights and
the edges. The restriction $\alpha>1$ gives finite mean weights, but
the results hold for every $\alpha>0$.

For $2\le h\le n/2$, retain the empirical measure notation of
\eqref{U-eq:smooth-statistic}, now with these degrees and $c=w_0$:
\begin{align}
 \nu_{n,h}=h^{-1}\sum_{i=1}^n
       \delta_{D_i(n)/\{w_0(n/h)^{1/\alpha}\}},\qquad
 S_{n,h}(f)=\sqrt h\{\nu_{n,h}f-\E\nu_{n,h}f\},\notag\\
 \nu(\dd z)=\alpha z^{-\alpha-1}\dd z\quad(z>0).
 \label{NR-eq:measure}
\end{align}
The empirical measures are defined on $[0,\infty)$, including degree zero.
At $h=a_n$, write $b_n=w_0\sqrt{a_n}$, $\nu_n=\nu_{n,a_n}$, and
$\mathbb G_n(x)=\sqrt{a_n}\{\nu_n(x,\infty)-\E\nu_n(x,\infty)\}$.
Use $\mathcal D_k$ and $\wh k$ from \eqref{eq:Dk}, with all $n$ degrees and
the stated convention for invalid ranks, including thresholds equal to
zero. Let $\mathbb W$ be standard Brownian motion, and set
$\mathbb G(x)=\mathbb W(x^{-\alpha})$ for $x>0$.

The next theorem gives the centered process at the MDSP scale and
quantifies the exclusion of ranks outside that scale.
\begin{theorem}\label{NR-thm:process-localization}
For the Norros--Reittu model specified in this section with fixed
$\alpha,w_0>0$, and every $\epsilon>0$,
\begin{equation}\label{NR-eq:process-limit}
 \mathbb G_n\wto\mathbb G
       \quad\hbox{in }\ell^\infty([\epsilon,\infty)).
\end{equation}
The limit has continuous paths and tends to zero at infinity.
There is $C<\infty$, depending only on $\alpha,w_0$, such that,
for all sufficiently large $n$, every $M\ge1$, and every $u>0$,
\begin{equation}\label{NR-eq:log-tail}
 \Prob\left\{\left|\int_M^\infty\mathbb G_n(z)
                         \frac{\dd z}{z}\right|>u\right\}
       \le C u^{-2}M^{-\alpha}.
\end{equation}
Jointly with \eqref{NR-eq:process-limit}, the integrals
$\int_x^\infty\mathbb G_n(z)\dd z/z$ converge in distribution to
$\mathbb L(x)=\int_x^\infty\mathbb G(z)\dd z/z$ in
$\ell^\infty([\epsilon,M_0])$ for every $0<\epsilon<M_0<\infty$.
For every sufficiently large $L$,
\begin{equation}\label{NR-eq:localization}
 \limsup_{n\to\infty}
 \Prob\{\wh k/a_n\notin[L^{-1},L]\}
       \le C\{L^{-1/3}+L^{-2/\alpha}\}.
\end{equation}
\end{theorem}
Thus the centered degree fluctuations have the iid Pareto covariance
at the selected scale. The bound \eqref{NR-eq:localization} places the
unrestricted minimizer at that scale, without assuming alignment of
the degree and weight order statistics. The proof first derives
concentration from the independent edges, then proves
\eqref{NR-eq:process-limit}--\eqref{NR-eq:log-tail} by comparison with
the weights. A quantitative normal approximation supplies the bound
needed to exclude small ranks.

\subsection{Concentration and the marginal degree law}

Each edge enters at most two degrees. Conditional on the weights,
this observation allows the generalized H\"older inequality
in~\cite[Lemma~16]{DuppalaEtAl2025} to compare an exponential moment
with a product over vertices. The exact conditional marginal in
\eqref{NR-eq:conditional-degree} then permits averaging that product
over the independent weights.

\begin{lemma}\label{NR-lem:concentration}
For a bounded measurable function $g$ of a weight and a degree, put
$Z_i=g(W_i,D_i(n))-\E g(W_i,D_i(n))$. For every real $\lambda$,
\begin{equation}\label{NR-eq:mgf}
 \E\exp\left(\lambda\sum_{i=1}^nZ_i\right)
       \le\{\E e^{2\lambda Z_1}\}^{n/2},\qquad
 \Var\left\{\sum_{i=1}^ng(W_i,D_i(n))\right\}
       \le2n\Var\{g(W_1,D_1(n))\}.
\end{equation}
In particular, every Borel set $A\subset[0,\infty)$ and every $u>0$ satisfy
\begin{equation}\label{NR-eq:Bernstein}
 \Prob\left\{\left|\sqrt h\{\nu_{n,h}(A)-\E\nu_{n,h}(A)\}\right|
       >2\sqrt{\E\nu_{n,h}(A)u}+2u/\sqrt h\right\}\le2e^{-u}.
\end{equation}
\end{lemma}
\begin{proof}
Given all the weights, the independent inputs are the variables
$A_{ij}$. Each appears in at most two factors $e^{\lambda Z_i}$,
so the cited inequality gives
\[
 \E\left[\exp\left(\lambda\sum_iZ_i\right)\,
                    \middle|\,W_1,\ldots,W_n\right]
 \le\prod_i\left\{\E[e^{2\lambda Z_i}\mid W_1,\ldots,W_n]\right\}^{1/2}.
\]
By \eqref{NR-eq:conditional-degree}, the $i$th conditional expectation
on the right is a function of $W_i$ alone. Independence of the weights
and Jensen's inequality therefore bound its expectation by
$\prod_i\{\E e^{2\lambda Z_i}\}^{1/2}$, proving the first inequality
in \eqref{NR-eq:mgf}. Comparing the second derivatives at zero gives
the variance inequality. With $b=w_0(n/h)^{1/\alpha}$, apply the exponential bound to
$g(w,j)=\mathbf1\{j/b\in A\}$ and use the bounded-variable Bernstein
bound~\cite[Chapter~2]{BoucheronLugosiMassart2013} to obtain
\eqref{NR-eq:Bernstein}. The factor two in the exponential argument
gives variance parameter $2n\Prob(D_1(n)/b\in A)$ and scale at most
$2$, which give the constants displayed there.
\end{proof}

The marginal law supplies both the centering at integer thresholds
and the size of individual atoms used to exclude large ranks.
Write $\pi_j=\Prob(D_1(n)=j)$ for $j\ge0$; this law does not depend on
$n$ by \eqref{NR-eq:conditional-degree}.
\begin{lemma}\label{NR-lem:marginal}
For integers $d>\alpha-1$,
\begin{equation}\label{NR-eq:gamma-tail}
 0\le w_0^\alpha\frac{\Gamma(d+1-\alpha)}{\Gamma(d+1)}
            -\Prob(D_1(n)>d)
 \le\frac{w_0^{d+1}}{(d+1-\alpha)d!}.
\end{equation}
Consequently,
\begin{equation}\label{NR-eq:tail-expansion}
 \Prob(D_1(n)>d)
 =w_0^\alpha d^{-\alpha}
       \left\{1+\frac{\alpha(\alpha-1)}{2d}+O(d^{-2})\right\}.
\end{equation}
There are positive constants $c_1,C_1$, depending only on
$\alpha,w_0$, such that, for every integer $j\ge1$,
\begin{equation}\label{NR-eq:atom-bounds}
 c_1j^{-\alpha-1}\le\pi_j\le C_1j^{-\alpha-1}.
\end{equation}
For every real $d\ge1$, the marginal pair $(W_1,D_1(n))$ satisfies
\begin{equation}\label{NR-eq:crossing}
 \E\left|\mathbf1\{D_1(n)>d\}-\mathbf1\{W_1>d\}\right|
       \le C_1d^{-\alpha-1/2}.
\end{equation}
\end{lemma}
\begin{proof}
Represent the marginal pair as $(W_1,N(W_1))$, where $N$ is a unit-rate
Poisson process independent of $W_1$. For integer $d$, let $G_{d+1}$
be its $(d+1)$st arrival time, with density $e^{-z}z^d/d!$.
Conditioning on this arrival time gives
\[
 \Prob(D_1(n)>d)=\E\min\{1,(G_{d+1}/w_0)^{-\alpha}\}.
\]
The full negative moment is the Gamma ratio in
\eqref{NR-eq:gamma-tail}. Its excess over the truncated moment is
nonnegative and at most
$w_0^\alpha\int_0^{w_0}z^{d-\alpha}\dd z/d!$, proving that bound.
The Gamma-ratio expansion~\cite[\S5.11(iii)]{DLMF} proves
\eqref{NR-eq:tail-expansion}. Direct mixing also gives
\[
 \pi_j=\frac{\alpha w_0^\alpha}{j!}
                \int_{w_0}^\infty e^{-w}w^{j-\alpha-1}\dd w.
\]
Extending the integral to zero adds at most
$\alpha w_0^j/\{(j-\alpha)j!\}$ for $j>\alpha$.
The Gamma-ratio expansion therefore yields
$\pi_j\sim\alpha w_0^\alpha j^{-\alpha-1}$.
Every remaining $\pi_j$ is positive, so adjusting the constants gives
\eqref{NR-eq:atom-bounds} for all $j\ge1$.

For integer $d\ge2w_0$, the probability in
\eqref{NR-eq:crossing} equals
$\E|\Prob(W_1>G_{d+1}\mid G_{d+1})-\Prob(W_1>d)|$.
On $\{G_{d+1}\ge d/2\}$, the derivative of the Pareto survival
function bounds the integrand by
$Cd^{-\alpha-1}|G_{d+1}-d|$, whose expectation is at most
$Cd^{-\alpha-1/2}$ because
$\E(G_{d+1}-d)^2=d+2$.
The complementary probability is at most
$e^{d/2}2^{-d-1}$, by Markov's inequality applied to $e^{-G_{d+1}}$.
For real $d$, replace the degree threshold by $\lfloor d\rfloor$.
The added probability that $W_1$ lies between these two thresholds
is at most $Cd^{-\alpha-1}$. Increasing $C_1$ covers the bounded
values of $d$ and proves \eqref{NR-eq:crossing}.
\end{proof}

\subsection{Proof of the centered process limit}

We prove the process and integral assertions of
Theorem~\ref{NR-thm:process-localization}. For each fixed $x>0$,
apply the variance inequality in \eqref{NR-eq:mgf} to
$g(w,j)=\mathbf1\{j>b_nx\}-\mathbf1\{w>b_nx\}$.
Equation~\eqref{NR-eq:crossing} gives
\begin{equation}\label{NR-eq:indicator-comparison}
 \E\left|\mathbb G_n(x)
       -\frac1{\sqrt{a_n}}\sum_{i=1}^n
          [\mathbf1\{W_i>b_nx\}-\Prob(W_i>b_nx)]\right|^2
       \le C_x b_n^{-1/2}.
\end{equation}
Each sum is centered by its own expectation. The terms from integer
rounding and the $d^{-1}$ correction in \eqref{NR-eq:tail-expansion}
have scale $\sqrt{a_n}/b_n=1/w_0$ in the uncentered process, so this
comparison retains the exact degree expectation.
The weights are independent, and the centered summands on the right
have absolute value at most $a_n^{-1/2}$.
For every finite linear combination, the Lindeberg condition holds,
and its covariance converges to $(x\vee y)^{-\alpha}$ because
\[
 \frac n{a_n}\Prob(W_1>b_n(x\vee y))=(x\vee y)^{-\alpha},\qquad
 \frac n{a_n}\Prob(W_1>b_nx)\Prob(W_1>b_ny)
       =\frac{a_n}{n}x^{-\alpha}y^{-\alpha}\longrightarrow0.
\]
These identities hold once $b_n\min(x,y)\ge w_0$.
The Lindeberg--Feller theorem~\cite[Proposition~2.27]{vanDerVaart1998}
and \eqref{NR-eq:indicator-comparison} prove convergence of the
finite-dimensional distributions to those of $\mathbb G$.

For tightness, put $\mu_n=\E\nu_n$.
Equations~\eqref{NR-eq:tail-expansion}--\eqref{NR-eq:atom-bounds}
give, for fixed $0<\epsilon<M<\infty$ and all sufficiently large $n$,
\begin{equation}\label{NR-eq:mean-bounds}
 \mu_n(x,\infty)\le Cx^{-\alpha}\quad(x\ge\epsilon),\qquad
 \mu_n((x,y])\le C_\epsilon(y-x+b_n^{-1})
                      \quad(\epsilon\le x\le y\le M).
\end{equation}
The second bound sums at most $b_n(y-x)+1$ atoms, each bounded by
$C_\epsilon/b_n$. The Bernstein inequality
\eqref{NR-eq:Bernstein} has precisely the form required by
\eqref{eq:nested-maximal}--\eqref{eq:nested-local-maximal}.
Partition $[\epsilon,M]$ into $N\le C_{\epsilon,M}/\delta$ intervals
of length at most $\delta$, with all but the last of length $\delta$.
Any interval with endpoints in $[\epsilon,M]$ at distance at most
$\delta$ meets at most three consecutive partition intervals.
Within each, its centered measure is the difference of two nested
tails. Applying \eqref{eq:nested-local-maximal} therefore gives
\begin{align}
 \E\sup_{\substack{x,y\in[\epsilon,M]\\|x-y|\le\delta}}
           |\mathbb G_n(x)-\mathbb G_n(y)|
 &\le C_{\epsilon,M}\left[
   \sqrt{(\delta+b_n^{-1})\log(2N)}
   +\frac{\log(a_n+2)\{\log(a_n+2)+\log(2N)\}}{\sqrt{a_n}}
                         \right].\label{NR-eq:modulus}
\end{align}
For $0<\delta<\min(1,M-\epsilon)$, the limit superior of this bound
tends to zero as $\delta\downarrow0$.
Applying \eqref{eq:nested-maximal} to $(M,\infty)$ also gives
\begin{equation}\label{NR-eq:tail-maximal}
 \E\sup_{x\ge M}|\mathbb G_n(x)|
       \le C\left\{M^{-\alpha/2}
                   +\frac{\log^2(a_n+2)}{\sqrt{a_n}}\right\},\qquad M\ge1.
\end{equation}
These bounds control increments over the entire index set and
the part beyond $M$. Together with finite-dimensional convergence,
they prove \eqref{NR-eq:process-limit} in
$\ell^\infty([\epsilon,\infty))$.

For the logarithmic integral, \eqref{NR-eq:mgf} and
\eqref{NR-eq:mean-bounds} yield
$\|\mathbb G_n(z)\|_2\le Cz^{-\alpha/2}$ for $z\ge1$.
Minkowski's integral inequality gives
\[
 \left\|\int_M^\infty\mathbb G_n(z)\frac{\dd z}{z}\right\|_2
       \le C\int_M^\infty z^{-\alpha/2}\frac{\dd z}{z}
       \le CM^{-\alpha/2}.
\]
Chebyshev's inequality proves \eqref{NR-eq:log-tail}.
For the limit, $\E|\mathbb G(z)|\le z^{-\alpha/2}$, so Tonelli's
theorem proves absolute convergence of $\mathbb L(x)$ almost surely,
simultaneously for all $x>0$ by first taking rational $x$.
Integration over $[x,M]$ is continuous in the uniform norm for
$x$ in a compact subinterval of $(0,M)$. The continuous mapping
theorem, \eqref{NR-eq:log-tail}, and the corresponding bound for
$\mathbb G$ give the asserted joint convergence of the integrals
in $\ell^\infty([\epsilon,M_0])$ for each $0<\epsilon<M_0<\infty$.

The same model verifies all the hypotheses of
Theorem~\ref{thm:general-lattice}, with
\begin{equation}\label{NR-eq:lattice-constants}
 c=w_0,\qquad Q=(1-\alpha)/2,\qquad
 B=\frac1{2(\alpha+1)},\qquad\eta=1,\qquad
 \mathbb G(x)=\mathbb W(x^{-\alpha}).
\end{equation}
To apply that theorem to positive integers, replace $D_i(n)$ by
$\max(D_i(n),1)$ and add $\pi_0$ to $\pi_1$.
Equation~\eqref{NR-eq:tail-expansion} gives the required reference
tail, and $e_j(n)=0$ for every $j$. The replacement leaves
$\mathbb G_n$ on $[\epsilon,\infty)$ unchanged for large $n$.
It also leaves all fitted ranks in any fixed compact range of
$k/a_n$ unchanged with probability tending to one: by
\eqref{NR-eq:Bernstein}, the number of degrees greater than one
exceeds that range except with probability at most $Ce^{-c_1n}$.
Thus \eqref{eq:lattice-uniform-criterion} applies to the original
degree fits with \eqref{NR-eq:lattice-constants}.

This last conclusion supplies a useful quantitative comparison at
one admissible rank:
\begin{equation}\label{NR-eq:benchmark}
 \limsup_n\Prob\{\sqrt{a_n}\mathcal D_{\lfloor a_n\rfloor}>M\}
       \le CM^{-2},\qquad M\ge1.
\end{equation}
Indeed, \eqref{eq:lattice-uniform-criterion} bounds the limiting
criterion there by
$\sup_{(t,\theta)\in[1/2,2]\times[0,1]}\mathcal D_{\rm lat}(t,\theta)$.
The threshold coordinate tends to one by
\eqref{NR-eq:mean-bounds}, \eqref{NR-eq:tail-expansion}, and
\eqref{NR-eq:Bernstein} at fixed thresholds on either side of one.
In \arxiveqref{M-eq:pa-Z}--\arxiveqref{M-eq:pa-D}, the coefficients of $\mathbb G$
and $\mathbb L$ are bounded on this range, uniformly over $y\ge1$,
and the deterministic terms are bounded.
Doob's $L^2$ inequality~\cite{RevuzYor1999} for Brownian motion and
$\int_{2^{-1/\alpha}}^\infty\|\mathbb G(z)\|_2\dd z/z<\infty$
give a finite second moment for this supremum. Markov's inequality,
with a threshold $M/2$ to absorb the remainder tending to zero in
probability, proves \eqref{NR-eq:benchmark}.

\subsection{A quantitative normal approximation}

Excluding small ranks needs a normal approximation whose errors
are summable over dyadic rank intervals. For smooth functions, the
comparison with weights is sharper than
\eqref{NR-eq:indicator-comparison} and gives such a bound.
\begin{proposition}\label{NR-prop:normal}
For every nonzero $f\in C_c^1((0,\infty))$, there is $C_f<\infty$
such that, for every integer $2\le h\le n/2$,
\begin{equation}\label{NR-eq:normal}
 \sup_{z\in\R}\left|
       \Prob\{S_{n,h}(f)\le z\}-\Phi\{z/\sqrt{\nu(f^2)}\}\right|
 \le C_f\left\{h^{-1/2}+(h/n)^{1/(3\alpha)}+h/n\right\}.
\end{equation}
\end{proposition}
\begin{proof}
Write $b=w_0(n/h)^{1/\alpha}$ and compare $S_{n,h}(f)$ with
$S^W_{n,h}(f)=h^{-1/2}\sum_i[f(W_i/b)-\E f(W_i/b)]$.
We first prove
\begin{equation}\label{NR-eq:smooth-comparison}
 \E|f(D_1(n)/b)-f(W_1/b)|^2\le C_fb^{-\alpha-1},\qquad
 \E|S_{n,h}(f)-S^W_{n,h}(f)|^2\le C_f/b.
\end{equation}
Choose $0<\ell<u$ with $\operatorname{supp}f\subset[\ell,u]$,
and extend $f$ by zero to $[0,\infty)$.
On $\ell b/2\le W_1\le2ub$, the Lipschitz bound for $f$ and
\eqref{NR-eq:conditional-degree} give conditional mean square at most
$\|f'\|_\infty^2 W_1/b^2\le C_f/b$.
The probability of this weight interval is at most $C_fb^{-\alpha}$.
For $W_1<\ell b/2$, the difference can be nonzero only if
$D_1(n)\ge\ell b$, whose conditional probability is at most
$e^{-b/C_f}$ by the Poisson Chernoff bound~\cite[Chapter~2]{BoucheronLugosiMassart2013}.
For $W_1>2ub$, it requires $D_1(n)\le ub$, whose conditional
probability is at most $e^{-W_1/C_f}$.
Both contributions are at most $C_fb^{-\alpha-1}$ for large $b$.
Increasing $C_f$ handles bounded $b\ge w_0 2^{1/\alpha}$ and proves
the first inequality in \eqref{NR-eq:smooth-comparison}.
Apply \eqref{NR-eq:mgf} to $g(w,j)=f(j/b)-f(w/b)$ and use
$(n/h)b^{-\alpha}=w_0^{-\alpha}$ to obtain the second.

For large $n/h$, the Pareto density gives the exact moment identities
$(n/h)\E f(W_1/b)=\nu f$ and
$(n/h)\E|f(W_1/b)|^p=\nu(|f|^p)$ for $p=2,3$.
Hence
\[
 \Var\{S^W_{n,h}(f)\}=\nu(f^2)-(h/n)(\nu f)^2,
 \qquad
 \frac n{h^{3/2}}\E|f(W_1/b)-\E f(W_1/b)|^3\le C_fh^{-1/2}.
\]
The variance tends to $\nu(f^2)>0$.
The Berry--Esseen inequality, in the form stated
in~\cite[Theorem~0.1]{Vershynin2026}, gives Kolmogorov error
$C_fh^{-1/2}$ after normalizing by this variance.
Replacing the normal variance by $\nu(f^2)$ adds at most $C_fh/n$.
For every $\varepsilon>0$, \eqref{NR-eq:smooth-comparison} bounds
$\Prob(|S_{n,h}(f)-S^W_{n,h}(f)|>\varepsilon)$ by
$C_f/(b\varepsilon^2)$. The normal density is bounded, so
\[
 \sup_z\left|\Prob\{S_{n,h}(f)\le z\}
                -\Phi\{z/\sqrt{\nu(f^2)}\}\right|
 \le C_f\{h^{-1/2}+h/n+\varepsilon+(b\varepsilon^2)^{-1}\}.
\]
Take $\varepsilon=b^{-1/3}$. For bounded $n/h$, enlarge $C_f$
so that the right side of \eqref{NR-eq:normal} exceeds one.
This proves the uniform bound for every stated $n,h$.
\end{proof}

\subsection{Exclusion of small and large ranks}

To exclude small ranks, we use a smooth function whose integrals
against Pareto laws vanish to second order at $\alpha$. This removes
the first-order effect of fitting the exponent in $\mathcal D_k$.
The deterministic reduction in
\eqref{U-eq:fitted-tail}--\eqref{U-eq:contrast-reduction} applies
once the four degree probabilities used there are concentrated.

\begin{lemma}\label{NR-lem:small-discrepancy}
Put $\kappa=\min\{1/4,1/(3\alpha)\}$.
There is $C<\infty$, depending only on $\alpha,w_0$, such that,
for every integer $2\le h\le a_n$ with $2h<n$ and every $0<u\le1$,
\begin{equation}\label{NR-eq:small-discrepancy}
 \Prob\left\{\min_{h\le k\le2h}\sqrt h\,\mathcal D_k\le u\right\}
       \le Cu+C\{h^{-\kappa}+(h/n)^\kappa\}.
\end{equation}
\end{lemma}
\begin{proof}
Use the nonzero function $f(z)=z^\alpha\eta''(\log z)$ from
\eqref{U-eq:orthogonal-function}, with
$\eta\in C_c^\infty((\log5,\log6))$.
For $\nu_\vartheta(\dd z)=\vartheta z^{-\vartheta-1}\dd z$,
\eqref{U-eq:orthogonal-bound} gives
$|\nu_\vartheta f|\le C(\vartheta-\alpha)^2$ on
$[\alpha/2,2\alpha]$, and $\nu(f^2)>0$.
Write $b=w_0(n/h)^{1/\alpha}$.
Equation~\eqref{NR-eq:tail-expansion}, with integer rounding retained,
gives
\begin{equation}\label{NR-eq:mean-tail-rate}
 |\E\nu_{n,h}(z,\infty)-z^{-\alpha}|\le C_z/b
 \quad\hbox{for each fixed }z>0\hbox{ and all sufficiently large }n/h.
\end{equation}
For all large $h,n/h$, \eqref{NR-eq:Bernstein} gives the following
conditions with probability at least $1-Ch^{-3}$:
\begin{enumerate}
\item For every $h\le k\le2h$,
$4^{-1/\alpha}\le D_{(k+1)}(n)/b\le2$, so all these thresholds
are less than $3b$. The limiting means of
$\nu_{n,h}(4^{-1/\alpha},\infty)$ and $\nu_{n,h}(2,\infty)$ are
$4$ and $2^{-\alpha}<1$, respectively. Their separation from $2$
and $1$ and \eqref{NR-eq:Bernstein} give failure probability $Ce^{-c_1h}$.
\item At $z=3,4$,
$|\nu_{n,h}(z,\infty)-z^{-\alpha}|
\le C\{b^{-1}+\sqrt{\log h/h}\}$, which locates the fitted
exponent. Equation~\eqref{NR-eq:Bernstein} with parameter $3\log h$,
together with \eqref{NR-eq:mean-tail-rate}, gives failure probability
$Ch^{-3}$.
\end{enumerate}
If these conditions hold and $\mathcal D_k\le u/\sqrt h$ at a rank
$k\in[h,2h]$, then \eqref{eq:Dk} gives
\eqref{U-eq:fitted-tail} for every $z\ge3$.
Taking the ratio at $z=3,4$ proves \eqref{U-eq:fitted-exponent},
because the logarithm has bounded derivative on fixed neighborhoods
of $3^{-\alpha}$ and $4^{-\alpha}$.
Integrating against $f'$, supported in $(5,6)$, and using
\eqref{U-eq:orthogonal-bound} gives
\[
 \sqrt h\,|\nu_{n,h}f|
 \le C\{u+\sqrt h/b^2+\log h/\sqrt h\}
 \le C\{u+\log h/\sqrt h\}.
\]
Here $h\le a_n$ implies $b\ge w_0\sqrt h$.
Thus $S_{n,h}(f)$ lies in an interval centered at
$-\sqrt h\,\E\nu_{n,h}f$ with length at most
$C\{u+\log h/\sqrt h\}$.
Equation~\eqref{NR-eq:normal} bounds this probability by a constant
times the interval length plus its displayed normal approximation
error, uniformly over the interval center.
The choice of $\kappa$ absorbs $h^{-3}$, $\log h/\sqrt h$, and
every term in \eqref{NR-eq:normal}, proving
\eqref{NR-eq:small-discrepancy}. Increasing $C$ covers bounded
$h$ or $n/h$.
\end{proof}

We finish the proof of Theorem~\ref{NR-thm:process-localization} by
comparing each excluded rank with $\lfloor a_n\rfloor$.
At a valid rank, the empirical survival function in \eqref{eq:Dk}
has a jump of size at least $1/k$ at a ratio strictly greater than
one, and the fitted survival function is continuous. Hence
$\mathcal D_k\ge1/(2k)$, also at invalid ranks.
On $\{\sqrt{a_n}\mathcal D_{\lfloor a_n\rfloor}\le M\}$, minimality therefore
gives $\wh k\ge\sqrt{a_n}/(2M)$.
Cover $[\sqrt{a_n}/(2M),a_n/L]$ by intervals $[h,2h]$ with dyadic
integer $h$, starting between $\sqrt{a_n}/(4M)$ and
$\sqrt{a_n}/(2M)$, and ending with the first interval that covers
$a_n/L$. In \eqref{NR-eq:small-discrepancy}, take
$u=M\sqrt{h/a_n}\le1$ when $L\ge M^2$.
Summing its three geometric sequences gives, for fixed $M,L$ and
all sufficiently large $n$,
\begin{align}
 \Prob\{\wh k<a_n/L\}
 &\le\Prob\{\sqrt{a_n}\mathcal D_{\lfloor a_n\rfloor}>M\}
       +CM/\sqrt L
       +C_Ma_n^{-\kappa/2}+C(a_n/n)^\kappa.
 \label{NR-eq:lower-localization}
\end{align}
After taking the limit superior in $n$, choose $M=L^{1/6}$ in
\eqref{NR-eq:benchmark} and \eqref{NR-eq:lower-localization} to obtain
$\limsup_n\Prob(\wh k<a_n/L)\le CL^{-1/3}$.

For large ranks, use $N_{n,\ge d}=N_{n,>d}+N_{n,=d}$, with these
counts defined in Section~\ref{sec:pa-interpretation}.
Equation~\eqref{NR-eq:atom-bounds} gives, uniformly over integers
$1\le d\le b_n$,
\[
 \E N_{n,=d+1}\ge c_1n(d+1)^{-\alpha-1},\qquad
 c_1nd^{-\alpha}\le\E N_{n,\ge d}\le C_1nd^{-\alpha}.
\]
The first lower bound is at least a constant times $\sqrt{a_n}$.
Apply \eqref{NR-eq:Bernstein} to these counts, or equivalently the
indicator case of \eqref{NR-eq:mgf}, with deviations equal to fixed
fractions of their means. A union bound over $d$ shows that
\begin{equation}\label{NR-eq:large-rank-atoms}
 N_{n,=d+1}\ge c_1n(d+1)^{-\alpha-1},\qquad
 N_{n,\ge d}\le C_1nd^{-\alpha}\quad(1\le d\le b_n)
\end{equation}
holds except with probability $Cb_n\exp(-c_1\sqrt{a_n})$, after
adjusting the constants.
At a valid rank with threshold $d$ in this range, $k\le N_{n,\ge d}$.
The jump at $(d+1)/d$ in \eqref{eq:Dk} therefore gives
\begin{equation}\label{NR-eq:large-rank-distance}
 \mathcal D_k\ge\frac{N_{n,=d+1}}{2k}\ge\frac{c_1}{d+1}.
\end{equation}

For fixed large $L$, set
$d_{n,L}=\lceil b_n(2/L)^{1/\alpha}\rceil$.
Equation~\eqref{NR-eq:tail-expansion} gives
$\E N_{n,>d_{n,L}}/a_n\to L/2$.
For sufficiently large $n$, \eqref{NR-eq:Bernstein} bounds
$\Prob(N_{n,>d_{n,L}}\ge La_n)$ by $2e^{-c_1La_n}$.
Outside this event, every valid $k\ge La_n$ has a threshold between
$1$ and $d_{n,L}\le b_n$. Equations
\eqref{NR-eq:large-rank-atoms}--\eqref{NR-eq:large-rank-distance}
then imply $\sqrt{a_n}\mathcal D_k\ge c_1L^{1/\alpha}$.
Invalid candidates have infinite criterion. Consequently,
\begin{align}
 \Prob\{\wh k>La_n\}
 &\le\Prob\{\sqrt{a_n}\mathcal D_{\lfloor a_n\rfloor}
                                      \ge c_1L^{1/\alpha}\}
       +Cb_n e^{-c_1\sqrt{a_n}}+2e^{-c_1La_n}.
 \label{NR-eq:upper-localization}
\end{align}
Equation~\eqref{NR-eq:benchmark}, with half the threshold if needed,
gives $\limsup_n\Prob(\wh k>La_n)\le CL^{-2/\alpha}$.
Combining this with \eqref{NR-eq:lower-localization} proves
\eqref{NR-eq:localization}.\hfill$\square$

\subsection{Uniqueness and the joint limit at the selected rank}
\label{NR-sec:selected-law}

The localization in Theorem~\ref{NR-thm:process-localization} reduces
the unrestricted search to compact ranges of $k/a_n$.
To identify the selected pair, use $\mathcal A$ and
$\mathcal D_{\rm lat}$ in \arxiveqref{M-eq:pa-A}--\arxiveqref{M-eq:pa-D} with
the constants in \eqref{NR-eq:lattice-constants} and
$\mathbb G(x)=\mathbb W(x^{-\alpha})$.
The tie fraction $\theta_{n,k}$ is given by
\eqref{eq:pa-tie-fraction}, now for the Norros--Reittu degrees.
Substituting $\wh\alpha_{n,k}$ for $\alpha$ in
$B=1/\{2(\alpha+1)\}$ defines $\wh B_{n,k}$ at every rank.
At a valid rank, subtract the resulting term from the Hill statistic:
\begin{equation}\label{NR-eq:corrected-estimator}
 \wh B_{n,k}=\frac1{2(1+\wh\alpha_{n,k})},\qquad
 \wt H_{n,k}=H_{n,k}
       -\frac{\wh B_{n,k}-\theta_{n,k}}{D_{(k+1)}(n)},\qquad
 \wt\alpha_{n,k}=\wt H_{n,k}^{-1}.
\end{equation}
Set $\wt H_{n,k}=\wt\alpha_{n,k}=0$ at invalid ranks, and
$\wt\alpha_{n,k}=0$ at valid ranks with $\wt H_{n,k}\le0$.
The selector $\wh k$ continues to minimize the
criterion \eqref{eq:Dk} formed with the uncorrected exponent.
In the joint vector in \eqref{NR-eq:selected-joint}, take the
criterion coordinate to be zero if every rank is invalid, so that
the vector is real-valued on every sample.

The next theorem gives a unique limiting pair and the joint law of
the quantities reported at the selected rank.
\begin{theorem}\label{NR-thm:selected-joint}
For the Norros--Reittu model specified in Sec.~\ref{app:nr},
with fixed $\alpha,w_0>0$, almost surely
$\mathcal D_{\rm lat}$ has exactly one global minimizer
$(T,\Theta)\in(0,\infty)\times[0,1]$.
With the correction \eqref{NR-eq:corrected-estimator},
\begin{equation}\label{NR-eq:selected-joint}
 \begin{pmatrix}
       \wh k/a_n\\
       \theta_{n,\wh k}\\
       D_{(\wh k+1)}(n)/b_n\\
       \sqrt{a_n}(\wh\alpha_{n,\wh k}-\alpha)\\
       \sqrt{a_n}\mathcal D_{\wh k}\\
       \sqrt{\wh k}(\wt\alpha_{n,\wh k}-\alpha)
 \end{pmatrix}
 \wto
 \begin{pmatrix}
       T\\
       \Theta\\
       T^{-1/\alpha}\\
       -\alpha^2\mathcal A(T,\Theta)\\
       \mathcal D_{\rm lat}(T,\Theta)\\
       -\dfrac\alpha{\sqrt T}
          \left\{\displaystyle\int_0^T\frac{\mathbb W(v)}v\dd v-\mathbb W(T)\right\}
 \end{pmatrix}
       \quad\hbox{in }\R^6.
\end{equation}
In particular, $\wh k/n\xrightarrow{p}0$ and
$(\wh\alpha_{n,\wh k},\wt\alpha_{n,\wh k})
\xrightarrow{p}(\alpha,\alpha)$ in $\R$ and $\R^2$, respectively.
The probability that the selected rank is valid and
$\wt H_{n,\wh k}>0$ tends to one.
\end{theorem}
The uncorrected exponent retains the shift
$-\alpha^2(B-\Theta)T^{1/\alpha}/w_0$ through
$\mathcal A(T,\Theta)$. The correction removes this term.
Its remaining fluctuation evaluates the Brownian path at the time
$T$ chosen from that same path.

\begin{proof}
We first verify attainment and the two Gaussian hypotheses of
Proposition~\ref{U-prop:unique-global}. We then use the uniform
criterion expansion to identify every subsequential limit of the
selected pair, and subtract the correction in
\eqref{NR-eq:corrected-estimator}.

The criterion tends to infinity at both ends of the time range.
At large times, \eqref{eq:pa-large-t}, which follows by evaluating
the residual at $y=1$, gives
$\mathcal D_{\rm lat}(t,\theta)\ge\alpha t^{1/\alpha}/(2w_0)$.
For small times, choose $1=r_1<r_2<r_3$ and coefficients $e_j$
satisfying \eqref{eq:pa-contrast-coefficients}, with $e_1=1$.
The variance
\[
 \sigma_e^2=\int_0^\infty
       \left(\sum_{j=1}^3e_j\mathbf1\{z>r_j\}\right)^2\nu(\dd z)
       >0
\]
is positive because the integrand is one on $(1,r_2)$.
For $b\le t\le2b$, take $z_j=b^{-1/\alpha}r_j$ and
$y=z_jt^{1/\alpha}$ in \arxiveqref{M-eq:pa-R}, with $\psi=0$.
The two identities in \eqref{eq:pa-contrast-coefficients} cancel
the coefficients of $z_j^{-\alpha}$ and
$z_j^{-\alpha}\log z_j$. Consequently,
\[
 \left|\sum_{j=1}^3e_j
       \left\{\mathbb G(z_j)-\frac{\alpha Q}{w_0}z_j^{-\alpha-1}\right\}\right|
 \le2b\,\mathcal D_{\rm lat}(t,\theta)\sum_{j=1}^3|e_j|.
\]
The random variable inside the absolute value is Gaussian with
variance $b\sigma_e^2$ and does not depend on $(t,\theta)$.
Its density is bounded by $(2\pi b\sigma_e^2)^{-1/2}$ regardless
of its mean, so
\begin{equation}\label{NR-eq:small-time}
 \Prob\left\{\inf_{\substack{b\le t\le2b\\0\le\theta\le1}}
                \mathcal D_{\rm lat}(t,\theta)\le z\right\}
       \le Cz\sqrt b,\qquad b,z>0.
\end{equation}
For each positive integer $z$, this bound is summable over $b=2^{-j}$.
The Borel--Cantelli lemma and a countable intersection over $z$
give divergence to infinity as $t\downarrow0$, uniformly in $\theta$.
Brownian continuity and the integrability defining $\mathbb L$
make the criterion continuous on compact time intervals: in
\arxiveqref{M-eq:pa-R}, the deterministic terms are bounded by
$Cy^{-\alpha}(1+\log y)$, and
$\sup_{x\ge M}|\mathbb G(x)|\to0$ almost surely as $M\to\infty$.
Thus the global minimum is attained in $(0,\infty)\times[0,1]$.

Here the process $Y(v)=\mathbb G(v^{-1/\alpha})$ is exactly
$\mathbb W(v)$. Its Gaussian increment moments and the Kolmogorov
continuity criterion place it in $E_\iota(T_0)$ almost surely for
every $T_0>0$ and $0<\iota<1/2$.
To verify nondegeneracy in that norm, let $\ell$ be a nonzero
continuous linear functional on $E_\iota(T_0)$.
The Riesz representation theorem gives a nonzero finite signed
measure $\mu$ on $(0,T_0]$ such that
$\ell(y)=\int_{(0,T_0]}v^{-\iota}y(v)\,\mu(dv)$.
Since $\int v^{1/2-\iota}|\mu|(dv)<\infty$, the Brownian isometry
commutes with this integral and gives
\begin{equation}\label{NR-eq:Brownian-nondegeneracy}
 \Var\ell(\mathbb W)
 =\int_0^{T_0}
       \left\{\int_{[s,T_0]}v^{-\iota}\mu(dv)\right\}^2\dd s>0.
\end{equation}
If the variance were zero, the inner integral would vanish for
almost every $s>0$. Its one-sided limits and differences determine
$v^{-\iota}\mu(dv)$ on every $[\epsilon,T_0]$, forcing $\mu=0$.
This proves the strict inequality.
Lemma~\ref{U-lem:Gaussian-roughness} gives the other required
property for Brownian motion, simultaneously at every $t>0$.
In \eqref{U-eq:shifted-path} the deterministic addition is
$\alpha^2v^{1+1/\alpha}/(2w_0)$, which belongs to $E_\iota(T_0)$
and is continuously differentiable for $v>0$.
All hypotheses of Proposition~\ref{U-prop:unique-global} hold,
so the global minimizing pair $(T,\Theta)$ is unique almost surely.
Its minimum is positive by \eqref{eq:pa-large-t}.

We next identify the selected pair jointly with the Gaussian path.
The verification following \eqref{NR-eq:lattice-constants} permits
using \eqref{eq:lattice-uniform-criterion} and
\eqref{eq:lattice-hill-expansion} for these degrees.
On each fixed compact range of $k/a_n$, they approximate the
criterion and Hill statistic uniformly at
$((D_{(k+1)}(n)/b_n)^{-\alpha},\theta_{n,k})$.
They also give
\begin{equation}\label{NR-eq:rank-threshold}
 \sup_{\lceil a_n/L\rceil\le k\le\lfloor La_n\rfloor}
       \left|\frac{k}{a_n}
          -\left(\frac{D_{(k+1)}(n)}{b_n}\right)^{-\alpha}\right|
       \xrightarrow{p}0\quad\hbox{in }\R
\end{equation}
for each fixed $L>1$, with the supremum set to zero if a threshold
in that range is zero, an event whose probability tends to zero.
The joint convergence of the process and its logarithmic integrals,
together with \eqref{NR-eq:localization}, gives joint tightness of
the selected rank and degree threshold, normalized as in
\eqref{NR-eq:selected-joint}, together with the tie fraction,
$\sqrt{a_n}\mathcal D_{\wh k}$, $\sqrt{a_n}(H_{n,\wh k}-1/\alpha)$,
and the processes on compact intervals.
In a joint subsequential limit, write $(T',\Theta')$ for the limit
of $(\wh k/a_n,\theta_{n,\wh k})$. Localization gives
$0<T'<\infty$ almost surely, and the uniform expansions identify
the criterion and Hill coordinates as
$\mathcal D_{\rm lat}(T',\Theta')$ and $\mathcal A(T',\Theta')$.

To prove that $(T',\Theta')$ minimizes globally, fix a rational pair
$(s,\phi)$ with $s>0$ and $\phi\in[0,1]$, and take the integer threshold
$d_n(s)=\lfloor b_ns^{-1/\alpha}\rfloor$.
Use the rank $k_n(s,\phi)$ defined in
\eqref{U-eq:rational-candidate} with the Norros--Reittu counts.
The atom expansion \eqref{eq:lattice-atoms}, the first line of
\eqref{eq:lattice-uniform-expansions}, and the uniform Hill expansion
give
\[
 \frac{N_{n,=d_n(s)}}{\sqrt{a_n}}
       \xrightarrow{p}\frac\alpha{w_0}s^{1+1/\alpha},\qquad
 \left(\frac{k_n(s,\phi)}{a_n},\theta_{n,k_n(s,\phi)},
              H_{n,k_n(s,\phi)}\right)
       \xrightarrow{p}(s,\phi,1/\alpha)
       \quad\hbox{in }\R\hbox{ and }\R^3,
\]
respectively, on events whose probabilities tend to one.
In particular, this rank is defined, lies between $2$ and $n-1$,
and is valid with probability tending to one. Use rank $2$ on the
complement, which preserves the displayed limits.
Global minimality of $\wh k$ gives
$\mathcal D_{\wh k}\le \mathcal D_{k_n(s,\phi)}$ whenever the latter rank is valid.
Passing to the joint limit yields
$\mathcal D_{\rm lat}(T',\Theta')
\le\mathcal D_{\rm lat}(s,\phi)$ almost surely.
Intersecting over rational pairs and using continuity extends this
inequality to all pairs. Uniqueness forces
$(T',\Theta')=(T,\Theta)$. This argument uses global uniqueness
and does not require uniqueness of the tie fraction at a boundary
of a compact search interval.

Every subsequential limit is now identified. In particular,
\eqref{NR-eq:rank-threshold} gives the selected degree limit, and
\eqref{eq:lattice-hill-expansion} gives
\begin{equation}\label{NR-eq:selected-Hill}
 \sqrt{a_n}(H_{n,\wh k}-1/\alpha)
       \wto\mathcal A(T,\Theta)\quad\hbox{in }\R,
\end{equation}
jointly with the selected pair, degree threshold, and criterion.
The compact expansions imply that the selected rank is valid with
probability tending to one. Inverting the Hill statistic at its
positive limit $1/\alpha$ gives
$\sqrt{a_n}(\wh\alpha_{n,\wh k}-\alpha)
\wto-\alpha^2\mathcal A(T,\Theta)$ in $\R$ in this same joint limit.

The correction follows by subtraction in
\eqref{NR-eq:selected-Hill}. Consistency of $\wh\alpha_{n,\wh k}$
gives $\wh B_{n,\wh k}\xrightarrow{p}B$ in $\R$.
Also $\sqrt{a_n}/D_{(\wh k+1)}(n)\wto T^{1/\alpha}/w_0$ in $\R$,
with the ratio taken as zero when its denominator vanishes.
Thus, jointly with all the preceding coordinates,
\begin{align}
 \sqrt{a_n}(\wt H_{n,\wh k}-1/\alpha)
 &\wto\mathcal A(T,\Theta)-\frac{(B-\Theta)T^{1/\alpha}}{w_0}\notag\\
 &=\frac1{\alpha T}
       \left\{\int_0^T\frac{\mathbb W(v)}v\dd v-\mathbb W(T)\right\}
       \quad\hbox{in }\R.
 \label{NR-eq:corrected-Hill-limit}
\end{align}
The equality uses the substitution $v=z^{-\alpha}$ in
$\mathbb L(T^{-1/\alpha})$.
Brownian membership in $E_\iota(T_0)$ makes the integral finite
for every finite $T$ almost surely, because
$\int_0^{T_0}v^{\iota-1}\dd v<\infty$.
Equation~\eqref{NR-eq:corrected-Hill-limit} implies
$\wt H_{n,\wh k}\xrightarrow{p}1/\alpha$ in $\R$, so its positivity
has probability tending to one. Inversion and the jointly converging
factor $\sqrt{\wh k/a_n}$ give the last coordinate of
\eqref{NR-eq:selected-joint}. All these operations use the joint
convergence already proved, and hence establish the full vector limit.
Finally, $a_n/n\to0$ and $a_n\to\infty$ give the stated consistency
limits. The convention for a sample with no valid rank changes
the vector only on an event whose probability tends to zero.
\end{proof}

\section{Observed networks}
\subsection{Degree sequences and threshold ties}\label{sec:unrestricted-corpus}
The analysis of observed degrees uses the saved degree histograms of the
corpus analyzed by Broido and Clauset~\cite{BroidoClauset2019,SFAnalysis}
and four SNAP edge lists~\cite{SNAPCAIDA,SNAPEnron,SNAPAstro,SNAPHepTh}.
Each histogram is expanded according to its recorded degree
multiplicities, and zero degrees are removed before fitting. The saved
exclusion file records $229$ corpus sequences with fewer than $50$
positive degrees, and each of the remaining $3440$ corpus sequences has
a valid unrestricted fit. Percentages over the corpus are proportions of
degree sequences unless stated otherwise. The $3440$ sequences carry
$1806$ distinct graph names in the corpus file names, because a graph
yields several sequences when it is directed (in-, out-, and total
degrees), multiplex (one sequence for each layer), or stored in several
simplifications. The proportion of an event per graph averages the
indicator of the event over the sequences of each graph and then over the
graphs, so that each graph has weight one.
Table~\ref{tab:corpus-weights} reports both weightings.

For each SNAP graph, the analysis preserves all observed vertex labels,
merges repeated unordered pairs, removes loops, and then fits the positive
degrees. For CAIDA it also ignores relationship direction and type.
The number of positive degrees is $n_+$, which replaces $n$ in
\eqref{eq:Dk}, so the search is over $2\le k<n_+$. 

For each corpus sequence, the ratio
\eqref{eq:lattice-ratio} at the selected rank,
$\rho_{n_+,\wh k}=\sqrt{\wh k}/D_{(\wh k+1)}$, compares the relative
integer spacing $1/D_{(\wh k+1)}$ with $\wh k^{-1/2}$. Exponent changes
are compared with $\wh\alpha_{n_+,\wh k}/\sqrt{\wh k}$, the asymptotic
standard deviation $\alpha/\sqrt k$ of $\wh\alpha_{n,k}$ for iid Pareto
observations in Section~\ref{sec:pa-correction}, evaluated at $\wh k$
and $\wh\alpha_{n_+,\wh k}$.

To measure the change of the exponent estimate across the selected tie,
write $d=D_{(\wh k+1)}$ and take $N_{n_+,>d}$ and $N_{n_+,=d}$ as in
Section~\ref{sec:pa-interpretation}. When the fits at $k=N_{n_+,>d}$ and
$k=N_{n_+,>d}+N_{n_+,=d}-1$ are valid and $N_{n_+,>d}\ge2$, this change
is
\begin{equation}\label{eq:observed-tie-change}
 \wh\alpha_{n_+,N_{n_+,>d}+N_{n_+,=d}-1}
 -\wh\alpha_{n_+,N_{n_+,>d}}.
\end{equation}
Every retained corpus sequence has a defined value of
\eqref{eq:observed-tie-change}, and the main article reports the
fraction of fits in which it exceeds $\wh\alpha_{n_+,\wh k}/\sqrt{\wh k}$.
The summaries also record $\theta_{n_+,\wh k}$, $\rho_{n_+,\wh k}$, and
$\wh k/n_+$. The fits to observed degree sequences use neither
\arxiveqref{M-eq:corrected-estimator} nor \arxiveqref{M-eq:unrestricted-interval}.

For each threshold $d$ such that some vertex has degree $d+1$, the
discrete fit compared in the main article fits the probabilities
proportional to $j^{-\alpha-1}$ on the integers $j>d$ to the degrees
greater than $d$. Their normalizing constant, a value of the Hurwitz zeta
function, is evaluated with \path|scipy.special.zeta|. The likelihood is
maximized over $\alpha\in\{0.01,0.02,\ldots,9.50\}$, the grid of the
function \texttt{pl} in the code of Ref.~\cite{SFAnalysis}, and the
selected $d$ is the smallest minimizer of the Kolmogorov--Smirnov distance
between the empirical and fitted distribution functions of the degrees
greater than $d$. For $60$ of the $3440$ corpus sequences, the selected
exponent is the largest grid value, $9.50$.

\begin{table}[htbp]
\centering
\caption{Proportions of the $3440$ corpus fits, over sequences and over
the $1806$ graphs from which the sequences come, with each graph weighted
equally. The ends of the selected tie are the estimates
$\widehat\alpha_{n_+,k}$ at $k=N_{n_+,>d}$ and
$k=N_{n_+,>d}+N_{n_+,=d}-1$ with $d=D_{(\widehat k+1)}$, the change across
the selected tie is \eqref{eq:observed-tie-change}, and the discrete
estimate is the fit of Ref.~\cite{ClausetEtAl2009} described in this
section. The last rows give, for each domain of the corpus, the
proportion of fits whose two ends lie on opposite sides of $2$.}
\label{tab:corpus-weights}
\small
\begin{tabular}{lrr}
\toprule
Event & Per sequence (\%) & Per graph (\%)\\
\midrule
$\rho_{n_+,\widehat k}>1$ & 38.0 & 35.9\\
Change across the selected tie exceeds
$\widehat\alpha_{n_+,\widehat k}/\sqrt{\widehat k}$ & 64.9 & 60.6\\
Ends of the selected tie on opposite sides of $2$ & 9.9 & 8.2\\
Discrete and continuous estimates on opposite sides of $2$ & 10.6 & 8.6\\
\midrule
\multicolumn{3}{l}{Ends of the selected tie on opposite sides of
$2$, by domain (sequences, graphs)}\\
Technological ($1951$, $933$) & 5.9 & 2.8\\
Biological ($818$, $473$) & 13.4 & 12.1\\
Social ($470$, $316$) & 18.7 & 15.9\\
Transportation ($171$, $67$) & 14.6 & 19.8\\
Informational ($26$, $16$) & 3.8 & 6.2\\
Economic ($4$, $1$) & 0.0 & 0.0\\
\bottomrule
\end{tabular}
\end{table}

\section{Numerical experiments}\label{app:unrestricted-experiments}
The experiments assess interval coverage, the PA covariance, and predictions of connectivity. Each design and its calibration are specified with its results.

\subsection{Coverage design and numerical calibration}\label{sec:unrestricted-design}
The coverage experiment uses $n\in\{10^4,10^5\}$, $m\in\{1,4\}$, and
$\alpha\in\{1.5,2,2.5,3\}$, with $\delta=m(\alpha-2)$. Each of the
$16$ parameter combinations uses $1000$ independent graphs generated by
\arxiveqref{M-eq:pa-model}. In each graph, $\wh k$ is the smallest minimizer of
$\mathcal D_k$ over all $2\le k<n$, as in \eqref{eq:Dk}, with the
parameter grid and replication counts fixed before simulation.
The two combinations with $m=4$ and $\alpha\in\{2.5,3\}$
were repeated at $n=10^6$, with $1000$ graphs each, the same seeding
scheme, and the calibration grid described next.

Calibration approximates the quantiles $\mathfrak q_{0.025}(\alpha,m)$
and $\mathfrak q_{0.975}(\alpha,m)$ in \arxiveqref{M-eq:selected-cdf}, which
$C_n(\upsilon)$ uses at $\upsilon=0.05$, for $m\in\{1,4\}$ and
$\alpha\in\{1.10,1.15,\ldots,4.50\}$. At each of these $138$ parameter
settings, $20{,}000$ independent Gaussian paths are generated in four
groups of $5000$. In the time coordinate, each path is evaluated at the
union of $641$ geometrically spaced times in $[10^{-6},64]$ and $321$
candidate times in $[1/64,64]$, set to zero at time zero, linearly
interpolated, and integrated against $\dd v/v$. For each path, the
criterion in \eqref{U-eq:criterion-transformed} is minimized over the
candidate times and over $0\le\theta\le1$, with the inner maximum
evaluated at the stored path times no greater than the candidate time,
and $\mathcal E_{\alpha,m}(t)$ in \eqref{eq:pa-error-time} is evaluated
at the minimizing time $t$. Empirical quantiles are
$\inf\{z:\wh F(z)\ge p\}$, as in \arxiveqref{M-eq:selected-cdf}, for the
empirical distribution function $\wh F$ of these values.

For each graph, the two quantiles at $\alpha=\wt\alpha_{n,\wh k}$ are
interpolated linearly between adjacent grid values of $\alpha$, or, in
the nine saved replications with $\wt\alpha_{n,\wh k}>1$ outside the
grid, simulated at $\wt\alpha_{n,\wh k}$ with the same settings. No
replication had $\wt\alpha_{n,\wh k}\le1$, so the convention
$C_n(\upsilon)=(1,\infty)$ of Section~\ref{sec:pa-unrestricted-inference}
was never used.

To check these settings, the calibration was repeated at the six
combinations of $\alpha\in\{1.5,2,3\}$ and $m\in\{1,4\}$ with as many
new paths, truncation $128$, $1281$ geometrically spaced times in
$[10^{-6},128]$, and $641$ candidate times in $[1/128,128]$. Table~\ref{tab:unrestricted-refinement} reports
the refined minus the original simulated quantiles. Because the refined
calibration also uses new paths, these changes include Monte Carlo
variation as well as the effect of the finer discretization. In the
original calibration, no minimizing time equaled an endpoint of
$[1/64,64]$.

\begin{table}[htbp]
\centering
\caption{Changes in the simulated quantiles $\mathfrak q_{0.025}(\alpha,m)$
and $\mathfrak q_{0.975}(\alpha,m)$ when the calibration is repeated with
truncation $128$ and finer grids. Entries are refined minus original
values.}
\label{tab:unrestricted-refinement}
\begin{tabular}{rrrr}
\toprule
$m$ & $\alpha$ & Change at $p=0.025$ & Change at $p=0.975$\\
\midrule
1 & 1.5 & $0.02120$ & $0.01639$\\
1 & 2.0 & $0.00986$ & $-0.01599$\\
1 & 3.0 & $0.00793$ & $0.02308$\\
4 & 1.5 & $-0.01015$ & $0.02827$\\
4 & 2.0 & $0.02971$ & $-0.01218$\\
4 & 3.0 & $0.03556$ & $0.00195$\\
\bottomrule
\end{tabular}
\end{table}

Two further intervals at the same selected rank serve as comparisons:
the normal interval
$\wt\alpha_{n,\wh k}\pm\Phi^{-1}(0.975)\sigma(\wt\alpha_{n,\wh k},m)/\sqrt{\wh k}$
and the Hill interval
$\wh\alpha_{n,\wh k}\pm\Phi^{-1}(0.975)\wh\alpha_{n,\wh k}/\sqrt{\wh k}$.
The second-order interval $C_n^{(2)}(\upsilon)$ substitutes
$\wt\alpha_{n,\wh k}^{(2)}$ from \arxiveqref{M-eq:second-order-estimator}
for $\wt\alpha_{n,\wh k}$ in \arxiveqref{M-eq:unrestricted-interval}, using
the same saved quantile grid. The summaries record, for each setting,
the coverage frequency of each interval (the fraction of graphs for
which it contains $\alpha$), its median length, and the empirical bias
and RMSE of each estimate, which replace the expectations in
\eqref{eq:mdsp-estimation-error} by averages over the graphs. All saved
intervals were finite. The coverage figure in the main article shows
95\% Wilson intervals~\cite{Wilson1927} for the coverage frequencies.

A further check replaces $\wt\alpha_{n,\wh k}$ in
\arxiveqref{M-eq:unrestricted-interval} by the solution $\alpha>1$, at
$k=\wh k$, of
\begin{equation}\label{eq:pa-exact-inversion}
 \frac{\sum_{j>d}\pi_j\log(j/d)}{\sum_{j>d}\pi_j+\theta\pi_d}=H_{n,k},
 \qquad d=D_{(k+1)}(n),\quad\theta=\theta_{n,k},
\end{equation}
whose left side is that of \arxiveqref{M-eq:hill-tie-expansion}, with
$\pi_j$ from \arxiveqref{M-eq:pa-stationary} at exponent $\alpha$ and
$\delta=m(\alpha-2)$. The coverage frequency of the resulting interval
differs from that of $C_n^{(2)}(\upsilon)$ by at most $1.6$ percentage
points in every setting.

\begin{table}[htbp]\centering\scriptsize
\caption{Coverage, median width, and RMSE at the unrestricted selected rank under the PA model, with $\upsilon=0.05$. Each row uses 1000 graphs and the rank $\wh k$ selected over $2\le k<n$. Coverage is the percentage of graphs for which $\alpha$ lies in $C_n(\upsilon)$ from Eq.~\arxiveqref{M-eq:unrestricted-interval}, in the second-order interval $C_n^{(2)}(\upsilon)$, centered at $\wt\alpha_{n,\wh k}^{(2)}$ from Eq.~\arxiveqref{M-eq:second-order-estimator}, in the normal interval $\wt\alpha_{n,\wh k}\pm1.96\,\sigma(\wt\alpha_{n,\wh k},m)\wh k^{-1/2}$, and in the Hill interval $\wh\alpha_{n,\wh k}\pm1.96\,\wh\alpha_{n,\wh k}\wh k^{-1/2}$. Monte Carlo standard errors of the coverage percentages are at most 1.6 percentage points. Widths are medians over graphs, and the RMSE columns compare $\wh\alpha_{n,\wh k}$ (Raw), $\wt\alpha_{n,\wh k}$, and $\wt\alpha_{n,\wh k}^{(2)}$ with $\alpha$.}
\label{tab:unrestricted}
\begin{tabular}{rrr rrrr rr rrr}
\toprule
 & & & \multicolumn{4}{c}{Coverage (\%)} & \multicolumn{2}{c}{Median width} & \multicolumn{3}{c}{RMSE}\\
$n$ & $m$ & $\alpha$ & $C_n(\upsilon)$ & $C_n^{(2)}(\upsilon)$ & Normal & Hill & $C_n(\upsilon)$ & Hill & Raw & $\wt\alpha_{n,\wh k}$ & $\wt\alpha_{n,\wh k}^{(2)}$\\
\hline
$10{,}000$ & 1 & 1.5 & 93.3 & 92.8 & 90.9 & 97.0 & 0.429 & 0.547 & 0.130 & 0.124 & 0.125\\
$10{,}000$ & 1 & 2.0 & 92.2 & 91.5 & 90.8 & 94.5 & 0.610 & 0.798 & 0.197 & 0.175 & 0.176\\
$10{,}000$ & 1 & 2.5 & 92.0 & 91.9 & 91.1 & 91.3 & 0.892 & 1.081 & 0.307 & 0.288 & 0.289\\
$10{,}000$ & 1 & 3.0 & 92.1 & 92.4 & 90.3 & 84.0 & 1.238 & 1.400 & 0.450 & 0.427 & 0.427\\
$10{,}000$ & 4 & 1.5 & 90.6 & 88.4 & 88.9 & 99.4 & 0.216 & 0.350 & 0.067 & 0.076 & 0.080\\
$10{,}000$ & 4 & 2.0 & 89.5 & 89.0 & 89.8 & 98.9 & 0.290 & 0.474 & 0.103 & 0.092 & 0.092\\
$10{,}000$ & 4 & 2.5 & 83.4 & 86.5 & 84.4 & 77.2 & 0.389 & 0.601 & 0.224 & 0.135 & 0.131\\
$10{,}000$ & 4 & 3.0 & 79.4 & 85.3 & 82.1 & 51.7 & 0.521 & 0.745 & 0.382 & 0.192 & 0.188\\
\hline
$100{,}000$ & 1 & 1.5 & 93.3 & 93.2 & 92.1 & 97.4 & 0.225 & 0.291 & 0.069 & 0.064 & 0.065\\
$100{,}000$ & 1 & 2.0 & 93.3 & 93.3 & 92.5 & 95.7 & 0.357 & 0.474 & 0.120 & 0.102 & 0.102\\
$100{,}000$ & 1 & 2.5 & 93.4 & 93.2 & 92.5 & 92.9 & 0.548 & 0.700 & 0.189 & 0.158 & 0.158\\
$100{,}000$ & 1 & 3.0 & 93.5 & 93.6 & 91.9 & 86.4 & 0.796 & 0.949 & 0.293 & 0.226 & 0.226\\
$100{,}000$ & 4 & 1.5 & 92.5 & 92.3 & 91.7 & 99.2 & 0.114 & 0.185 & 0.036 & 0.036 & 0.036\\
$100{,}000$ & 4 & 2.0 & 91.7 & 92.1 & 92.3 & 98.9 & 0.166 & 0.272 & 0.059 & 0.050 & 0.050\\
$100{,}000$ & 4 & 2.5 & 90.5 & 91.8 & 92.0 & 81.6 & 0.238 & 0.374 & 0.132 & 0.072 & 0.071\\
$100{,}000$ & 4 & 3.0 & 84.7 & 89.3 & 86.8 & 53.4 & 0.335 & 0.494 & 0.245 & 0.111 & 0.109\\
\hline
$1{,}000{,}000$ & 4 & 2.5 & 92.6 & 93.5 & 93.7 & 84.2 & 0.146 & 0.233 & 0.079 & 0.042 & 0.041\\
$1{,}000{,}000$ & 4 & 3.0 & 90.1 & 92.3 & 91.0 & 57.6 & 0.213 & 0.321 & 0.155 & 0.065 & 0.065\\
\bottomrule\end{tabular}\end{table}

\subsection{Matched iid comparison}\label{sec:variance-decisions}
For each combination of $n\in\{10^4,10^5\}$, $m\in\{1,4\}$,
and $\alpha\in\{1.8,2.2\}$, the comparison uses $500$ independent
PA graphs and $500$ independent samples from the limiting degree law.

For each graph, an independent sample of size $n$ is drawn from its
limiting PA degree law. Independently for each observation, draw a success
probability from $\operatorname{Beta}(m(\alpha-1),\alpha)$, and take the
degree $D$ to be $m$ plus the number of successes before the first failure
in conditionally independent Bernoulli trials. Integrating the geometric
probabilities against the beta density gives
\[
 \Prob(D=m+j)
 =\frac{\alpha\Gamma(r+\alpha)\Gamma(r+j)}
        {\Gamma(r)\Gamma(r+j+\alpha+1)},\qquad r=m(\alpha-1),
\]
which equals $\pi_{m+j}$ in \arxiveqref{M-eq:pa-stationary}.
Both samples are evaluated at deterministic
$k=\lfloor n^{2/(\alpha+2)}\rfloor$. The recorded quantity is
$\sqrt{k}(\wt\alpha_{n,k}-\alpha)/\alpha$, using the same correction
in \arxiveqref{M-eq:corrected-estimator}. Corollaries~\ref{cor:pa-corrected}
and~\ref{cor:iid-lattice} give limiting variances
$\sigma(\alpha,m)^2/\alpha^2$ for PA and $1$ for iid observations.
Every deterministic-rank estimate was valid. The figures use empirical
variances with divisor $499$ and compare their
ratio with $\sigma(\alpha,m)^2/\alpha^2$. No PA confidence interval is
assigned to an iid sample.

\subsection{Connectivity prediction in PA and Norros--Reittu graphs}
\label{sec:percolation-experiment}\label{sec:nr-percolation-experiment}

\subsubsection{Shared prediction design}
Bond percolation retains an edge when its independent uniform label
$U_e$ is at most $p$. With $L_{n,1}(p)$ the largest component size,
the outcome is
\begin{equation}\label{eq:percolation-crossing}
 \mathcal Q_{n,f}=\inf\{p\in[0,1]:L_{n,1}(p)/n\ge f\},
\end{equation}
set to one if the target is unreached. The primary fraction is
$f=0.10$, with $f=0.05,0.20$ prescribed as sensitivity checks.
This finite-fraction crossing differs from the asymptotic critical
probability for an extensive component~\cite{HazraEtAl2023}.
Inserting edges in label order gives the crossing without a grid in $p$.
Parallel edges receive separate labels; loops do not join components.
Every simulated test, calibration, and validation graph reaches all
three targets.

Both experiments use $n\in\{10^4,10^5\}$,
$\alpha\in\{1.8,2.2,2.6,3.0\}$, and known $m\in\{2,4\}$ for PA or known $w_0\in\{2,4\}$ for Pareto NR.
Each of the $16$ settings per model has $300$ independent training/test
pairs. Only the training degrees enter estimation. Every method
predicts the same independent test crossing within a pair, using a
shared calibration table independent of both graphs. The protocols
fixed settings, replication counts, methods, losses, sensitivity
checks, and handling of invalid fits before simulation.

Raw Hill uses the smallest minimizer of \eqref{eq:Dk} over all valid
ranks. Both PA corrections and the NR correction retain this rank.
The NR comparator uses raw Hill at
$\max\{2,\lfloor\wh k/2\rfloor\}$. Each experiment also evaluates
the calibration table at the true exponent. This reference retains
variation of the test outcome and finite calibration error.

\subsubsection{Graph laws and the PA likelihood}
The PA generator uses \arxiveqref{M-eq:pa-model} with
$\delta=m(\alpha-2)$ and two seed vertices joined by $m$ parallel
edges. Each new vertex attaches $m$ edges sequentially to older
vertices, updating target degrees after each edge. It becomes an
eligible target only after completing its attachments. Loops are
excluded and parallel edges retained.

The full PA likelihood uses all final degrees~\cite{GaoVanDerVaart2017}.
For this seed and attachment rule, its log likelihood, up to a term
independent of $\delta$, is
\begin{align}\label{eq:percolation-likelihood}
 \ell_n(\delta)
 ={}&\sum_{v=1}^n\{\log\Gamma(D_v(n)+\delta)-\log\Gamma(m+\delta)\}
 \notag\\
 &-\sum_{t=3}^n\sum_{i=1}^m
 \log\{(2m+\delta)(t-1)-2m+i-1\}.
\end{align}
The first sum contains the log attachment weights; the second contains
the log total weights before each edge.
The fit searches $\alpha\in[1.0001,20]$, solving
$\ell_n'(\delta)=0$ to absolute tolerance $10^{-10}$ when endpoint
derivatives have opposite signs, and otherwise comparing a bounded
maximizer with both endpoints. No optimum reaches a boundary.

NR uses the Pareto weights and independent conditional Poisson edge
multiplicities of Sec.~\ref{app:nr}. Degrees count a loop once,
and all vertices, including degree-zero vertices, enter the ordered
sample. Zero thresholds and zero Hill statistics make ranks invalid.
The estimator receives no latent weights. The correction uses
$\wh B_{n,\wh k}=1/\{2(1+\wh\alpha_{n,\wh k})\}$.
Poisson splitting generates nonloop edges with means
$W_iW_j/\sum_vW_v$, and independent loops with means
$W_i^2/\sum_vW_v$. Conditional on their number, the nonloop endpoint
pairs are iid. Assigning uniform order statistics to these pairs
therefore gives the same percolation law as sorting independent edge
labels. The implementation stops once all three targets are reached.

\subsubsection{Calibration, loss, and uncertainty}
For each $(n,m)$ or $(n,w_0)$, calibration uses $80$ independent graphs
at each of the $34$ exponents
$\{1.1,1.2,\ldots,4.1\}\cup\{5,6,8\}$.
Linear interpolation of mean crossings converts a fitted exponent
into a prediction, without smoothing. Positive finite estimates are
projected onto $[1.1,8]$; an invalid estimate predicts one and remains
in the analysis. No estimate requires projection or is invalid.

For prediction $\wh{\mathcal Q}_{n,f}^{(b,j)}$ by method $j$ and
independent test crossing $\mathcal Q_{n,f}^{(b)}$, the reported loss is
\begin{align}
 \operatorname{MSE}_j(f)&=\frac1{300}\sum_{b=1}^{300}
 (\wh{\mathcal Q}_{n,f}^{(b,j)}-\mathcal Q_{n,f}^{(b)})^2,
 \label{eq:percolation-loss}\\
 \text{reduction in MSE}&=
 100\{1-\operatorname{MSE}_j(f)/\operatorname{MSE}_{\rm Hill}(f)\}
 \ \text{percent}.\label{eq:percolation-reduction}
\end{align}
Pooled MSE averages the $16$ setting MSEs equally; its reduction uses
the ratio of pooled losses. RMSE is the square root of MSE.
The $2000$ paired bootstrap repetitions~\cite{Efron1979} independently
resample calibration graphs within each grid point and training/test
pairs within each setting. Each resampled calibration table is shared
across methods and settings with the same model parameters. All methods
and fractions remain paired, and predictions are recomputed. Intervals
are the $2.5$ and $97.5$ percentiles. Setting-specific intervals are
pointwise and concern prediction loss, not exponent coverage.

Validation uses $100$ new graphs at each
$\alpha\in\{1.85,2.25,2.65,3.05\}$ for every model and parameter pair.
Across fractions, the largest absolute calibration discrepancy is
$0.268$ percentage points for PA and $0.494$ for NR. Standardizing by
the combined calibration and validation standard error gives maxima
$2.28$ and $2.87$, respectively; $4$ and $5$ of the $48$ comparisons
exceed two in magnitude. Separate predictions from the first and last
$40$ calibration graphs give primary correction reductions of
$26.4$ and $28.7$ percent for PA, and $22.46$ and $22.45$ percent for NR.
These checks leave the prescribed calibration unchanged.

Table~\ref{tab:connectivity-settings} retains every primary correction
comparison, including the adverse PA setting. Table~\ref{tab:connectivity-pooled}
reports every original method and fraction. 
\begin{table}[tbp]
\centering
\small
\caption{Primary connectivity comparisons in every prescribed PA and Pareto NR setting. The second column gives $m$ for PA and $w_0$ for NR. Raw RMSE is in percentage points of edge survival probability. Correction entries are percentage MSE reductions relative to raw Hill, with pointwise 95\% bootstrap intervals. Each entry uses 300 independent training/test pairs and known $m$ or $w_0$.}
\label{tab:connectivity-settings}
\begin{tabular}{rrrrlrl}
\toprule
$n$ & $m,w_0$ & $\alpha$ & PA raw & PA correction & NR raw & NR correction\\
\midrule
$10^{4}$ & 2 & 1.8 & 1.043 & 8.3 [3.3, 12.6] & 2.999 & 21.9 [16.0, 28.0]\\
$10^{4}$ & 2 & 2.2 & 1.056 & 25.7 [18.6, 31.6] & 4.896 & 24.2 [19.4, 29.2]\\
$10^{4}$ & 2 & 2.6 & 1.083 & 32.9 [25.2, 39.5] & 5.424 & 20.5 [14.5, 26.3]\\
$10^{4}$ & 2 & 3.0 & 1.114 & 31.8 [22.4, 39.9] & 5.154 & 23.5 [17.3, 29.8]\\
$10^{4}$ & 4 & 1.8 & 0.379 & -3.8 [-7.3, -0.8] & 1.303 & 19.4 [14.6, 24.1]\\
$10^{4}$ & 4 & 2.2 & 0.420 & 26.5 [20.6, 31.6] & 1.643 & 22.0 [15.7, 27.6]\\
$10^{4}$ & 4 & 2.6 & 0.434 & 37.3 [28.7, 44.6] & 1.937 & 27.9 [21.5, 34.2]\\
$10^{4}$ & 4 & 3.0 & 0.472 & 28.7 [18.1, 37.5] & 2.087 & 15.0 [8.9, 21.3]\\
$10^{5}$ & 2 & 1.8 & 0.473 & 16.4 [11.8, 20.9] & 1.738 & 24.6 [19.8, 29.0]\\
$10^{5}$ & 2 & 2.2 & 0.461 & 29.4 [22.1, 35.4] & 2.739 & 18.5 [12.9, 24.1]\\
$10^{5}$ & 2 & 2.6 & 0.481 & 48.0 [41.6, 53.5] & 3.489 & 26.6 [21.5, 32.1]\\
$10^{5}$ & 2 & 3.0 & 0.556 & 57.0 [51.8, 61.8] & 4.014 & 21.8 [16.2, 27.4]\\
$10^{5}$ & 4 & 1.8 & 0.139 & 3.4 [-1.3, 8.0] & 0.690 & 19.1 [14.7, 23.2]\\
$10^{5}$ & 4 & 2.2 & 0.169 & 33.4 [26.2, 39.1] & 0.930 & 17.9 [11.7, 23.6]\\
$10^{5}$ & 4 & 2.6 & 0.196 & 58.0 [51.5, 63.4] & 1.187 & 25.4 [19.6, 30.8]\\
$10^{5}$ & 4 & 3.0 & 0.215 & 68.3 [62.9, 72.9] & 1.351 & 20.4 [13.3, 27.1]\\
\bottomrule
\end{tabular}
\end{table}

\begin{table}[tbp]
\centering
\small
\caption{Pooled connectivity results for all original methods and target fractions. The fraction columns give RMSE in percentage points. The final column gives percentage MSE reduction [95\% interval] relative to raw Hill at $f=0.10$. Pooling weights the 16 settings equally within each graph model. The true-exponent reference retains test variation and calibration error.}
\label{tab:connectivity-pooled}
\begin{tabular}{llrrrl}
\toprule
Model & Method & $f=0.05$ & $f=0.10$ & $f=0.20$ & Reduction (\%)\\
\midrule
PA & Raw Hill & 0.798 & 0.635 & 0.487 & Reference\\
PA & First correction & 0.702 & 0.540 & 0.407 & 27.6 [24.8, 30.2]\\
PA & Second correction & 0.701 & 0.539 & 0.406 & 27.9 [24.9, 30.7]\\
PA & Full likelihood & 0.589 & 0.425 & 0.316 & 55.3 [51.9, 58.0]\\
PA & True exponent & 0.585 & 0.422 & 0.314 & 55.9 [52.2, 58.6]\\
NR & Raw Hill & 3.215 & 3.009 & 2.769 & Reference\\
NR & Correction & 2.850 & 2.649 & 2.426 & 22.5 [20.3, 24.6]\\
NR & Half rank & 3.560 & 3.354 & 3.096 & -24.3 [-28.5, -20.3]\\
NR & True exponent & 0.917 & 0.687 & 0.500 & 94.8 [94.3, 95.1]\\
\bottomrule
\end{tabular}
\end{table}

\FloatBarrier
\subsection{Estimated scale and departures from Pareto weights}
\label{sec:nr-robustness-experiment}

\subsubsection{Weight laws and fixed comparisons}
The robustness experiment measures how the correction behaves when the
weight scale is estimated and the NR weights depart from an exact Pareto
law. The conditional Poisson edge law is retained, and the weight distribution changes. The Pareto correction is applied without adapting
its coefficient to the generating distribution.

Let $Z$ have survival function $\Prob(Z>z)=z^{-\alpha}$ for $z\ge1$.
The four weight laws are
\begin{equation}\label{eq:nr-robust-weight-laws}
 W\in\left\{w_0Z,\quad w_0(Z+0.2/Z),\quad
                 w_0(Z+0.2),\quad w_0(Z-0.2)\right\}.
\end{equation}
In each setting one of these formulas defines the weights, independently
across vertices. Every modified weight differs from $w_0Z$ by at most
$20$ percent. All four laws have survival exponent $\alpha$ and tail
constant $w_0^\alpha$. The additive shifts change the tail at relative
order $w^{-1}$, because
\[
 \Prob\{w_0(Z\pm0.2)>w\}
 =\left(\frac{w}{w_0}\mp0.2\right)^{-\alpha}
 =\left(\frac{w_0}{w}\right)^\alpha
       \left\{1\pm\frac{0.2\alpha w_0}{w}+O(w^{-2})\right\}.
\]
For $w_0(Z+0.2/Z)$, the inverse transformation at $w/w_0$ equals
$w/w_0-0.2w_0/w+O(w^{-3})$, so the relative tail change starts
at order $w^{-2}$. These choices distinguish a change concentrated
near the lower end from changes in the leading tail correction.

For every weight law, the design uses $n\in\{10^4,10^5\}$,
$w_0\in\{2,4\}$, and $\alpha\in\{1.8,2.2,2.6,3.0\}$,
with $300$ independent training/test pairs per setting.
Within each $(n,w_0,\alpha)$ and replication, the four laws use common
Pareto draws and matching graph and damage random-number streams.
The training and test streams are independent. Each method receives
only training degrees and predicts the same independent test crossing
within that law and pair. Degrees count loops once, and connectivity
ignores loops, as in Sec.~\ref{sec:nr-percolation-experiment}.

The primary outcome is the capped crossing $\mathcal Q_{n,0.10}$.
The primary comparison gives equal weight to the $48$ non-Pareto
settings and compares corrected and raw Hill under estimated scale.
The exact-Pareto controls, each individual weight law, and every
setting are reported separately. Fractions $f=0.05,0.20$, exponent
errors, scale errors, and predictions using the true $w_0$ are
prescribed secondary comparisons. The saved protocol fixes all
choices before simulation, with no outcome-dependent changes or exclusions.

\subsubsection{Exponent and scale estimation}
Raw Hill, its NR correction, and raw Hill at half the selected rank
use the definitions of Sec.~\ref{sec:nr-percolation-experiment}.
The correction retains the unrestricted MDSP rank and the coefficient
in \eqref{NR-eq:corrected-estimator}. The discrete comparator fits
a discrete power-law tail by the conditional likelihood of
Clauset et al.~\cite{ClausetEtAl2009}. At each observed positive
cutoff $x_{\min}$ retaining at least $50$ vertices, it maximizes
the sum of log probabilities
\[
 \log\left\{\frac{j^{-a-1}}
  {\sum_{\ell=x_{\min}}^\infty\ell^{-a-1}}\right\},
 \qquad j\ge x_{\min},\qquad 0.05\le a\le10,
\]
over the retained degrees. It then selects the cutoff with the
smallest KS distance, choosing the smallest cutoff at a tie.
Both sides of every empirical jump enter the distance.
The fit optimizes $a$ continuously and retains each threshold tie
whole. It uses a tail likelihood, without a likelihood for the full
network or a model for the degrees smaller than its cutoff.

For a Pareto weight law with $\alpha>1$,
$\E D_i(n)=\E W_i=w_0\alpha/(\alpha-1)$.
This identity motivates estimating the scale from the mean degree.
For each fitted exponent, first project it onto $[1.1,8]$ and denote
the result by $a$. The corresponding scale estimate is
\begin{equation}\label{eq:nr-estimated-scale}
 \wh w_0(a)=\frac{a-1}{a}\,\frac1n\sum_{i=1}^nD_i(n).
\end{equation}
Every method uses this same rule with its own exponent.
The true-exponent reference substitutes $a=\alpha$ and still
estimates scale from the training graph. It separates exponent
estimation from scale estimation and misspecification of the
Pareto prediction family. It is unavailable in practice and is
not a lower bound on prediction error.

Under the modified weight laws, \eqref{eq:nr-estimated-scale}
matches a Pareto mean to the observed mean degree. Even with the
true exponent, its population value divided by $w_0$ is
\[
 1+\frac{0.2(\alpha-1)}{\alpha+1}
 \quad\text{for }W=w_0(Z+0.2/Z),\qquad
 1\pm\frac{0.2(\alpha-1)}{\alpha}
 \quad\text{for }W=w_0(Z\pm0.2).
\]
These formulas follow from
$\E Z=\alpha/(\alpha-1)$ and $\E Z^{-1}=\alpha/(\alpha+1)$.
Thus the scale used by the Pareto predictor need not equal the
tail parameter $w_0$ under misspecification. Exponent estimation
and prediction of a finite component fraction consequently assess
different aspects of the fitted model.

The protocol projects scale estimates onto $[0.5,16]$ for prediction.
An invalid exponent or scale gives prediction $1$ and remains in
the analysis. Exponent errors use the values before projection.
Scale errors compare \eqref{eq:nr-estimated-scale} with the common
tail parameter $w_0$. No fit is invalid, no projection is needed,
and no discrete likelihood optimum reaches a boundary in this run.

\subsubsection{Calibration and paired uncertainty}
The predictor uses one new Pareto calibration shared by all four weight laws.
For each $n$, it contains $80$ independent graphs at each of the
$34$ exponents used in Sec.~\ref{sec:nr-percolation-experiment},
all generated at weight scale $16$. To evaluate a candidate scale
$0<w\le16$, transform each calibration crossing by
\begin{equation}\label{eq:nr-scale-crossing}
 \mathcal Q_{n,f}(w)
 \stackrel{\mathrm d}{=}
 \min\{1,16\mathcal Q_{n,f}(16)/w\}.
\end{equation}
Here the argument records the weight scale, and each crossing uses
the cap at one when the target is unreached.
To verify \eqref{eq:nr-scale-crossing}, condition on the Pareto draws.
Multiplying every weight by $w/16$ multiplies every Poisson edge
mean by $w/16$. Retaining only calibration edges whose uniform labels
are at most $w/16$ therefore produces the graph at scale $w$.
Rescaling these retained labels by $16/w$ yields independent uniform
labels, and gives the stated crossing, including the unreached case.

We average the transformed crossings and interpolate linearly in
the exponent between consecutive grid points. Scale requires no
interpolation. The calibration remains Pareto under every generating
weight law, so the experiment tests a misspecified predictor as well
as a modified degree distribution. No calibration graph is reused
as a training or test graph.

The $2000$ bootstrap repetitions resample the $300$ replication
indices within each $(n,w_0,\alpha)$, keeping the four laws, all
methods, both scale conventions, and all component fractions together.
They independently resample the $80$ calibration graphs at each
$(n,\alpha)$ grid point. Each resampled calibration is shared across
all methods, weight laws, and prediction settings with that $n$.
Predictions and losses are recomputed after both resampling steps.
The $2.5$ and $97.5$ percentiles give the reported intervals, with
pointwise intervals for individual settings. Losses and pooled
percentage reductions use \eqref{eq:percolation-loss} and
\eqref{eq:percolation-reduction}.

Independent calibration checks use $100$ Pareto graphs at each
$\alpha\in\{1.85,2.25,2.65,3.05\}$ and
$w_0\in\{1,3,6\}$, for both sizes. Across these settings and
fractions, the largest absolute discrepancy between an interpolated
calibration mean and a validation mean is $0.531$ percentage points.
The largest discrepancy divided by its combined calibration and
validation standard error is $2.46$, and one of the $72$ comparisons
exceeds two in magnitude. Using the first and last $40$ calibration
graphs separately gives primary correction reductions of $11.21$
and $11.69$ percent. The prescribed calibration remains unchanged.
Every calibration, validation, and test graph reaches all three targets.

\subsubsection{Results under the modified weights}
Tables~\ref{tab:nr-settings-1} and \ref{tab:nr-settings-2} retain all
$64$ primary correction comparisons. Table~\ref{tab:nr-robust-pooled}
reports every method under estimated and known scale, together with
exponent and scale RMSE. Table~\ref{tab:nr-robust-fractions} gives the
prescribed secondary fractions. Full setting-level results for all
methods remain in \path|summary.csv| and \path|exponent_summary.csv|.
The correction lowers pooled exponent MSE under each weight law.
Under estimated scale, its prediction MSE falls in $16$, $3$, $4$,
and $16$ of the $16$ settings for the four laws in
\eqref{eq:nr-robust-weight-laws}, respectively. The adverse settings
remain in every relevant pooled comparison.

\subsubsection{Crossing fitted exponents and scales}
The reversal between exponent error and prediction error motivates
changing one fitted parameter at a time.
This exploratory analysis reuses all saved training/test pairs and
calibration samples. It crosses the five exponent sources with the
five scale estimates and the known $w_0$, holding each scale column
fixed across exponent rows. The analysis plan was saved after the
original robustness outcomes were known and before computing the
crossed predictions. It supplies a decomposition of the existing
comparison, without independent validation.

Table~\ref{tab:nr-crossed} reports the four combinations of corrected
and discrete fits for the primary component fraction. Replacing only
the exponent or only the scale raises prediction error, and replacing
both lowers it. The isolated replacements move the mean prediction
error in opposite directions. The fitted pair therefore benefits
from compensation between exponent and scale in the Pareto predictor.
The result does not assign the improvement uniquely to either
parameter, or separate systematic bias compensation from covariance
between fitted estimates.
\begin{table}[H]
\centering
\small
\caption{Crossed exponent and scale estimates for the ten-percent component target, pooling the $48$ non-Pareto settings equally. All entries use the same saved test outcomes and Pareto calibration. RMSE and mean prediction error are in percentage points of edge survival probability.}
\label{tab:nr-crossed}
\begin{tabular}{llrr}
\toprule
Exponent & Scale & RMSE & Mean error\\
\midrule
Corrected & Corrected & 1.7457 & -0.1350\\
Discrete & Corrected & 3.3319 & 2.4539\\
Corrected & Discrete & 2.8916 & -1.0887\\
Discrete & Discrete & 1.6636 & 1.1433\\
\bottomrule
\end{tabular}
\end{table}

The bootstrap retains every crossed combination together with the
coupling across weight laws and the resampling of calibration graphs
specified in this section. The saved matrices include each setting,
weight law, and component fraction. Under both weight changes for
which the correction worsens prediction, changing the raw exponent
to the corrected one at the fixed raw scale already
increases prediction error. Updating the scale partly offsets that
increase. The adverse result consequently cannot be attributed
solely to recomputing scale after correcting the exponent.

\begin{table}[tbp]
\centering
\small
\caption{All primary correction comparisons for the indicated NR weight laws, with estimated scale and Pareto calibration. Raw RMSE is in percentage points. Correction entries give percentage MSE reductions [pointwise 95\% intervals] at $f=0.10$, using 300 pairs per setting. Negative entries retain adverse results.}
\label{tab:nr-settings-1}
\begin{tabular}{rrrlrl}
\toprule
$w_0$ & $\alpha$ & Raw, $10^4$ & Correction, $10^4$ & Raw, $10^5$ & Correction, $10^5$\\
\midrule
\multicolumn{6}{l}{$W/w_0=Z$}\\
2 & 1.8 & 1.975 & 13.0 [8.5, 17.5] & 1.034 & 15.2 [10.1, 20.3]\\
2 & 2.2 & 2.923 & 17.0 [12.3, 21.5] & 1.553 & 15.5 [9.9, 21.1]\\
2 & 2.6 & 2.972 & 11.3 [5.2, 17.4] & 1.768 & 16.3 [9.9, 22.8]\\
2 & 3.0 & 3.117 & 11.4 [5.7, 17.3] & 1.942 & 15.7 [9.5, 22.7]\\
4 & 1.8 & 0.903 & 14.9 [10.9, 18.7] & 0.436 & 13.8 [8.7, 19.2]\\
4 & 2.2 & 1.189 & 13.3 [7.6, 18.2] & 0.586 & 15.2 [9.3, 20.9]\\
4 & 2.6 & 1.208 & 8.6 [2.6, 14.2] & 0.706 & 13.3 [6.1, 19.5]\\
4 & 3.0 & 1.157 & 9.6 [3.3, 15.7] & 0.659 & 9.5 [2.2, 16.4]\\
\addlinespace
\multicolumn{6}{l}{$W/w_0=Z+0.2/Z$}\\
2 & 1.8 & 1.896 & 4.4 [-1.1, 9.7] & 1.034 & -8.2 [-14.6, -2.7]\\
2 & 2.2 & 2.700 & 5.5 [0.0, 10.7] & 1.456 & -15.3 [-21.8, -9.3]\\
2 & 2.6 & 2.639 & -10.8 [-19.1, -3.7] & 1.676 & -28.0 [-36.4, -20.6]\\
2 & 3.0 & 2.532 & -7.6 [-14.2, -2.0] & 1.805 & -25.8 [-33.7, -19.3]\\
4 & 1.8 & 0.856 & 3.9 [-0.7, 8.1] & 0.461 & -9.3 [-13.9, -5.0]\\
4 & 2.2 & 1.059 & -2.8 [-9.8, 3.4] & 0.621 & -21.9 [-27.1, -17.4]\\
4 & 2.6 & 1.097 & -14.9 [-22.6, -8.0] & 0.759 & -23.2 [-30.2, -17.1]\\
4 & 3.0 & 1.013 & -13.6 [-20.7, -7.3] & 0.730 & -37.7 [-46.1, -30.7]\\
\addlinespace
\bottomrule
\end{tabular}
\end{table}

\begin{table}[tbp]
\centering
\small
\caption{All primary correction comparisons for the indicated NR weight laws, with estimated scale and Pareto calibration. Raw RMSE is in percentage points. Correction entries give percentage MSE reductions [pointwise 95\% intervals] at $f=0.10$, using 300 pairs per setting. Negative entries retain adverse results.}
\label{tab:nr-settings-2}
\begin{tabular}{rrrlrl}
\toprule
$w_0$ & $\alpha$ & Raw, $10^4$ & Correction, $10^4$ & Raw, $10^5$ & Correction, $10^5$\\
\midrule
\multicolumn{6}{l}{$W/w_0=Z+0.2$}\\
2 & 1.8 & 1.806 & 2.6 [-2.6, 7.9] & 1.027 & -8.3 [-14.1, -3.1]\\
2 & 2.2 & 2.447 & 4.1 [-2.1, 9.4] & 1.385 & -11.8 [-18.9, -5.7]\\
2 & 2.6 & 2.609 & -3.0 [-9.5, 3.1] & 1.582 & -18.2 [-27.4, -10.3]\\
2 & 3.0 & 2.439 & -3.2 [-10.1, 3.4] & 1.733 & -18.4 [-25.9, -11.9]\\
4 & 1.8 & 0.831 & 2.1 [-2.6, 6.1] & 0.443 & -10.4 [-14.9, -6.3]\\
4 & 2.2 & 1.069 & 2.6 [-3.4, 8.0] & 0.593 & -18.3 [-24.1, -13.3]\\
4 & 2.6 & 1.023 & -6.6 [-13.8, -0.1] & 0.689 & -20.0 [-28.3, -13.6]\\
4 & 3.0 & 0.980 & -8.8 [-15.2, -2.4] & 0.592 & -33.2 [-43.2, -24.9]\\
\addlinespace
\multicolumn{6}{l}{$W/w_0=Z-0.2$}\\
2 & 1.8 & 2.392 & 20.5 [16.6, 24.4] & 1.450 & 28.8 [25.5, 32.3]\\
2 & 2.2 & 3.730 & 26.7 [23.2, 30.1] & 2.334 & 33.3 [30.1, 36.7]\\
2 & 2.6 & 4.034 & 23.3 [19.1, 27.6] & 2.911 & 34.2 [30.9, 37.4]\\
2 & 3.0 & 4.024 & 25.4 [20.8, 29.8] & 3.010 & 30.9 [27.4, 34.4]\\
4 & 1.8 & 1.076 & 20.5 [17.2, 23.9] & 0.632 & 25.1 [22.1, 28.2]\\
4 & 2.2 & 1.524 & 23.6 [19.7, 27.4] & 1.002 & 27.8 [24.7, 31.0]\\
4 & 2.6 & 1.634 & 22.8 [18.5, 27.0] & 1.170 & 30.2 [26.9, 33.5]\\
4 & 3.0 & 1.639 & 23.8 [19.0, 28.4] & 1.190 & 33.0 [30.1, 36.2]\\
\addlinespace
\bottomrule
\end{tabular}
\end{table}

\begin{table}[tbp]
\centering
\small
\caption{All NR methods pooled within each weight law and across the 48 non-Pareto settings. At $f=0.10$, prediction RMSE is in percentage points; reduction compares estimated-scale prediction MSE with raw Hill and gives its 95\% interval. The final columns give exponent RMSE and scale RMSE relative to the true $\alpha$ and $w_0$. Known-scale prediction uses the true $w_0$.}
\label{tab:nr-robust-pooled}
\begin{tabular}{llrlrrr}
\toprule
$W/w_0$ & Method & Est. scale & Reduction (\%) & Known scale & Exponent & Scale\\
\midrule
$Z$ & Raw & 1.730 & 0.0 [0.0, 0.0] & 3.007 & 0.325 & 0.200\\
 & Corrected & 1.609 & 13.5 [11.5, 15.6] & 2.646 & 0.256 & 0.177\\
 & Half rank & 1.903 & -21.0 [-25.6, -16.5] & 3.361 & 0.381 & 0.230\\
 & Discrete & 1.662 & 7.7 [0.8, 13.4] & 2.983 & 0.304 & 0.199\\
 & True exponent & 0.722 & 82.6 [80.8, 83.8] & 0.697 & 0.000 & 0.052\\
\addlinespace
$Z+0.2/Z$ & Raw & 1.567 & 0.0 [0.0, 0.0] & 3.275 & 0.327 & 0.393\\
 & Corrected & 1.635 & -8.9 [-11.2, -6.6] & 2.718 & 0.256 & 0.331\\
 & Half rank & 1.753 & -25.1 [-29.9, -20.3] & 3.599 & 0.377 & 0.406\\
 & Discrete & 1.311 & 30.1 [23.8, 35.3] & 4.551 & 0.448 & 0.552\\
 & True exponent & 0.999 & 59.3 [54.8, 62.8] & 0.972 & 0.000 & 0.264\\
\addlinespace
$Z+0.2$ & Raw & 1.494 & 0.0 [0.0, 0.0] & 3.733 & 0.345 & 0.521\\
 & Corrected & 1.531 & -5.0 [-7.3, -2.9] & 3.049 & 0.265 & 0.454\\
 & Half rank & 1.647 & -21.6 [-26.0, -17.1] & 3.964 & 0.393 & 0.527\\
 & Discrete & 1.331 & 20.6 [13.5, 26.8] & 5.253 & 0.468 & 0.702\\
 & True exponent & 0.980 & 57.0 [52.2, 60.8] & 1.375 & 0.000 & 0.367\\
\addlinespace
$Z-0.2$ & Raw & 2.374 & 0.0 [0.0, 0.0] & 2.939 & 0.301 & 0.390\\
 & Corrected & 2.031 & 26.8 [25.5, 28.2] & 3.115 & 0.249 & 0.439\\
 & Half rank & 2.534 & -13.8 [-16.9, -10.7] & 3.283 & 0.347 & 0.403\\
 & Discrete & 2.194 & 14.6 [10.7, 18.4] & 1.251 & 0.157 & 0.336\\
 & True exponent & 1.437 & 63.4 [60.0, 66.5] & 1.488 & 0.000 & 0.367\\
\addlinespace
Non-Pareto & Raw & 1.855 & 0.0 [0.0, 0.0] & 3.332 & 0.325 & 0.439\\
 & Corrected & 1.746 & 11.5 [9.9, 13.0] & 2.966 & 0.257 & 0.412\\
 & Half rank & 2.017 & -18.2 [-20.3, -16.0] & 3.626 & 0.373 & 0.449\\
 & Discrete & 1.664 & 19.6 [15.8, 23.2] & 4.077 & 0.385 & 0.551\\
 & True exponent & 1.158 & 61.0 [59.1, 62.4] & 1.298 & 0.000 & 0.336\\
\addlinespace
\bottomrule
\end{tabular}
\end{table}

\begin{table}[tbp]
\centering
\small
\caption{Secondary component fractions, pooled equally over all 48 non-Pareto settings. Entries give prediction RMSE in percentage points and percentage MSE reduction [95\% interval] relative to raw Hill under the same scale convention.}
\label{tab:nr-robust-fractions}
\begin{tabular}{llcc}
\toprule
Scale & Method & $f=0.05$ & $f=0.20$\\
\midrule
Estimated & Raw & 2.177; 0.0 [0.0, 0.0] & 1.398; 0.0 [0.0, 0.0]\\
Estimated & Corrected & 2.059; 10.6 [8.9, 12.2] & 1.314; 11.6 [10.2, 13.1]\\
Estimated & Half rank & 2.367; -18.1 [-20.3, -15.8] & 1.509; -16.4 [-18.3, -14.3]\\
Estimated & Discrete & 1.976; 17.6 [13.6, 21.5] & 1.235; 22.0 [18.5, 25.3]\\
Estimated & True exponent & 1.365; 60.7 [58.6, 62.2] & 0.912; 57.5 [55.6, 58.9]\\
Known & Raw & 3.414; 0.0 [0.0, 0.0] & 3.352; 0.0 [0.0, 0.0]\\
Known & Corrected & 3.038; 20.8 [19.3, 22.3] & 3.025; 18.6 [17.4, 19.7]\\
Known & Half rank & 3.727; -19.1 [-21.5, -16.9] & 3.606; -15.7 [-17.8, -13.7]\\
Known & Discrete & 4.123; -45.8 [-52.0, -40.3] & 4.142; -52.7 [-57.9, -47.9]\\
Known & True exponent & 1.273; 86.1 [85.2, 86.7] & 1.624; 76.5 [75.7, 77.3]\\
\bottomrule
\end{tabular}
\end{table}

\FloatBarrier
\section{Estimating the finite-degree coefficient from degrees}
\label{sec:degree-bias-estimation}

\subsection{A consistent estimate at smaller thresholds}
The correction in \arxiveqref{M-eq:corrected-estimator} requires the
coefficient $B$ of inverse degree in the population Hill ratio.
At thresholds of the selected order in PA and Pareto-weight NR,
$B/d$ and the sampling fluctuation have comparable magnitudes.
A smaller threshold supplies enough observations to estimate $B$
consistently.

Let $V_{n,1},\ldots,V_{n,n}$ be nonnegative integer observations and
let $(\pi_j)_{j\ge0}$ be a probability mass function whose tail
$\overline F(d)=\sum_{j>d}\pi_j$ satisfies
\eqref{eq:lattice-reference-tail}. The coefficient remains
$B=(Q+\alpha/2)/(\alpha+1)$. For an integer $d\ge1$, define the
Hill ratio over strict exceedances
\begin{equation}\label{eq:degree-bias-strict-hill}
 H_n(d)=\frac{\sum_{v:V_{n,v}>d}\log(V_{n,v}/d)}{N_{n,>d}},
 \qquad N_{n,>d}=\sum_v\mathbf1\{V_{n,v}>d\},
\end{equation}
and set $H_n(d)=0$ when the denominator is zero. Put
$u_n=\lceil\log n\rceil$ and, for $j=0,\ldots,8$, let
$d_{n,j}$ be the nearest integer to $u_n2^{j/8}$, choosing the larger
integer at a tie. These nine thresholds are distinct for all
sufficiently large $n$. Write $x_{n,j}=1/d_{n,j}$ and
$\bar x_n=9^{-1}\sum_{j=0}^8x_{n,j}$. The least-squares slope of
$H_n(d_{n,j})$ against $x_{n,j}$ is
\begin{equation}\label{eq:degree-bias-slope}
 \check B_n=\sum_{j=0}^8w_{n,j}H_n(d_{n,j}),\qquad
 w_{n,j}=\frac{x_{n,j}-\bar x_n}
                  {\sum_{\ell=0}^8(x_{n,\ell}-\bar x_n)^2}.
\end{equation}
The population version of
\eqref{eq:degree-bias-strict-hill} is
$1/\alpha+B/d+O(d^{-1-\min(\eta,1)})$, by
\eqref{eq:lattice-log-expansion} and \eqref{eq:lattice-reference-tail}. The zero sum of the slope weights removes the constant $1/\alpha$.

The following sufficient conditions separate the tail expansion from
the control of sampling error. For every integer $d\in[u_n,2u_n]$,
suppose that a constant $C$ independent of $n,d$ satisfies
\begin{align}
 &|\E N_{n,>d}-n\overline F(d)|
 +\left|\E\sum_{v:V_{n,v}>d}\log(V_{n,v}/d)
           -n\sum_{j>d}\pi_j\log(j/d)\right|\le C\log n,
 \label{eq:degree-bias-means}\\
 &\Var(N_{n,>d})+
 \Var\left\{\sum_{v:V_{n,v}>d}\log(V_{n,v}/d)\right\}
 \le C\{n\overline F(d)+(\log n)^2\}.
 \label{eq:degree-bias-variances}
\end{align}
These conditions yield consistency of the slope without requiring
independence among the observations.
\begin{proposition}\label{prop:degree-bias-consistency}
Under \eqref{eq:lattice-reference-tail},
\eqref{eq:degree-bias-means}, and \eqref{eq:degree-bias-variances},
$\check B_n\xrightarrow{p}B$ in $\R$.
The conditions hold for homogeneous PA with fixed $m\ge1$ and
$\delta>-m$, for the Pareto-weight NR model with fixed
$\alpha,w_0>0$, and for iid observations with the stated tail expansion.
For PA and Pareto-weight NR, replacing $\wh B_{n,k}$ by the single
statistic $\check B_n$ in the correction leaves the selected limits
in Corollary~\ref{cor:pa-corrected-selected} and
Theorem~\ref{NR-thm:selected-joint}, respectively, unchanged.
\end{proposition}
The coefficient can therefore be estimated from degrees in these
models. Interval calibration in PA uses its selected law and the known value of $m$.

\subsection{Proof of Proposition~\ref{prop:degree-bias-consistency}}
\begin{proof}
We bound the difference between each empirical Hill ratio and its
population value, apply the slope
weights, and verify the bounds for the three stated sampling laws.
The weights satisfy
\begin{equation}\label{eq:degree-bias-weights}
 \sum_jw_{n,j}=0,\qquad \sum_jw_{n,j}/d_{n,j}=1,
 \qquad \sum_j|w_{n,j}|\le C u_n.
\end{equation}
The last inequality follows by dividing the numerator and denominator
in \eqref{eq:degree-bias-slope} by $u_n^{-1}$ and $u_n^{-2}$,
respectively. The limiting denominator is the positive sum of squared
deviations of $2^{-j/8}$, $j=0,\ldots,8$.
Applying \eqref{eq:degree-bias-weights} to the population expansion
shows that its fitted slope differs from $B$ by at most
$Cu_n^{-\min(\eta,1)}$.

The empirical ratios approach their population values faster than
$u_n^{-1}$. Fix $\epsilon>0$. Chebyshev's inequality and a union bound
over the eighteen centered counts and logarithmic sums show that the
probability that any exceeds $\epsilon n\overline F(u_n)/u_n$ in
absolute value is at most
\begin{equation}\label{eq:degree-bias-probability}
 \frac{C u_n^2\{n\overline F(u_n)+(\log n)^2\}}
      {\epsilon^2\{n\overline F(u_n)\}^2}\longrightarrow0.
\end{equation}
Here $n\overline F(u_n)\sim nc^\alpha(\log n)^{-\alpha}$,
and $\overline F(d_{n,j})/\overline F(u_n)$ is bounded away from zero.
On the complementary event, \eqref{eq:degree-bias-means} bounds the
error of every ratio in \eqref{eq:degree-bias-strict-hill} by
$C\epsilon/u_n+C\log n/\{n\overline F(u_n)\}$ for all sufficiently
large $n$. Its denominators are then positive. Multiplication by
\eqref{eq:degree-bias-weights}, followed by $n\to\infty$ and
$\epsilon\downarrow0$, proves consistency.

For PA, \eqref{eq:pa-profile-error} bounds the count error by a
constant. Jensen's inequality applied to the increasing concave
function $x\mapsto\{1+\log(1+x)\}^2$ and the measure with masses
$|e_j(n)|$ bounds its logarithmic second moment by $C(\log n)^2$.
The first logarithmic moment is therefore at most $C\log n$.
Differentiating
\eqref{eq:pa-product-mgf} twice at zero bounds the variance of each
sum by twice the sum of its marginal variances. Truncation extends
this inequality to the logarithms. The tail bound
$\overline F(d)\le C d^{-\alpha}$ and integration of the tail give
\begin{equation}\label{eq:degree-bias-log-moment}
 \sum_{j>d}\pi_j\log^2(j/d)\le C\overline F(d).
\end{equation}
The discrepancy between the corresponding finite-graph second
moment and its stationary value is at most $C(\log n)^2$ by
\eqref{eq:pa-profile-error}. This proves
\eqref{eq:degree-bias-variances}. In NR the marginal degree law is
exactly $(\pi_j)$, so both errors in
\eqref{eq:degree-bias-means} vanish. Lemma~\ref{NR-lem:concentration}
and \eqref{eq:degree-bias-log-moment} give
\eqref{eq:degree-bias-variances}; truncation is justified by the finite
logarithmic second moment. For iid observations the mean errors also
vanish and independence gives the variance bound directly.

Finally, at the unrestricted selected rank in PA or Pareto-weight NR,
$\sqrt{\wh k}/D_{(\wh k+1)}$ converges in distribution to
$T^{1/2+1/\alpha}/c$ in $\R$, and the model-based coefficient
$\wh B_{n,\wh k}$ converges in probability to $B$ in $\R$.
Consequently the difference of the two corrected Hill statistics,
multiplied by $\sqrt{\wh k}$, converges in probability to zero in
$\R$. Both Hill statistics converge to $1/\alpha>0$, so taking
reciprocals and applying Slutsky's theorem proves the selected-limit
claim.
\end{proof}

\subsection{A fixed rule for finite samples}
The finite-sample rule uses the same slope, with
$u_n=\max\{\lceil\log n\rceil,2\min_{v:V_{n,v}>0}V_{n,v}\}$
and duplicate thresholds removed. This protects the fit against
using thresholds close to the smallest positive degree. It agrees
with the preceding construction eventually in PA and with
probability tending to one in NR and iid Zipf samples. In PA the
minimum is $m$. In NR and Zipf, the proportion of degree-one
observations converges in probability to a positive constant, by
the count variance bound. With duplicates removed, the averages and
sums defining the slope use the remaining thresholds.

To measure instability of the slope, define the empirical second
moment
\begin{equation}\label{eq:degree-bias-proxy}
 \check v_n=2\sum_{v=1}^n
 \left[\sum_{j:V_{n,v}>d_{n,j}}\frac{w_{n,j}}
                          {N_{n,>d_{n,j}}}
   \{\log(V_{n,v}/d_{n,j})-H_n(d_{n,j})\}\right]^2.
\end{equation}
A summand with zero denominator is set to zero. The factor two
matches the variance bounds of Lemmas~\ref{lem:pa-exponential} and
\ref{NR-lem:concentration} for sums with deterministic coefficients.
The empirical denominators and weights in
\eqref{eq:degree-bias-proxy} do not give a finite-sample variance
bound or a calibrated confidence interval.

At the selected rank, abbreviate $d=D_{(\wh k+1)}$ and
$h=H_{n,\wh k}+\theta_{n,\wh k}/d$ within this subsection.
The rule reports $(h-\check B_n/d)^{-1}$ only when all the following
conditions hold: at least four thresholds are distinct,
$2u_n\le d$, $N_{n,>2u_n}\ge50$, $h-\check B_n/d>0$,
$|\check B_n|/d\le h/2$, and
$\sqrt{\check v_n}/d\le H_{n,\wh k}/(2\sqrt{\wh k})$.
Otherwise it retains raw Hill and records why the correction was
withheld. These constants were fixed before the validation run.

Under PA and Pareto-weight NR, the probability of withholding tends
to zero. The selected degree divided by $\log n$ and the number of
observations exceeding $2u_n$ tend to infinity in probability, and
$h\to1/\alpha$ in probability in $\R$.
For the last condition, Cauchy--Schwarz in
\eqref{eq:degree-bias-proxy}, the count bounds in
\eqref{eq:degree-bias-probability}, and the logarithmic second moments
used in the proof imply $\check v_n\to0$ in probability in $\R$.
Explicitly, outside an event of probability at most $\epsilon+o(1)$,
Markov's inequality bounds the sum of the nine empirical logarithmic
second moments by
$C\{n\overline F(u_n)+(\log n)^2\}/\epsilon$, and every denominator
is at least $C^{-1}n\overline F(u_n)$. Thus
\[
 \check v_n\le
 \frac{C u_n^2\{n\overline F(u_n)+(\log n)^2\}}
      {\epsilon\{n\overline F(u_n)\}^2}\longrightarrow0.
\]
Together with the selected limit of $\sqrt{\wh k}/d$, this proves
that the last condition holds with probability tending to one.
The other conditions follow from consistency of $\check B_n$.
Retaining raw Hill on the withheld samples therefore leaves the
selected limits unchanged in these two graph models.

\subsection{Fresh validation and comparison methods}
The validation fixed the code, seeds, withholding conditions, and
reporting criteria before generating any production samples.
The primary design uses $n\in\{10^4,10^5\}$ and
$\alpha\in\{1.8,2.2,3\}$. PA uses $m\in\{1,4\}$.
NR uses $w_0\in\{2,4\}$ and each of the four weight laws in
\eqref{eq:nr-robust-weight-laws}. The iid comparator has
$\Prob(V=j)=j^{-\alpha-1}/\zeta(\alpha+1)$ for $j\ge1$.
There are $400$ fresh independent samples in each of these $66$
settings. Six secondary settings at $n=10^6$ use $200$ samples each:
PA at $(m,\alpha)=(1,3),(4,3)$, Pareto NR at
$(w_0,\alpha)=(2,2.2),(2,3)$, and Zipf at $\alpha=2.2,3$.
The complete run contains $27{,}600$ samples.

Every method receives the same degrees within a replication. The
comparisons are raw Hill, the tie adjustment, the degree-only rule,
the correction using the known generating family with its coefficient
evaluated at raw Hill, the correction using the true coefficient,
and the discrete likelihood with its own selected cutoff from
Sec.~\ref{sec:nr-robustness-experiment}. For Zipf the true
coefficient is $1/2$. For the shifted NR weights $W=w_0Z+s$ it is
$(1/2-s)/(\alpha+1)$, with $s=\pm0.2w_0$. Indeed,
$\Prob(W>w)=w_0^\alpha w^{-\alpha}\{1+\alpha s/w+O(w^{-2})\}$.
The mixed-Poisson identity
$\Prob(D>d)=\E\Prob(W>G_{d+1}\mid G_{d+1})$, where
$G_{d+1}$ has a Gamma distribution with shape $d+1$ and rate one,
and the Gamma-ratio expansion used in
\eqref{NR-eq:tail-expansion} give
$Q=(1-\alpha)/2-s$. The exponentially small probability of
$G_{d+1}$ lying in a fixed bounded interval does not affect this
expansion. The law $w_0(Z+0.2/Z)$ has $s=0$ at this order.
Thus the known-family comparator uses information about each
modified weight law that the degree-only rule does not receive.

We also compare a direct adaptation of the extended Pareto estimator
of Beirlant, Joossens, and Segers~\cite{BeirlantJoossensSegers2009},
using their equations (3.8)--(3.9) with second-order exponent $-1$.
For strict exceedances of an integer threshold $u$, let $h_n^+(u)$
be the mean of $\log\{V_{n,v}/(u+1/2)\}$ and let $r_n^+(u)$ be the
mean of $(u+1/2)/V_{n,v}$ over $V_{n,v}>u$. Substituting $-h_n^+(u)$ for the second-order parameter in the cited formulas
gives the reciprocal of
\[
 h_n^+(u)+\frac{\{1+2h_n^+(u)\}\{1+h_n^+(u)\}^2}{h_n^+(u)^2}
             \left\{r_n^+(u)-\frac1{1+h_n^+(u)}\right\}.
\]
It is evaluated at the selected threshold and at $u_n$ separately,
requiring at least three strict exceedances and a positive reciprocal.
This specifies the numerical adaptation; the cited asymptotic
results assume independent continuous observations.

Squared exponent errors use the generating $\alpha$, as in
\eqref{eq:mdsp-estimation-error}. Withheld corrections and invalid
comparator estimates retain raw Hill in the unconditional loss.
Invalid means nonfinite or nonpositive. The records also report
errors conditional on a valid estimate and, for the degree-only
rule, conditional on offering the correction. Settings receive
equal weight within each family. Overall pooling first gives PA,
NR, and Zipf equal weight. Paired percentile intervals use $2000$
bootstrap repetitions, resampling replication indices independently
within each setting and retaining all methods together.

The fixed criterion for applying the new correction to ICON required
that the overall MSE reduction have a positive lower $95$ percent
bootstrap endpoint, that corrections be offered in at least one
quarter of samples in each family, and that each family's upper
$95$ percent endpoint for the MSE ratio to raw Hill be at most
$1.10$. Each NR weight law is also reported separately. The criterion
failed because the Zipf offer fraction did not reach the required
minimum. No new corpus correction or accuracy claim is made from this
experiment.

\begin{table}[tbp]
\centering
\small
\caption{Primary exponent RMSE at $n=10^4,10^5$, pooled as specified in Sec.~\ref{sec:degree-bias-estimation}. Withheld or invalid estimates retain raw Hill. The degree-only rule uses no generating parameters. Known-family corrections use the appropriate coefficient for each weight law. Positive MSE reductions favor the degree-only rule over raw Hill.}
\label{tab:degree-bias-pooled}
\begin{tabular}{lrrr}
\toprule
Method & PA & NR & Zipf\\
\midrule
Raw Hill & 0.2307 & 0.3171 & 0.4646\\
Tie adjustment & 0.3077 & 0.2429 & 0.4517\\
Known family & 0.1832 & 0.2463 & 0.5009\\
True coefficient & 0.1641 & 0.2508 & 0.5009\\
Degree-only rule & 0.2170 & 0.3172 & 0.4647\\
Discrete likelihood & 0.3651 & 0.3567 & 0.0245\\
Extended Pareto, selected & 0.6088 & 11.1875 & 3349.9900\\
Extended Pareto, pilot & 0.6011 & 2.3115 & 21.9116\\
\midrule
\multicolumn{4}{l}{Degree-only rule: MSE reduction, percent [95\% interval]}\\
Equal-family pool & \multicolumn{3}{r}{1.60 [1.39, 1.82]}\\
PA & \multicolumn{3}{r}{11.47 [10.30, 12.64]}\\
NR & \multicolumn{3}{r}{-0.09 [-0.55, 0.32]}\\
Zipf & \multicolumn{3}{r}{-0.04 [-0.10, 0.02]}\\
NR, Pareto & \multicolumn{3}{r}{-0.69 [-1.86, 0.27]}\\
NR, body change & \multicolumn{3}{r}{0.31 [-0.53, 1.21]}\\
NR, positive shift & \multicolumn{3}{r}{1.03 [0.41, 1.70]}\\
NR, negative shift & \multicolumn{3}{r}{-1.39 [-2.18, -0.71]}\\
\midrule
Correction offered (percent) & 31.08 & 27.85 & 0.33\\
\bottomrule
\end{tabular}
\end{table}

\begin{table}[tbp]
\centering
\small
\caption{Secondary exponent RMSE at $n=10^6$, with $200$ samples per setting. The parameter after PA is $m$ and after NR is $w_0$; both NR settings use Pareto weights. The last column is the percentage of samples offered the degree-only correction.}
\label{tab:degree-bias-secondary}
\begin{tabular}{lrrrrrr}
\toprule
Family & $\alpha$ & Raw & Degree-only & Known family & Discrete & Offered\\
\midrule
PA, 1 & 3.0 & 0.2026 & 0.1958 & 0.1558 & 0.3203 & 33.5\\
PA, 4 & 3.0 & 0.1463 & 0.0543 & 0.0592 & 0.2687 & 100.0\\
NR, 2 & 2.2 & 0.1022 & 0.1019 & 0.0921 & 0.1913 & 91.0\\
NR, 2 & 3.0 & 0.2872 & 0.2822 & 0.2300 & 0.5083 & 12.0\\
Zipf & 2.2 & 0.1914 & 0.1928 & 0.1862 & 0.0026 & 3.5\\
Zipf & 3.0 & 0.4343 & 0.4343 & 0.4182 & 0.0042 & 0.0\\
\bottomrule
\end{tabular}
\end{table}

Table~\ref{tab:degree-bias-pooled} reports the primary results.
The degree-only rule improves PA estimation and has little effect
on NR or Zipf. Table~\ref{tab:degree-bias-secondary} reports the
larger samples without combining them with the primary design.
The known-family and discrete comparisons reverse their ordering
across families. The direct extended Pareto adaptations can produce
large positive exponents when the corrected reciprocal approaches
zero; all such values remain in the loss calculation.

Invalid estimates before fallback number $0$, $43$, and $152$
for the selected-threshold extended Pareto adaptation in PA, NR,
and Zipf, respectively, and $0$, $3$, and $419$ at the pilot threshold.
The unrestricted slope correction has $2$, $9$, and $564$ invalid
estimates. The other comparison methods have none. Returning to raw
Hill after withholding prevents invalid output from the degree-only
rule, without establishing accuracy of the unrestricted correction.
The complete setting-level estimates, coefficient errors, conditional
losses, and withholding reasons remain in \path|summary.csv| and the
individual records.

The saved secondary PA interval analysis uses the independent
$20{,}000$-path calibration, issuing an interval only when the
correction is offered and its exponent lies in $[1.1,4.5]$.
At $n=10^5$, $m=4$, $\alpha=3$, it issues $381$ of $400$ intervals,
of which $334$ cover, giving conditional coverage $87.7$ percent.
These conditional frequencies do not establish unconditional coverage.
All issuance and containment counts, including settings with no issued
interval, remain in \path|summary.csv|.

\bibliographystyle{apsrev4-2}
\makeatletter
\immediate\write\@auxout{\string\citation{apsrev42Control}}
\makeatother
\bibliography{mainNotes,prx-references}